\documentclass[11pt,reqno]{amsart}
\usepackage[T1]{fontenc}
\usepackage[utf8]{inputenc}
\usepackage{lmodern}
\usepackage{amsmath,amssymb,amsthm,mathtools,bm}
\usepackage{microtype}
\usepackage[margin=1.02in]{geometry}
\usepackage{booktabs,array,enumitem,tikz,fancyvrb}
\usetikzlibrary{arrows.meta}
\usepackage[hidelinks,pagebackref]{hyperref}

\usepackage{embedfile}
\newtheorem{theorem}{Theorem}[section]
\newtheorem{proposition}[theorem]{Proposition}
\newtheorem{lemma}[theorem]{Lemma}
\newtheorem{corollary}[theorem]{Corollary}
\theoremstyle{definition}
\newtheorem{definition}[theorem]{Definition}

\theoremstyle{remark}
\newtheorem{remark}[theorem]{Remark}
\numberwithin{equation}{section}
\allowdisplaybreaks[2]
\title[Falconer's conjecture through descent and chaining]
{Falconer's conjecture through descent and chaining}
\author{Will Burstein, Alex Iosevich, and Alexia Yavicoli}
\thanks{Alex Iosevich was partially supported by NSF grant
DMS-2506858 and a Simons Fellowship in Mathematics. A.~Y. was supported in part by the Natural Sciences and Engineering Research Council of Canada (NSERC), grants GR030571 and GR030540.}
\date{\today}
\begin{document}
\begin{abstract}
OpenAI proved that a compact set $E\subset{\mathbb R}^d$ of Hausdorff
dimension greater than $\frac{d}{2}$ determines a set of Euclidean
distances of positive length. We give a complete proof of this theorem
which follows the overall architecture of the OpenAI argument, from the
construction of the initial pair measures to the final Fourier estimate,
but which modifies that argument substantially at the points where the
geometric costs are paid and where the estimates are reassembled. The
reductions are organized by the descent and chaining viewpoint of
Burstein, Iosevich, and Krause: one potential, built from two separated
depths together with a penalty for insufficient separation, pays for
every reduction on every interval of scales, the terminal intervals
included, and a finite-record chaining lemma keeps the fixed power losses
attached to the few reductions which incur them. In odd dimensions the
radial improvement is obtained directly at each finite scale from a
quantitative contact interval, chosen by minimizing a mass curve against
lines of two slopes, with explicit margins in place of a compactness
argument, and the conditional estimates on the three parts of that
interval are multiplied by a lemma on full source fibers. A random sign
identity represents the coarse pair filter by products of endpoint
multipliers and reduces the cost of reassembling the spatial pieces
inside the shell estimate. The argument gives more than positive length.
Since the mass first excluded near heavy hyperplanes at each late scale
is summable with a power rate, all the Fourier estimates compare with one
fixed measure, and the distance set supports a probability measure with
an $L^p$ density for some $p>1$ and with power Fourier decay, so that
$\Delta(E)$ has positive Fourier dimension. We also show that
a planar $N$-point set in the unit square with averaged $s$-energy at
most $A$, for some $1<s<2$, determines at least $c_{s,A,\varepsilon}N^{1-\varepsilon}$
distances separated by $N^{-1}$, for every $\varepsilon>0$, which generalizes a result of
Guth, Solomon, and Wang.
The introduction closes with a discussion of further questions suggested
by Falconer's sharpness construction. For the $l^d$ metric on
${\mathbb R}^d$, $d\ge3$, the elementary lattice count gives an
obstruction at dimension one and leads to the open question of the
density of sums of $d$ positive $d$th powers, and we also discuss the
unresolved sharp thresholds for distance sets on compact Riemannian
manifolds and for dot product sets.
\end{abstract}
\maketitle
\tableofcontents
\section{The distance set and the question behind the proof}\label{sec:introduction}

Let $E$ be a compact subset of ${\mathbb R}^d$. Its distance set is
$$ \Delta(E)=\{|x-y|:x,y\in E\}. $$ 

Falconer's conjecture asks how large $E$ must be to guarantee that
$\Delta(E)$ has positive length. OpenAI proved the conjecture in all
dimensions in the manuscript released on October 6, 2026
\cite{OpenAI2026}. The proof given in this paper keeps the architecture
of that argument, from the construction of the initial pair measures to
the last Fourier estimate, and it keeps the three mechanisms which drive
it, namely the directional tests, the graph reductions, and the payment
of geometric costs by two separated depths. It departs from that argument
in the way the costs are accounted for, in the way the radial improvement
is obtained in odd dimensions, and in the way the spatial pieces are
reassembled, and these departures are described in
section~\ref{sec:changes} once the objects involved have been introduced.
At the last stage we obtain information about the distribution of the
distances themselves.

The sharpness construction also guides the discussion of what remains
open in section~\ref{sec:future}, where the Euclidean count is compared
with the count for the $l^d$ metric, which leads to a question about sums
of powers, and where we then turn to distance sets on compact Riemannian
manifolds and to dot product sets, all of which remain open at the
thresholds suggested by the lattice examples.

Here and throughout the paper, $d$ denotes the dimension of the Euclidean
space under discussion and $|x-y|$ denotes the Euclidean distance between
$x$ and $y$. We write $\dim_H$ for the Hausdorff dimension and ${\mathcal L}^n$
for the Lebesgue measure on ${\mathbb R}^n$, so that ${\mathcal L}^1$ measures
length. Applied to a finite set, $|\cdot|$ denotes instead the cardinality of
that set, which we also denote by $\#$. The context will always make it clear
which meaning is intended.

\vskip.125in 

How does the exponent $\frac{d}{2}$ come about? Let $q\ge2$ be a positive integer, $s \in (0,d)$, and let $E_{q,s}$ denote the closed $q^{-\frac{d}{s}}$ neighborhood of $\frac{1}{q}\{{\mathbb Z}^d\cap{[0,q]}^d\}$, where ${\mathbb Z}^d$ is the integer lattice. Let $q_1=2$ and take $q_{i+1}>q_i^i$. A result due to Falconer~\cite{Falconer1985} shows that $E_s=\bigcap_{i=1}^{\infty} E_{q_i,s}$ has Hausdorff dimension $s$; see also~\cite[section~1.2]{IosevichRudnev2005}. It is not difficult to see that, for a constant $C>0$ depending only on $d$, 
$$ {\mathcal L}^1(\Delta(E_{q_i,s})) \leq Cq_i^{-\frac{d}{s}} \cdot q_i^2,$$ since the quadratic form $x_1^2+\dots+x_d^2$ takes on at most $dq_i^2+1$ values on ${\mathbb Z}^d \cap {[0,q_i]}^d$. It follows that the Lebesgue measure of $\Delta(E_s)$ is $0$ if $s<\frac{d}{2}$. With a bit of number theoretic input, namely that $x_1^2+x_2^2$ takes on at most $C\frac{q^2}{\sqrt{\log(q)}}$ values on ${\mathbb Z}^2 \cap {[0,q]}^2$~\cite{Landau1908}, one eliminates the endpoint in two dimensions. 

\begin{figure}[!htbp]
\centering
\begin{tikzpicture}[x=1cm,y=1cm,>=Stealth,font=\small]
\node[align=center,text width=14.5cm] at (6.9,6.25)
  {Write $\delta=q^{-\frac{d}{s}}$. The drawing on the left uses
   $d=2$, $q=4$, $s=1$, so $\delta=\frac{1}{16}$.};
\node[align=center,text width=5.4cm] at (2,5.3)
  {{(a) The thickened lattice}\\
   Grid centres $a,b$ and points $x,y$ in their balls};
\node[align=center,text width=7.2cm] at (10.15,5.3)
  {{(b) One covering interval}\\
   The distance axis is enlarged};
\draw[gray!60,thin] (0,0) rectangle (4,4);
\node[anchor=east,inner sep=1pt] at (-.3,.17) {$(0,0)$};
\node[anchor=west,inner sep=1pt] at (4.18,4.22) {$(1,1)$};
\foreach \u in {0,...,4} {
  \foreach \v in {0,...,4} {
    \filldraw[fill=gray!15,draw=gray!55,thin] (\u,\v) circle[radius=.25];
    \fill (\u,\v) circle[radius=.022];
  }
}
\draw[dashed,thick] (1,1) -- (3,2);
\draw[thick] (.82,1.10) -- (3.18,1.89);
\fill (1,1) circle[radius=.035] node[below right,inner sep=2pt] {$a$};
\fill (3,2) circle[radius=.035] node[above left,inner sep=2pt] {$b$};
\fill (.82,1.10) circle[radius=.035] node[above left,inner sep=2pt] {$x$};
\fill (3.18,1.89) circle[radius=.035] node[below right,inner sep=2pt] {$y$};
\draw[->,thin] (4,3) -- (4.25,3);
\draw[thin] (4.22,3.03) -- (4.58,3.36);
\node[anchor=west,inner sep=1pt] at (4.58,3.36) {$\delta$};
\draw[thin] (0,-.32) -- (0,-.63);
\draw[thin] (1,-.32) -- (1,-.63);
\draw[<->,thin] (0,-.51) -- node[below] {$\frac{1}{q}$} (1,-.51);
\node[anchor=west] at (1.55,-.51) {grid spacing};
\node[align=center,text width=7.5cm] at (10.15,3.68)
  {$||x-y|-|a-b||\le |x-a|+|y-b|\le2\delta$};
\fill[gray!20] (7.65,1.8) rectangle (12.65,2.2);
\draw[->] (7.05,2) -- (13.3,2);
\draw[thick] (7.65,1.77) -- (7.65,2.23);
\draw[thick] (12.65,1.77) -- (12.65,2.23);
\draw (10.15,1.82) -- (10.15,2.18);
\node[above] at (7.65,2.28) {$|a-b|-2\delta$};
\node[above] at (12.65,2.28) {$|a-b|+2\delta$};
\node[below] at (10.15,1.71) {$|a-b|$};
\fill (11.42,2) circle[radius=.045];
\draw[thin] (11.42,1.91) -- (11.72,1.21);
\node[below,inner sep=1pt] at (11.72,1.21) {$|x-y|$};
\draw[thin] (7.65,.24) -- (7.65,.56);
\draw[thin] (12.65,.24) -- (12.65,.56);
\draw[<->,thin] (7.65,.4) -- node[below] {length $4\delta$} (12.65,.4);
\node[align=center,text width=14.5cm] at (6.9,-1.48)
  {At most $dq^2+1$ grid distances, each covered by an interval of length $4\delta$\\[4pt]
   $\displaystyle{\mathcal L}^1(\Delta(E_{q,s}))\le4\delta(dq^2+1)
      \le Cq^{2-\frac{d}{s}}$.};
\end{tikzpicture}
\caption{The sharpness example. Moving each endpoint by at most $\delta$
changes the distance by at most $2\delta$. For an integer $n$ between
$0$ and $dq^2$, each grid distance has the form $\frac{\sqrt n}{q}$.
The covering intervals may overlap; their total length still gives an
upper bound. In the plane, Landau's estimate reduces the number of
intervals to at most $C\frac{q^2}{\sqrt{\log(q)}}$. At $s=1$, their total
length is at most $\frac{C}{\sqrt{\log(q)}}$, which tends to zero along
the sequence $q=q_i$.}
\label{fig:lattice-sharpness}
\end{figure}
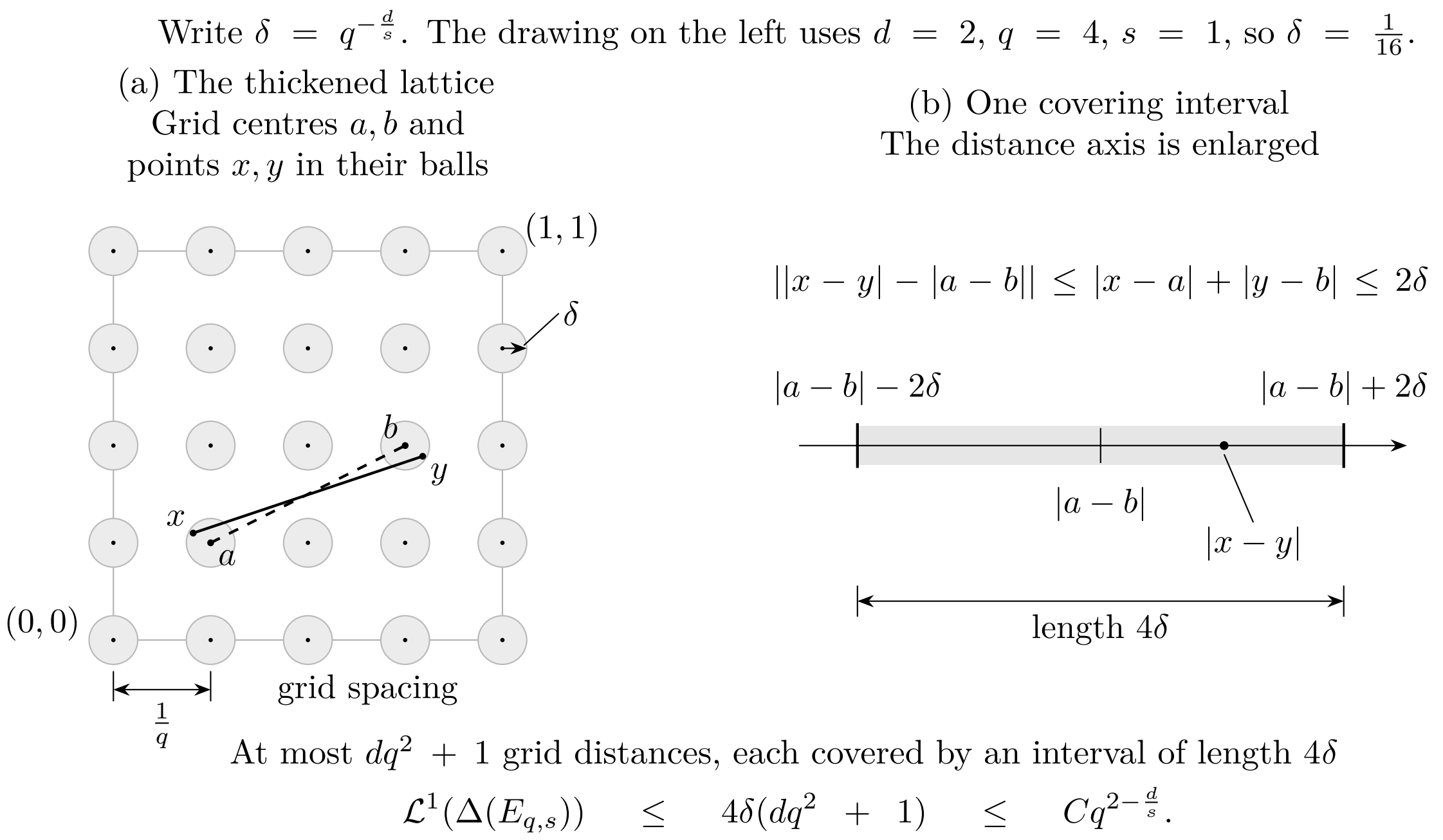

The strict inequality in the following theorem of OpenAI
\cite{OpenAI2026} agrees with the threshold in the lattice example.

\begin{theorem}[OpenAI]\label{thm:main}
Let $d\ge 2$ and let $E\subset {\mathbb R}^d$ be compact. If
$$
\dim_H(E)>\frac{d}{2},
$$
then $\Delta(E)$ has positive Lebesgue measure.
\end{theorem}

For a finite measure $\lambda$ on $\mathbb R^n$ we use the convention
$$
 \widehat\lambda(\xi)=\int e^{-ix\cdot\xi}d\lambda(x).
$$
The Fourier dimension $\dim_F F$ of a compact set $F\subset\mathbb R$
is the supremum of the exponents $u\in[0,1]$ for which $F$ supports a
probability $\nu$ satisfying
$|\widehat\nu(r)|\le C(1+|r|)^{-\frac{u}{2}}$ for every $r\in\mathbb R$.
Sets of positive length can have Fourier dimension zero
\cite[examples~7 and~18]{EkstromPerssonSchmeling2015}, so a positive
Fourier dimension for $\Delta(E)$ does not follow from
theorem~\ref{thm:main} itself, and to obtain it we return to the
particular measure built in the proof. Throughout the
paper $f\lesssim g$ means that $f\le Cg$ for a constant $C$ which does
not depend on the variables displayed, and the dependence of $C$ on fixed
parameters is specified when it matters. One more function space is
needed to state the second theorem.
Fix a smooth function $\vartheta$ on the line which equals one on
$[-1,1]$ and vanishes outside $[-2,2]$, let $Q_N$ be the Fourier
multiplier with symbol $\vartheta(2^{-N}r)$, and put $P_N=Q_N-Q_{N-1}$,
so that $P_N$ is a smooth cutoff to the frequency band
$2^{N-1}\le|r|\le2^{N+1}$. For $0<\delta<1$, the Besov space
$B^\delta_{1,\infty}({\mathbb R})$ consists of the integrable functions
$f$ for which
$$
\|Q_0f\|_{L^1}+\sup_{N\ge1}2^{\delta N}\|P_Nf\|_{L^1}<\infty,
$$
and this condition is equivalent to the integral H\"older condition
$\|f(\cdot+h)-f\|_{L^1}\le C|h|^\delta$ for every real $h$. The space
does not depend on the choice of $\vartheta$
\cite[chapter~2]{BahouriCheminDanchin2011}.

\begin{theorem}[A distance measure with power Fourier decay]
\label{thm:fourier-distance}
Let $E\subset\mathbb R^d$ be compact, where $d\ge2$ and
$\dim_H E>\frac{d}{2}$. Then $\Delta(E)$ supports a probability measure
$\nu=fdt$ and there is a $\delta>0$ such that
$$
 f\in B^\delta_{1,\infty}(\mathbb R) \text{ and } |\widehat\nu(r)|\le C(1+|r|)^{-\delta}
 \text{ for every } r\in\mathbb R.
$$
The density belongs to $L^p(\mathbb R)$ for every $1\le p<\infty$
with $1-\frac{1}{p}<\delta$. In particular
$\dim_F\Delta(E)\ge\min\{1,2\delta\}>0$.
\end{theorem}

The Besov condition says that the piece of the density living at
frequencies comparable to $2^N$ has $L^1$ norm $O(2^{-\delta N})$, and
this is the form in which it is proved in \eqref{asm:band-l1}. The
quantitative estimate which makes it possible is the bound for the mass
first excluded near heavy hyperplanes at a late scale, which is
\cite[equation~(5.16)]{OpenAI2026} in the source manuscript; summing
that bound over the late scales gives the rate at which the retained pair
measures approach their limit, and this rate brings every Fourier estimate back
to one fixed measure.

\subsection{What is kept from the OpenAI argument and what is changed}
\label{sec:changes}

The reader will want to know at the outset which parts of the argument
below are the ones in \cite{OpenAI2026} and which are not. The overall
plan is the same. One prepares two separated pieces of $E$ carrying
measures whose radial projections satisfy a cap estimate with the
critical exponent $\frac{d}{2}$, which in odd dimensions requires an
angular improvement through a planar incidence theorem; one divides each
measure into regular pieces on which the mass of a cube is governed by a
mass curve; one tests directions at every scale; one converts a count of
nearly equal distances into spherical Fourier energies by stationary
phase and back; and one passes to smaller cells by a graph estimate whose
cost is an excess of the mass curve, paying that cost with the values of
the mass curve at two separated depths. The directional tests, the graph
reductions, and the idea of paying with two separated depths are
OpenAI's, and so is the heavy-plate construction of the pair filters.

The argument differs from the source in four respects, each of which
changes the proof rather than its presentation, and the first of them is
the accounting. A depth counts dyadic subdivisions, and OpenAI's calculation
pays for a reduction with the minimum of a sum of two excess values at
depths separated by a fixed amount in the current interval, so that the
short intervals on which the iteration stops have to be handled
separately. In the descent viewpoint of Burstein, Iosevich, and Krause
\cite{BursteinIosevichKrause2026} one asks instead for one quantity which
pays for every reduction through its change at that reduction, and the
potential $\mathcal V_h$ of section~\ref{ledger:section}, which adds to
the two excess values a penalty proportional to the amount by which their
separation falls short, is such a quantity on every interval, the
terminal ones included. The same normalization is then used throughout,
the whole iteration becomes a single instance of the weighted descent and
chaining lemma of section~\ref{bik:section}, and that lemma keeps the
fixed power losses attached to the few reductions which incur them, so
that the losses at the many low-frequency reductions remain $R^{o(1)}$.

The second is the radial improvement in odd dimensions, where, writing
$d=2k-1$, the exponent $k-1$ supplied by finite energy has to be raised
to $k-\frac{1}{2}$ on the middle portion of a depth interval, and the
interval has to be long enough for the angular gain and short enough for
the radial map to be linearized on it. We choose this interval at each
finite scale by two minimizations of the mass curve, first against a line
of the larger slope and then against a line of the smaller one, and
lemma~\ref{prep:contact} gives explicit margins for its length and its
position which do not depend on the particular mass curve. The one-scale
improvement, lemma~\ref{prep:one-scale}, therefore works directly with
the finite mass curves, with every loss chosen before the scale, and no
compactness argument for families of mass curves is needed. The three
conditional estimates on the three parts of the interval are multiplied
by lemma~\ref{reg:pruning-multiplication}, which records the full
conditional source estimate that both radial arguments use.

The third difference lies in the reassembly of the spatial pieces, since
in odd dimensions the retained pair set at scale $2^N$ is a union of
rectangles from two finite spatial partitions and the Fourier estimate
for the retained distance measure has to be proved with this filter in
place. The identity \eqref{asm:signed-filter} represents the filter as an
average of products of random signed multipliers at the two endpoints,
so that the estimate with arbitrary bounded endpoint multipliers, which
the shell conversion supplies, can be applied directly, at a cost of
$\sqrt{3M}$ in the Fourier norm, where $M$ bounds the number of cells in
either partition, in place of the cost of a term by term summation over
the rectangles.

The fourth difference is in what is proved, since the power tail of the
discarded pair mass, which is already present in the heavy-plate
construction of the source, is used here to compare every shell
approximant with one fixed measure in total variation, and lemma~\ref{asm:fixed-shell} turns
the two estimates for the approximants into the Besov regularity and the
Fourier decay of theorem~\ref{thm:fourier-distance}. The discrete
consequences of section~\ref{discrete:section} are new as well, and they
use a quantitative planar length bound of Cheng, Fang, Li, Qu, and
Xiao~\cite{ChengFangLiQuXiao2026} rather than the odd-dimensional
construction.

All the proofs are given in full, including the proofs of the
intermediate estimates which are taken over from the source, since the
descent accounting requires the support conditions and the normalizations
of those estimates to be stated exactly.

\subsection{Start with the conclusion we want}

Let us work backwards. Suppose that we could place a nonzero measure
$\nu$ on $\Delta(E)$ and prove that it has a density $f$ with respect to
the Lebesgue measure. Then
$$
0<\nu({\mathbb R})=\int_{\Delta(E)}f(t)dt
$$
would imply that $\Delta(E)$ has positive length. The calculation is
trivial, but it tells us what to build, namely a nonzero measure on the
distance set which possesses a density.

It is convenient to use the word law as another name for a finite positive
Borel measure, with no normalization implied by the word itself; a
probability law has total mass one, and a subprobability has total mass at
most one. The operations we shall perform on such measures are the familiar
ones. For finite measures $\lambda,\eta$ and Borel sets $A,B$, the product
$\lambda\otimes\eta$ assigns to the rectangle $A\times B$ the mass
$\lambda(A)\eta(B)$, and if $F$ is a measurable map, the pushforward
$F_\#\lambda$ is given by $F_\#\lambda(B)=\lambda(F^{-1}(B))$. We write
$\lambda|_A$ for the restriction of $\lambda$ to a Borel set $A$, and
$1_A$ for the function equal to one on $A$ and zero elsewhere, while
$g\lambda$ denotes more generally the measure with density $g$ relative to
$\lambda$. None of these operations is normalized, and therefore a restriction
becomes a probability only when we explicitly divide by its mass. For an
integer $j\ge1$, $\lambda^{\otimes j}$ is the product of $j$ copies of
$\lambda$. The marginals of a measure on pairs are its pushforwards under
the two coordinate projections; for a product of two probabilities these
are the two original measures, but restricting to a set of pairs need not
preserve the product structure. In an ordered pair $(x,y)$ we call $x$ the
pin, from which the source point $y$ is viewed, and the pin law and the
source law are the measures on these two coordinates before any
pair-dependent deletion has taken place. Fixing one coordinate in a set of
pairs leaves a fiber in the other coordinate, and reversing the orientation
interchanges the pin and the source. Conditioning a measure
$\lambda$ on a set $Q$ of positive mass means passing to the probability
$\frac{\lambda|_Q}{\lambda(Q)}$.

There is a natural candidate for the measure we are after. Choose a
probability measure $\mu$ on $E$, take two independent points with this
distribution, and record their distance. In formulas, with $D(x,y)=|x-y|$,
this measure is
$$
\nu=D_\#(\mu\otimes\mu),
$$
where $D_\#$ denotes the pushforward, so that for a Borel set
$B\subset {\mathbb R}$,
$$
\nu(B)=(\mu\otimes\mu)\{(x,y):|x-y|\in B\}.
$$

The model measure $\nu$ accounts for all pairs, whereas we shall use its
variant on two separated pieces and keep only the pairs for which the
Fourier estimates can be proved, so that the resulting density reflects
the distribution of the retained distances.

A few more conventions will be needed. A decreasing sequence of sets is one
in which each set contains the next. Absolute continuity is written
$\lambda\ll\eta$ and means that every $\eta$-null set is also
$\lambda$-null, and for a positive measure we sometimes write
$\|\lambda\|$ for the total mass. For positive measures, $\lambda\le\eta$
means that $\lambda(B)\le\eta(B)$ for every Borel set $B$, so that a
decreasing sequence of measures loses mass without moving the mass that
remains. The index $N$ below is a positive integer which will tend to
infinity. We say that two sets are separated if the distance between them
has a positive lower bound, and a pair filter is a measurable set which
specifies which ordered pairs are kept. We shall work with two separated
pieces $E_1,E_2$ of $E$, carrying probabilities $\mu_1,\mu_2$ supported on
them, the point of the separation being that the distance function is
smooth on the product. Write $\lambda=\mu_1\otimes\mu_2$. The retained
pair sets $\Gamma_N$ decrease as $N$ increases, and their intersection
$\Gamma_\infty=\bigcap_N\Gamma_N$ satisfies $\lambda(\Gamma_\infty)>0$.
On this product the weight $W(x,y)=|x-y|^{-\frac{d-1}{2}}$ is bounded
above and below by positive constants, and its role is to cancel the
leading spherical stationary phase factor. Define
$$
 \alpha_N=D_\#(W1_{\Gamma_N}\lambda) \text{ and } \alpha_\infty=D_\#(W1_{\Gamma_\infty}\lambda).
$$
In even dimensions the prepared pair set is already the full product,
while in odd dimensions \eqref{prep:coarse-tail} bounds the mass remaining
to be discarded by $\lambda(\Gamma_N\setminus\Gamma_\infty)
\lesssim 2^{-c_\Gamma N}$ for some $c_\Gamma>0$, the first scale in this
bound being allowed to depend on the measures.

At frequency $2^N$ a further angular selection gives a positive measure
$\tau_N\le\alpha_N$ whose discarded mass is bounded by
$C2^{-\eta N}$ with $\eta>0$. Recall that the total variation norm
of a signed measure $\rho$ is $\sup\sum_j|\rho(B_j)|$, the supremum
being over finite Borel partitions, and that for a positive measure it is
the mass. Thus the loss in passing from $\alpha_N$ to $\alpha_\infty$ and
the loss in passing from $\alpha_N$ to $\tau_N$ give, together with the
Fourier estimate,
$$
 \|\alpha_\infty-\tau_N\|_{\mathrm{TV}}\lesssim 2^{-eN}
 \text{ and }
 \int_{2^{N-1}\le |r|\le2^{N+1}}|\widehat{\tau_N}(r)|^2dr
 \lesssim 2^{-\gamma N}
$$
for some $e,\gamma>0$. The first estimate compares every approximant
to the same nonzero measure, whereas the second is needed only on the
shell for which $\tau_N$ was chosen.

The convolution kernel of the band cutoff $P_N$ defined above has
bounded $L^1$ norm, so that the first estimate also bounds
$\|P_N(\alpha_\infty-\tau_N)\|_1$, while on a fixed interval containing
all the distances, Cauchy--Schwarz and Plancherel give
$\|P_N\tau_N\|_{L^1}\lesssim 2^{-\frac{\gamma N}{2}}$, with the norm
taken on that interval. Outside a unit enlargement of the distance
interval the kernel estimate for $P_N\alpha_\infty$ gives rapid decay,
and consequently
$$
 \|P_N\alpha_\infty\|_1\lesssim 2^{-\delta N} \text{ with } \delta=\min\left\{e,\frac{\gamma}{2}\right\}.
$$
Summing these bands constructs the density, and since at any fixed large
frequency only a bounded number of bands contribute, the same inequality
gives power decay of its Fourier transform. The details are in
lemma~\ref{asm:fixed-shell}.

\subsection{What should the descent remember?}

The measure remains fixed during the descent; what changes is the estimate
we are trying to prove and the interval of scales still available to prove
it. A depth is an integer $n$ between $0$ and $N$.

A dyadic cube in the usual grid is a Cartesian product of half-open
intervals $[j2^{-n},(j+1)2^{-n})$ with integer $j$, the half-open endpoints
ensuring that the cubes at a fixed depth form a partition; translated grids
will be specified when they are needed. A dyadic cube of depth $n$ has side
length $2^{-n}$, and we call such a cube occupied if it has positive mass.
A regular piece is a normalized restriction on which all the occupied cubes
at a given depth have comparable masses, and these pieces are constructed
in section~\ref{reg:section}. In this paragraph $m:[0,1]\to[0,d]$ denotes
their nondecreasing mass curve, with $m(0)=0$ and Lipschitz constant at
most $d$, so that $|m(u)-m(v)|\le d|u-v|$. Fix a frequency $R=2^N$ and
suppose that, on a regular piece of the measure, an occupied cube of side
length $2^{-n}$ has mass approximately
$$
M(n)=R^{-m(\frac{n}{N})},
$$
the approximation permitting a controlled small power of $R$. The function
$m$ records the decay of the mass as the cubes become smaller. Set
$$
S=\frac{d}{2} \text{ and } A_n=m\left(\frac{n}{N}\right)-\frac{Sn}{N}.
$$
Here $S$ is the critical dimension and $A_n$ is the excess at
depth $n$. The ratio $\frac{n}{N}$ is a normalized depth in
$[0,1]$, whereas $n$ itself is an integer depth, and the minimum in the
potential below is taken over integer depths, not over the continuum. The
number $A_n$ measures how far the decay of the mass exceeds the critical
decay $2^{-Sn}$, and we call $m$ the mass curve and $A$ its excess. The
excess need not increase with $n$, and the iteration must therefore
remember more than its current scale. Let us check the sign, since we
shall rely on it repeatedly. The formula for the model mass can be written
equivalently as
$$
M(n)=2^{-Sn}R^{-A_n},
$$
and therefore a positive excess means that there is less mass in a cell than the critical
model predicts, which is a saving available to an estimate, whereas a
negative excess means that the critical model has underestimated the mass
of the cell. The restrictions performed later need not make these savings
increase from one depth to the next.

We write $a$ and $t$ for integer depths with $0\le a\le t\le N$, these
being the starting and the ending depths of the current estimate, and we
write $a'$ and $t'$ for the corresponding depths of a child estimate, the
word child referring to a term on the right side of a reduction inequality,
as in the tree of estimates in section~\ref{bik:section}. By a spherical
Fourier energy we mean an average of the squared modulus of a Fourier
integral over unit directions, and the positive angular energy used in the
descent combines two such averages at the same radial frequency. The proof
moves back and forth between two quantities, one of which counts pairs of
distances that differ by at most $2^{-t}$, while the other is a positive
product of spherical Fourier energies. Let $X$ denote either quantity and
let $[a,t]$ be its current interval of depths. The analytic estimates
produce children $X'$ whose starting depth is $a'\ge a$, and, apart from
controlled errors, such a child is multiplied by
$$
R^{-2(A_{a'}-A_a)+c},
$$
where $c\ge0$ is the cost of that reduction. The first exponent is merely a
change of normalization; the cost $c$ is the part that must be paid.

Choose a small $\kappa>0$, put $q_0=\lceil\kappa N\rceil$ and let
$h=10q_0$. The bound $|A_{n+1}-A_n|\le\frac{d}{N}$ suggests charging
$\frac{d}{N}$ for each unit by which a pair of depths falls short
of separation $h$. With $u_+=\max(u,0)$, put
$$
 \mathcal V_h([a,t])=
 \min_{a\le i\le j\le t}
 \left\{A_i+A_j+\frac{d}{N}(h-(j-i))_+\right\}.
$$
All indices here are integers. If $t-a\ge h$, any pair with gap
less than $h$ can be moved apart one unit at a time inside $[a,t]$.
At each move the penalty drops by $\frac{d}{N}$ and the two excess
values together rise by at most that amount. Thus on a long interval
$\mathcal V_h$ equals OpenAI's minimum of $A_i+A_j$ over pairs
separated by at least $h$. On every interval the choice $i=j=a$ gives
$\mathcal V_h([a,t])-2A_a\le\frac{dh}{N}$.

Now let $J=[a,t]$ be a parent and $J'=[a',t']$ a child. Restricting an
interval removes candidates from the minimum, so
$\mathcal V_h(J')\ge\mathcal V_h(J)$. This pays for the transitions of
zero cost. For each transition with positive cost $c$, the geometric
choice in lemma~\ref{ledger:payment} gives
$$
 c\le \mathcal V_h(J')-\mathcal V_h(J).
$$
This is a descent inequality for the nonnegative quantity
$\Phi(J)=4d-\mathcal V_h(J)$, since $\kappa$ is small and $N$ is
sufficiently large. To see how it enters an estimate, normalize the
parent and child by
$$
 Y=R^{\mathcal V_h(J)-2A_a}X \text{ and } Y'=R^{\mathcal V_h(J')-2A_{a'}}X'.
$$
In the normalized estimate the child term becomes
$$
 R^{\mathcal V_h(J)-2A_a}R^{2A_a-2A_{a'}+c}X'
   = R^{c+\mathcal V_h(J)-\mathcal V_h(J')}Y'\le Y'.
$$
The geometrically determined cost has disappeared, and only the small
analytic losses remain; lemma~\ref{bik:tree} keeps track of which
reductions incur them.

One minimum of $A$ cannot pay for all the reductions, since a low value at
the left endpoint can remain the minimum even when a reduction has a
positive cost on the right. The remedy is to minimize the sum of two values
at separated depths, for then one of them can be moved to pay for the
reduction while the other preserves the required separation. This is proved
explicitly for every type of reduction in section~\ref{ledger:section}. We
call an initial portion of a depth interval a prefix and a final portion a
suffix, so that a suffix reduction is one which removes a final portion and
keeps the interval to its left. A numerical example shows why the second
depth is needed. For this example only, replace the separation $h$ by $2$,
take a number $\varepsilon>0$, and consider the seven excess values
$$
(A_0,A_1,A_2,A_3,A_4,A_5,A_6)
=(0,\varepsilon,2\varepsilon,2\varepsilon,
\varepsilon,\varepsilon,0).
$$
Use the minimum of $A_i+A_j$ with $j-i\ge 2$ on the two intervals
below, which both have width at least $2$. Denote this minimum by $\mathcal V_2$.
The
minimum of $A$ is zero on both $[0,6]$ and $[0,4]$, and it cannot pay the
suffix cost $A_4-\min_{4\le n\le6}A_n=\varepsilon$. The pair $(0,6)$,
however, gives $\mathcal V_2([0,6])=0$, whereas $(0,4)$ gives
$\mathcal V_2([0,4])=\varepsilon$, and no admissible pair in $[0,4]$ has a
smaller sum. The increase at two depths pays exactly the cost that the
single minimum misses. This is an illustration of the algebra and not an
instance of the full choice rule, whose intervals are much longer.

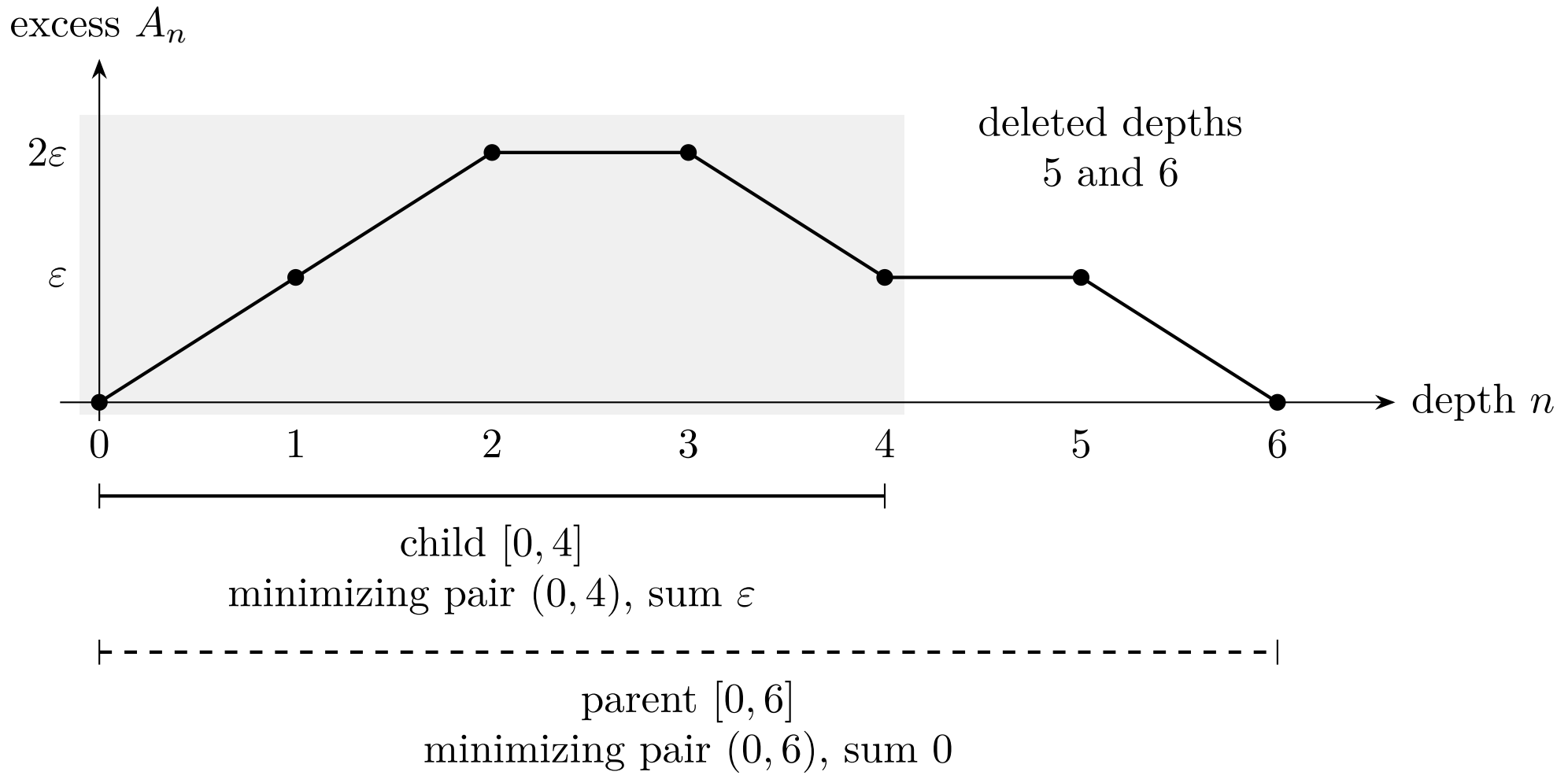
\begin{figure}[htbp]
\centering
\begin{tikzpicture}[x=1.65cm,y=1.05cm,>=Stealth,font=\small]
\fill[gray!12] (-.1,-.1) rectangle (4.1,2.3);
\draw[->] (-.2,0)--(6.6,0) node[right] {depth $n$};
\draw[->] (0,-.15)--(0,2.75) node[above] {excess $A_n$};
\draw[thick] (0,0)--(1,1)--(2,2)--(3,2)--(4,1)--(5,1)--(6,0);
\foreach \x/\y in {0/0,1/1,2/2,3/2,4/1,5/1,6/0}
{\fill (\x,\y) circle (2pt);\node[below] at (\x,-.1) {$\x$};}
\node[left] at (-.08,1) {$\varepsilon$};
\node[left] at (-.08,2) {$2\varepsilon$};
\draw[thick] (0,-.75)--(4,-.75);
\draw (0,-.65)--(0,-.85) (4,-.65)--(4,-.85);
\node[below,align=center] at (2,-.85)
{child $[0,4]$\\minimizing pair $(0,4)$, sum $\varepsilon$};
\draw[thick,dashed] (0,-2)--(6,-2);
\draw (0,-1.9)--(0,-2.1) (6,-1.9)--(6,-2.1);
\node[below,align=center] at (3,-2.1)
{parent $[0,6]$\\minimizing pair $(0,6)$, sum $0$};
\node[align=center] at (5.15,2.05)
{deleted depths\\$5$ and $6$};
\end{tikzpicture}
\caption{The shaded interval is the child. The depth $0$ remains a
single minimizer after the deletion, so its value records no cost.
The separated pair must change, and its sum increases by
$\varepsilon$. The required separation is $2$ in this example.}
\label{fig:two-depth-example}
\end{figure}
\subsection{How the pieces fit together}

A handful of terms will be used repeatedly, and it is convenient to record
them before the argument begins. By a cap we mean a ball on the unit
sphere, and a cap estimate is a bound for the measure of such a ball in
terms of its radius. The directional tests below ask how many cells lie
near a line and how often scalar projections bring two points close
together, and we refer to the latter events as projection collisions. By a
mask we mean a function with values in $[0,1]$, smooth in the direction
variable, which equals one on a good set and vanishes outside a slightly
larger good set, the space between the two leaving room for smoothing; we
call this space the buffer. The degree of a vertex of a graph is, as usual,
the number of edges incident to it.

The descent and chaining lemma is proved in section~\ref{bik:section},
before any of the geometric details appear. Sections~\ref{sec:conventions}
and~\ref{reg:section} introduce the measures, the regular pieces, and the
mass estimates, while sections~\ref{ang:section} and~\ref{prep:section}
prepare directions with the required cap estimates. In the plane, and more
generally in even dimensions, the integer cap estimate is enough, and the
extra angular improvement is needed only in odd dimensions.
Section~\ref{test:section} turns the directional bounds into tests which
remain valid during the descent, and sections~\ref{osc:section}
and~\ref{graph:section} relate distance frequency bands to angular energies
and then to smaller cells. The resulting reductions, together with the
conditions they preserve, are recorded in section~\ref{rec:section}.
Section~\ref{ledger:section} uses the same potential on every interval
to pay the cost of each reduction. In section~\ref{asm:section} the original
class weights are restored, an average of products of signed endpoint
multipliers recombines the coarse spatial pieces inside the shell estimate,
and the geometric tail of the discarded pair mass allows us to compare
each shell approximant with one fixed measure, whose dyadic bands are
summable in $L^1$ and whose rate of decay gives positive Fourier dimension
for the distance set. Finally, section~\ref{discrete:section} uses the
quantitative planar estimate of Cheng, Fang, Li, Qu, and
Xiao~\cite{ChengFangLiQuXiao2026} to pass from bounded discrete energy
to separated distances.

The reader who is interested first in the connection with descent can read
sections~\ref{bik:section}, \ref{rec:section}, and
\ref{ledger:section}, treating the analytic reductions as inputs on that
first reading; all of their proofs are retained in the preceding sections.

The contact argument in lemma~\ref{prep:contact} makes one part of this
route shorter. We subtract a line from a mass curve and choose the place
where the difference is smallest, and then repeat with a second slope.
These two choices produce an interval which is long enough for the angular
gain and short enough for radial linearization, and lemma~\ref{prep:one-scale}
uses the resulting numerical margins directly at each finite scale. The
idea is related to the lower barriers and the finite witnesses of Iosevich and
Yavicoli \cite{IosevichYavicoli2026}, though the proof we give here is
elementary and self-contained. This step belongs to the additional
preparation needed in odd dimensions.

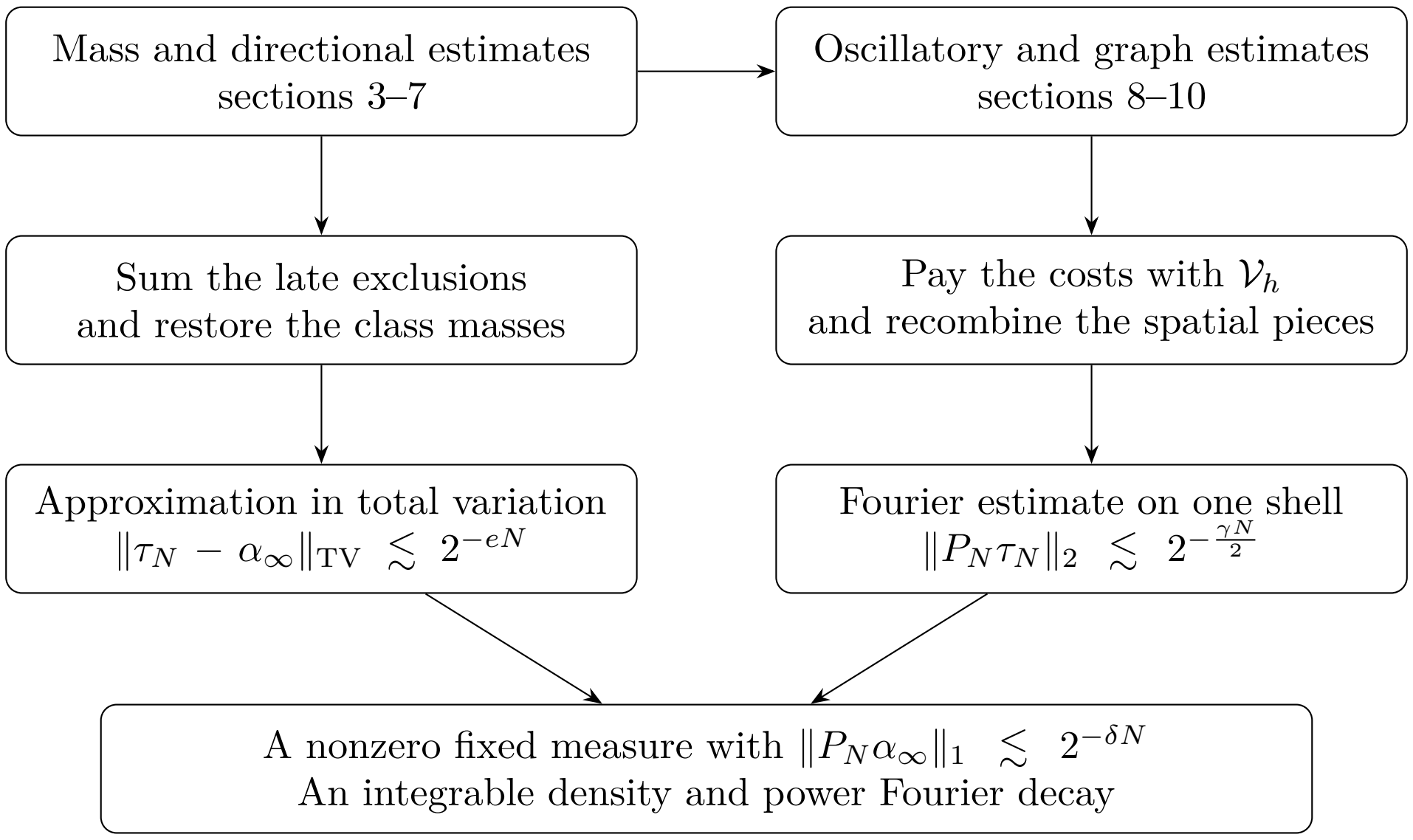
\begin{figure}[htbp]
\centering
\begin{tikzpicture}[>=Stealth,font=\small,
box/.style={draw,rounded corners,align=center,text width=5.3cm,
minimum height=1.15cm,inner sep=5pt}]
\node[box] (geom) at (-3.45,0)
{Mass and directional estimates\\sections~\ref{sec:conventions}--\ref{test:section}};
\node[box] (analytic) at (3.45,0)
{Oscillatory and graph estimates\\sections~\ref{osc:section}--\ref{rec:section}};
\node[box] (retain) at (-3.45,-2.05)
{Sum the late exclusions\\and restore the class masses};
\node[box] (descent) at (3.45,-2.05)
{Pay the costs with $\mathcal V_h$\\and recombine the spatial pieces};
\node[box] (error) at (-3.45,-4.1)
{Approximation in total variation\\
$\|\tau_N-\alpha_\infty\|_{\mathrm{TV}}\lesssim 2^{-eN}$};
\node[box] (shell) at (3.45,-4.1)
{Fourier estimate on one shell\\
$\|P_N\tau_N\|_2\lesssim 2^{-\frac{\gamma N}{2}}$};
\node[box,text width=10.5cm] (limit) at (0,-6.25)
{A nonzero fixed measure with $\|P_N\alpha_\infty\|_1\lesssim2^{-\delta N}$\\
An integrable density and power Fourier decay};
\draw[->] (geom)--(analytic);
\draw[->] (geom)--(retain);
\draw[->] (analytic)--(descent);
\draw[->] (retain)--(error);
\draw[->] (descent)--(shell);
\draw[->] (error)--(limit);
\draw[->] (shell)--(limit);
\end{tikzpicture}
\caption{The two estimates for the shell approximant $\tau_N$.
The geometric tail supplies the fixed target $\alpha_\infty$ in
\eqref{asm:fixed-tv}. The descent supplies the Fourier estimate in
proposition~\ref{asm:shell-decay}. Lemma~\ref{asm:fixed-shell} combines
them on a fixed interval containing the distance set.}
\label{fig:proof-map}
\end{figure}

\begin{table}[htbp]
\centering\small
\renewcommand{\arraystretch}{1.35}
\begin{tabular}{p{3.4cm}p{4.5cm}p{5.4cm}}
\toprule
Part of the proof & What it supplies & Where it is used\\
\midrule
Regularization, lemma~\ref{reg:regularization}
& Comparable cell masses and a Lipschitz mass curve
& Cell counts, conditional ball bounds, and the excess $A_n$.\\
Contact interval, lemma~\ref{prep:contact}
& Two minima with uniform scale margins
& The direct one-scale radial improvement in odd dimensions.\\
Pair preparation, proposition~\ref{prep:starting-pairs}
& Positive limiting mass, a power tail, and cap exponent $S$
& Comparison with a fixed distance measure and control of exceptional directions.\\
Directional tests, section~\ref{test:section}
& Tube counts, projection collisions, and smooth buffered masks
& Graph degrees, suffix reduction, and support preservation.\\
Shell conversion, proposition~\ref{osc:shell-conversion}
& A bound for a distance-frequency band by angular energies
& High reductions and the outer shell estimate.\\
Graph estimate, proposition~\ref{graph:estimate}
& Passage to smaller cells, with explicit excess cost
& High and direct reductions, and the common entry depth.\\
Descent, theorem~\ref{ledger:bound}
& A bound for each single-mass-curve angular energy
& The two-class shell estimate after Cauchy--Schwarz.\\
Signed recombination, \eqref{asm:signed-filter}
& A coarse pair filter represented by products of endpoint tests
& The Fourier shell estimate for $\tau_N$.\\
Fixed measure, lemma~\ref{asm:fixed-shell}
& Decay of the dyadic $L^1$ norms in \eqref{asm:band-l1}
& Positive length and positive Fourier dimension of $\Delta(E)$.\\
Uniform planar bound, lemma~\ref{discrete:uniform-length}
& A lower bound for length in terms of the planar energy
& Thickening and interval covering in theorem~\ref{discrete:main}.\\
\bottomrule
\end{tabular}
\vskip.125in 
\caption{A guide to the dependencies. The complete descent calculation can be
read first with the analytic estimates treated as inputs, and then checked
against their proofs.}
\label{tab:dependencies}
\end{table}

\vskip.125in 

\subsection{What the reader should know before beginning}

The argument combines several subjects, but it does not require a course in
each of them. A reader who knows graduate real analysis and the basic
properties of the Fourier transform should be able to work through the
paper with the preparation described below. In going through it, it is
useful to distinguish between a theorem whose conclusion we use and a
technique whose proof we develop here, and the references are meant as a
guide to that distinction.

We use finite Borel measures, restriction and normalization, product
measures, pushforward, the theorems of Tonelli and Fubini, and the
Radon--Nikodym theorem, and our geometric starting point is Frostman's
lemma. We write $B(x,r)=\{y:|y-x|<r\}$ for the open Euclidean ball, and we
recall that a probability is a positive measure of total mass one, whereas a
subprobability has total mass at most one. If a compact set has Hausdorff
dimension greater than $s$, it supports a probability $\mu$ for which
$$
\mu(B(x,r))\le Cr^s,
$$
and a probability satisfying this bound for every $0<r\le1$ is called an
$s$-Frostman probability. The Riesz potential of $\mu$ at $x$ is
$U_s^\mu(x)=\int |x-y|^{-s}d\mu(y)$, and a bounded-potential set is one on
which this function has a finite uniform upper bound. The reader needs this
statement and the meaning of the Hausdorff dimension; the proof of
Frostman's lemma is not reproduced, and Mattila
\cite[theorem 8.8]{Mattila1995} is a convenient reference. The further
facts we need about finite energy, restriction to bounded-potential sets,
and typical orthogonal projections are proved in
section~\ref{sec:conventions}, and the dyadic regularization is proved in
section~\ref{reg:section}.

\vskip.125in 

From Fourier analysis we use inversion, Plancherel, convolution with smooth
cutoffs, and the rapid decay of the Fourier transforms of smooth compactly
supported functions. For a finite measure $\mu$ on ${\mathbb R}^n$ and a
real number $s>0$, the Riesz $s$-energy of $\mu$ is
$$
I_s(\mu)=\iint |x-y|^{-s}d\mu(x)d\mu(y),
$$
which may be infinite. Finiteness says that pairs of very close points do
not carry too much mass in the aggregate. The useful link with geometric
size is the Riesz-energy identity, valid for $0<s<n$,
$$
I_s(\mu)=c_{n,s}\int_{\mathbb R^n}
|\widehat\mu(\xi)|^2|\xi|^{s-n}d\xi,
$$
where $n$ is the dimension of the underlying space and $c_{n,s}>0$ is a
normalization constant. One should understand how this identity converts an
energy hypothesis into a weighted $L^2$ bound; see Mattila
\cite[chapter 12]{Mattila1995}. The constants depend on the Fourier
normalization and do not affect the exponents. The calculation with
Gaussian kernels and the annular estimates used later are given at their
point of use. In lemma~\ref{asm:fixed-shell}, smooth frequency cutoffs,
Cauchy--Schwarz on a bounded interval, and the decay of the cutoff kernels
give the $L^1$ estimate for each band. Their integrable sum is the density
of the fixed measure.

\vskip.125in 

For a point $x$ outside the source support, the radial projection is
$$
\pi_x(y)=\frac{y-x}{|y-x|}.
$$
The phenomenon that concerns us is that sufficient energy of the source and
the pin measures forces the radial projections to have $L^p$ densities for
some $p>1$. Orponen's work supplies the finite-energy method behind our
integer directional preparation
\cite[section 3 and the proof of theorem 1.13]{Orponen2019}, and the
precise version we use is lemma~\ref{prep:radial-lp}. We denote the unit
sphere in ${\mathbb R}^{h+1}$ by $S^h=\{w\in{\mathbb R}^{h+1}:|w|=1\}$ and
its surface measure by $\sigma$, and in the display that follows the
quotient of measures is a Radon--Nikodym density, $h\ge1$ is an integer and
$0<\epsilon<\min(h,1)$. Separated probabilities $\alpha,\beta$ in
$\mathbb R^{h+1}$ with finite $I_{h+\epsilon}$ energies have radial
densities satisfying
$$
\int\int_{S^h}
\left|\frac{d(\pi_x)_\#\alpha}{d\sigma}(w)\right|^p
d\sigma(w)d\beta(x)<\infty
$$
for some $p>1$, and the reverse orientation has the same property.
Truncating the two densities at height $T$ costs at most a constant
times $T^{1-p}$ in product mass. The remaining fibers have the required
pointwise cap bounds. For the integrated estimate it is enough to be
comfortable with $L^p$ duality, the inequalities of H\"older and Minkowski,
finite energy, and approximation of a measure by convolution.

\vskip.125in 

The preparation in odd dimensions uses a discretized incidence theorem of
Orponen and Shmerkin, in the tube formulation of Orponen, Shmerkin, and
Wang \cite[theorem 1.3]{OrponenShmerkin2023}
\cite[definitions 2.5--2.6 and theorem 2.7]{OrponenShmerkinWang2024}. This
is a deep external input; we state it in full as
theorem~\ref{ang:furstenberg} and do not reproduce its proof. At scale
$\delta>0$, a covering number is the least number of balls of radius
$\delta$ needed to cover a set, a set is $\delta$-separated if its distinct
points are at least $\delta$ apart, and a nonconcentration bound limits the
mass, or the relative number of points, in every small ball. An incidence
is a pair satisfying a specified relation, such as a point lying in a tube,
and its incidence mass is the product measure of the set of such pairs.
What the reader needs here is a covering number, a separated set at scale
$\delta$, and a nonconcentration bound for points or for line parameters.
The theorem says that a sufficiently spread-out collection of points, each
supporting a sufficiently spread-out family of incident tubes, forces a
fixed-power improvement in the total tube count. We call directions
transverse when they are linearly independent, and a quantitative
transversality bound places a stated positive lower bound on the
determinant or the volume which measures that independence. The passage
from measures to finite sets and the use of the theorem after slicing are
proved here, and at the endpoint we use three dimensional slices together
with a separate transverse tube argument, which is also proved here.

\vskip.125in 

The slicing argument uses averaging over rotations and orthogonal
subspaces, singular values of a matrix, and determinants as measures of
transversality, where Haar measure is simply the rotation invariant
probability law on the relevant family of subspaces. The band-probability
estimates we require, and all the restrictions of the slices, are proved in
lemmas~\ref{ang:haar} and~\ref{ang:slice-selection}. For vectors
$u_1,\ldots,u_j$, the exterior product $u_1\wedge\cdots\wedge u_j$ has norm
equal to the $j$-dimensional volume of the parallelepiped they span, or
equivalently the square of that norm is the determinant of their Gram
matrix $(u_i\cdot u_l)_{i,l=1}^j$, and it vanishes exactly when the vectors
are linearly dependent. Thus $u\wedge v$ records the area of the
parallelogram spanned by $u$ and $v$, and the three-vector version records
volume. The tube-counting inequality of Loomis and Whitney
\cite{LoomisWhitney1949} is proved in lemma~\ref{ang:loomis-whitney},
including the finite-entropy calculation used there.

\vskip.125in 

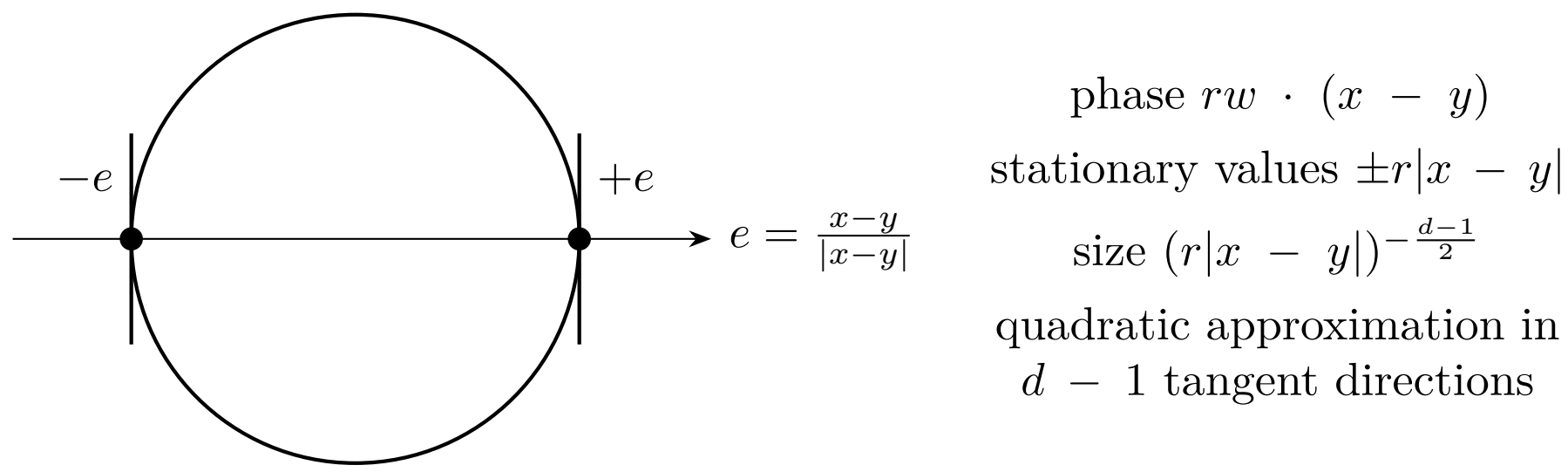
\begin{figure}[htbp]
\centering
\begin{tikzpicture}[>=Stealth,font=\small]
\draw[thick] (0,0) circle (1.7);
\draw[->] (-2.6,0)--(2.7,0) node[right] {$e=\frac{x-y}{|x-y|}$};
\fill (1.7,0) circle (2.5pt);
\fill (-1.7,0) circle (2.5pt);
\node[above right] at (1.7,.2) {$+e$};
\node[above left] at (-1.7,.2) {$-e$};
\draw[thick] (1.7,-.8)--(1.7,.8);
\draw[thick] (-1.7,-.8)--(-1.7,.8);
\node[align=center,text width=7cm] at (7,0)
{phase $r w\cdot(x-y)$\\[4pt]
stationary values $\pm r|x-y|$\\[4pt]
size $(r|x-y|)^{-\frac{d-1}{2}}$\\[4pt]
quadratic approximation in $d-1$ tangent directions};
\end{tikzpicture}
\caption{Stationary phase on the sphere, shown in a two-dimensional
cross-section. At $\pm e$ the displacement $x-y$ is normal to the
tangent space, so the first tangential derivative vanishes, and these two
points produce the positive and negative distance phases, whereas away
from them a tangential derivative of the phase is nonzero and permits
integration by parts. The distance weight in the shell conversion
is chosen to match the displayed stationary-phase factor.}
\label{fig:sphere-stationary}
\end{figure}

In an oscillatory integral, the phase is the real function in the
exponential and the amplitude is the factor multiplying that exponential,
while a stationary point is a point at which the derivatives of the phase
in the variables of integration vanish. The reader should know Taylor's
theorem and integration by parts, and should be familiar with the idea that
rapid oscillation produces cancellation when a derivative of the phase does
not vanish. Section~\ref{osc:section} derives the stationary-phase
expansion needed on the sphere, its inverse form, its derivative bounds,
and the separation of the remaining spatial coefficient, so that a general
stationary-phase theorem with parameters is not a separate prerequisite.
The source of the factor $(r|x-y|)^{-\frac{d-1}{2}}$ is explained in
figure~\ref{fig:sphere-stationary}.

\vskip.125in 

Both the elementary entropy-decrease calculation and the weighted descent
lemma used here are proved in section~\ref{bik:section}, and no restricted-isometry or random-selection theorem from Burstein--Iosevich--Krause is
assumed. Elsewhere we use compactness of a uniformly bounded Lipschitz
family and elementary weak convergence of finite measures, where weak
convergence means convergence of the integrals against every bounded
continuous function, and the contact argument in lemma~\ref{prep:contact}
needs only the existence of a minimum for a continuous function on a closed
interval. The relevant references are the real analysis text of Folland
\cite{Folland1999} and the measure theory preliminaries of Simon
\cite{Simon2014}. The parameter choices, the count of the reduction steps,
and the limiting measure argument are all carried out explicitly.

For a first reading, one can take Frostman's lemma, the Riesz energy
identity, and the stated planar Furstenberg theorem as known, and work
through the other specialized estimates in the order above. The next
subsection turns this list of prerequisites into a reading program, and
the historical references in section~\ref{sec:history} are meant to give
context rather than to add further prerequisites.

\subsection{A reading guide for a first-year graduate student}
\label{sec:reading-guide}

The preceding subsection describes what is assumed, and this one says what to
read, and in what order, so that a student at the beginning of graduate
study can follow the proofs in detail rather than take the inputs on
faith. The books named below are standard, and in each case only a few
chapters are needed.

The foundation is measure theory, and Folland \cite{Folland1999} is the
reference used throughout the paper. Chapters 1 through 3 contain
everything needed about measures, integration, product measures, the
theorems of Tonelli and Fubini and the Radon--Nikodym theorem, chapter 6
treats the $L^p$ spaces and their duality, chapter 7 proves the Riesz
representation theorem and the density of continuous functions, and
chapter 8 introduces the Fourier transform, the Schwartz space,
convolution and the Plancherel theorem. A student who has worked through
these chapters can read sections~\ref{sec:conventions}
and~\ref{reg:section} immediately, and the weak convergence and
disintegration arguments in section~\ref{prep:section} use nothing beyond
chapters 6 and 7.

Hausdorff dimension and the measures which detect it come next. Chapters
2 and 4 of Falconer's book \cite{Falconer2014} introduce Hausdorff measure
and dimension together with the mass distribution principle, which is the
elementary half of the Frostman correspondence, and chapter 8 of Mattila
\cite{Mattila1995} proves Frostman's lemma in the form used in
lemma~\ref{conv:energy} and introduces the Riesz energies. Mattila's
second book \cite{Mattila2015} is written for exactly the reader we have
in mind. Its chapter 2 reviews the measure theoretic preliminaries,
chapter 3 proves the energy identity which converts a Frostman bound into
a weighted $L^2$ estimate for the Fourier transform, and chapter 4 proves
Marstrand's projection theorem and Falconer's threshold $\frac{d+1}{2}$
for the distance problem. With these chapters in hand, the classical
threshold in lemma~\ref{prep:classical-plane} and the projection
statements in lemma~\ref{conv:energy} become routine, while the Haar
measure on the Grassmannian used in section~\ref{ang:section} is
constructed in chapter 3 of \cite{Mattila1995}.

The oscillatory estimates of section~\ref{osc:section} rest on the method
of stationary phase, and Wolff's lectures \cite{Wolff2003} are the best
first encounter with it. Chapters 1 through 3 of those lectures set up the
Fourier transform at the level of this paper, chapter 6 proves the
stationary phase expansion with the error terms made explicit, chapter 8
develops Hausdorff measures and Frostman's lemma from the harmonic
analyst's point of view, and chapter 9 presents Mattila's spherical
average formulation of the distance problem, which is the ancestor of the
Fourier estimates in sections~\ref{osc:section} and~\ref{asm:section}.
The fuller treatment of oscillatory integrals in chapter VIII of Stein
\cite{Stein1993} is needed for the sphere expansions of
lemma~\ref{osc:stationary}, and chapters 14 and 15 of \cite{Mattila2015}
explain the decay of the Fourier transform of surface measure and the way
in which spherical averages control distance sets, which is the
calculation behind lemma~\ref{test:half-dimensional-fourier} and
proposition~\ref{osc:shell-conversion}.

The descent calculation of section~\ref{bik:section} and the
Loomis--Whitney argument of lemma~\ref{ang:loomis-whitney} use entropy.
Chapter 2 of Cover and Thomas \cite{CoverThomas2006} contains the
definitions of entropy and relative entropy, the chain rule, and the fact
that conditioning does not increase entropy, which are the three facts
used in the proof of lemma~\ref{ang:loomis-whitney}, and chapter 17 of the
same book proves Han's inequality. The multiplicative weights update and
its relative entropy analysis are explained in sections 2.1 and 2.2 of
Arora, Hazan, and Kale \cite{AroraHazanKale2012}, and the student who then
reads section 4.1 of Frieze and Kannan \cite{FriezeKannan1999} and the
introduction of \cite{BursteinIosevichKrause2026} will see where the
potential of section~\ref{ledger:section} comes from. The probabilistic
inequalities used in the sampling arguments, namely Hoeffding's
exponential moment bound, the Chernoff tail bound for binomial counts and
Khintchine's inequality, are all treated in chapter 2 of Vershynin
\cite{Vershynin2018}, and they enter in section~\ref{bik:section}, in
lemma~\ref{ang:planar-measures} and in the signed recombination
\eqref{asm:signed-filter}.

The multiscale and incidence ideas have no textbook treatment yet, and
here one must read papers. Section 3 of Keleti and Shmerkin
\cite{KeletiShmerkin2019} is the model for the regularization of
lemma~\ref{reg:regularization} and for the use of mass curves, section 3
of Orponen \cite{Orponen2019} contains the finite energy method behind the
radial estimate of lemma~\ref{prep:radial-lp}, and the introductions of
Orponen and Shmerkin \cite{OrponenShmerkin2023} and of Orponen, Shmerkin,
and Wang \cite{OrponenShmerkinWang2024} explain what the discretized
Furstenberg theorem says and why it bears on radial projections. The
theorem itself is used here as a black box, in the form of
theorem~\ref{ang:furstenberg}, and a first reading of
section~\ref{ang:section} should accept it as such.

With this preparation a first pass through the paper can proceed as
follows. Read sections~\ref{sec:conventions} and~\ref{reg:section} for the
conventions and the regular classes, then section~\ref{bik:section} for
the descent lemma, then sections~\ref{rec:section}
and~\ref{ledger:section}, taking the reductions of
proposition~\ref{rec:reductions} as given, to see how the potential pays
for every step, and finally section~\ref{asm:section}, where the pieces
are assembled. The analytic inputs in sections~\ref{test:section}
through~\ref{graph:section} can then be read in order,
section~\ref{prep:section} after them, and the angular gain of
section~\ref{ang:section} last, since it is used only in odd dimensions
and is the most technical part of the argument.
Section~\ref{discrete:section} depends only on section~\ref{sec:conventions}
and on the planar length bound of Cheng, Fang, Li, Qu, and Xiao which it
quotes, and it can be read at any time.

\subsection{Relation to earlier work}
\label{sec:history}

The historical discussion below includes the dimension estimates which
preceded the public release of the source manuscript.

The problem originates in the work of Falconer~\cite{Falconer1985}, who
established the sufficient threshold $\frac{d+1}{2}$. The lattice
construction above rules out every threshold below $\frac{d}{2}$.
In the plane, Landau's count of sums of two squares makes the covering
length tend to zero even at dimension one. Mattila's spherical-average
criterion~\cite{Mattila1987} made Fourier $L^2$ estimates a central
approach to the problem, and Bourgain connected improvements in dimensions
two and three with Fourier restriction phenomena on
spheres~\cite{Bourgain1994}. Among the subsequent milestones were Wolff's
planar threshold $\frac{4}{3}$~\cite{Wolff1999} and Erdo\u{g}an's threshold
$\frac{d}{2}+\frac{1}{3}$ for $d\ge3$~\cite[theorem~2]{Erdogan2005}, where
once again the positive measure conclusions require strict inequalities
between dimensions. For the parallel problem for finite point sets in
${\mathbb R}^d$, $d\ge3$, see the distinct distances theorem
in~\cite[theorem~1.1]{OpenAIDistinctDistances2026}.

Weighted restriction and decoupling led to further improvements. Du, Guth,
Ou, Wang, Wilson and Zhang obtained the threshold $\frac{9}{5}$ in
dimension three, together with improvements in higher dimensions
\cite[theorem~1.2]{DuGuthOuWangWilsonZhang2021}, and Du and Zhang obtained
the threshold $\frac{d}{2}+\frac{1}{4}+\frac{1}{8d-4}$, which improves the
preceding bounds for $d\ge4$ \cite[theorem~2.6]{DuZhang2019}. These
arguments draw on refined restriction estimates and on the decoupling
theory of Bourgain and Demeter~\cite{BourgainDemeter2015}. For a pin $y$ we
write $\Delta_y(E)=\{|x-y|:x\in E\}$. Liu's $L^2$ identity connects pinned
distance measures with spherical extension estimates and transfers
spherical-average bounds to pinned
conclusions~\cite[theorems~1.4 and~1.9]{Liu2019}, and the two restriction
advances just mentioned also give ${\mathcal L}^1(\Delta_y(E))>0$ for some
$y\in E$ under the corresponding strict dimension inequalities; see
\cite[corollary~1.11]{DuGuthOuWangWilsonZhang2021}
and \cite[theorem~2.7]{DuZhang2019}. For compact $E\subset{\mathbb R}^2$
with $\dim_HE>\frac{5}{4}$, Guth, Iosevich, Ou and Wang proved the same
pinned conclusion \cite[theorem~1.2]{GuthIosevichOuWang2020}, their
argument combining radial-projection control of the discarded wave packets
with refined decoupling for the part that is retained. Du, Iosevich, Ou,
Wang and Zhang obtained the analogous pinned conclusion in even dimensions
$d\ge4$ when
$\dim_HE>\frac{d}{2}+\frac{1}{4}$~\cite[theorem~1.2]{DuIosevichOuWangZhang2021},
and for compact $E\subset{\mathbb R}^d$, $d\ge3$, Du, Ou, Ren and Zhang
obtained the same existential pinned conclusion when
$$
 \dim_HE>\frac{d}{2}+\frac{1}{4}-\frac{1}{8d+4}
$$
\cite[theorem~1.2]{DuOuRenZhang2024}. At the critical dimension itself,
Shmerkin and Wang proved that $\dim_H\Delta(E)=1$ for Borel
$E\subset{\mathbb R}^d$, $d\ge2$, whose Hausdorff and packing dimensions
both equal $\frac{d}{2}$~\cite[theorem~1.4]{ShmerkinWang2025}.

For planar sets of dimension just above one there are substantial results
which do not assume that the Hausdorff and packing dimensions agree.
Recall that an analytic set is a continuous image of a complete separable
metric space, and that every Borel set is analytic. Let
$E\subset{\mathbb R}^2$ be analytic and put $t=\dim_HE>1$. Stull proved
that
$$
 \dim_H\Delta_y(E)\ge\frac{1}{2}+\frac{t}{4}
$$
for every $y\in{\mathbb R}^2$ outside a set of Hausdorff dimension at most
one~\cite[theorem~1]{Stull2022}, and since that exceptional set cannot
contain $E$, the same lower bound holds for $\dim_H\Delta(E)$. His argument
uses effective dimension, which measures the information needed to describe
a point at increasingly fine precisions, together with radial projection
estimates. Fiedler and Stull subsequently obtained, in particular,
$$
 \sup_{y\in E}\dim_H\Delta_y(E)
 \ge\min\left\{1,\frac{t(4-t)}{5-t}\right\}
$$
\cite[corollary~2]{FiedlerStull2023}, while their more general theorem also
uses the packing dimension of $E$ to strengthen the conclusion
\cite[theorem~1]{FiedlerStull2023}. Both of the displayed bounds tend to
$\frac{3}{4}$ as $t$ decreases to one.

Iosevich and Yavicoli improved these estimates in a preprint submitted on 4
October 2026, two days before the public release of the OpenAI manuscript
on 6 October 2026~\cite{IosevichYavicoli2026,OpenAI2026}. They proved that
there is an absolute constant $0<\eta_0<\frac{3}{40}$ for which
$\dim_H\Delta(E)\ge\frac{4}{5}+\eta_0$ whenever $t>1$, and therefore the
improvement remains positive even as the dimension of $E$ decreases to one.
They also obtained the explicit estimate
$$
 \dim_H\Delta(E)\ge\frac{t^2-2t+5}{2t^2-7t+10}
$$
in the range $1<t\le\frac{11-\sqrt{41}}{4}$, which reaches $\frac{7}{8}$ at
the right endpoint, along with the bound $\dim_H\Delta(E)\ge\frac{7}{8}$
for larger $t$~\cite[theorem~1.2]{IosevichYavicoli2026}. For pinned
distances they raised the uniform lower bound from $\frac{3}{4}$ to an
absolute constant greater than $0.7813$, again outside a set of pins of
Hausdorff dimension at most one~\cite[theorem~1.1]{IosevichYavicoli2026}.
Their argument builds on the complexity estimates of Stull and of Fiedler
and Stull, retaining information from both endpoints and combining the
resulting scale estimates with incidence bounds. These advances concern
general planar sets above dimension one, whereas the theorem of Shmerkin
and Wang quoted above reaches full dimension at the critical dimension
under its additional packing-dimension hypothesis.

Write $\dim_P$ for packing dimension. For Borel $E,F\subset{\mathbb R}^2$,
Liu proved that $\dim_HE>1$, $\dim_HE+\dim_HF>2$, and $\dim_HF=\dim_PF$
imply ${\mathcal L}^1(\Delta_y(E))>0$ for some $y\in F$
\cite[theorem~1.1]{Liu2026}. Taking $F=E$ gives the conjectured planar
conclusion for sets whose Hausdorff and packing dimensions are equal and
greater than one, the regularity of the pin set being a hypothesis of this
result rather than a consequence of its Hausdorff dimension. Borges, Ou and
Pasquariello study positive measure for pinned $k$-star distance sets and
for admissible distance graph
configurations~\cite[theorem~A and corollary~1.11]{BorgesOuPasquariello2026};
their threshold in the case $k=1$ lies above $\frac{d}{2}$, and they record
the sharper earlier one-edge bounds
in~\cite[remark~1.2]{BorgesOuPasquariello2026}. None of these results gives
theorem~\ref{thm:main}. The result proved here is unpinned, and no
corresponding claim about a single pin is asserted at the threshold
$\frac{d}{2}$.

The multiscale analysis of irregular measures has important precedents.
Keleti and Shmerkin decompose measures into dyadically regular pieces,
control bad projections across scale intervals, and optimize variations of
Lipschitz mass curves in their work on the dimensions of planar distance
sets~\cite[corollary~3.5, section~3.2 and section~5]{KeletiShmerkin2019}.
These are the antecedents of the regular classes, the scale-dependent
tests, and the mass curve accounting used below, though their estimates
concern Hausdorff dimension, and their optimization controls a sum of
interval drops, whereas section~\ref{ledger:section} uses a potential built
from two separated depths to pay for the costs of the present recursion.
Shmerkin and Wang later combined scale mass curves with radial projection
estimates and slicing~\cite[sections~5.5 and~6.1--6.3]{ShmerkinWang2025}.

There is a further connection with the recent paper of Iosevich and
Yavicoli, Improved dimension bounds for Falconer's distance problem in the
plane \cite{IosevichYavicoli2026}. Their reciprocal estimates keep
information from both endpoints, including the excesses at intermediate
scales, until those excesses can be used together, and the common entry
depth in lemma~\ref{asm:entry} serves a related purpose here; the
discussions preceding lemma~\ref{reg:pruning-multiplication} and
the equation~\eqref{asm:net-saving} explain the full source comparisons and the
common entry calculation. Their dimension estimates and the Fourier
summability required here remain different conclusions. The contact
construction in lemma~\ref{prep:contact} also follows their practice of
retaining a definite margin at selected scales, except that here two minima
give explicit margins, which remove the separate mass-curve compactness
argument from lemma~\ref{prep:one-scale}. No result from that paper is
added to the hypotheses of the proof. Our angular preparation is related to
the multiscale incidence approach of Orponen and
Shmerkin~\cite[section~1.5]{OrponenShmerkin2023}, to the radial thin-tube
bootstrap of Orponen, Shmerkin and
Wang~\cite[section~1.3]{OrponenShmerkinWang2024}, and to the higher-dimensional radial and plate-concentration methods of Ren~\cite{Ren2023}
and of Du, Ou, Ren and Zhang~\cite{DuOuRenZhang2024}. The incidence input
we retain is the theorem of Orponen and Shmerkin
\cite[theorem~1.3 and remark~1.4]{OrponenShmerkin2023}, in the tube
formulation of
\cite[definitions~2.5--2.6 and theorem~2.7]{OrponenShmerkinWang2024}, while
the auxiliary radial estimate is proved here by Orponen's finite-energy
method~\cite[section~3 and the proof of theorem~1.13]{Orponen2019}, and it
is applied only in a dimension where its energy hypotheses hold.

Auxiliary projections used to prepare radial information, while keeping the
final distance map in the original space, also occur in
\cite[sections~2--3]{DuIosevichOuWangZhang2021}, and the separation between
approximation of a distance measure and Fourier control of its retained
part has a precedent in
\cite[propositions~2.1--2.2]{GuthIosevichOuWang2020}. Liu's regular-pin
theorem uses a multiscale good and bad decomposition together with weighted
extension estimates~\cite[section~1.1]{Liu2026}. These comparisons
describe relationships between methods, and the relationship with the
source manuscript \cite{OpenAI2026} itself is the one described in
section~\ref{sec:changes}: the angular tests, the graph reductions and
the payment by two separated depths are taken from it, whereas the
accounting, the odd-dimensional radial improvement, the reassembly of the
spatial pieces, and the passage to a fixed measure with a regular density
are not. Cheng, Fang, Li, Qu, and Xiao also give a fixed-target
reconstruction in the plane and a length bound uniform in the Frostman
data \cite[lemma~2.1 and corollary~1.3]{ChengFangLiQuXiao2026}, and our
discrete application uses their quantitative length bound after a single
restriction to a set of bounded potential.

\subsection{The future of the Falconer distance problem}
\label{sec:future}

The sharpness construction is a natural place to begin a discussion of
what remains open, since the arithmetic which produces the exponent
$\frac{d}{2}$ changes as soon as the metric changes. On the integer points
of a box of side length $q$ the squared Euclidean distance takes at most
$dq^2+1$ values, and once the box has been rescaled and its points
thickened by $q^{-\frac{d}{s}}$, the covering intervals have total length
at most $C_dq^{2-\frac{d}{s}}$, where here and in the rest of this
subsection $C_d$ denotes a constant depending only on $d$. This count is
the source of the exponent $\frac{d}{2}$ in Falconer's construction
\cite{Falconer1985,IosevichRudnev2005}, and the way in which distance set
estimates depend on the geometry of the unit ball is studied in
\cite{IosevichLaba2005,IosevichRudnev2005}.

Consider the $l^d$ norm on ${\mathbb R}^d$, where $d\ge3$, and the
corresponding distance set of a compact set $E$,
$$
 \|x\|_{l^d}=\left(\sum_{j=1}^d|x_j|^d\right)^{\frac{1}{d}}
 \text{ and }
 \Delta_{l^d}(E)=\{\|x-y\|_{l^d}:x,y\in E\}.
$$
Let $G_q=\frac{1}{q}({\mathbb Z}^d\cap[0,q]^d)$ be the rescaled lattice
whose thickening is the set $E_{q,s}$ of the introduction, and let
$$
 \mathcal S_d(q)=\{n_1^d+\cdots+n_d^d:n_j\in\{0,\ldots,q\}\}
$$
be the set of sums of $d$ nonnegative $d$th powers of integers not
exceeding $q$. The $l^d$ distances determined by $G_q$ are exactly
$$
 \Delta_{l^d}(G_q)
 =\left\{\frac{m^{\frac{1}{d}}}{q}:m\in\mathcal S_d(q)\right\},
$$
and since $\mathcal S_d(q)\subset[0,dq^d]\cap{\mathbb Z}$ there are at
most $dq^d+1$ of them. Moving each endpoint by at most $q^{-\frac{d}{s}}$
changes an $l^d$ distance by at most $C_dq^{-\frac{d}{s}}$, by the
triangle inequality, so that
\begin{equation}\label{future:ld-cover}
 {\mathcal L}^1(\Delta_{l^d}(E_{q,s}))
 \le C_dq^{-\frac{d}{s}}\#\mathcal S_d(q)
 \le C_dq^{d-\frac{d}{s}},
\end{equation}
and the right side tends to zero when $s<1$. Falconer's limiting
construction therefore produces compact sets of every Hausdorff dimension
$0<s<1$ whose $l^d$ distance sets have measure zero, and the elementary
count now points to the exponent one rather than to $\frac{d}{2}$.

How many of the integers in $[0,dq^d]$ actually belong to
$\mathcal S_d(q)$? This question leads directly to an open problem in
number theory, since for $d\ge3$ it is not known whether a positive
proportion of the positive integers can be written as a sum of $d$
positive $d$th powers. The case of three cubes is discussed in Maynard's
recent survey~\cite[section~1]{Maynard2026}, and the general problem,
together with the related conjectures on the second moments of
representation counts, is discussed by Vaughan and
Wooley~\cite[section~17]{VaughanWooley2002}. A positive proportion among
the integers up to $q^d$ would give $\#\mathcal S_d(q)\ge c_dq^d$ for all
sufficiently large $q$, and allowing zero summands, as the definition of
$\mathcal S_d(q)$ does, adds at most $d(q+1)^{d-1}$ values and so leaves
the density question unchanged, while the absolute values in the norm
keep the summands nonnegative even when $d$ is odd.

Two different questions are in play here. The obstruction in
\eqref{future:ld-cover} for $s<1$ is unconditional, whereas the number
theoretic question concerns the actual size of the lattice distance set,
and positive density would be stronger than the estimate
$\#\mathcal S_d(q)=q^{d-o(1)}$, which would already rule out any fixed
power saving in this count. Neither statement would prove that
$\dim_HE>1$ implies ${\mathcal L}^1(\Delta_{l^d}(E))>0$ for an arbitrary
compact set, and the correct threshold for the $l^d$ metric is a separate
problem, which the lattice suggests may have an answer quite different
from the Euclidean one.

There are also basic questions in which the metric varies from point to
point. Let $(M,g)$ be a compact connected Riemannian manifold of dimension
$d\ge2$ without boundary, write $d_g(x,y)$ for its geodesic distance, and
for a compact set $E\subset M$ put
$$
 \Delta_g(E)=\{d_g(x,y):x,y\in E\}.
$$
The question is whether $\dim_HE>\frac{d}{2}$ guarantees that
$\Delta_g(E)$ has positive length. On suitably small separated coordinate
patches the distance is a smooth function of both endpoints, but the
coefficients of the metric depend on position, so that this is a variable
coefficient version of the distance problem, and generalized Radon
transforms and Fourier integral operators provide a framework for such
questions \cite{EswarathasanIosevichTaylor2011,IosevichTaylorUriarteTuero2021}.
Iosevich, Liu, and Xi proved the sufficient threshold $\frac{5}{4}$ on
Riemannian surfaces~\cite[theorem~1.1]{IosevichLiuXi2022}, and Jian, Liu,
and Xi~\cite{JianLiuXi2026} established positive length on general
Riemannian manifolds under the condition
$$
 \dim_HE>\frac{d}{2}+\frac{1}{4}+\frac{1-(-1)^d}{8d},
$$
whereas the threshold $\frac{d}{2}$ remains open for general metrics.

The sharpness construction itself deserves further attention in this
setting. A small Euclidean ball in a flat torus reproduces the Euclidean
lattice example, so that a theorem valid for every compact Riemannian
manifold cannot have a threshold below $\frac{d}{2}$, but for an
individual curved metric in dimension at least three the same sharpness
conclusion does not follow by transferring a lattice into local
coordinates, because the equal distances which make the Euclidean example
work need not remain equal. Local comparison of the two metrics controls
the sizes of the distances without preserving the arithmetic coincidences
on which the covering argument depends, and this issue is discussed in the
introduction of~\cite{IosevichLiuXi2022}.

For a compact set $E\subset{\mathbb R}^d$, another basic example is the
dot product set
$$
 \Pi(E)=\{x\cdot y:x,y\in E\}.
$$
For $x=\frac{a}{q}$ and $y=\frac{b}{q}$ in $G_q$ the product
$x\cdot y=\frac{a_1b_1+\cdots+a_db_d}{q^2}$ takes at most $dq^2+1$ values,
and since the dot product is Lipschitz on a fixed bounded box, the same
thickening argument gives
$$
 {\mathcal L}^1(\Pi(E_{q,s}))\le C_dq^{2-\frac{d}{s}},
$$
so that the lattice obstruction is again $\frac{d}{2}$. Whether
$\dim_HE>\frac{d}{2}$ forces $\Pi(E)$ to have positive length remains
open, as does the stronger pinned question, which asks for an $x\in E$
such that $\{x\cdot y:y\in E\}$ has positive length and which is
discussed, together with recent partial results, by Bright, Marshall, and
Senger~\cite[section~1]{BrightMarshallSenger2026}. The sufficient
condition $\dim_HE>\frac{d+1}{2}$ follows from the general results of
Iosevich, Taylor, and Uriarte-Tuero~\cite{IosevichTaylorUriarteTuero2021}.

The size of the image set and a uniform bound for the mass of every short
interval are different things, and the dot product illustrates the
distinction. Iosevich and Senger constructed measures with finite
$s$-energy, for any $0<s<\frac{d+1}{2}$, for which the estimate
$$
 (\mu\otimes\mu)\{(x,y):1\le x\cdot y\le1+\varepsilon\}\le C\varepsilon
$$
fails for every constant $C$ independent of $\varepsilon$
\cite{IosevichSenger2020}. This rules out that particular estimate below
$\frac{d+1}{2}$ without settling the positive measure question for
$\Pi(E)$, and the gap between the lattice obstruction and the known
positive measure thresholds remains to be closed.

\section{One descent budget for the whole argument}\label{bik:section}

Before returning to distances, let us isolate the calculation we are going
to need. It is useful whenever an estimate is improved through a succession
of choices, each of which may depend on everything that has happened before
it, provided that one and the same potential continues to measure the cost
of all of them.

\subsection{The elementary calculation}

Work in a real Hilbert space, with $\langle\cdot,\cdot\rangle$ denoting the
inner product and $\|\cdot\|$ the associated norm, and let $\varepsilon>0$
be the accuracy at which a test detects an error. Let $f$ be a target and
let $g$ be its current approximation, suppose that $\|u\|\le1$ and
$\langle f-g,u\rangle\ge\varepsilon$, and replace $g$ by
$g'=g+\varepsilon u$. Expanding the square gives
$$
\|f-g'\|^2
=\|f-g\|^2-2\varepsilon\langle f-g,u\rangle
+\varepsilon^2\|u\|^2
\le\|f-g\|^2-\varepsilon^2.
$$
Write $g_0$ for the initial approximation and $\varepsilon_j$ for the
accuracy at the step $j$. If the accuracies are allowed to change from one
step to another, exactly the same calculation gives
$$
\sum_j\varepsilon_j^2\le\|f-g_0\|^2,
$$
and no factor counting the number of accuracy levels appears,
because the approximation was never restarted when the accuracy changed.
The same squared-error calculation underlies the matrix decomposition
argument of Frieze and Kannan \cite[section~4.1]{FriezeKannan1999}.

In the entropic version used by Burstein--Iosevich--Krause the potential is
a relative entropy and the update changes the weights in a convex average,
the conclusion being a bound for the weighted number of corrections. Here
the quantities being corrected are estimates for distance counts and
Fourier energies, and we use a different potential, adapted to those
estimates, while retaining the same principle.

\subsection{The entropy calculation behind the method}

The same calculation can be carried out for probabilities, and this is the
quickest way to see what is borrowed from the Burstein--Iosevich--Krause
argument and what is not, since their sampling theorem plays no role in it.

All unlabelled logarithms are natural logarithms, $\log_b x$ denotes the
logarithm to the base $b>0$, $b\ne1$, and $m\ge1$ is the number of
coordinates in the two probability vectors. Let $p=(p_1,\ldots,p_m)$ be a
fixed target probability and let $q=(q_1,\ldots,q_m)$ be a probability with
all $q_l>0$. Their relative entropy is
$$
\mathcal E(p,q)=\sum_{l:p_l>0}p_l\log\frac{p_l}{q_l}.
$$
This quantity is nonnegative. Indeed, concavity of the logarithm gives
$$
\sum_{l:p_l>0}p_l\log\frac{q_l}{p_l}
\le\log\left(\sum_{l:p_l>0}q_l\right)\le0.
$$
Suppose that a real test $h=(h_l)$, with $|h_l|\le1$, detects an error
$$
\sum_l(p_l-q_l)h_l\ge\epsilon>0,
$$
and increase the weights where $h_l$ is large by setting
$$
q'_l=\frac{q_l e^{\epsilon h_l}}{\sum_jq_j e^{\epsilon h_j}}.
$$
To estimate the change in entropy, define
$$
F(t)=\log\left(\sum_lq_l e^{th_l}\right),
$$
and, for each real $t$, introduce the tilted probability with weights
$$
\omega_l(t)=\frac{q_l e^{th_l}}{\sum_j q_j e^{th_j}}.
$$
Differentiating the finite sum gives
$$
F'(t)=\sum_l\omega_l(t)h_l
\text{ and }
F''(t)=\sum_l\omega_l(t)h_l^2-
\left(\sum_l\omega_l(t)h_l\right)^2\le1,
$$
where the last inequality holds because the first term is at most one and
the subtracted term is nonnegative. Thus $F(0)=0$,
$F'(0)=\sum_lq_lh_l$, and $F''(t)$ is the variance of $h$ under the tilted
probability at time $t$, which is at most one since $|h_l|\le1$, so that Taylor's
formula with integral remainder gives Hoeffding's exponential moment
estimate \cite[equation~(4.16)]{Hoeffding1963},
$$
F(\epsilon)\le\epsilon\sum_lq_lh_l+\frac{\epsilon^2}{2},
$$
and substitution in the entropy difference yields
$$
\begin{aligned}
\mathcal E(p,q')-\mathcal E(p,q)
&=-\epsilon\sum_lp_lh_l+F(\epsilon)\\
&\le-\epsilon\sum_l(p_l-q_l)h_l+\frac{\epsilon^2}{2}\\
&\le-\frac{\epsilon^2}{2}.
\end{aligned}
$$
A detected error thus pays for its correction by a definite decrease of the
entropy, and since the target $p$ stays fixed while the tests and the
accuracies change from one step to the next, the same entropy serves
throughout. Hence
$$
\sum_j\epsilon_j^2\le2\mathcal E(p,q_{\mathrm{initial}}).
$$

This exponential change of weights is the classical multiplicative weights
update, whose relative entropy analysis is explained in
\cite[sections~2.1--2.2]{AroraHazanKale2012}. The feature we use below is
that the test and its accuracy may change at every step while the target
stays fixed.

In the comparison that follows, $r$ is the sparsity parameter of the
amplitude prediction problem, sparsity meaning here an upper bound for the
number of nonzero coefficients, and $\eta>0$ is the prescribed accuracy,
while the numbers $a_j>0$ are amplitude levels, which group the test values
according to their size, $d_j$ counts the corrections at the level $a_j$,
and $\Phi_{\mathrm{initial}}$ is the initial relative entropy potential.
The constant $C_\eta$ may depend on $\eta$. All these letters are local to
the comparison, and in particular this $r$ is not a Fourier frequency. In
the amplitude predictor of Burstein--Iosevich--Krause the margin at the
amplitude $a_j$ is proportional to $\frac{\eta a_j}{\sqrt r}$, and if $d_j$
corrections occur at that amplitude, their entropy calculation gives
$$
\sum_j a_j^2d_j\le C_\eta r\Phi_{\mathrm{initial}}.
$$
The corrections are recorded by their positions and labels, and counting
these records, while keeping their amplitude weights, produces the finite
families used in their probability estimate
\cite{BursteinIosevichKrause2026}.

For the distance problem the record consists of the intermediate depths
and the types of reduction. The interval potential $\mathcal V_h$ keeps
the accumulated cost when the scale changes or a distance count gives way
to an angular energy. The next lemma includes the change of normalization
and the error terms.

\begin{table}[htbp]
\centering\small
\renewcommand{\arraystretch}{1.3}
\begin{tabular}{p{3.1cm}p{4.8cm}p{5.2cm}}
\toprule
Role & Amplitude descent & Distance descent\\
\midrule
Quantity retained & Predictor distribution across amplitude levels & Current estimate and its tests across depth intervals\\
Potential & Relative entropy from a fixed target & $\Phi(J)=4d-\mathcal V_h(J)$\\
Cost of one step & Squared error margin & Excess cost in a geometric reduction\\
Finite record & Correction positions and labels & Reduction types and depth choices\\
Final summation & Weighted count and probability estimates & Uniform estimate over polynomially many finite records\\
\bottomrule
\end{tabular}
\vskip.125in 
\caption{The correspondence with the Burstein--Iosevich--Krause
method. The geometric payment inequality for the right-hand column
must still be proved; the entropy calculation does not supply it.}
\label{tab:descent-comparison}
\end{table}

\subsection{A finite tree of estimates}

It is convenient to record the successive substitutions in an inequality by
means of a finite rooted tree, in which a vertex is one estimate, the
children of a vertex are the estimates to which it reduces, and a leaf is
an estimate that can already be bounded directly.

Fix $R=2^N$, associate to every vertex $v$ a nonnegative number $X_v$ and
two real numbers $u_v,V_v$, and let $c_{vw}\ge0$ be the cost of the edge
from $v$ to $w$. Suppose that
$$
c_{vw}\le V_w-V_v.
$$
If $V_v\le B$ on the tree, then $\Phi_v=B-V_v$ is nonnegative and
$$
c_{vw}\le\Phi_v-\Phi_w,
$$
so that along a branch $v_0,v_1,\ldots,v_l$,
$$
\sum_{j=0}^{l-1}c_{v_jv_{j+1}}
\le\Phi_{v_0}-\Phi_{v_l}\le\Phi_{v_0}.
$$
This is the descent budget in its simplest form, and we now retain in
addition the changes of normalization and the small errors, both of which
are needed in the distance problem.

By the depth of the tree we mean the largest number of edges in a path from
the root to a leaf. The real number $u_v$ records the normalization at the
vertex $v$, while $V_v$ records the increasing quantity whose complement is
the descent potential, and $b_v\ge0$ is a contribution which can already be
bounded at that vertex, without further reduction. The marked edges are
those which are allowed to carry the additional fixed loss. We
write $o(1)$ for a quantity tending to zero as $N\to\infty$, with all the
stated parameters held fixed, and we call a factor bounded by $R^{o(1)}$
subpower.

\begin{lemma}[Weighted descent and finite-record chaining]\label{bik:tree}
Let $L,K,C,B_0,\varepsilon,\tau\ge0$ be fixed independently of $N$,
with $L,K$ integers, and consider a rooted tree of depth at most $L$ in
which every vertex has at most $N^C$ children and every branch contains
at most $K$ marked edges.

Let $|V_v-u_v|\le B_0$, and let $T>0$ be the exponent in the additive error bound. Suppose that each nonleaf satisfies
$$
X_v\le b_v+
\sum_{w\text{ child of }v}
R^{u_v-u_w+c_{vw}+e_{vw}+\delta_N}X_w+R^{-T},
$$
where $\delta_N\ge0$ tends to zero, $c_{vw}\le V_w-V_v$,
$e_{vw}=0$ on unmarked edges, and $0\le e_{vw}\le\varepsilon$
on marked edges. Assume
$$
b_v\le R^{u_v-V_v+\tau+\delta_N}
$$
and the same bound for $X_v$ at every leaf. If
$T>B_0+K\varepsilon+1$, then
$$
X_{v_0}\le
R^{u_{v_0}-V_{v_0}+\tau+K\varepsilon+o(1)}.
$$
The estimate is uniform over all trees and quantities with the stated
fixed bounds.
\end{lemma}

\begin{proof}
Set
$$
Y_v=R^{V_v-u_v}X_v.
$$
The coefficient of $Y_w$ in the normalized inequality is
$$
R^{c_{vw}+V_v-V_w+e_{vw}+\delta_N}
\le R^{e_{vw}+\delta_N},
$$
and therefore
$$
Y_v\le R^{\tau+\delta_N}
+\sum_{w\text{ child of }v}R^{e_{vw}+\delta_N}Y_w
+R^{B_0-T}.
$$
Substitute now the corresponding inequality at every child, stopping at the
leaves. A term is specified by the finite list of child choices which leads
to it, followed by the choice to stop or to use the error term, so that the
number of such lists is at most a fixed multiple of $(1+N^C)^{L+1}$, and
its logarithm divided by $\log R=N\log2$ tends to zero.

Along each list the sum of the $e_{vw}$ is at most $K\varepsilon$ and the
sum of the $\delta_N$ terms is at most $(L+1)\delta_N$, and hence each
stopping or leaf term is at most
$$
R^{\tau+K\varepsilon+(L+1)\delta_N},
$$
while each propagated error is at most
$$
R^{B_0-T+K\varepsilon+L\delta_N},
$$
the latter exponent being negative for all sufficiently large $N$. Summing
the finite lists changes the exponent by $o(1)$, which proves the asserted
bound for $Y_{v_0}$, and hence for $X_{v_0}$.
\end{proof}

For example, a fixed number $L$ of stages with at most $N^C$ choices at
each stage gives at most $N^{CL}$ complete records, whose base-$R$
logarithm is $\frac{CL\log N}{N\log2}$ and therefore tends to zero. The
order of the limits matters here, since one must first fix the small
parameters which determine $L$ and only then let $N$ increase; nothing is
asserted if the number of stages is allowed to grow with $N$. Note also
that $K$ may be much smaller than $L$, which will matter below, since some
of the reductions contribute a fixed power loss $R^{C_d\sigma}$, where
$C_d$ depends only on $d$ and $\sigma>0$ is a fixed small loss exponent,
whereas the low frequency reductions contribute only $R^{o(1)}$, and
charging the fixed power loss at every step would cost us the gain we are
trying to obtain.

The word chaining refers here to this deterministic succession of
approximating estimates, together with the count of their finite records,
and no probability estimate is concealed in the lemma. In the sampling
arguments of Burstein--Iosevich--Krause the finite records are instead
combined with probability estimates, the common feature being that the cost
attached to a given scale remains attached to it throughout the summation.

\subsection{What still has to be proved}

The lemma does not choose the reductions, nor does it produce their
geometric estimates, and what it requires is a valid inequality at each
step together with a potential which pays every cost. A small total cost
does not by itself prove that a branch terminates either, since a branch
may contain a great many steps of zero cost, and termination will be proved
separately by following the endpoints of the scale interval and the test
depths already passed.

This separation is convenient when adaptive choices are involved. Such a
choice may depend on the current interval, on the mass curve, and on all
the tests still present, and it remains legitimate provided that the
resulting child satisfies the same cost inequality, the same support
conditions, and the same termination bound. No adaptive rule can substitute
for these verifications.

\section{Conventions and preliminary facts}\label{sec:conventions}

The proof will repeatedly restrict a measure, divide by its mass, and look
at what is left inside a smaller cube. These operations are elementary, but
a missing normalization would change the exponent at the end, and so we fix
the conventions before anything else is done. In particular, restricting a
measure does not silently turn it into a probability.

\begin{table}[htbp]
\centering\small
\renewcommand{\arraystretch}{1.25}
\begin{tabular}{p{3.2cm}p{10.3cm}}
\toprule
Notation & Meaning\\
\midrule
$R=2^N$, $S=\frac{d}{2}$ & Reciprocal final spatial scale and critical exponent.\\
$a,g,p,t,v$ & Integer dyadic depths; increasing depth means smaller spatial cells.\\
$m$, $M(n)$ & Mass curve of one normalized class and model cell mass $R^{-m(\frac{n}{N})}$.\\
$A_n$ & Excess $m(\frac{n}{N})-\frac{Sn}{N}$ above critical mass decay.\\
$k(g,p)$ & Rise $A_p-\min_{g\le n\le p}A_n$, always nonnegative.\\
$q_0=\lceil\kappa N\rceil$ & Block length in the choice of intermediate depths.\\
$i$, $J_0$ & Current smoothing level and largest test index. The test levels are $0,\ldots,J_0$.\\
$\mathcal B$ & Retained base interval of an angular energy; it determines the test list.\\
$\mathcal V_h$, $\Phi$ & Increasing interval quantity and decreasing potential $\Phi=4d-\mathcal V_h$.\\
$\alpha_\infty$, $\alpha_N$, $\tau_N$ & Fixed distance measure, its coarse approximant, and the retained measure tested at shell $N$.\\
\bottomrule
\end{tabular}
\vskip.125in 
\caption{Notation used across the proof. Auxiliary letters inside
the angular-gain lemmas are local to those lemmas and are defined there.}
\label{tab:notation}
\end{table}

Only the angular masks used to construct $\tau_N$ are smooth, and absolute
continuity of the limiting distance measure is a conclusion of the proof
rather than an assumption about these approximants.

Some elementary notation will be used in the local lemmas. The support
$\operatorname{supp}\lambda$ of a measure is the smallest closed set outside which it
vanishes, $\operatorname{diam} A$ is the diameter of a set $A$, and
$\operatorname{dist}(x,A)=\inf_{y\in A}|x-y|$, and $\operatorname{Lip}(f)$ denotes the least
Lipschitz constant of $f$. We write $f\asymp g$ when
$C^{-1}g\le f\le Cg$ for a positive fixed constant, $f\lesssim g$ when
$f\le Cg$ for such a constant, and $f\gtrsim g$ when $g\lesssim f$, a
subscript on a comparison constant listing the permitted dependencies. As
usual, $O(g)$ stands for a quantity of absolute value at most $Cg$, and
$o(g)$ for one whose ratio to $g$ tends to zero in the specified limit.
Finally, an equality $x=R^{b\pm e}$ is shorthand for the two inequalities
$R^{b-e}\le x\le R^{b+e}$, and not for a choice of sign.

We write $C^j$ for the class of functions with continuous derivatives
through the order $j$ and $C^\infty$ for those with derivatives of every
order, a subscript $c$ requiring in addition compact support, and we recall
that a Schwartz function and all its derivatives decay faster than every
inverse power at infinity. The convolution of a function $f$ with a finite
measure $\lambda$ is $(f*\lambda)(x)=\int f(x-y)d\lambda(y)$, and by an
approximate identity we mean a family of integrable kernels of integral
one, with uniformly bounded $L^1$ norms, whose absolute integrals outside
every fixed neighborhood of the origin tend to zero. A majorant of a
function is a pointwise upper bound for it, and, for $r>0$, an $r$-net is a
collection of representatives such that every point of the set under
consideration lies within the distance $r$ of one of them, this distance
being measured in the parameter space when lines or tubes are involved.

\subsection{Measures, scales and directions}

All our measures are finite Borel measures, a subprobability being one of
mass at most one, and restrictions are left unnormalized unless we say
otherwise. We write
$$
 I_s(\mu)=\iint |x-y|^{-s}d\mu(x)d\mu(y)
 \text{ and }
  \widehat\mu(\xi)=\int e^{-i x\cdot\xi}d\mu(x),
$$
and we use the same Fourier convention for measures on the line, while
$\|\cdot\|_{\mathrm{TV}}$ is the total variation norm and $f\lesssim g$
means $f\le Cg$ with the dependence of $C$ specified by the surrounding
discussion.

We say that a measure has an $s$-ball bound with loss $K$, down to the
scale $\delta$, if
$$
 \mu(B(z,r))\le K r^s
  \text{ for every } z \text{ and every } \delta\le r\le1,
$$
and a cap bound on a sphere has the analogous meaning for the Euclidean
chordal metric there. Passing between that metric, angular distance and
fixed smooth charts changes only constants, so that we may work with dyadic
radii and then enlarge the constant to recover all the radii in the
specified interval.

For $x\ne y$ we put
$$
 \pi_x(y)=\frac{y-x}{|y-x|} \text{ and } D(x,y)=|x-y|.
$$
A tube is a neighborhood of an affine line and a plate is a neighborhood of
an affine plane of the dimension indicated, while directions are oriented
unless an antipodal pair of caps is explicitly used. All the configurations
we consider lie in bounded sets, so that fixed rescalings and a finite
number of coordinate charts cost us nothing.

The original ambient dimension is always denoted by $d$, with
$S=\frac{d}{2}$, and $R=2^N$ denotes the reciprocal of the final scale. In
odd dimension we write $d=2k-1$, in which case the linear angular gain is
used in ${\mathbb R}^{k+1}$, while $S=k-\frac{1}{2}$ continues to refer to
the original dimension. Letters denoting dimensions inside a preliminary
lemma are local to that lemma.

For a sequence of numbers $z=(z_j)$ we write
$\|z\|_{\ell^2}=(\sum_j|z_j|^2)^{\frac{1}{2}}$ for its square-sum norm,
which, unlike $L^2$, is taken with respect to counting measure on an index
set; it will be used for the masses of distance bins.

\subsection{Finite energy and projection}

We record here the elementary consequences of the usual Frostman and
projection principles which will be needed later, including the restriction
operation. The three assertions of the lemma below allow us to move between
geometric size and integrable kernels: a ball bound gives an energy
estimate with a slightly smaller exponent, finite energy gives back a ball
bound after a small restriction, and the last assertion carries this
information through an auxiliary projection. Such a projection will be used
to study directions, and it does not replace the distance map in the final
theorem.

For an $m$-dimensional linear subspace $V\subset{\mathbb R}^d$ we write
$P_V$ for the orthogonal projection onto $V$ and $V^\perp$ for the
orthogonal complement, the set of all such $V$ being the Grassmannian
$G(d,m)$, whose Haar probability measure, written $dV$, is the probability
measure invariant under orthogonal changes of coordinates. Thus a
Haar-random plane is one chosen without favoring any direction.

\begin{lemma}\label{conv:energy}
Let $\mu$ be supported in a bounded subset of ${\mathbb R}^d$.
\begin{enumerate}[label=(\roman*)]
\item If $\mu(B(z,r))\le C r^s$ for $0<r\le1$, then
$I_t(\mu)<\infty$ for every $0<t<s$.
\item If $I_t(\mu)<\infty$, then for every $\varepsilon>0$ there is a
compact restriction of mass at least $\mu({\mathbb R}^d)-\varepsilon$ that satisfies
a $t$-ball bound with a finite constant.
\item If $0<t<m\le d$, $I_t(\mu)<\infty$, and $P_V$ is
orthogonal projection to a Haar-random $m$-plane, then
$I_t((P_V)_\#\mu)<\infty$ for almost every $V$.
\end{enumerate}
\end{lemma}

\begin{proof}
For (i), fix $x$ and divide the region $0<|x-y|\le1$ into the annuli
$2^{-j-1}<|x-y|\le2^{-j}$, with $j\ge0$. On the $j$th annulus the kernel is
at most $2^{(j+1)t}$ while the mass is at most $C2^{-js}$, so that
$$
\int_{0<|x-y|\le1}|x-y|^{-t}d\mu(y)
\le C2^t\sum_{j\ge0}2^{-j(s-t)}<\infty,
$$
the series converging because $s-t>0$. Since the ball bound with $s>0$ also
implies that $\mu$ has no atoms, the diagonal carries no product mass, and
on $|x-y|>1$ the kernel is at most one; thus integrating the resulting
uniform bound in $x$ gives $I_t(\mu)<\infty$.

For (ii), let $U_t^\mu(x)=\int |x-y|^{-t}d\mu(y)$, which is finite for
$\mu$-almost every $x$, so that a sufficiently large level $M$ discards
arbitrarily little mass. Choose a compact subset
$K\subset\{U_t^\mu\le M\}$ losing at most the remaining prescribed mass. If
$K\cap B(z,r)\ne\varnothing$, choose $x$ in this intersection, and then
$$
 \mu|_K(B(z,r))
 \le\mu(B(x,2r))\le (2r)^t U_t^\mu(x)\le M(2r)^t,
$$
while the same bound is immediate when the intersection is empty.

For (iii), rotational invariance together with the integrability of the
negative $t$-moment of the length of an $m$-dimensional projection gives
$$
 \int_{G(d,m)}|P_V z|^{-t}d V\le C_{d,m,t}|z|^{-t}
  \text{ for } z\ne0.
$$
For $m<d$ this follows by integrating the density of the first $m$
coordinates of a uniform unit vector near the origin, the radial integral
there being bounded by a constant times
$\int_0^1 r^{m-1-t}d r<\infty$, and for $m=d$ the identity projection makes
the assertion immediate. Integrating in $(x,y)$ and applying Tonelli proves
the claim.
\end{proof}

We use Frostman's lemma in its standard form
\cite[theorem~8.8]{Mattila1995}: if a compact set $E$ has $\dim_H E>s$,
then it supports a nonzero finite measure satisfying an $s$-ball bound at
all scales, and dividing by its mass gives a probability. The positive
dimension margin permits a slight decrease of $s$ whenever a finite energy
is needed in place of a ball bound.

Dyadic cubes are half open, and therefore they form genuine partitions, and one
may choose a translate of the dyadic grid whose boundaries have zero mass
for the finitely many starting measures. For each coordinate hyperplane
this follows by integrating its mass over the translation, and the set of
grid hyperplanes is countable.

\subsection{Subpower losses and fixed gains}

For a sequence $R\to\infty$, we call a positive factor at most subpower if
it is bounded by $R^{o(1)}$, and a retained mass at least subpower if it is
bounded below by $R^{-o(1)}$. These are statements about the sequence, and
using them does not amount to asserting that some fixed positive exponent
vanishes; constants and powers of $\log R$ are subpower in this sense. For
$\delta\to0$ the corresponding notation uses $R=\delta^{-1}$.

We shall use the following quantifier principle, which makes precise those
compactness passages in which the input losses tend to zero while the
argument produces a fixed positive improvement. Two different limits occur
in it, since the scale tends to zero while the permitted losses also tend
to zero along a hypothetical sequence of counterexamples, and the
conclusion supplies one fixed loss and one fixed gain which work at every
sufficiently fine scale. It asserts nothing about what happens to the gain
when a dimension gap closes.

\begin{lemma}[Sequential uniformity principle]\label{conv:uniformity}
For fixed background parameters, let $\mathcal C(\varepsilon,R)$ be classes
of data, increasing in the allowed loss $\varepsilon>0$.
Let $\mathcal G(c,R)$ be conclusions, weakening as the gain
$c>0$ decreases. Suppose that for every sequence
$$
 \varepsilon_j\downarrow0, R_j\to\infty \text{ and } X_j\in\mathcal C(\varepsilon_j,R_j),
$$
there are a subsequence and a fixed $c>0$ on which
$X_j$ satisfies $\mathcal G(c,R_j)$ for all sufficiently large $j$.
Then there exist $\varepsilon_0,c_0>0$ and $R_0<\infty$ such that
every $X\in\mathcal C(\varepsilon_0,R)$, $R\ge R_0$, satisfies
$\mathcal G(c_0,R)$.
\end{lemma}

\begin{proof}
If the conclusion fails, then for each $j$ we may choose $R_j\ge j$ and
$X_j\in\mathcal C(\frac{1}{j},R_j)$ failing $\mathcal G(\frac{1}{j},R_j)$.
The assumed subsequence satisfies $\mathcal G(c,R_j)$ for some fixed $c>0$,
and for large $j$ we have $\frac{1}{j}<c$, so that monotonicity gives
$\mathcal G(\frac{1}{j},R_j)$, a contradiction.
\end{proof}

A bound is called power-small if it is $O(R^{-c})$, or equivalently
$O(\delta^c)$ when $R=\delta^{-1}$, for a fixed $c>0$. Fixed constants in a
power-small output may be absorbed by decreasing the output exponent and
increasing the scale threshold, and when several losses occur in an input,
their maximum serves as $\varepsilon$ in lemma~\ref{conv:uniformity}.
Compactness of uniformly Lipschitz normalized mass curves is the usual
Arzel\`a--Ascoli compactness \cite[theorem~4.43]{Folland1999}. All the
background parameters, including any strict dimension or transversality
exponent, are fixed before this principle is applied, and the arguments
below specify the ordering of the losses that each of them requires.

For measurable selection at a fixed scale it is enough to work with dyadic
grids, finite tube nets and a finite number of charts, and strict
inequalities can first be imposed with a fixed enlargement of widths. Unions
of exceptional configurations can therefore be chosen Borel, and no
selection from an uncountable unparameterized family is ever needed.

\section{Regularization and multiplication of fibre bounds}\label{reg:section}

An occupied cube always means a cube of positive mass for the measure
currently under discussion, even when a larger, earlier measure occupied
many more cubes. A Frostman bound controls how much mass a ball can
contain, but it says nothing about whether occupied cells at the same depth
have comparable masses, and so we begin by dividing the measure into pieces
on which that comparison is available at every relevant depth
\cite[lemma~3.1]{OpenAI2026}. The
discarded mass has to be kept small, since we shall eventually return to
the original measure, and the second lemma of this section explains how
estimates inside the conditional pieces multiply when the pieces are put
back together.

We first record two facts about finite trees, neither of which requires
lower density estimates for the original measure. Dyadic cubes are taken
half open, so that they form genuine measurable partitions, and a fixed
translation and rescaling places each measure under consideration in one
unit cube.
Related decompositions into dyadically regular measures occur in
\cite[corollary~3.5]{KeletiShmerkin2019}; we prove the version needed here.
Their decomposition also retains almost all of the mass and gives each
normalized piece controlled branching at successive scales, and what we add
is that the mass of every piece is kept as well, since those masses will be
the coefficients when the estimates are assembled in
section~\ref{asm:section}.

The root of the dyadic tree is its depth-zero cube, its descendants are the
smaller cubes contained in it, and its leaves are the cubes at the final
depth, while a dyadic bin is an interval $[2^j,2^{j+1})$ with $j$ an
integer, and therefore two numbers in one bin differ by a factor of less than
two. For a uniform measure on a dyadic tree with $b_j$ children at the
generation $j$, every occupied cell at depth $n$ has mass
$$
\frac{1}{b_1\cdots b_n},
$$
and taking logarithms turns this product of branching numbers into a sum;
the function $m$ in the lemma below is the continuous interpolation of that
sum, divided by $N$. An arbitrary measure of course need not have uniform
branching, but we recover the same description, up to a small power of $R$,
by putting leaves with comparable masses and comparable branching histories
into the same class.

\begin{lemma}[Regular classes]\label{reg:regularization}
Let $n\ge1$, let $\rho$ be a probability on a unit dyadic cube in
${\mathbb R}^n$, and put $R=2^N$.  For every fixed $0<\lambda<1$ and all
sufficiently large $N$, there are pairwise disjoint sets $E_\alpha$,
each a union of depth-$N$ cubes, with the following properties.
\begin{enumerate}
\item Their number is at most $N^{\frac{C_n}{\lambda}}$, each has
      $\rho(E_\alpha)\ge R^{-\lambda}$, and
      $$
      \rho\left({\mathbb R}^n\setminus\bigcup_\alpha E_\alpha\right)
      \le R^{-2}+N^{\frac{C_n}{\lambda}}R^{-\lambda}
      =R^{-\lambda+o(1)}.
      $$
\item For $\rho_\alpha=\frac{\rho|_{E_\alpha}}{\rho(E_\alpha)}$ there is a
      nondecreasing $n$-Lipschitz function $m_\alpha:[0,1]\to[0,n]$,
      with $m_\alpha(0)=0$, such that every occupied depth-$j$ cube
      satisfies
\begin{equation}\label{reg:cell-mass}
R^{-m_\alpha(\frac{j}{N})-C_n\lambda}
      \le\rho_\alpha(Q)\le
      R^{-m_\alpha(\frac{j}{N})+C_n\lambda} \text{ for } 0\le j\le N.
\end{equation}
\item If $\rho(B(x,r))\le Kr^s$ for $R^{-1}\le r\le1$, with
      $0\le s\le n$ and $K\ge1$, then
\begin{equation}\label{reg:profile-lower-bound}
m_\alpha(u)\ge su-C_n\lambda-\log_R K
       \text{ for } 0\le u\le1.
\end{equation}
\end{enumerate}
Here and below an occupied cube means a cube of positive mass for the
specified class.  In particular, for occupied $P\subset Q$ at depths
$b\ge a$,
\begin{equation}\label{reg:conditional-mass}
\frac{\rho_\alpha(P)}{\rho_\alpha(Q)}
 =R^{-[m_\alpha(\frac{b}{N})-m_\alpha(\frac{a}{N})]\pm C_n\lambda}.
\end{equation}
Constants depend only on $n$; the threshold for $N$ can depend on
$\lambda$.
\end{lemma}

\begin{proof}
The construction has four stages. We choose a finite set of depths, label
the leaves by mass and by branching, read the mass curve off the labels,
and finally discard the classes whose original weights are too small.
Whenever a class is normalized, its original mass is kept separately for
the later recombination.

Discard first the depth-$N$ cubes of mass less than $R^{-n-2}$. There are
at most $R^n$ of them, so that this costs at most $R^{-2}$. Choose depths
$$
0=j_0<j_1<\cdots<j_J=N
 \text{ with }
 j_i-j_{i-1}\le\lambda N+1
 \text{ and }
 J\le \frac{C}{\lambda},
$$
and give each remaining leaf its dyadic mass bin as an initial label, of
which there are $O_n(N)$ possibilities. Now work upwards through the
selected depths. If $\ell$ is a label at depth $j_i$, count, within each
depth-$j_{i-1}$ parent, the depth-$j_i$ children carrying $\ell$, and
extend $\ell$ by the dyadic bin of this positive count, the extension being
given to all the leaves of label $\ell$ in that parent. A parent may carry
several labels, but the sets of leaves having different labels are
disjoint, and at every extension there are at most $O_n(N)$ possible count
bins, and hence the number of final labels is at most
$N^{\frac{C_n}{\lambda}}$ after the constants have been adjusted for
sufficiently large $N$.

Fix a final label. For each $i$ there is a positive integer $b_i$, a power
of two, such that every occupied depth-$j_{i-1}$ parent has between $b_i$
and $2b_i$ occupied depth-$j_i$ children in this label. To see that these
comparisons survive the later extensions, fix a partial label and a node at
the current selected depth. All the leaves with that label in this node
receive the same extension, which is determined by its parent and the old
label, and at each subsequent higher step the extension is again constant
on the entire group. Backward induction through the selected depths
therefore shows that a final label retains or removes each earlier child
group as a whole and never splits it, so that the branching comparisons
hold simultaneously. We also have
$$
1\le b_i\le 2^{n(j_i-j_{i-1})}.
$$
All the leaf masses lie in one dyadic bin, so that counting the leaves
below a selected node and then normalizing the final class gives
$$
\rho_\alpha(Q)=2^{\pm C J}
                  \left(\prod_{h\le i}b_h\right)^{-1}
 \text{ at depth } j_i.
$$
Define
$$
m_\alpha\left(\frac{j_i}{N}\right)=\frac{1}{N}\sum_{h\le i}\log_2b_h
$$
and interpolate linearly. On the interval between two consecutive selected
depths its slope is
$$
\frac{m_\alpha\left(\frac{j_i}{N}\right)-m_\alpha\left(\frac{j_{i-1}}{N}\right)}{\frac{j_i-j_{i-1}}{N}}
=\frac{\log_2 b_i}{j_i-j_{i-1}},
$$
and the bound on $b_i$ places this number in $[0,n]$; thus the Lipschitz
bound is exact even though the cell masses are only approximately equal. At
an intermediate depth an occupied cube is contained in its preceding
selected-depth parent and contains an occupied descendant at the next
selected depth, and both the possible mass variation and the variation of
the mass curve across this gap are at most
$R^{n\lambda+O_n(\frac{1}{N})}$. The remaining comparison error is
$2^{O(J)}$, whose base-$R$ logarithm is $O(\frac{1}{\lambda N})$, so that
for fixed $\lambda$ and sufficiently large $N$ these errors give
\eqref{reg:cell-mass}.

It remains to discard the final classes of original mass less than
$R^{-\lambda}$, whose total mass is bounded by the number of labels times
$R^{-\lambda}$, as asserted. For a retained class, normalization multiplies
an original upper ball bound by at most $R^\lambda$, and covering a
depth-$j$ cube by a bounded number of balls of radius $2^{-j}$ gives
$$
\rho_\alpha(Q)\le C_nKR^\lambda2^{-sj}.
$$
Comparison with the lower bound in \eqref{reg:cell-mass} proves
\eqref{reg:profile-lower-bound} at integer depths, with a larger
dimensional constant, and linear interpolation gives it at every $u$, while
dividing two instances of \eqref{reg:cell-mass} proves
\eqref{reg:conditional-mass}.
\end{proof}

For use in the scale-compactness arguments, the proof also gives a precise
interpretation of subpower regularity. If $\lambda=\lambda_N\to0$
sufficiently slowly that $\frac{\lambda_N^2N}{\log N}\to\infty$, then the
number of classes, all the mass comparison losses, and the normalization
losses of the retained classes are $R^{o(1)}$, while the discarded mass
tends to zero and the mass curves remain exactly $n$-Lipschitz. This by
itself asserts nothing about a fixed power of discarded mass; fixed
positive deletion exponents are obtained by choosing a fixed small
regularization parameter after the relevant compactness argument has been
carried out.

By a witness we mean a retained pair at which an estimate is required.
Once a pin and a target cap are fixed, a child cube contributes only if it
contains a witness whose direction lies in the cap, and since the local
test in the parent counts the full masses of all such children, its bounds
can be multiplied over successive depths.

\begin{figure}[htbp]
\centering
\begin{tikzpicture}[>=Stealth,font=\small,
 box/.style={draw,align=center,text width=3.7cm,minimum height=1.1cm}]
\node[box] (original) at (0,0) {Original probability $\rho$\\total mass $1$};
\node[box] (class) at (-3,-2) {Restriction $\rho|_{E_\alpha}$\\mass $a_\alpha=\rho(E_\alpha)$};
\node[box] (other) at (3,-2) {Other classes and\\discarded leaves};
\node[box] (normal) at (-3,-4) {Class probability\\$\rho_\alpha=\frac{\rho|_{E_\alpha}}{a_\alpha}$};
\node[box] (cell) at (3,-4) {Occupied depth-$j$ cell $Q$\\$\rho_\alpha(Q)=R^{-m_\alpha(\frac{j}{N})\pm C_n\lambda}$};
\draw[->] (original)--(class);
\draw[->] (original)--(other);
\draw[->] (class)--node[left] {divide by $a_\alpha$}(normal);
\draw[->] (normal)--node[above,align=center] {cell\\masses}(cell);
\end{tikzpicture}
\caption{Regularization separates the original weight $a_\alpha$ from
the geometry of a normalized class. The mass curve describes
$\rho_\alpha$ rather than $\rho|_{E_\alpha}$, and in every final sum the
class probability is multiplied by $a_\alpha$ again, so that discarding
a class of small weight is charged in the original measure.}
\label{fig:class-normalization}
\end{figure}
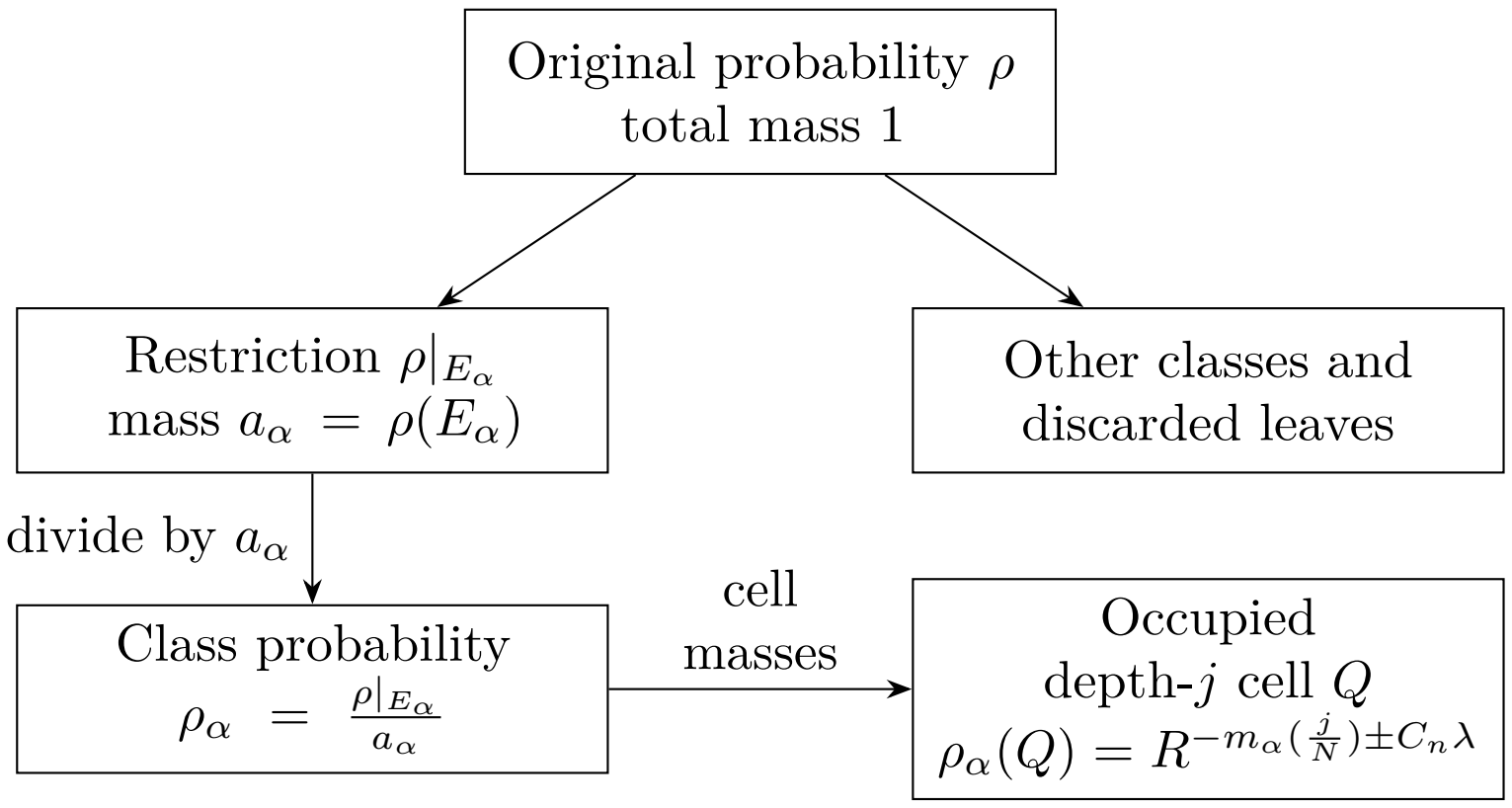

For any set $G$ of pairs $(x,z)$, its source fiber at $x$ is
$G_x=\{z:(x,z)\in G\}$ and its pin fiber at $z$ is
$G^z=\{x:(x,z)\in G\}$, and we use the same notation for pair sets denoted
by other letters. Integrating a source-fiber mass against the pin measure
recovers the product measure of the pair set, by Tonelli's theorem.

A local radial test refers to two source points, one of which, $z$, is
a witness in a pair that survives the deletions, while the other, $z'$,
ranges over a whole cell $Q$ of the source measure. If $\rho(Q)>0$, the
probability in that test is
$$
 \rho_Q=\frac{\rho|_Q}{\rho(Q)},
$$
whose denominator is the original source mass of the cell, before any
condition involving the pin, and this allows us to multiply successive
local bounds by summing the full masses of the cells they cover. The same distinction between a retained witness and a comparison
measure occurs in the conditional source estimates of Iosevich and
Yavicoli~\cite{IosevichYavicoli2026}.

\begin{lemma}[Multiplication of full source bounds]\label{reg:pruning-multiplication}
Let $\rho$ be a probability on a bounded subset of $\mathbb R^n$, let
$G$ be a measurable set of pairs $(x,z)$ with $z\in\operatorname{supp}\rho$, and let
$0\le a_0<\cdots<a_J\le N$ be dyadic depths. Fix $C_1>0$. Suppose
that the pins and the full source support lie in bounded regions
separated by a fixed positive distance. For each $j<J$, suppose that
every $(x,z)\in G$, with $z$ in a depth-$a_j$ cube $Q$, satisfies
\begin{equation}\label{reg:local-fibre}
 \rho\{z'\in Q:|\pi_x(z')-\pi_x(z)|\le C2^{-a_{j+1}}\}
 \le f_j\rho(Q),
\end{equation}
where $f_j\ge0$ and $C$ is sufficiently large in terms of $C_1$ and
the fixed geometry. Then for each pin $x$, depth-$a_0$ cube $Q_0$
and spherical ball $B$ of radius $C_12^{-a_J}$,
\begin{equation}\label{reg:multiplied-fibre}
 \rho\{z\in Q_0\cap G_x:\pi_x(z)\in B\}
 \le \rho(Q_0)\prod_{j=0}^{J-1} f_j .
\end{equation}
For bounded pin and source regions the same conclusion holds even when
they touch, provided $G$ omits the diagonal and the local test on $z'$
in \eqref{reg:local-fibre} is replaced by
$\operatorname{dist}(z',x+\mathbb R(z-x))\le C2^{-a_{j+1}}$.
\end{lemma}

\begin{proof}
Fix $x,Q_0,B$, and call a cube inside $Q_0$ hit if it contains a witness
$z\in G_x$ with $\pi_x(z)\in B$. In each hit parent $Q$ at depth $a_j$
choose such a witness $z_Q$, and if $P\subset Q$ is hit at depth
$a_{j+1}$, choose a witness $z_P$ in $P$, so that both directions lie in
the final cap. Since $\pi_x$ is uniformly Lipschitz on the source support
by the separation assumption, every $z'\in P\cap\operatorname{supp}\rho$
satisfies
$$
 |\pi_x(z')-\pi_x(z_Q)|
 \le C_2 2^{-a_{j+1}}+2C_1 2^{-a_J}
 \le C2^{-a_{j+1}}.
$$
The hit children are disjoint and the local test at $z_Q$ counts all of
their source mass, and therefore
$$
 \sum_{\substack{P\subset Q\\P\text{ hit at }a_{j+1}}}\rho(P)
 \le f_j\rho(Q).
$$
Iterating this inequality proves \eqref{reg:multiplied-fibre}, since the
final hit cubes cover the fiber in that equation. For the physical tube
test, the bounded pin--source distance and the final cap condition put
$z_P$ within $O(2^{-a_J})$ of the line through $x$ and $z_Q$, every point
in $P$ lies within a further $\sqrt{n}2^{-a_{j+1}}$ of that line, and the
same hit-child estimate and iteration apply.
\end{proof}

All the selections above can be made measurably, since the cube partitions
and label sets are finite and the fibre integrals are measurable, and
enlarging the fixed test widths allows the use of closed caps or closed
tubes without altering any exponent.

\section{An angular projection gain}\label{ang:section}

The following question arises in odd dimensions. If a measure on
directions already avoids excessive concentration in small caps, can one
improve that estimate using the fact that those directions come from pairs
of points? The extra half dimension needed later is obtained by reducing an
incidence configuration to the plane, and the finite-scale estimate proved
in this section bounds the incidences whose fiber at each line direction
is covered by a small number of parallel tubes, the only deep input being
the planar discretized Furstenberg theorem. We state the main estimate
first and then prove it, following
\cite[theorem~4.1 and lemmas~4.3--4.5]{OpenAI2026} and keeping track
throughout of the mass lost in slicing, since an improvement which retained only a negligible amount of
mass would be of no use in the distance problem. A reader interested only
in the even dimensional application can pass directly to the integer cap
construction in the next section.

All sets in this section lie in fixed bounded regions, and directions are
unit vectors, with antipodal vectors identified when lines are under
discussion. A projective chart chooses one representative from each nearby
antipodal pair; in ${\mathbb R}^n$, after a coordinate axis has been
chosen, it has the form $w(p)=\frac{(1,p)}{\sqrt{1+|p|^2}}$ with
$p\in{\mathbb R}^{n-1}$ in a bounded set, so that directions can be
described by ordinary Euclidean coordinates without identifying distinct
nearby points, and a finite number of such charts removes the distinction
altogether. For a bounded set $A$ we write $|A|_r$ for its covering number
by balls of radius $r$, and a tube has bounded length and the stated
perpendicular width, its length being enlarged by a fixed constant whenever
that is convenient. We use the ball-bound and subpower conventions of
section~\ref{sec:conventions}, bounds for radii in a fixed larger bounded
interval following by a change of constants; in particular, incidence mass
at least subpower means a lower bound $\delta^{o(1)}$.

\subsection{The gain and its two slicing regimes}

\begin{theorem}[Linear projection gain]\label{ang:linear-gain}
Fix $k\ge2$ and
$$
 k-1\le q\le k-\frac{1}{2}<t<k.
$$
Fix a tube-width constant $C$. In the case $q=k-1$, also fix
$\zeta>0$. There exist $\epsilon,c>0$ such that, for all sufficiently
small $\delta$, the following holds. Let $\mu$ be a bounded point
subprobability in ${\mathbb R}^{k+1}$ and $\nu$ a direction subprobability on $S^k$
with respective $t$- and $q$-ball bounds of loss $\delta^{-\epsilon}$
down to $\delta$. At $q=k-1$ assume in addition
\begin{equation}\label{ang:endpoint-plate}
\nu\{w:\operatorname{dist}(w,H)\le r\}\le\delta^{-\epsilon}r^\zeta
        \text{ for } \delta\le r\le1
\end{equation}
for every linear hyperplane $H$. Let $E\subset{\mathbb R}^{k+1}\times S^k$ be a
measurable incidence relation, and write $E_w=\{x:(x,w)\in E\}$.
If each direction fiber $E_w$ is covered by at most
$\delta^{-q-\epsilon}$ tubes parallel to $w$ of width $C\delta$, then
$$
 (\mu\otimes\nu)(E)\le\delta^c.
$$
All constants are uniform over the measures and the incidence relations
with these fixed parameters.
\end{theorem}

A cover by $\delta^{-q-\epsilon}$ tubes means that the perpendicular
projection of a direction fiber is covered by that many balls of radius
comparable to $\delta$, and the theorem says that such covers cannot
capture more than a power-small part of the incidence between the two
nonconcentrated laws.

The two cases in its proof are distinguished by what is left after the
slicing. If $q>k-1$, then slicing by codimension $h=k-1$ leaves a plane,
with exponents
$$
 b=q-h\in\left(0,\frac{1}{2}\right] \text{ and } f=t-h\in\left(\frac{1}{2},1\right),
$$
so that $b<f$, and the planar incidence estimate applies in this range. At $q=k-1$ the
same slicing would leave direction exponent zero, so we use instead
$h=k-2$, which leaves three-dimensional slices with direction exponent one
and point exponent greater than $\frac{3}{2}$. The endpoint argument then
finds three transverse directions in one such slice, while a weak radial
estimate with exponent $\ell>k-2$ controls the triples which are almost
dependent.

In the finite sampling arguments $\Pr$ and $\mathbb P$ both denote
probability, $\mathbb E Y$ denotes the expectation of a random variable
$Y$, and a subscript, as in $\mathbb E_V$, specifies the random object
being averaged, and complements are denoted by a superscript $c$.

\subsection{The planar incidence input and its measure formulation}

As the metric on bounded affine line parameters we may take the sum of the
distance between the associated orthogonal projection operators and the
distance between the closest-to-origin points on the lines, the operator
distance being measured by the Euclidean operator norm
$\|T\|=\sup_{|z|=1}|Tz|$ of a linear map $T$; other smooth bounded-chart
metrics are equivalent to this one. A nonempty bounded set $P$ of points or
line parameters is a $(\delta,s,K)$-set if
$$
 |P\cap B(a,r)|_\delta\le K r^s |P|_\delta
        \text{ for } \delta\le r\le1,
$$
the normalization by $|P|_\delta$ being part of the defining bound.

The planar theorem below is an external input rather than a consequence of
the descent argument, and its conclusion compares two counts: a
sufficiently spread-out family of directions through sufficiently many
points forces more than $\delta^{-2b}$ distinct tubes, whereas an assumed
small cover of each direction fiber will later give at most
$\delta^{-2b-o(1)}$ tubes. The fixed excess $\epsilon_0$ in the lower bound
makes these two counts incompatible.

\begin{theorem}[Planar discretized Furstenberg estimate]
\label{ang:furstenberg}
Let $0<b<1$ and $b<f<2$. There is $\epsilon_0>0$ such that the
following holds for sufficiently small $\delta$. Suppose that $P$ is
a nonempty bounded $(\delta,f,\delta^{-\epsilon_0})$-set and that
each $x\in P$ has a nonempty family $\mathcal T_x$ of incident
$\delta$-tubes whose line parameters form a
$(\delta,b,\delta^{-\epsilon_0})$-set. Then
$$
 \left|\bigcup_{x\in P}\mathcal T_x\right|_\delta
       \ge\delta^{-2b-\epsilon_0}.
$$
The exponent may be decreased to absorb fixed bounded charts and fixed
width multiples.
\end{theorem}

This is the tube formulation in
\cite[definitions~2.5--2.6 and theorem~2.7]{OrponenShmerkinWang2024}, following the
discretized theorem of \cite[theorem~1.3]{OrponenShmerkin2023}. In the application
below, $f<1$; that range is included in the theorem.

A label is an element of the underlying measure space, and its geometric
coordinate is the image under the specified coordinate map. The measure
formulation below keeps the labels being sampled separate from their
geometric coordinates, since several labels may share a coordinate and
still carry different incidence fibers.

The planar theorem is stated for separated finite sets, whereas our data
are measures, possibly with repeated geometric coordinates. We pass from
one formulation to the other by sampling labels, retaining many incidences,
and then selecting separated coordinates, and throughout this operation a
direction label must keep its own tube cover, since two labels at the same
direction may have different fibers.

\begin{lemma}[Planar measure incidence obstruction]
\label{ang:planar-measures}
Fix $0<b<1$ and $b<f<2$. There is no sequence $\delta\to0$ with the
following data. Let $\mu,\nu$ be subprobabilities on Borel label sets
$X,W$ in Euclidean spaces, and let
$$
 p:X\longrightarrow{\mathbb R}^2 \text{ and } v:W\longrightarrow S^1
$$
be measurable coordinate maps with fixed bounded images. Suppose that
$p_\#\mu$ and $v_\#\nu$ have respective $f$- and $b$-ball bounds
of loss $\delta^{-o(1)}$ down to $\delta$. Let $E\subset X\times W$ be a measurable
incidence relation with
$$
 (\mu\otimes\nu)(E)\ge\delta^{o(1)}.
$$
For each original direction label $\omega\in W$, put
$E_\omega=\{\xi\in X:(\xi,\omega)\in E\}$. Assume that the set of
point coordinates $p(E_\omega)$ is covered by at most
$\delta^{-b-o(1)}$ strips of width $C\delta$ parallel to $v(\omega)$.
\end{lemma}

\begin{proof}
Restrict the coordinate images to fixed bounded charts and normalize the
two laws, which costs only subpower factors because their masses are at
least the incidence mass. Let $I$ be the resulting incidence mass, sample
$n=\lceil\delta^{-b}\rceil$ independent direction labels
$\omega_1,\ldots,\omega_n$ from $\nu$, and write $v_i=v(\omega_i)$. There
is a subpower quantity $T\to\infty$ such that, with failure probability
smaller than any prescribed power of $\delta$, the sample satisfies
\begin{equation}\label{ang:sample-caps}
\#\{i:v_i\in B(v_0,r)\}\le Tnr^b
        \text{ for } \delta\le r\le1.
\end{equation}
To verify this, use a fixed collection of dyadic grids in the direction
chart, of which only polynomially many cells are relevant. Their expected
counts are at most $Knr^b$, where $K=\delta^{-o(1)}$ and $nr^b\ge1$, so
choose $T$ larger than a fixed large multiple of
$K+\log(\frac{1}{\delta})$. For a binomial count $Y$, the bound
$$
 \Pr(Y\ge j)\le \left(e\frac{\mathbb E Y}{j}\right)^j
$$
and a union bound prove \eqref{ang:sample-caps}, including
arbitrary balls by fixed grid comparison. Increasing $T$ by another
subpower factor makes the failure probability smaller than $\frac{I}{4}$.

The expected empirical incidence mass equals $I$, and since it is bounded
by one, there is a sample satisfying \eqref{ang:sample-caps} whose
empirical incidence mass is at least $\frac{I}{2}$. The set $X_0$ of point
labels incident to at least $\frac{In}{4}$ sampled direction labels then
has $\mu$-mass at least $\frac{I}{4}$, its complement contributing at most
$\frac{I}{4}$ to the empirical incidence integral. Normalize $\mu|_{X_0}$
and sample $m=\lceil\delta^{-f}\rceil$ independent point labels; if their
coordinates are $x_1,\ldots,x_m$, the same argument gives a sample with
$$
 \#\{j:x_j\in B(x,r)\}\le T'mr^f
  \text{ for } \delta\le r\le1, \text{ where } T'=\delta^{-o(1)}.
$$
In particular the multiplicities in $\delta$-balls are subpower. Color a
spatial $\delta$-grid with finitely many colors in such a way that cells of
the same color are $\delta$-separated, choose a color with the most
occupied cells, and take one sampled label from each of them. Their point
coordinates form a separated set $P$ with $|P|_\delta\ge m\delta^{o(1)}$,
and the displayed local count makes $P$ a $(\delta,f,\delta^{-o(1)})$-set,
each $x\in P$ retaining the list of at least $\frac{In}{4}$ good sampled
direction labels which belong to its chosen point label.

For $x\in P$, take the lines exactly through $x$ in the coordinates
$v_i$ of those listed labels, and let $\mathcal T_x$ be their
$\delta$-tubes. The estimate~\eqref{ang:sample-caps} bounds the multiplicity
in a $\delta$-ball of directions by $CT$, and the count in an $r$-ball
by $Tnr^b$. Consequently
$$
 |\mathcal T_x|_\delta\ge \frac{cIn}{T}
  \text{ and }
  |\mathcal T_x\cap B(\ell,r)|_\delta\le CTnr^b.
$$
Here the direction component of the line metric supplies the lower
distance comparison; the reverse comparison holds because $x$ is bounded.
The normalized loss is therefore at most $\frac{CT^2}{I}=\delta^{-o(1)}$.
In particular, $\mathcal T_x$ is a nonempty
$(\delta,b,\delta^{-o(1)})$-set of line parameters even when different
sampled labels have the same direction coordinate.

For each sampled label $\omega_i$, its prescribed strips contain the
coordinates of all the sampled point labels incident to it, and moving the
central line of such a strip so as to pass through one of these points
changes its intercept by $O(\delta)$. Each original strip therefore
accounts for only a bounded number of $\delta$-cells in line-parameter
space, so that altogether the union of our line families has covering count
at most
$$
 n\delta^{-b-o(1)}=\delta^{-2b-o(1)}.
$$
This contradicts theorem~\ref{ang:furstenberg} for small $\delta$. The
count treats each original sampled direction label separately, and it
never asks for one common fiber cover after coincident coordinates have
been identified.
\end{proof}

\subsection{Slicing with controlled losses}

Before applying them, we state the estimates we need for Haar-random
subspaces, Haar measure being normalized to have total mass one.

A band of width $\delta$ about a plane consists of the points within
distance $\delta$ of that plane, and a covering tube is called active for a
point band when it meets that band in the bounded region under
consideration. Since a codimension-$h$ plane imposes $h$ small transverse
coordinates on one vector, a band of width $\delta$ has probability
comparable to $\delta^h$, while for several vectors the constraints can be
dependent, and it is their singular values which measure that dependence.
By the singular values of vectors $w_1,\ldots,w_j$ we mean here the
nonnegative square roots of the eigenvalues of the Gram matrix
$(w_i\cdot w_l)_{i,l=1}^j$, listed in decreasing order as
$d_1,\ldots,d_j$, or equivalently the singular values of the matrix whose
columns are those vectors. A singular value below $\delta$ creates no
further effective restriction, which accounts for the regularized
factor $(d_i+\delta)^{-h}$ below.

\begin{lemma}[Haar band estimates]\label{ang:haar}
Let $0\le h\le n-1$, and let $V$ be a Haar-random codimension-$h$
linear subspace of ${\mathbb R}^n$.
For a unit vector $w$,
$$
 \Pr\{\operatorname{dist}(w,V)\le\delta\}\asymp\delta^h.
$$
If $j\le n-h$ unit vectors $w_1,\ldots,w_j$ have singular values
$d_1\ge\cdots\ge d_j\ge0$, then
\begin{equation}\label{ang:haar-multiple}
\delta^{-jh}\Pr\{\operatorname{dist}(w_i,V)\le\delta\text{ for all }i\}
       \le C_{n,j}\prod_{i=1}^j(d_i+\delta)^{-h}.
\end{equation}
When $h=0$, the probability and the factors with exponent $h$ are one.
\end{lemma}

\begin{proof}
The one-vector estimate follows from spherical coordinates. For the
second estimate, singular-value decomposition implies that the band
conditions constrain the projections of orthonormal singular directions
$u_i$ by $\frac{C_j\delta}{d_i}$. Keep only the constraints whose right-hand
sides are less than a sufficiently small dimensional constant, and
order these right-hand sides increasingly. By rotational invariance,
one can fix $V^\perp$ and reveal successive vectors of a random
orthonormal frame. Conditional on the earlier constrained vectors,
projection onto $V^\perp$, restricted to their orthogonal complement,
has rank $h$ and nonzero singular values bounded below, because the sum
of squared lengths of the earlier projections is less than $\frac{1}{2}$
when the dimensional constant is chosen small enough. To see the
lower bound explicitly, suppose that $u_1,\ldots,u_\ell$ have already
been chosen, and let $U$ be their span. For $z\in V^\perp$,
$$
\begin{aligned}
 |P_{U^\perp}z|^2
 &=|z|^2-\sum_{i=1}^\ell|z\cdot u_i|^2\\
 &\ge\left(1-\sum_{i=1}^\ell
                  |P_{V^\perp}u_i|^2\right)|z|^2
 \ge\frac{1}{2}|z|^2.
\end{aligned}
$$
The maps $P_{U^\perp}|_{V^\perp}$ and
$P_{V^\perp}|_{U^\perp}$ are adjoints, so their nonzero singular
values agree and are at least $\frac{1}{\sqrt{2}}$. The next vector is
uniform on the unit sphere of that complement, and its probability of
having projection of length at most $r$ is at most $C_n r^h$, the
ambient dimension at this last conditional step being at least
$n-j+1\ge h+1$, which includes the limiting case of a one-dimensional
kernel. Multiplying the conditional estimates proves
$$
 \Pr\{w_i\text{ all in the band}\}
       \le C_{n,j}\prod_i\min\left(1,\frac{\delta}{d_i}\right)^h,
$$
which is equivalent to \eqref{ang:haar-multiple}.
\end{proof}

Shmerkin and Wang \cite[sections~6.1--6.3]{ShmerkinWang2025} use normalized bands
to define slices and combine sliced with projected thin-tube bounds.
Here we keep bands at width $\delta$ and control the retained incidence
mass, ball bounds, and active tubes for each original direction label.

Slicing is required to do three things at once here, namely to retain
incidence mass, to preserve ball bounds, and to reduce the number of active
tubes by the expected factor, and none of them follows merely from
selecting a slice carrying a large amount of mass. We therefore average
first, prove second-moment estimates, and only then remove the exceptional
slices and points, in a specified order, the normalizations $\delta^{-h}$
compensating for the probability of lying in a band. The proof also keeps
the direction restriction independent of the translation $z$ of the point
slice, which will matter at the endpoint; independence means here that the
restriction is chosen without using $z$, rather than probabilistic
independence of two events. For fixed $V$, the integral in $z$ of the total
point-band mass is uniformly bounded, and therefore a set of removed directions
can be charged once against that integral, which it could not be if the
direction restriction were allowed to depend on $z$ in an arbitrary way.

\begin{lemma}[Selection of slices]\label{ang:slice-selection}
Let $0\le h\le n-2$. Let point and direction subprobabilities $\mu,\nu$ in ${\mathbb R}^n$ have
$t$- and $q$-ball bounds with subpower losses, where $t,q>h$.
Assume the directions lie in a fixed small projective chart. Let an
incidence relation $E$ have mass at least subpower, and suppose each
direction fiber is covered by at most $\delta^{-q-o(1)}$ tubes of
width $C\delta$.

For Haar-random codimension-$h$ subspaces $V$ and
$z\in V^\perp$, set
$$
\begin{aligned}
 d\nu_V(w)&=\delta^{-h}1_{\{|P_{V^\perp}w|\le\delta\}}d\nu(w),\\
 d\mu_{V,z}(x)&=\delta^{-h}
       1_{\{|P_{V^\perp}x-z|\le\delta\}}d\mu(x).
\end{aligned}
$$
Integrate $z$ over a fixed sufficiently large ball. There are families of restrictions
$\widetilde\nu_V\le\nu_V$, $\widetilde\mu_{V,z}\le\mu_{V,z}$,
and relations $\widetilde E_{V,z}\subset E$, measurable in the displayed
parameters and restricted to the two bands defining $\mu_{V,z}$ and
$\nu_V$, such that
$$
 \mathbb E_V\int
 (\widetilde\mu_{V,z}\otimes\widetilde\nu_V)(\widetilde E_{V,z})dz
 \ge\delta^{o(1)},
$$
and the following conclusions hold.
\begin{enumerate}
\item The restricted measures have masses and ball losses at most
$\delta^{-o(1)}$, with exponents $q-h$ and $t-h$, respectively.
\item Each remaining direction fiber is covered by at most
$\delta^{-(q-h)-o(1)}$ width-$C'\delta$ tubes.
\item Put
$$
 p_{V,z}(x)=P_Vx \text{ and } v_V(w)=\frac{P_Vw}{|P_Vw|}
  \text{ for } |P_{V^\perp}w|\le\delta.
$$
The first map is projection to $V+z$ followed by translation to $V$;
the second is well-defined on the direction band for small $\delta$.
The coordinate pushforwards of the restricted measures under these maps
satisfy the ball bounds in item~1, with fixed constant changes. For every
original direction label $w$ in that band, the coordinates
$p_{V,z}(\{x:(x,w)\in\widetilde E_{V,z}\})$ have the tube cover in
item~2 with direction $v_V(w)$. This cover remains attached to $w$ even
when two projected direction coordinates coincide.
\end{enumerate}

The notation $\widetilde\nu_V\le\nu_V$ means a restriction of the
direction band measure, before any normalization. Likewise the
retained point measure is a submeasure of its point band measure.
The direction restriction $\widetilde\nu_V$ depends on $V$ only, and its
mass has an a priori subpower upper bound valid for every $V$.
\end{lemma}

\begin{proof}
We keep the averaging calculation separate from the restrictions. The
averaging calculation concerns the full band measures, and every
subsequent restriction decreases the relevant integrals; when we come to
estimate the incidence cost of a restriction, however, we must also
control the mass of the other measure in the product, and this dictates
the order of the four truncations below.

We first choose the fiber covers measurably, enlarging their width by a
fixed dimensional factor while preserving each actual direction $w$, and we
let $M_\delta=\delta^{-q-o(1)}$ be their common size bound. The choice has
to be made before slicing, so that a direction label has the same list of
tubes in every slice, and it can be made with a finite set of spatial
anchors while leaving the actual direction unchanged. Use a finite spatial
$\delta$-grid in a fixed bounded region containing one point on every axis
whose tube meets the source region; rounding such a point to the grid and
keeping the direction $w$ enlarges the width by only $C_n\delta$. The
candidate covers are then the finitely many lists of at most
$\lfloor M_\delta\rfloor$ grid anchors, and for each candidate $J$ the
function
$$
 w\longmapsto
 \mu\left(E_w\setminus\bigcup_{a\in J}T(w,a)\right)
$$
is measurable by Fubini, where $T(w,a)$ is the enlarged tube parallel
to $w$ through $a$. Choose the first candidate for which this function
is zero. Coverage up to a $\mu$-null set suffices, since every point
slice $\mu_{V,z}$ is absolutely continuous with respect to $\mu$.
For a relation measurable only in the completed product sigma-algebra,
first choose a Borel subset agreeing with the relation modulo a
product-null set, and discard the $\nu$-null set of
directions on which its fiber differs by positive $\mu$-mass. Every
direction slice is absolutely continuous with respect to $\nu$, so
this replacement also changes no sliced incidence.
Writing $J(w)$ for the selected candidate, now replace $E$ by
$$
 E\cap\left\{(x,w):x\in\bigcup_{a\in J(w)}T(w,a)\right\},
$$
with the exceptional null directions removed. This loses zero original
or sliced incidence and makes the selected covers contain every fiber
of the restricted relation.

Intersecting the relation with the two bands does not change its sliced
incidence, and the expected sliced incidence is comparable from below to
the original one. Indeed, integrating the point-band indicator in $z$ gives
a constant times $\delta^h$, while the one-vector Haar estimate supplies
another factor comparable to $\delta^h$, and these cancel the two slice
normalizations.

Write $M(V)=\|\nu_V\|$ and $L(V,z)=\|\mu_{V,z}\|$.
Uniformly in $V$,
\begin{equation}\label{ang:point-slice-total}
\int L(V,z)dz\le C.
\end{equation}
The two-vector estimate in a projective chart gives the direction-pair
kernel
\begin{equation}\label{ang:direction-pair-kernel}
\delta^{-2h}\Pr\{w,w'\text{ both within }\delta\text{ of }V\}
       \le C(|w-w'|+\delta)^{-h}.
\end{equation}
For points, integrate both band indicators in $z$ first. This yields
at most $C\delta^h$ times the indicator
$|P_{V^\perp}(x-x')|\le2\delta$. Haar averaging therefore gives
the analogous normalized kernel $C(|x-x'|+\delta)^{-h}$.

The strict inequality in the exponents comes from the following summation.
Suppose that a measure $\zeta$ has an $s$-ball bound with constant $K$ down
to $\delta$, and that $s>h$. For $\delta\le r\le1$, divide the ball of
radius $r$ about $z$ into annuli with radii decreasing by factors of two,
stopping at $\delta$. On an annulus of radius $u$ the integral of
$(|z-z'|+\delta)^{-h}$ is at most $CKu^{s-h}$, whereas the innermost ball
contributes at most $K\delta^{s-h}$, and hence the geometric series gives
$$
\int_{|z-z'|\le r}(|z-z'|+\delta)^{-h}d\zeta(z')
\le C_{s,h}Kr^{s-h}.
$$
We apply this with $s=q$ to directions and with $s=t$ to points, which
accounts both for the bounded second moments and for the exponents $q-h$
and $t-h$ in the close-pair bounds. At $s=h$ the series would instead have
a logarithmic loss and would give no positive close-pair exponent, and this
is why the endpoint requires a different slicing argument.
Dyadic annulus summation, using $q,t>h$, proves
\begin{equation}\label{ang:slice-second}
\mathbb E_V M(V)^2+\mathbb E_V\int L(V,z)^2dz
       \le\delta^{-o(1)},
\end{equation}
together with
\begin{equation}\label{ang:direction-close}
\mathbb E_V\iint_{|w-w'|\le r}d\nu_V(w)d\nu_V(w')
       \le\delta^{-o(1)}r^{q-h}
\end{equation}
and
\begin{equation}\label{ang:point-close}
\mathbb E_V\int\iint_{|x-x'|\le r}
             d\mu_{V,z}(x)d\mu_{V,z}(x')dz
       \le\delta^{-o(1)}r^{t-h}
\end{equation}
for $\delta\le r\le1$. In \eqref{ang:direction-close}, for example, the
shells below $r$ form a geometric sum with successive exponents $q-h>0$,
and the finest shell uses the original ball estimate at radius $\delta$.

In the next table $K\ge1$ is a subpower threshold, chosen large enough to
absorb the input losses and the reciprocal incidence mass, its precise
choice being made immediately after the table.
\begin{table}[htbp]
\centering\small
\renewcommand{\arraystretch}{1.25}
\begin{tabular}{p{4.2cm}p{4.1cm}p{5cm}}
\toprule
Restriction & Quantity removed & Incidence cost\\
\midrule
Direction-slice mass $M>K^3$ & Weighted mass at most $CK^{-2}$ & Integrate the point mass in $z$ first.\\
Point ball bound fails at threshold $K^6$ & Expected point mass at most $CK^{-5}$ & Multiply by $M\le K^3$ to obtain $CK^{-2}$.\\
Direction ball bound fails at threshold $K^3$ & Expected direction mass at most $CK^{-2}$ & The integrated point mass is bounded independently of $V$.\\
Too many active tubes & Bad $z$-volume at most $K^{-11}$ per $(V,w)$ & Point mass $O(K^6)$ and direction mass $K^3$ give $CK^{-2}$.\\
\bottomrule
\end{tabular}
\caption{The order of restrictions in slice selection. The initial
integrated incidence is at least $K^{-1}$, so four losses of order
$K^{-2}$ leave positive incidence for sufficiently large $K$.}
\label{tab:slice-losses}
\end{table}

The truncations below first bound the direction-slice mass, then enforce
both ball bounds, and finally bound the active tube count. Take a subpower
number $K\to\infty$ absorbing all the constants, all the input losses, the
dyadic radius sums in the equations~\eqref{ang:direction-close} and
\eqref{ang:point-close}, and the inverse original incidence mass.
Enlarge it so the sliced incidence is at least $K^{-1}$ and the
right-hand side of \eqref{ang:slice-second} is at most $K$.
First remove planes with $M(V)>K^3$. By
the equations~\eqref{ang:point-slice-total} and \eqref{ang:slice-second},
their incidence contribution is at most $CK^{-2}$.

In the remaining point slices, remove points at which
$$
 \mu_{V,z}(B(x,Cr))>K^6r^{t-h}
$$
at some dyadic radius $r\in[\delta,1]$. Markov's inequality together with
\eqref{ang:point-close} bounds the expected deleted point mass by
$CK^{-5}$, and since $M(V)\le K^3$, its incidence cost is at most
$CK^{-2}$. The restricted point measure then has the asserted ball bound at
every center, since any ball meeting the retained set can be enlarged by a
fixed factor around a retained point, and it has mass $O(K^6)$ by covering
the bounded support with a fixed number of unit balls.

Remove similarly the direction-slice points for which
$\nu_V(B(w,Cr))>K^3r^{q-h}$ at some dyadic radius. Their expected mass is
at most $CK^{-2}$ by \eqref{ang:direction-close}, and the
estimate~\eqref{ang:point-slice-total} bounds the corresponding integrated
incidence loss by the same quantity, no independence between direction mass
and point mass being used at this step.

Finally, for a direction $w$ in the band about $V$, let $J(w)$ be its
selected list of grid anchors. A bounded tube through $a$, parallel to $w$,
can meet the source region in the band about $V+z$ only if
$|P_{V^\perp}a-z|\le C'\delta$, and its transverse projection within the
bounded source region has diameter $O(\delta)$, so that the measurable
upper count
$$
 N(V,z,w)=\sum_{a\in J(w)}
       1_{\{|P_{V^\perp}a-z|\le C'\delta\}}
$$
dominates the number of active covering tubes and satisfies
$$
 \int N(V,z,w)dz
       \le C M_\delta\delta^h\le K\delta^{h-q}.
$$
Delete the direction and $z$ pairs for which
$N(V,z,w)>K^{12}\delta^{h-q}$, whose $z$-volume for each $V,w$ is at most
$K^{-11}$. Since the truncated point mass is $O(K^6)$ and the truncated
direction mass is at most $K^3$, this costs at most $CK^{-2}$ in integrated
incidence. Delete the corresponding incidences while keeping the marginal
restrictions already constructed; at large $K$ the surviving integrated
incidence is still at least $(2K)^{-1}$.

Projection changes each direction by $O(\delta)$, and the points of one
band whose coordinates $p_{V,z}(x)$ lie in a ball of radius $r\ge\delta$
lie in a ball of radius $Cr$ in the original space, the same being true of
directions after unit normalization. These preimage bounds give the
asserted ball bounds on the coordinate pushforwards, and for each original
label $w$ every active tube in its own selected list becomes a tube of
width $C'\delta$ parallel to $v_V(w)$, which proves item~3 for each
original direction label.
\end{proof}

\begin{figure}[htbp]
\centering
\begin{tikzpicture}[>=Stealth,font=\small]
\fill[gray!15] (0,0)--(8,0)--(10,2)--(2,2)--cycle;
\draw[thick] (0,0)--(8,0)--(10,2)--(2,2)--cycle;
\draw[dashed] (0,.5)--(8,.5)--(10,2.5)--(2,2.5)--cycle;
\draw[dashed] (0,-.5)--(8,-.5)--(10,1.5);
\draw[very thick,->] (2,.75)--(7.5,1.0) node[right] {$w$};
\draw[->] (5,1)--(5,3.4) node[above] {$V^\perp$};
\fill (6,1.5) circle (2pt) node[above right] {$x$};
\draw[dotted,thick] (6,1.5)--(6,1);
\node[below] at (4,-.55) {point band $|P_{V^\perp}x-z|\le\delta$};
\node[align=left] at (11,1) {$V+z$\\point slice};
\node[align=center] at (5,-1.6)
{direction band $|P_{V^\perp}w|\le\delta$\\
each band measure is multiplied by $\delta^{-h}$};
\end{tikzpicture}
\caption{A thickened affine slice, drawn in three dimensions, the
actual codimension being $h$. A bounded tube in direction $w$ has
transverse displacement $O(\delta)$ when $w$ lies in the direction
band, so it can meet the point band only for a set of
translations $z$ of volume $O(\delta^h)$, and this reduces the average
number of active tubes. The direction band depends on
$V$ alone, while the point band depends on $z$ as well.}
\label{fig:bands}
\end{figure}
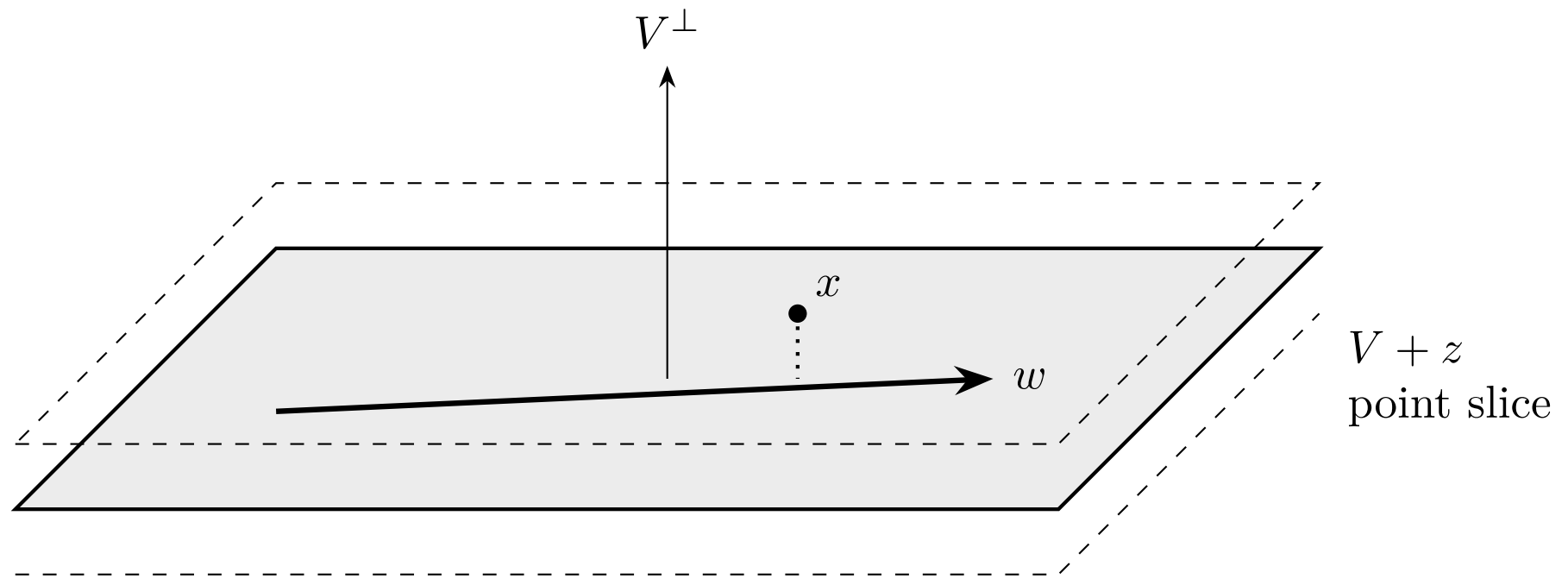

\begin{table}[htbp]
\centering\small
\renewcommand{\arraystretch}{1.35}
\begin{tabular}{p{3cm}p{4.7cm}p{5.5cm}}
\toprule
Original exponent & Slice and remaining exponents & Estimate used in the slice\\
\midrule
$q>k-1$ & $h=k-1$; dimension $2$; directions $q-h>0$; points $t-h>q-h$
& Planar measure incidence obstruction, lemma~\ref{ang:planar-measures}.\\
$q=k-1$ & $h=k-2$; dimension $3$; directions $1$; points $t-h>\frac{3}{2}$
& Three transverse directions, weak radial control of nearly dependent triples, and Loomis--Whitney.\\
\bottomrule
\end{tabular}
\caption{Why the endpoint requires a different slice. A planar slice
at $q=k-1$ would leave direction exponent zero, outside the planar
incidence theorem used here.}
\label{tab:slicing-cases}
\end{table}

\subsection{Endpoint preparation by transverse tubes and weak radial estimates}

We begin with an elementary counting fact which is used twice below. In
the plane a horizontal strip and a vertical strip meet in a single small
rectangle, and the higher-dimensional version of this uses all the
coordinate projections at once. One first changes coordinates so that the
given transverse directions become the coordinate axes, the determinant
controlling the distortion of that change of variables, and the discrete
Loomis--Whitney inequality then bounds how many small cubes can satisfy all
the tube-cover restrictions simultaneously.

\begin{lemma}[Transverse tube counting]\label{ang:loomis-whitney}
Let $v_1,\ldots,v_m$ be unit vectors in ${\mathbb R}^m$, $m\ge2$, and suppose
$|\det(v_1,\ldots,v_m)|\ge\Delta>0$. If a bounded set $A$ is covered,
for each $i$, by $M_i$ tubes of width $C\delta$ parallel to $v_i$, then
$$
 |A|_{C_m\delta}\le C_{m,C}\Delta^{-m}
       \prod_{i=1}^m M_i^{\frac{1}{m-1}}.
$$
Changing the fixed constant in the covering radius changes only the
constant in this inequality.
\end{lemma}

\begin{proof}
In this proof $\det(v_1,\ldots,v_m)$ denotes the determinant of the matrix
with those columns, matrix norms are Euclidean operator norms,
$\operatorname{pr}_{\widehat i}$ denotes the coordinate projection which
deletes the $i$th coordinate, and $X_{\widehat i}$ denotes the random
vector with the coordinate $i$ omitted. Let $L$ be the matrix with columns
$v_i$, so that $\|L\|\le C_m$ and $\|L^{-1}\|\le C_m\Delta^{-1}$, and let
$S\subset{\mathbb Z}^m$ index the $\delta$-grid cubes meeting $L^{-1}A$.
The image of such a cube has diameter at most $C_m\delta$, while each
original tube parallel to $v_i$ becomes a tube parallel to the $i$th
coordinate axis, of width at most $C_{m,C}\delta\Delta^{-1}$. Therefore
$$
 |\operatorname{pr}_{\widehat i}S|
       \le C_{m,C}\Delta^{-(m-1)}M_i.
$$
The discrete Loomis--Whitney inequality \cite{LoomisWhitney1949}
$|S|^{m-1}\le\prod_i|\operatorname{pr}_{\widehat i}S|$ now proves the
claim, and we include a proof of that inequality. For a finite random
variable $Y$, define
$$
H(Y)=-\sum_y\mathbb P(Y=y)\log\mathbb P(Y=y),
$$
with $0\log0=0$. If $Y$ takes at most $M$ values, then comparison with the
uniform law, by nonnegativity of the relative entropy, gives
$H(Y)\le\log M$. Define the conditional entropy by averaging the entropy of
the conditional laws, and factor a joint probability into a marginal
probability and a conditional probability, which gives the chain rule
$H(Y,Z)=H(Z)+H(Y\mid Z)$. Conditioning on more variables cannot increase
the conditional entropy, since each law before the extra conditioning is a
convex combination of the laws after it, and $-x\log x$ is concave.

Now take $X=(X_1,\ldots,X_m)$ uniform on $S$. For each $i$,
apply the chain rule to $X_{\widehat i}$ in the original coordinate
order. Conditioning additionally on the omitted earlier coordinate,
when it is present, can only decrease each term. Thus
$$
 H(X_{\widehat i})\ge
 \sum_{j\ne i}H(X_j\mid X_1,\ldots,X_{j-1}).
$$
Summing in $i$ yields $(m-1)\log|S|\le
\sum_i\log|\operatorname{pr}_{\widehat i}S|$, as required.
\end{proof}

The entropy inequality in this proof is Han's entropy form of
Loomis--Whitney \cite{Han1978}, a special case of Shearer's inequality
\cite{ChungGrahamFranklShearer1986}. Its relation
to projection counts and the proof by the chain rule are discussed in
\cite[sections~3--4]{BalisterBollobas2012}. The factor $m-1$ has a
simple explanation, since each coordinate survives in exactly $m-1$ of the
projections, and no independence of the coordinates is assumed anywhere.

For $x\ne z$, put $\pi_x(z)=\frac{z-x}{|z-x|}$. A restriction of the
product measure to a pair set is not normalized in the next statement.

At the endpoint of the angular gain we must exclude almost dependent
triples of directions, and the following radial estimate supplies that
exclusion. Its exponent is deliberately below the final angular
exponent, since only the strict inequality $\ell>k-2$ will be used there.
The proof gives an estimate at one scale first, and then sums the
exceptional sets over scales, using a slightly stronger intermediate
exponent so that this summation still loses only a small power of the final
scale. For an affine hyperplane $H$, we write
$H^{(r)}=\{x:\operatorname{dist}(x,H)\le r\}$ for its closed $r$-neighborhood.

\begin{lemma}[Weak radial estimate]\label{ang:weak-radial}
Fix $k\ge2$, $\zeta>0$, and
$$
 k-2<\ell<\frac{(k-1)^2}{k}.
$$
For every $e>0$ there exist $\epsilon,c>0$ with the following property
for sufficiently small $\delta$. Suppose that a probability $\lambda$
in a bounded subset of ${\mathbb R}^k$ satisfies
$$
 \lambda(B(x,r))\le\delta^{-\epsilon}r^{k-1}
 \text{ and }
 \lambda(H^{(r)})\le\delta^{-\epsilon}r^\zeta
  \text{ for } \delta\le r\le1
$$
and for every affine hyperplane $H$.
There is a measurable pair set $G$, omitting the diagonal, such that
$$
 (\lambda\otimes\lambda)(G^c)\lesssim\delta^c
 \text{ and }
  \lambda\{z:(x,z)\in G,\ \pi_x(z)\in B(\theta,r)\}
       \le\delta^{-e}r^\ell
$$
for every $x$, every spherical cap center $\theta$, and
$\delta\le r\le1$. The conclusion
also holds for balls in projective direction space, with a fixed change
of constant, by taking the two antipodal caps.
\end{lemma}

\begin{proof}
The contradiction at one scale compares the same common source set in two
different ways. Regularity gives a lower bound for the number of its small
cells, while several transverse tube covers give an upper bound through
Loomis--Whitney, the transverse directions being obtained by choosing
several pins with a large common source fiber and then excluding the nearly
dependent choices by means of the plate bound.

We first prove a single-scale assertion at $u=2^{-N}$, with reciprocal
scale $R=2^N$. For any fixed $\ell_1<\frac{(k-1)^2}{k}$, sufficiently small
fixed losses in the hypotheses permit deletion of $O(R^{-c_1})$ pair
mass so that every remaining fiber has cap mass at most $R^{-\ell_1}$
at radius $u$. All choices made below are fixed independently of $R$.

Apply lemma~\ref{reg:regularization} to the source measure with a small
parameter $\tau>0$. It is enough to treat one normalized source class
$\rho$, retaining the original pin law $\lambda$, since the discarded
source mass is a fixed negative power of $R$ and the classes will be
recombined with their original weights. Write $m$ for the class mass curve;
its comparison errors are $R^{C_k\tau}$, it is $k$-Lipschitz, and
$$
 m(1)\ge k-1-C_k\tau-\epsilon.
$$
The constant absorbs the normalization of a class of mass at least
$R^{-\tau}$. At depths $a<b$, an occupied depth-$b$ child of an
occupied depth-$a$ cell $Q$ has conditional $\rho$-mass at most
\begin{equation}\label{ang:conditional-child}
R^{C_k\tau}D^{-1} \text{ with } D=R^{m(\frac{b}{N})-m(\frac{a}{N})}.
\end{equation}
Divide the depth interval into $J$ nearly equal blocks, with $J$ a
large fixed integer, and ignore the first two blocks. In every later
block $[a,b]$, remove the pairs whose pin is within $R^{-\frac{1}{4J}}$ of
the center $z_Q$ of the source cell $Q$ at depth $a$. Summing with the
source-cell weights shows that this costs at most
$C_JR^{-\frac{k-1}{4J}+\epsilon}$ pair mass. For the other pairs the center
direction
$$
 v(x,Q)=\frac{z_Q-x}{|z_Q-x|}
$$
differs from the direction to any point of $Q$ by
$O(2^{-a}R^{\frac{1}{4J}})$. In particular the transverse displacement inside
$Q$ caused by freezing this direction is
\begin{equation}\label{ang:freezing-error}
O\left(2^{-2a}R^{\frac{1}{4J}}\right)=o(2^{-b});
\end{equation}
here $a\ge\frac{2N}{J}+O(1)$ and $b-a=\frac{N}{J}+O(1)$.

Fix a small $e'>0$. A pair $(x,z)$, $z\in Q$, is troublesome if the
full conditional source mass in the width-$C2^{-b}$ tube about the
line through $x,z$ exceeds
\begin{equation}\label{ang:local-fraction}
F=R^{e'}D^{-\frac{k-1}{k}}.
\end{equation}
We say that a family of sets has bounded overlap if each point belongs to
at most a fixed number of its members, and the fixed width $C$ may be taken
as large as any subsequent application requires. If $F\ge1$ there are no
troublesome pairs, whereas otherwise the estimate~\eqref{ang:freezing-error}
and a grid in $v(x,Q)^\perp$ show that all the troublesome sources for this
pin lie in at most $\frac{C}{F}$ parallel tubes of width $C'2^{-b}$.
Indeed, each relevant projected grid cell has a fixed enlargement of
conditional mass greater than $F$, and these enlargements have bounded
overlap.

For this fixed $Q$, let $B_x$ be its troublesome source set when the
pin survives the near-pin deletion, and let $B_x=\varnothing$ otherwise.
Put $\rho_Q=\frac{\rho|_Q}{\rho(Q)}$ and
$$
 I_Q=\int\rho_Q(B_x)d\lambda(x).
$$
To see why a large $I_Q$ gives a common source set for many
pins, put $f_Q(z)=\lambda\{x:z\in B_x\}$. Then
$\int f_Qd\rho_Q=I_Q$, and Fubini identifies $\int f_Q^k d\rho_Q$
with the average mass of the intersection of $k$ pin fibers.
Since $\rho_Q$ is a probability, Jensen gives
$\int f_Q^kd\rho_Q\ge I_Q^k$.
This is the common-neighborhood averaging calculation used in dependent
random choice, written for measures
\cite[section~2]{FoxSudakov2011}.

We claim that $I_Q\le R^{-c_2}$, uniformly in $Q$. Suppose instead that
$I_Q\ge R^{-c_2}$. Jensen's inequality, applied to the pin-incidence
mass at each source, gives
\begin{equation}\label{ang:common-k-incidence}
\int \rho_Q\left(\bigcap_{i=1}^k B_{x_i}\right)
           d\lambda(x_1)\cdots d\lambda(x_k)
      \ge I_Q^k.
\end{equation}
Tuples containing a removed pin contribute zero to this integral. For a
tuple of surviving pins, if
$|\det(v(x_1,Q),\ldots,v(x_k,Q))|<R^{-\eta}$, one of the
successive Gram--Schmidt distances is less than
$R^{-\frac{\eta}{k-1}}$. Conditional on the preceding pins, the corresponding
pin lies in a $C R^{-\frac{\eta}{k-1}}$-neighborhood of an affine hyperplane
through $z_Q$. Thus the product measure of these tuples is at most
$$
 C_kR^{\epsilon-\frac{\eta\zeta}{k-1}}.
$$
Choose $\eta$ small compared with $e'$, and then choose
$c_2,\epsilon$ so small that this is $o(R^{-kc_2})$. There is a tuple
of determinant at least $R^{-\eta}$ whose common source set has
conditional mass at least $\frac{1}{2}R^{-kc_2}$. By \eqref{ang:conditional-child}, it needs at least
$$
 c D R^{-kc_2-C_k\tau}
$$
balls of radius comparable to $2^{-b}$. On the other hand,
lemma~\ref{ang:loomis-whitney}, applied to its $k$ tube covers, bounds
this number by
$$
 C D R^{k\eta-\frac{e'k}{k-1}}.
$$
The two estimates are contradictory if
$$
 k\eta+kc_2+C_k\tau<\frac{ke'}{k-1}.
$$
These parameter requirements are compatible, since we may choose $\eta$
first, then $c_2$, and finally $\tau,\epsilon$, which proves the claim.
Deleting the troublesome pairs and summing over $Q$ costs at most
$R^{-c_2}$ per block, because the weights $\rho(Q)$ sum to one.

Apply the physical tube form of lemma~\ref{reg:pruning-multiplication}
to the surviving pairs. At a retained witness $(x,z)$ the local test
in \eqref{ang:local-fraction} bounds the full conditional source mass
in the tube about the line through $x$ and $z$, and for any final cap
the hit children are contained in the enlarged tube about one such
witness from their parent, so that their full masses sum to at
most $F\rho(Q)$. This comparison uses only the boundedness of the
supports, and it holds also for pins close to the source.

Multiplying these inequalities telescopes the mass curve increments. Sum
the resulting bounds over the occupied initial cells at the end of the two
ignored blocks; their masses sum to one, and therefore no factor equal to the
number of cells is introduced, and the ignored growth is at most
$\frac{2k}{J}+o(1)$. Consequently every final cap has normalized class mass
at most
$$
 R^{Je'+o(1)}
 R^{-\frac{k-1}{k}(k-1-\frac{2k}{J}-C_k\tau-\epsilon)}.
$$
Choose $J$ large depending on $\frac{(k-1)^2}{k}-\ell_1$, then $e'$
small, and then the remaining losses as above, which proves the
single-scale assertion with all deletion exponents positive fixed numbers.
Recombining the classes with their original weights preserves the cap
estimate and the power deletion, the negligible discarded class mass being
included in the deletion.

Now choose $\ell<\ell_1<\frac{(k-1)^2}{k}$, and note that it is enough
to consider $0<e<1$. Apply the single-scale assertion at dyadic $u$
between $\delta$ and $\delta^{\frac{e}{2k}}$; at these scales the input
loss $\delta^{-\epsilon}$ is at most $u^{-\frac{2k\epsilon}{e}}$, so that
choosing $\epsilon$ sufficiently small makes all the single-scale
hypotheses uniformly valid. The sum of the deletions is at most
$$
 C\log\left(\frac{1}{\delta}\right)\delta^{\frac{ec_1}{2k}},
$$
which is a positive power of $\delta$ after its exponent has been
decreased. Intersect the surviving pair sets; dyadic comparison then gives
the desired cap estimate for $r\le\delta^{\frac{e}{2k}}$, the slack between
$\ell$ and $\ell_1$ absorbing the fixed constants, while for larger $r$ we
have $\delta^{-e}r^\ell\ge1$ at small scales and the estimate is trivial.
\end{proof}

\subsection{Proof of the linear gain}

\begin{proof}[Proof of theorem~\ref{ang:linear-gain}]
We rule out any sequence with subpower input losses, subpower incidence
mass from below, and cover size $\delta^{-q-o(1)}$, which implies the
stated fixed-power result. Indeed, if no positive $\epsilon,c$ worked, we
could choose counterexamples with $\epsilon=c=\frac{1}{j}$ at scales
$\delta_j<\frac{1}{j}$ tending to zero, and they would form precisely such
a sequence. The same observation justifies all the later uses of the gain
with sufficiently small fixed losses.

By choosing one of finitely many small projective charts, we may assume
that all the relevant directions lie in a single chart and that the
incidence still has subpower mass. If $q>k-1$, apply
lemma~\ref{ang:slice-selection} with $h=k-1$, select a slice with subpower
incidence, and normalize its marginal restrictions, whose masses are
bounded above by subpower quantities and below by the incidence divided by
the other marginal mass. The coordinate maps then take values in a plane,
with
$$
 b=q-h\in\left(0,\frac{1}{2}\right] \text{ and } f=t-h\in\left(\frac{1}{2},1\right),
$$
so that $b<f$, and with tube count at most $\delta^{-b-o(1)}$. Apply
lemma~\ref{ang:planar-measures} with the original sliced points and
directions as labels, and with coordinate maps $p_{V,z}$ and $v_V$.
Item~3 of lemma~\ref{ang:slice-selection} gives the ball bounds on their
coordinate pushforwards and the separate cover for each original direction
label. The labeled incidence relation still has subpower mass, so these
data contradict the planar obstruction.

At $q=k-1$, slicing to the plane would leave direction exponent zero,
outside the range of lemma~\ref{ang:planar-measures}. We therefore put
$h=k-2$, so the selected slices have dimension three, and let
$$
 g=t-h-\frac{3}{2}>0.
$$
Choose once and for all $0<a<\frac{g}{4}$. We will find in one slice three
directions of determinant at least $c\delta^a$ with a common incident
point set of mass at least subpower. The resulting loss in
lemma~\ref{ang:loomis-whitney} is $O(\delta^{-3a})$, which the
strict gap $g$ will absorb.

Parameterize our chart by $p\in{\mathbb R}^k$, with unit direction
$w(p)=\frac{(1,p)}{\sqrt{1+|p|^2}}$. Restricting and normalizing the chart
costs only a subpower factor; denote the resulting law in $p$ by
$\lambda$. In the rest of the endpoint proof the direction law supplied
to slicing is $w_\#\lambda$, and the incidence relation is restricted to
this chart. The slice law of the chart labels is
$$
 d\lambda_V(p)=\delta^{-h}
 1_{\{|P_{V^\perp}w(p)|\le\delta\}}d\lambda(p)
 \text{ and }
  \nu_V=w_\#\lambda_V.
$$
The estimate~\eqref{ang:endpoint-plate}
becomes the affine-hyperplane bound in the $p$ coordinates, because
an affine hyperplane there is the intersection of a linear hyperplane
with the chart. The ball exponent remains $k-1$. Choose
$$
 h<\ell<\frac{(k-1)^2}{k}=h+\frac{1}{k}.
$$
Apply lemma~\ref{ang:weak-radial} to this chart law, with a small
fixed output loss $\delta^{-e}$. Its hypotheses hold at sufficiently
small scales for every fixed input-loss threshold. This deletes a
pair set of original product mass $O(\delta^{c_e})$ and gives, for
every anchor $p$, an $\ell$-cap bound $\delta^{-e}r^\ell$ on
the undeleted displacement directions. Also delete pairs at distance
less than $\delta^{b_0}$, where $b_0>0$ will be small. Denote the union
of these two exceptional pair relations by $\mathcal E_{\mathrm{pair}}$,
and call an anchored pair retained when it lies outside this union.

We first record how these deletions transfer to normalized slice-pair
mass. The estimate~\eqref{ang:direction-pair-kernel}
shows that for every measurable chart-pair relation $\mathcal E$ with
$(\lambda\otimes\lambda)(\mathcal E)\le\delta^c$, and every fixed
$\theta>0$,
\begin{equation}\label{ang:exceptional-pair-transfer}
\mathbb E_V(\lambda_V\otimes\lambda_V)(\mathcal E)
 \le C\delta^{c-\theta h}
       +\delta^{\theta(q-h)-o(1)}.
\end{equation}
For $0<\theta\le1$, the first term treats distances greater than
$\delta^\theta$, and the second treats smaller distances by the ball
bound and annular summation. For $\theta\ge1$, the kernel is at most
$C\delta^{-h}$, so the whole expectation is at most
$C\delta^{c-h}\le C\delta^{c-\theta h}$. The same annulus calculation
bounds the expected mass of the close pairs by
$\delta^{b_0(q-h)-o(1)}$. We leave $\theta$ free until $e$, and hence
the weak-radial deletion exponent $c_e$, has been fixed.

We next bound triples that are close to linearly dependent. For two
retained anchored pairs $(p,p'),(p,p'')$, set
$$
 \alpha=\left|
 \frac{p'-p}{|p'-p|}\wedge\frac{p''-p}{|p''-p|}\right|.
$$
In a bounded chart,
\begin{equation}\label{ang:chart-wedge}
|w(p)\wedge w(p')\wedge w(p'')|
       \asymp |p'-p||p''-p|\alpha.
\end{equation}
To see this, put $u=p'-p$ and $v=p''-p$ and subtract the first
unnormalized column from the other two. The orthogonal decomposition
of exterior powers gives the identity
$$
 |(1,p)\wedge(0,u)\wedge(0,v)|^2
       =|u\wedge v|^2+|p\wedge u\wedge v|^2.
$$
It supplies the lower bound, while the boundedness of $p$ supplies the
upper bound, and normalizing the three columns changes only fixed
constants. A wedge smaller than $2\delta^a$ therefore implies
$\alpha\le C\delta^{a-2b_0}$.

For $(w_1,w_2,w_3)=(w(p),w(p'),w(p''))$, the largest singular value is
comparable to one, the first two columns span area at least
$c\delta^{b_0}$, so that the second singular value is at least
$c\delta^{b_0}$, and by \eqref{ang:chart-wedge} the third is at least
$c\delta^{2b_0}\alpha$. Lemma~\ref{ang:haar} therefore bounds
their normalized expected slicing weight by
\begin{equation}\label{ang:triple-weight}
\delta^{-3h}\Pr_V\{\max_{1\le i\le3}|P_{V^\perp}w_i|\le\delta\}
 \le C\delta^{-3hb_0}(\alpha+\delta)^{-h}.
\end{equation}
This includes singular values below $\delta$, since the kernel is
truncated there. At fixed $p,p'$, the weak radial estimate and
the two antipodal caps give
$$
 \lambda\{p'':(p,p'')\text{ is retained},\ \alpha\le r\}
       \le C\delta^{-e}r^\ell
        \text{ for } \delta\le r\le1.
$$
Let $\mathcal T_\delta$ be the chart triples whose two anchored pairs
are retained and whose original wedge has size less than $2\delta^a$.
Summing dyadic annuli in \eqref{ang:triple-weight}, up to
$C\delta^{a-2b_0}$, gives
\begin{equation}\label{ang:bad-triples}
\mathbb E_V\lambda_V^{\otimes3}(\mathcal T_\delta)
 \le C\delta^{(a-2b_0)(\ell-h)-e-3hb_0}.
\end{equation}
Since $\ell>h$, first choose $b_0>0$ sufficiently small depending
on $a,\ell,h$, and then $e>0$ sufficiently small, so that this
exponent is positive. These choices also ensure $a-2b_0>0$.
With $e$ fixed, the weak radial estimate fixes $c_e>0$. Now choose
$\theta>0$ small enough that $c_e-\theta h>0$. Since $q-h=1$, the
second exponent in \eqref{ang:exceptional-pair-transfer} is
also positive. Together with the close-pair bound, this proves that
$\mathcal E_{\mathrm{pair}}$ has power-small expected slice-pair mass.

Perform the restrictions in lemma~\ref{ang:slice-selection}, and
write $\widetilde\lambda_V=(w^{-1})_\#\widetilde\nu_V$ for the restricted
law of chart labels. It is dominated by $\lambda_V$ and has the same
uniform subpower mass bound as $\widetilde\nu_V$. Define
$$
 F(V,z,x)=\widetilde\lambda_V
 \{p:(x,w(p))\in\widetilde E_{V,z}\}.
$$
Its integrated first moment is the retained sliced incidence and is at
least subpower, whereas the total integrated point mass is bounded by a
fixed constant.
Regard $dVdzd\widetilde\mu_{V,z}(x)$ as one finite measure of mass at
most $C$. If its first moment of $F$ is $I$, H\"older gives
$$
I\le\left(\int F^3\right)^{\frac{1}{3}}C^{\frac{2}{3}},
 \text{ so that }
 \int F^3\ge\frac{I^3}{C^2},
$$
the integrals in this calculation including all three
variables. Since $I$ is at least subpower, this lower bound is
also at least subpower, and thus
\begin{equation}\label{ang:triple-incidence}
\mathbb E_V\int\int F(V,z,x)^3
                 d\widetilde\mu_{V,z}(x)dz
       \ge\delta^{o(1)}.
\end{equation}
Dropping the point-incidence indicators bounds the contributions of
bad triples and exceptional anchored pairs by the left-hand sides below.
The estimate~\eqref{ang:point-slice-total} and
$\widetilde\lambda_V\le\lambda_V$ give
$$
\begin{aligned}
 \mathbb E_V\int \|\widetilde\mu_{V,z}\|
      \widetilde\lambda_V^{\otimes3}(\mathcal T_\delta)dz
 &\le C\mathbb E_V\lambda_V^{\otimes3}(\mathcal T_\delta),\\
 \mathbb E_V\int \|\widetilde\mu_{V,z}\|\|\widetilde\lambda_V\|
      (\widetilde\lambda_V\otimes\widetilde\lambda_V)
          (\mathcal E_{\mathrm{pair}})dz
 &\le C\sup_V\|\widetilde\lambda_V\|
      \mathbb E_V(\lambda_V\otimes\lambda_V)
          (\mathcal E_{\mathrm{pair}}).
\end{aligned}
$$
Both quantities are power-small by \eqref{ang:bad-triples} and by the
preceding bound for $\mathcal E_{\mathrm{pair}}$, the latter covering the
two possible anchored pairs as well after a change of constant, and no
independence between the point and direction slices is used. Hence a
subpower amount of the triple incidence in \eqref{ang:triple-incidence} has
original wedge at least $2\delta^a$, and since projection to $V$ and unit
normalization change each vector by $O(\delta)$, their determinant in $V$
is at least $\delta^a$ for small $\delta$.

The total parameter mass for integrating $(V,z)$ and three restricted
chart labels is at most subpower, since $\|\widetilde\lambda_V\|$ has a
uniform subpower bound and the $z$-domain is bounded.
It follows that some slice and some such triple have a common
incident point set of sliced mass at least subpower. Let $A$ be its
projected point coordinates. The pushforward of its restricted point
measure has the same mass and inherits the point-coordinate ball bound,
which implies
$$
 |A|_\delta\ge\delta^{-(t-h)+o(1)}.
$$
For each of the three original direction labels, the projected common
point set has its own cover by
$\delta^{-(q-h)-o(1)}=\delta^{-1-o(1)}$ tubes in that label's projected
direction. By
lemma~\ref{ang:loomis-whitney},
$$
 |A|_\delta\le\delta^{-\frac{3}{2}-3a-o(1)}.
$$
The exponents now give the contradiction directly, since the lower bound is
$\delta^{-\frac{3}{2}-g+o(1)}$ while the upper bound is
$\delta^{-\frac{3}{2}-3a-o(1)}$, and their ratio tends to infinity because
$g-3a>0$. When $k=2$ we have $h=0$, so that all the slicing weights equal
one and the same proof applies, with the $z$-integral interpreted as
integration over the singleton zero-dimensional space. This proves the
endpoint, and with it the theorem.
\end{proof}

The gain just proved concerns an incidence relation, and the next corollary
packages it as a statement about a set of bad directions, independent of
any particular witness point, which is the form needed when a radial
direction is moved to the center of a small spatial cell. By the outer
$\nu$-measure of a set we mean the infimum of $\nu(B)$ over all Borel sets
$B$ containing it, which allows us to bound the geometric bad set without
first proving that this set is Borel.

\begin{corollary}[Small mass of bad directions]\label{ang:bad-directions}
Fix the parameters of theorem~\ref{ang:linear-gain}. There are
$c_1,\epsilon_1>0$ such that the following holds at sufficiently small
$\delta$. Let $\mu$ be a bounded point probability satisfying the
$t$-ball bound with loss $\delta^{-\epsilon_1}$. Let $B_\delta$
be the set of directions for which at most
$\delta^{-q-c_1}$ width-$C\delta$ tubes capture $\mu$-mass at
least $\delta^{c_1}$. There is a Borel set
$B_\delta^+\supset B_\delta$, depending only on the point law and this
geometric test, such that every direction subprobability $\nu$
satisfying the $q$-ball bound with loss $\delta^{-\epsilon_1}$,
and, when $q=k-1$, the estimate~\eqref{ang:endpoint-plate} with
$\epsilon_1$ in place of $\epsilon$, satisfies
$$
 \nu(B_\delta^+)\le\delta^{c_1}.
$$
In particular this bounds the outer $\nu$-measure of the original
geometric bad set. The direction law need not be normalized on either
set or on any individual incidence fiber. All predetermined fixed
widths are allowed.
\end{corollary}

\begin{proof}
Choose a fixed enlarged width $C'>C$ large enough for the spatial
anchor rounding used in lemma~\ref{ang:slice-selection}. Let
$\epsilon,c$ be supplied by theorem~\ref{ang:linear-gain} for width
$C'\delta$, and choose $\epsilon_1\le\frac{\epsilon}{2}$ and
$0<c_1<\min(\frac{\epsilon}{2},\frac{c}{3})$. In a fixed bounded spatial
$\delta$-grid consider all candidate lists $J$ of at most
$\lfloor\delta^{-q-c_1}\rfloor$ anchors. Keep the actual direction
$w$ in every tube. Write $T_{C'}(w,a)$ for the tube of width
$C'\delta$ about the line through the anchor $a$ parallel to
$w$, restricted to the fixed bounded region. Define
$$
 B_\delta^+
 =\left\{w:\max_J
       \mu\left(\bigcup_{a\in J}T_{C'}(w,a)\right)
                    \ge\delta^{c_1}\right\}.
$$
The candidate masses are measurable, so this is a Borel set, and
rounding a point on each axis in an original bad cover gives a
candidate at the enlarged width with at least the same captured mass,
so that $B_\delta\subset B_\delta^+$, equality not being needed.
For each $w\in B_\delta^+$ choose the first candidate attaining the
threshold. The resulting measurable incidence relation has mass at least
$\delta^{c_1}\nu(B_\delta^+)$ and uses at most
$\delta^{-q-c_1}\le\delta^{-q-\epsilon}$ tubes per fiber, so
theorem~\ref{ang:linear-gain} bounds its incidence by $\delta^c$, and
hence $\nu(B_\delta^+)\le\delta^{c-c_1}\le\delta^{c_1}$.
\end{proof}

\section{Preparation of the pair measures}\label{prep:section}

We now return from abstract directions to pairs of points in the original
set. What we want is that, after a small exceptional set at scale $2^{-N}$
has been excluded, each endpoint should see a measure on directions with
cap exponent $S=\frac{d}{2}$, and that the construction should leave a
fixed positive amount of pair mass in the limit. These conclusions are
collected in proposition~\ref{prep:starting-pairs}, and the arguments
leading to it explain why integer and half-integer values of $S$ require
different preparation.

Proposition~\ref{prep:starting-pairs} provides decreasing pair filters of
positive limiting mass. The mass remaining to be removed after scale
$2^{-N}$ is at most a negative power of $2^N$. At that scale the filter
is a union of rectangles from two finite spatial partitions, and a further
power-small deletion gives cap bounds of exponent $S$ with an arbitrarily
small power loss. The auxiliary projections estimate
directions only, all the distances remaining those of the original pairs.
We prove first the improvement needed in odd dimensions and then construct
the measures to which it applies.

\subsection{From integer to half-integer directional bounds}

For separated supports, write $\pi_x(z)=\frac{z-x}{|z-x|}$. If $G$ is a
set of pairs with pin probability $\nu$ and source probability $\rho$, its directional
measures are
$$
 (\pi_x)_\#(\rho|_{G_x}) \text{ and } (\pi_z)_\#(\nu|_{G^z}).
$$
These are subprobabilities, and they are never normalized separately.
A directional plate means $\{w\in S^k:\operatorname{dist}(w,V)\le u\}$,
where $V$ is a linear hyperplane in ${\mathbb R}^{k+1}$.
In the cap and plate bounds of this section, $\lambda$ denotes
either of the two directional measures, and $V$ ranges over all
linear hyperplanes in the indicated direction space.

\begin{definition}\label{prep:availability}
Let $\rho,\nu$ be fixed separated probabilities in ${\mathbb R}^{k+1}$, and
let $G_N$ be Borel pair sets. An exponent $q$ is available for this
family if, for every $e>0$, there are $c_e>0$, $C_e<\infty$, and Borel
$B_N\subset G_N$, for all sufficiently large $N$, such that
$$
 (\nu\otimes\rho)(B_N)\le C_e2^{-c_eN},
$$
and both directional measures of $G_N\setminus B_N$ give every spherical
cap of radius $u\in[2^{-N},1]$ mass at most $2^{eN}u^q$.
The estimates hold at every pin in the retained supports, after removal of
null fibers if necessary.
\end{definition}

We know that the mass curve lies approximately above a line of slope $s$,
but its slope on a particular short interval can be smaller than $s$, and
we need one interval on which the conditional mass has slope at
least $t>q$, with slope-$q$ estimates before and after it. Subtracting the
line of slope $t$ makes the first choice simple, because at a minimum of
the difference the graph cannot subsequently fall below that line, and a
second minimum, after subtracting the slope $q$, supplies the last
interval. We keep track of the distances between the chosen depths, since
an interval which disappeared as the scale changed would be of little use
to us.

In this local argument $A$ and $B$ are real numbers in $(0,1)$ marking
normalized depths, and they are unrelated to the excess values $A_n$. A
physical cell at depth approximately $AN$ has diameter comparable to
$R^{-A}$, where $N$ is the final depth and $R=2^N$, and on separated
supports the second-order Taylor error of the radial map on that cell is
bounded by a constant times $R^{-2A}$. The inequality $B<2A$ makes this
error smaller than the tube width $R^{-B}$, while a positive lower bound
for $2A-B$ makes the comparison uniform.

\begin{lemma}[A quantitative contact interval]\label{prep:contact}
Let $L$ be a positive Lipschitz bound and suppose
$0\le q<t<s\le L$. Let $\varepsilon$ be a nonnegative error satisfying
$\varepsilon\le\frac{s-t}{4}$. Suppose $m:[0,1]\to{\mathbb R}$ is nondecreasing
and $L$-Lipschitz, $m(0)=0$, and
$$
m(u)\ge su-\varepsilon \text{ for } 0\le u\le1.
$$
Define the positive numbers
$$
a_* =\frac{s-t}{4L}
\text{ and }
h_* =\frac{a_*(t-q)}{2L}.
$$
There exist $0<A<B<1$ such that
\begin{equation}\label{prep:contact-first}
m(w)\ge qw-\varepsilon \text{ for } 0\le w\le A,
\end{equation}
\begin{equation}\label{prep:contact-middle}
m(w)-m(A)\ge t(w-A) \text{ for } A\le w\le B,
\end{equation}
and
\begin{equation}\label{prep:contact-last}
m(w)-m(B)\ge q(w-B) \text{ for } B\le w\le1.
\end{equation}
The choices have the uniform margins
\begin{equation}\label{prep:contact-margins}
\begin{aligned}
 A&\ge a_*, & h_*\le B-A&\le\frac{a_*}{2},\\
 2A-B&\ge\frac{a_*}{2}, & 1-B&\ge\frac{3a_*}{2}.
\end{aligned}
\end{equation}
In particular $B<2A$. These constants depend only on $L,s,t,q$,
not on the particular mass curve or on its error $\varepsilon$.
\end{lemma}
\begin{proof}
Put
$$
H(u)=m(u)-tu \text{ and } G(u)=m(u)-qu.
$$
Both functions are $L$-Lipschitz, since for $u<v$ the increment
$m(v)-m(u)$ lies between $0$ and $L(v-u)$, and both slopes $t,q$
lie in $[0,L]$.

Choose $A$ where $H$ attains its minimum on $[a_*,1]$.
Then $H(w)\ge H(A)$ for every $w\ge A$, or equivalently
$$
m(w)-m(A)\ge t(w-A) \text{ for } A\le w\le1.
$$
This is stronger than \eqref{prep:contact-middle}, and there is also a
definite distance left before the endpoint $1$. Since $H(0)=0$, we have
$$
H(A)\le H(a_*)\le La_* =\frac{s-t}{4}
\text{ and }
H(1)\ge s-t-\varepsilon\ge\frac{3(s-t)}{4}.
$$
Consequently
$$
L(1-A)\ge H(1)-H(A)\ge\frac{s-t}{2}=2La_*,
$$
so $A\ge a_*$ and $1-A\ge2a_*$.

Now choose $B$ where $G$ attains its minimum on $[A+h_*,1]$.
This interval is nonempty, since $h_*\le\frac{a_*}{2}$.
The minimum property gives $G(w)\ge G(B)$ for $w\ge B$,
which is the estimate~\eqref{prep:contact-last}.
The first minimum also gives
$$
G(w)-G(A)\ge(t-q)(w-A) \text{ for } w\ge A.
$$
Comparing the second minimum with the left endpoint of its interval, we find that
$$
\begin{aligned}
(t-q)(B-A)
&\le G(B)-G(A)\\
&\le G(A+h_*)-G(A)\\
&\le Lh_* =\frac{a_*(t-q)}{2}.
\end{aligned}
$$
Thus $h_*\le B-A\le\frac{a_*}{2}$, and it follows that
$$
2A-B=A-(B-A)\ge\frac{a_*}{2}
\text{ and }
1-B=(1-A)-(B-A)\ge\frac{3a_*}{2}.
$$
Finally, $m(w)\ge sw-\varepsilon\ge qw-\varepsilon$ proves \eqref{prep:contact-first}.
\end{proof}

\begin{figure}[htbp]
\centering
\begin{tikzpicture}[>=Stealth,font=\small,
box/.style={draw,align=center,text width=5.35cm,minimum height=1.35cm,
inner sep=6pt}]
\node[box] (first) at (-3.15,0)
{Minimize $H(u)=m(u)-tu$\\on $[a_*,1]$ to choose $A$};
\node[box] (firstout) at (3.15,0)
{For $w\ge A$\\$m(w)-m(A)\ge t(w-A)$\\
$A\ge a_*$ and $1-A\ge2a_*$};
\node[box] (second) at (-3.15,-2.3)
{Minimize $G(u)=m(u)-qu$\\on $[A+h_*,1]$ to choose $B$};
\node[box] (secondout) at (3.15,-2.3)
{For $w\ge B$\\$m(w)-m(B)\ge q(w-B)$\\
$h_*\le B-A\le\frac{a_*}{2}$};
\draw[->] (first)--(firstout);
\draw[->] (first)--node[left,align=right]
{the interval for\\the second minimum exists}(second);
\draw[->] (second)--(secondout);
\draw[->] (firstout)--node[right,align=left]
{the slope gap $t-q$\\bounds $B-A$}(secondout);
\end{tikzpicture}
\caption{The two choices and what each supplies. The first minimum
provides the slope-$t$ estimate to its right, and the second provides
the slope-$q$ estimate after $B$. Comparing the second minimum
with $G(A+h_*)$ bounds the length of the middle interval from above,
while its domain supplies the lower bound.}
\label{fig:two-contact-minima}
\end{figure}
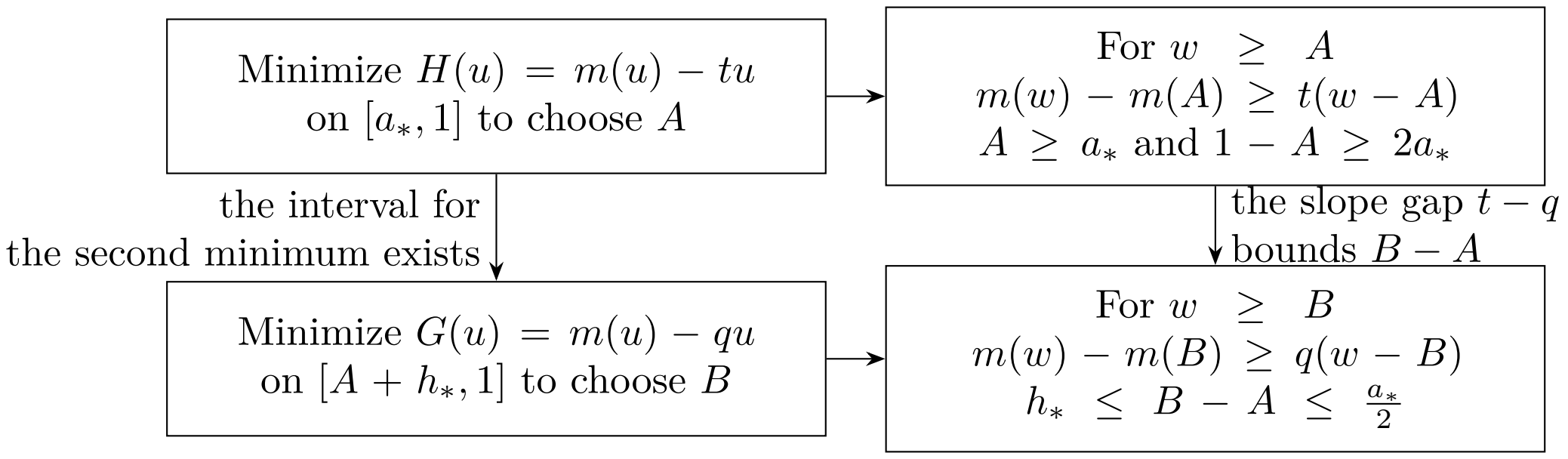

This construction is an adaptation of the lower-barrier and
finite-witness viewpoint in Iosevich and Yavicoli
\cite{IosevichYavicoli2026}. Their lemma on uniform finite witnesses
for the partition gain preserves a strict margin when the functions and the
scales vary, whereas here the two minimum choices give the required margins
explicitly, and we shall use them before taking any limit of mass curves.
The angular estimate itself is still corollary~\ref{ang:bad-directions}.
The argument below has precedents in the thin-tube bootstrap of Orponen,
Shmerkin and Wang~\cite[section~1.3 and lemma~2.8]{OrponenShmerkinWang2024}
and in Shmerkin and Wang's use of scale mass curves to improve radial
bounds~\cite[section~5.5]{ShmerkinWang2025}, the contact interval being
what allows us to combine three local estimates while keeping a power bound
for the deleted pair mass.

It is at this point that the angular gain improves the radial information.
We use the old exponent $q$ on the first and the last portions of the depth
interval and the new incidence estimate on the middle portion, and
multiplying the three conditional estimates gives a gain in the exponent
for the whole interval. All the comparisons involve the full conditional
source measure, since we restrict the set of witnesses rather than
renormalizing the law separately for each pin.

\begin{figure}[htbp]
\centering
\begin{tikzpicture}[>=Stealth,font=\small]
\draw[->] (0,0)--(12.6,0) node[right] {depth $u$};
\draw[line width=3pt,gray!55] (0,0)--(4.7,0);
\draw[line width=3pt] (4.7,0)--(7.3,0);
\draw[line width=3pt,gray!55] (7.3,0)--(12,0);
\foreach \x/\v in {0/0,4.7/A,7.3/B,12/1}{
\draw (\x,-.15)--(\x,.15);\node[below] at (\x,-.2) {$\v$};}
\node[align=center] at (2.3,1) {old radial bound\\exponent $q$};
\node[align=center] at (6,1) {angular gain\\conditional slope $t>q$};
\node[align=center] at (9.8,1) {old radial bound\\exponent $q$};
\node[align=center] at (6,-1.4)
{cell radius $R^{-A}$ and final tube width $R^{-B}$\\
$B-A\ge h_*$ and $2A-B\ge\frac{a_*}{2}$};
\end{tikzpicture}
\caption{The three intervals in the one-scale improvement, drawn
schematically. The first and last estimates use the available exponent
$q$, while the middle interval supplies a strict gain, and measuring the
full mass in each child allows the three conditional bounds to multiply.
The condition below the axis controls the quadratic Taylor error, so that
the middle radial problem may be replaced by a linear tube problem.}
\label{fig:contact-interval}
\end{figure}
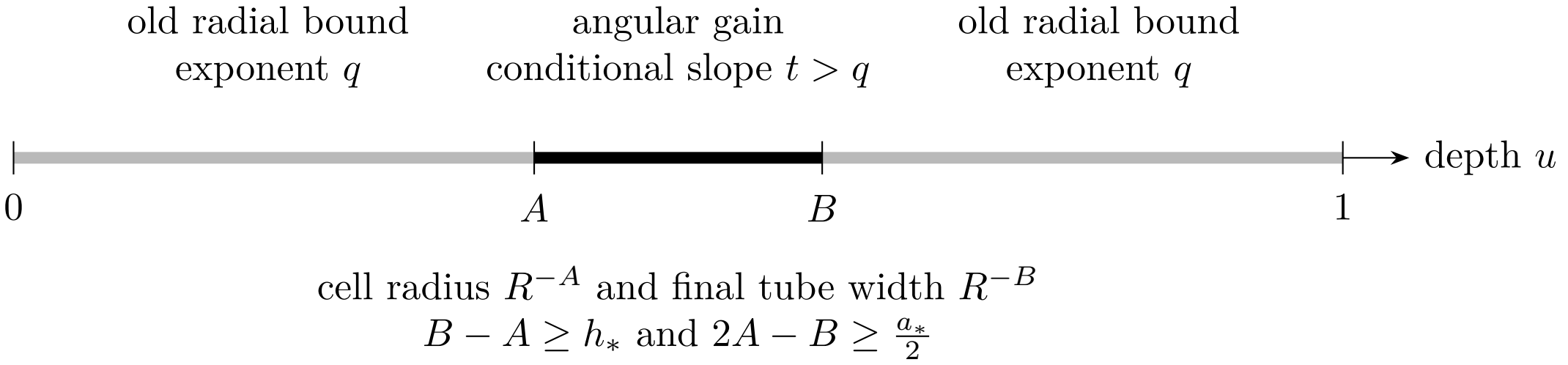

\begin{lemma}[One-scale improvement]\label{prep:one-scale}
Fix $k\ge2$, $k-\frac{1}{2}<s<k$, $k-1\le q<k-\frac{1}{2}$,
and $\zeta>0$. Fix bounded support, separation, and $s$-Frostman
constants, as well as a plate-bound constant $C$. There are $\varepsilon_0,c>0$
such that the following holds for
all sufficiently small dyadic $\delta$. Let $\rho,\nu$ be probabilities
with the fixed bounds, and let $G$ be a Borel pair set whose two
directional measures satisfy, for $\delta\le u\le1$,
$$
 \lambda(B(w,u))\le\delta^{-\varepsilon_0}u^q
 \text{ and }
  \lambda\{w:\operatorname{dist}(w,V)\le u\}\le C u^\zeta.
$$
Then one can delete pair mass at most $C\delta^c$ so that both remaining
directional measures give every cap of radius $\delta$ mass at most
$\delta^{q+c}$.
The constants are uniform over the measures and pair sets with these bounds.
\end{lemma}

\begin{proof}
The outer intervals use a first-moment estimate followed by Markov's
inequality, while the middle interval uses the angular incidence theorem,
and the quantitative contact lemma supplies the same lower bound for the
length of that middle interval in every regular class. We can therefore
choose all the losses before the scale, and work directly with the finite
mass curves.

It is enough to prove the conclusion in one orientation and intersect the
two resulting pair sets. Write $R=\delta^{-1}=2^M$, where $M$ is the
final dyadic depth, fix $t\in(k-\frac{1}{2},s)$, take $L=k+1$, the
Lipschitz bound for mass curves of measures in ${\mathbb R}^{k+1}$, and
let $a_*,h_*$ be the constants in lemma~\ref{prep:contact} for these
parameters.

Choose the fixed tube widths large enough for the geometric comparisons
below, and let $c_1,\epsilon_1>0$ be supplied by
corollary~\ref{ang:bad-directions} for those widths, so that $c_1$ is
the angular gain and $\epsilon_1$ is the permitted loss in its ball
bounds. Choose the positive number
$$
\theta=\frac{c_1h_*}{40}
$$
for the two outer Markov thresholds, and then a positive comparison
tolerance $\varepsilon$ so small that
\begin{equation}\label{prep:finite-loss-choice}
\varepsilon<\min\left\{\frac{s-t}{4},
                 \frac{\epsilon_1h_*}{4},
                 \frac{\theta}{8}\right\}.
\end{equation}
Regularize only the source law by lemma~\ref{reg:regularization}, with a
fixed parameter $\lambda>0$ small enough that its cell and conditional cell
comparison exponents are at most $\varepsilon$, choosing it smaller if
necessary so that the lower bound for each mass curve is
$m(u)\ge su-\varepsilon$ for all sufficiently large $R$. This is possible
because, if $K$ denotes the fixed original Frostman constant, its
contribution $\log_R K$ tends to zero as $R$ increases. Finally take the
input loss $\varepsilon_0\le\varepsilon$. None of these choices depends on
the regular class or on $M$.

During the estimate for one class, $\rho$ denotes its normalized
probability and $m$ its mass curve. Restrict $G$ to pairs whose
source belongs to that class, and discard the null pairs outside
$\operatorname{supp}\nu\times\operatorname{supp}\rho$. Deleted pair mass in this calculation
is measured by $\nu\otimes\rho$. Only the reverse laws
$(\pi_z)_\#(\nu|_{G^z})$ are used as direction measures below.
They remain unnormalized restrictions of the original pin law,
so normalizing a source class does not alter their cap or plate constants.

Choose $A,B$ from lemma~\ref{prep:contact} for this finite mass
curve, and use the four depths $0$, $a=\lfloor AM\rfloor$,
$b=\lfloor BM\rfloor$ and $M$. The choices $A,B$ may vary from class to class, whereas their margins in
\eqref{prep:contact-margins} do not. Rounding changes a normalized depth by
at most $\frac{1}{M}$, and therefore the Lipschitz bound for $m$ changes a cell
comparison by at most a fixed factor, and such factors, together with the
at most $M+1$ annuli in the calculation below, are written as $R^{o(1)}$;
these $o(1)$ terms tend to zero uniformly over the classes once all the
preceding parameters have been fixed.

We establish local radial
fiber estimates on the three intervals. Here $u,v$ are normalized
depths, so the initial integer depth is $uM$, rounded as above,
and the final spatial scale is $R^{-v}$. For an occupied source cell
$Q$ at that initial depth, write
$\rho_Q=\frac{\rho|_Q}{\rho(Q)}$.
Let $C_{\mathrm{rad}}$ be a sufficiently large fixed width constant,
included in the geometric choices made before applying
corollary~\ref{ang:bad-directions}. Define
$$
 \mathcal R_{Q,v}(x,z)
 =\rho_Q\{z'\in Q:
       |\pi_x(z')-\pi_x(z)|\le C_{\mathrm{rad}}R^{-v}\}
  \text{ for } x\in\operatorname{supp}\nu \text{ and } z\in Q\cap\operatorname{supp}\rho.
$$
The comparison source $z'$ always ranges over the full conditional
class law in $Q$, and we shall bound this quantity only at retained witness
pairs $(x,z)$, no individual pair fiber being normalized. The function
$\mathcal R_{Q,v}$ is Borel by parameterized integration.

First consider $[u,v]=[0,A]$ or $[B,1]$, and a source cell $Q$ at
depth $uM$.
For every $u\le w\le v$, the mass curve comparisons and
the equations~\eqref{prep:contact-first} and~\eqref{prep:contact-last} give
\begin{equation}\label{prep:local-ball}
\rho_Q(B(z,CR^{-w}))\le R^{-q(w-u)+2\varepsilon+o(1)}.
\end{equation}
For the first interval one $\varepsilon$ comes from the lower bound for $m$
and one from the cell comparison, whereas for the last interval the contact
inequality has no error, so that the displayed bound leaves more room than
we need. A ball meets a bounded number of cubes at the corresponding dyadic
depth, which accounts for a fixed factor only.

At a fixed witness $z\in Q$, integrate pins only on $G^z$.
Bounded separation implies
$$
 |\pi_x(z)-\pi_x(z')|\le C_{\mathrm{rad}}R^{-v}
  \Longrightarrow \pi_z(x)\text{ lies in two caps of radius }
 \frac{C'R^{-v}}{|z-z'|}.
$$
For $\frac{1}{2}R^{-w}<|z-z'|\le R^{-w}$, the product of the cap bound and
\eqref{prep:local-ball} is at most $R^{-q(v-u)+3\varepsilon+o(1)}$, the
third $\varepsilon$ bounding the input cap loss. Summing over the
intermediate dyadic distance annuli, and noting that the innermost ball
$|z-z'|\lesssim R^{-v}$ obeys the same bound by \eqref{prep:local-ball},
we obtain
\begin{equation}\label{prep:local-mean}
\int_{G^z}\mathcal R_{Q,v}(x,z)d\nu(x)
 \le R^{-q(v-u)+3\varepsilon+o(1)}.
\end{equation}
Fubini and Markov delete at most
$R^{-\theta+3\varepsilon+o(1)}$ pair mass and leave
$$
 \mathcal R_{Q,v}(x,z)\le R^{-q(v-u)+\theta}
$$
at every retained witness $(x,z)$ in this cell. More precisely,
the contribution of $Q$ to the deleted product mass is bounded
by this power times $\rho(Q)$. These weights sum to one at the
chosen depth. By \eqref{prep:finite-loss-choice}, the
total deletion in either outer interval is at most
$R^{-\frac{\theta}{2}}$ for sufficiently large $R$.

For the middle interval put $D=R^{-(B-A)}$. After rescaling an
$a$-cell $Q$ to bounded size, its conditional probability has a
$t$-ball bound down to scale $D$, with loss $R^{\varepsilon+o(1)}$, by \eqref{prep:contact-middle}. To apply the angular theorem,
we must express this loss in the scale $D$. For large $R$,
$$
R^{\varepsilon+o(1)}\le R^{2\varepsilon}
=D^{-\frac{2\varepsilon}{B-A}}
\le D^{-\frac{2\varepsilon}{h_*}}
\le D^{-\frac{\epsilon_1}{2}}.
$$
This is where the lower bound for the interval length is used.
The unchanged reverse laws have cap loss at most $R^\varepsilon$, which
satisfies the same comparison. Their required radii are at least
$D\ge R^{-1}$, and hence are covered by the hypothesis. At $q=k-1$,
the fixed plate constant is at most $D^{-\epsilon_1}$ for all
sufficiently large $R$. All scale thresholds here are uniform,
since $D\le R^{-h_*}$.

Corollary~\ref{ang:bad-directions}, with the fixed widths chosen
above, supplies a common Borel set $\mathcal B_Q$ containing all
directions for which a union of
at most $D^{-q-c_1}$ parallel tubes of width $C'D$ carries at least
$D^{c_1}$ conditional source mass. For every $z\in Q$,
\begin{equation}\label{prep:witness-bad}
((\pi_z)_\#(\nu|_{G^z}))(\mathcal B_Q)\le D^{c_1}.
\end{equation}
The directional law in this statement is a subprobability and has not been
renormalized, and the set $\mathcal B_Q$ depends on $Q$ and on its
conditional source law, not on $z$, while the same fixed exponent $c_1$
already bounds the exceptional mass in that corollary, whose statement also
includes the fixed-width enlargements and the gridding. Robustness at a
direction means that no union of at most $D^{-q-c_1}$ tubes parallel to
that direction, of the permitted width in the rescaled cube, captures
conditional mass at least $D^{c_1}$.

Delete the pairs in \eqref{prep:witness-bad}, at cost
$O(D^{c_1})$ by Fubini. For a pin still represented in $Q$, at least
one retained witness has a good direction. Let $c_Q$ be the cell center
and $w_{x,Q}=\pi_x(c_Q)$. All witness directions differ, up to sign,
from $w_{x,Q}$ by $O(R^{-A})$. Since
\begin{equation}\label{prep:linearization-margin}
\frac{R^{-A}}{D}=\frac{R^{-2A}}{R^{-B}}
=R^{-(2A-B)}\le R^{-\frac{a_*}{2}},
\end{equation}
a tube cover at $w_{x,Q}$ with a slightly smaller fixed width would
also be a forbidden cover at that witness direction. Thus the
conditional source law is robust at the center direction, uniformly
over every represented pin. The directions may be regarded as
unoriented for this assertion.

For $z,z'\in Q$, Taylor expansion of $y\mapsto\pi_x(y)$ on the
separated supports gives
$$
 |\pi_x(z)-\pi_x(z')|\le C_{\mathrm{rad}}R^{-B}
  \Longrightarrow |P_{w_{x,Q}^{\perp}}(z-z')|
 \le C''(R^{-B}+R^{-2A})\le C'''R^{-B}.
$$
Grid the perpendicular projection at this scale. For each pin $x$, let
$$
 \mathcal H_{x,Q}
 =\left\{z\in Q\cap\operatorname{supp}\rho:\mathcal R_{Q,B}(x,z)>D^{q+\frac{c_1}{2}}\right\}.
$$
If $z\in\mathcal H_{x,Q}$, then one of a bounded number of neighboring
projected grid cells has full $\rho_Q$-mass at least a fixed multiple of
$D^{q+\frac{c_1}{2}}$, and there are at most $C D^{-q-\frac{c_1}{2}}$ such
heavy cells. Enlarging each of them by its bounded number of neighboring
grid cells, we see that $\mathcal H_{x,Q}$ lies in at most that many tubes,
of a fixed width permitted by the preceding robustness assertion, and for
small $D$ this is less than $D^{-q-c_1}$. For every represented pin,
robustness therefore gives $\rho_Q(\mathcal H_{x,Q})\le D^{c_1}$.

From the pair relation surviving the bad-direction deletion, delete the
pairs $(x,z)$ with $z\in\mathcal H_{x,Q}$. If a pin has no surviving
witness in $Q$ this deletion has zero fiber mass, and otherwise the
preceding bound applies, and hence integrating over the full pin law and then
summing with the cell weights $\rho(Q)$ costs $O(D^{c_1})$. This deletion
intersects the first surviving Borel relation with the threshold event for
$\mathcal R_{Q,B}$, and is therefore Borel and needs no measurable choice
of the witness used in the robustness argument. Every remaining witness
satisfies $\mathcal R_{Q,B}(x,z)\le D^{q+\frac{c_1}{2}}$.

Denote the three local factors, in interval order, by
$$
 f_0=R^{-qA+\theta},
 f_1=R^{-(q+\frac{c_1}{2})(B-A)}
 \text{ and }
 f_2=R^{-q(1-B)+\theta}.
$$
Let $G'$ be the common Borel pair relation left in the class-restricted
$G$ after the two outer Markov deletions and the two middle deletions. At each $(x,z)\in G'$ in the initial
cell $Q$ of the corresponding interval $[u,v]$, the local estimate is
$$
 \rho\{z'\in Q:
       |\pi_x(z')-\pi_x(z)|\le C_{\mathrm{rad}}R^{-v}\}
 =\rho(Q)\mathcal R_{Q,v}(x,z)\le f_j\rho(Q).
$$
The comparison integrates the full source law in $Q$, so that for a
fixed final cap the enlarged comparison around a retained witness in a
hit parent covers the full source mass in each hit child, and
lemma~\ref{reg:pruning-multiplication} applies at the four endpoint
depths, the fixed width $C_{\mathrm{rad}}$ including the enlargement
needed for the rounded depths. The final cap mass on $G'$ is at most
$$
 f_0 f_1 f_2
 =R^{-q-\frac{c_1}{2}(B-A)+2\theta}
 \le R^{-q-\frac{9c_1h_*}{20}},
$$
since $B-A\ge h_*$ and $\theta=\frac{c_1h_*}{40}$. The two outer
deletions cost at most $2R^{-\frac{\theta}{2}}$ and the middle ones cost
$O(D^{c_1})=O(R^{-c_1h_*})$. Both the gain in the cap exponent
and the power saving in deleted mass are uniform over the classes.

Restore now the original source weight of each class. These weights sum to
at most one, so that both the cap bounds and the deletion bounds sum
without a factor equal to the number of classes, while the classes omitted
by regularization have original mass $R^{-\lambda+o(1)}$, which for
sufficiently large $R$ is at most $R^{-\frac{\lambda}{2}}$. One may
therefore choose the final gain and deletion exponent $c>0$ smaller than
each of $\frac{c_1h_*}{8}$, $\frac{\theta}{4}$ and $\frac{\lambda}{4}$.
Finally, repeat the argument with the two original measures interchanged
and intersect the outputs; decreasing $c$ absorbs the fixed constants and
proves the lemma.
\end{proof}

The preceding gain is local both in exponent and in scale. To reach
$k-\frac{1}{2}$, consider the supremum of the exponents already available
in the precise sense of definition~\ref{prep:availability}; if that
supremum were smaller than the desired exponent, the one-scale
improvement, followed by a sum over scales, would increase it. It is
essential to keep the deletion exponents in this argument, since an
improvement of the exponent with no estimate for the discarded mass would
not give the later approximation in total variation.

\begin{proposition}[Thin-tube upgrade]\label{prep:thin-upgrade}
Let $k\ge2$ and $k-\frac{1}{2}<s<k$. Let $\rho,\nu$ be bounded,
separated $s$-Frostman probabilities in ${\mathbb R}^{k+1}$. Suppose that the
two directional measures of $G_N$ have plate bounds $Cu^\zeta$,
uniformly for $2^{-N}\le u\le1$, with fixed $\zeta>0$.
If $k-1$ is available, then $k-\frac{1}{2}$ is available.
\end{proposition}
\begin{proof}
Availability is downward closed, and it is also closed under increasing
limits inside $[k-1,k-\frac{1}{2}]$. For a desired loss $e$, choose an
available $q$ within $\frac{e}{2}$ of the supremum, apply its estimate with
loss $\frac{e}{2}$, and use $u\ge2^{-N}$.
Writing $q_*$ for that supremum, the numerical comparison is
$$
2^{\frac{eN}{2}}u^q
=2^{\frac{eN}{2}}u^{q_*}u^{-(q_*-q)}
\le 2^{\frac{eN}{2}}2^{N(q_*-q)}u^{q_*}
\le 2^{eN}u^{q_*}.
$$
The deletion estimate is the one belonging to this fixed available $q$, and
its exponent need not be uniform as $q$ approaches $q_*$, since $q$ is
chosen once, after $e$ and before $N$.

Suppose that an available exponent $q<k-\frac{1}{2}$ is given. Let
$c,\varepsilon_0$ be the constants in lemma~\ref{prep:one-scale}, and choose
$q'>q$ with $q'\le q+\frac{c}{2}$ and $q'\le k-\frac{1}{2}$. For any desired
loss $e>0$, choose $e_1>0$ so small that $2e_1q'<e$.
First use availability of $q$ with input loss
$\varepsilon<\frac{e_1\varepsilon_0}{2}$. For every integer $M\in[e_1N,N]$, that
input loss is at most $2^{\varepsilon_0M}$, so
lemma~\ref{prep:one-scale} applies at $\delta=2^{-M}$ to the same
remaining pair set. Intersect the outputs over these $M$.
Their combined deletion is at most
$$
 C\sum_{M\ge e_1N}2^{-cM}=O(2^{-ce_1N}),
$$
in addition to the original power-small deletion. At these radii the output
is stronger than the requested $2^{eN}u^{q'}$ bound, while at larger radii
that bound is trivial from mass at most one, for all large $N$, by the
choice of $e_1$, dyadic enlargement being absorbed by the remaining loss.
Hence $q'$ is available, and applying this to the supremum of the available
exponents rules out a supremum below $k-\frac{1}{2}$, while closure gives
the endpoint itself.
\end{proof}

\subsection{Starting measures and radial densities}

We record explicitly the classical exceptional case, which allows the
remaining construction to assume that no relevant affine plane is charged.
For this case we use the Fourier estimate underlying the classical
threshold $\frac{n+1}{2}$, which gives a power saving on each frequency
shell for a fixed nonzero distance measure.

\begin{lemma}[The classical distance threshold]\label{prep:classical-plane}
Let $\rho$ be a compactly supported probability in ${\mathbb R}^n$, $n\ge2$,
with an $s$-ball bound for some $s>\frac{n+1}{2}$. Then
$\Delta(\operatorname{supp}\rho)$ supports a nonzero finite positive measure $\tau$
satisfying, for every $R\ge 2$,
\begin{equation}\label{prep:classical-shell}
 \int_{\frac{R}{2}\le |r|\le 2R}|\widehat\tau(r)|^2dr
 \lesssim R^{\frac{n+1}{2}-s}+R^{-1}.
\end{equation}
In particular $\tau$ has an $L^2$ density and
$\Delta(\operatorname{supp}\rho)$ has positive Lebesgue measure.
\end{lemma}
\begin{proof}
Choose positive normalized restrictions $\rho_1,\rho_2$ to sufficiently
small separated compact patches. Put $a=\frac{n-1}{2}$ and
$$
 \tau=D_\#[D^{-a}(\rho_1\otimes\rho_2)]
 \text{ with }
  D(x,y)=|x-y|,
$$
which is a nonzero finite positive measure on the distance line.
Write $\sigma_{S^{n-1}}$ for surface measure on $S^{n-1}$.
The elementary sphere-kernel estimate
$$
 |\widehat{\sigma_{S^{n-1}}}(\xi)|\lesssim(1+|\xi|)^{-a}
$$
follows by a partition into the two stationary patches and their
complement, using stationary phase on the former and integration by
parts on the latter. Consequently
\begin{equation}\label{prep:classical-spherical}
\int_{S^{n-1}}|\widehat{\rho_1}(rw)|^2d\sigma(w)
 \lesssim r^{-a}I_a(\rho_1) \text{ for } r\ge1.
\end{equation}
The energy is finite since $a<s$. The ball bound and a smooth Fourier
ball majorant give
\begin{equation}\label{prep:classical-annulus}
\int_{r\asymp R}\int_{S^{n-1}}
       |\widehat{\rho_2}(rw)|^2d\sigma(w)\frac{dr}{R}
 \lesssim R^{-n}\int_{|\xi|\lesssim R}|\widehat{\rho_2}(\xi)|^2d\xi
 \lesssim R^{-s}.
\end{equation}
For the last inequality, the inverse transform of the majorant is
$O_L(R^n(1+R|x-y|)^{-L})$; integrating by dyadic annuli against
the $s$-ball bound proves the claim.

Choose a smooth angular cutoff $\psi$ equal to one on all the cross
directions $\frac{x-y}{|x-y|}$ and zero near their antipodes. Let $c_n$ be
the nonzero leading constant in the spherical stationary-phase
expansion for this choice of Fourier sign. The leading spherical stationary-phase
term, uniformly on the separated patches, gives
$$
 \widehat\tau(r)=c_n^{-1}r^a
 \int_{S^{n-1}}\psi(w)\widehat{\rho_1}(rw)
                 \overline{\widehat{\rho_2}(rw)}d\sigma(w)
 +O(r^{-1}),
$$
the error before multiplication by $r^a$ being $O(r^{-a-1})$
because distances stay in a compact subset of $(0,\infty)$.
Cauchy--Schwarz, the equations~\eqref{prep:classical-spherical}
and~\eqref{prep:classical-annulus}, and $2a+1=n$ imply
$$
 \int_{r\asymp R}|\widehat\tau(r)|^2dr
 \lesssim R^{n-a-s}+R^{-1}
 =R^{\frac{n+1}{2}-s}+R^{-1}.
$$
The dyadic sum converges, the negative frequencies are conjugates, and the
low frequencies are bounded, so that Plancherel gives an $L^2$ density for
$\tau$, which proves the conclusion. This is the classical threshold
argument of Falconer~\cite{Falconer1985}, in the Fourier formulation
developed by Mattila~\cite{Mattila1987}.
\end{proof}

We shall need a finite collection of projections, with pointwise ball
bounds on their images, simultaneously for both starting measures, and
almost-everywhere finiteness of the projected energy is not yet such a
statement. We therefore choose projections with finite energy first, and
then restrict the original measures to preimages of compact sets on which
the projected potentials are bounded, a finite intersection still keeping
positive mass. By a generic orthogonal projection we mean here one outside
a set of Haar measure zero in the appropriate Grassmannian.

\begin{lemma}[Simultaneous projected Frostman bounds]\label{prep:projections}
Let $E\subset{\mathbb R}^d$ be compact and $\dim_HE>S=\frac{d}{2}$.
There is a probability on $E$ with a ball exponent strictly greater
than $S$. Finitely many successive generic orthogonal projections to
dimensions strictly greater than $S$ preserve finite energies with
exponents strictly greater than $S$. After a positive compact restriction
of the original probability, all these projected laws, as well as the
original law, have a common ball exponent $S+\chi$, with $\chi>0$.

If $d=2k-1$ and the original measure gives zero mass to every affine
$k$-plane, the projection to ${\mathbb R}^{k+1}$ may additionally be chosen
to give zero mass to every affine hyperplane.
\end{lemma}
\begin{proof}
Choose $S<s_0<\dim_HE$ and use Frostman's lemma to obtain an $s_0$-ball
probability, and fix $u>S$ below $s_0$ and below all the target dimensions.
Lemma~\ref{conv:energy} yields finite $u$-energy for almost every successive
projection, and in order to make all the projected ball bounds hold on the
same restriction we use the potential proof of that lemma: for each
projected law $\eta$, retain a set on which $U_u^\eta$ is bounded and pull
this set back to the original space. Each discarded preimage can be made of
arbitrarily small mass; thus the intersection of these finitely many
restrictions has positive mass, and an inner compact restriction has the
desired ball bound in every space. Decreasing the exponent slightly if
necessary produces one common $S+\chi$.

For the last assertion, almost every $(k+2)$-tuple of original points is
affinely independent, as one sees by sampling successively and using zero
mass on every span of dimension at most $k$. For each fixed such tuple,
almost every projection to ${\mathbb R}^{k+1}$ preserves its independence,
and Fubini makes this true for almost every tuple for almost every
projection, simultaneously with the preservation of energy. If a projected
hyperplane had mass $a>0$, the dependent tuples would have mass at least
$a^{k+2}$, which is a contradiction. Subsequent restrictions preserve this
property.
\end{proof}

The following lemma adapts the finite-energy proof of
Orponen~\cite[section~3 and the proof of theorem~1.13]{Orponen2019}.
The use of this radial estimate to control exceptional directions already
appears in \cite[propositions~3.11--3.12]{KeletiShmerkin2019}, and
Guth, Iosevich, Ou and Wang use it to bound the contribution of bad tubes
to pinned distance measures~\cite[section~3]{GuthIosevichOuWang2020}, the
common step being to integrate over the pins before removing the
exceptional pairs. To construct the integer cap bound we estimate the
$p$th power of a radial density jointly over pins and directions, and its
tail will bound the product mass removed from the set of pairs.

\begin{lemma}[Radial $L^p$ estimate]\label{prep:radial-lp}
Let $h\ge1$ and $0<\epsilon<\min(h,1)$. Let $\alpha,\beta$ be
bounded, separated probabilities in ${\mathbb R}^{h+1}$ with finite
$I_{h+\epsilon}$ energies. For some $p>1$ there exists a nonnegative
jointly Borel function $F(x,w)$ which is a density of $(\pi_x)_\#\alpha$
for $\beta$-almost every $x$ and satisfies
$$
 \int\int_{S^h}F(x,w)^pd\sigma(w)d\beta(x)<\infty.
$$
The same assertion holds in the reverse orientation.
\end{lemma}
\begin{proof}
Put
$$
 u=h-\epsilon, p=1+\frac{\epsilon}{4h}
 \text{ and }
  t=h-\frac{\epsilon}{2}.
$$
The two strict margins used below can be checked directly from these
choices, and they are
$$
t-up=\frac{\epsilon}{4}+\frac{\epsilon^2}{4h}>0
\text{ and }
h(2-p)-t=\frac{\epsilon}{4}>0.
$$
The first of these permits the energy comparison in H\"older's inequality,
while the second makes the weighted Fourier integral integrable after the
projection average, and both are fixed positive margins. Then $1<p<2$ and
$up<t<h(2-p)$. Suppose initially that the source is smooth, meaning that it
has a smooth density, which in the Lebesgue integrals below is also denoted
by $\alpha$. For $w\in S^h$, let $P_w$ be the projection to $w^\perp$, let
$f_w$ be the density of $(P_w)_\#\alpha$, and put
$\eta_w=(P_w)_\#\beta$. The radial density is
$$
 F(x,w)=\int_0^\infty\alpha(x+rw)r^hdr
 \lesssim f_w(P_wx),
$$
because the supports lie in a fixed bounded region.

If $g\ge0$ and $\|g\|_{L^{\frac{p}{p-1}}(\eta_w)}=1$, H\"older on the
product space gives
$$
 I_u(g\eta_w)\le I_{up}(\eta_w)^{\frac{1}{p}}
 \lesssim I_t(\eta_w)^{\frac{1}{p}}.
$$
In the first inequality, apply H\"older to the product measure
$\eta_w\otimes\eta_w$, with exponents $p$ and
$p'=\frac{p}{p-1}$. The kernel contributes
$I_{up}(\eta_w)^{\frac{1}{p}}$, whereas the other factor is
$\|g\|_{L^{p'}(\eta_w)}^2=1$.
The second inequality uses bounded support and $up<t$.
The energy identity
$I_u(\eta)=c_{h,u}\int|\widehat\eta(\xi)|^2|\xi|^{u-h}d\xi$
follows by writing the Riesz kernel as an integral of Gaussian
kernels; see \cite[lemma~12.12]{Mattila1995} for this identity with
our Fourier convention. Fourier duality and Cauchy--Schwarz with reciprocal
weights therefore imply
\begin{equation}\label{prep:radial-duality}
\|f_w\|_{L^p(\eta_w)}
 \lesssim I_t(\eta_w)^{\frac{1}{2p}}
 \left(\int_{w^\perp}|\widehat\alpha(\xi)|^2
                       |\xi|^{h-u}d\xi\right)^{\frac{1}{2}}.
\end{equation}
The pairing can first be justified with bounded smooth truncations,
then passed to the limit by the same energy bound.

Write $A(w)$ for the Fourier integral in \eqref{prep:radial-duality} and $B(w)=I_t(\eta_w)$.
Polar coordinates and rotational invariance give
$$
 \int_{S^h}A(w)d\sigma(w)
 =C\int_{{\mathbb R}^{h+1}}|\widehat\alpha(\xi)|^2
                      |\xi|^{h-u-1}d\xi
 =C'I_{h+\epsilon}(\alpha).
$$
Set $r_0=\frac{1}{2-p}$. The one-vector estimate
$\sigma\{w:|P_wz|\le a|z|\}\lesssim a^h$, together with
$tr_0<h$, shows that
$$
 \||P_wz|^{-t}\|_{L^{r_0}(d\sigma)}\lesssim|z|^{-t}.
$$
Minkowski yields $\|B\|_{L^{r_0}}\lesssim I_t(\beta)$.
Raise \eqref{prep:radial-duality} to the $p$th power and
apply H\"older with exponents $\frac{2}{p}$ and $\frac{2}{2-p}$ to obtain
\begin{equation}\label{prep:radial-integrated}
\int\int F(x,w)^pd\sigma(w)d\beta(x)
 \lesssim I_{h+\epsilon}(\alpha)^{\frac{p}{2}}I_t(\beta)^{\frac{1}{2}}.
\end{equation}
For a general source, let $\alpha_\varepsilon$ be its convolution with
a nonnegative smooth approximate identity of scale $\varepsilon$ and
sufficiently small support. Separation persists, and the
Fourier energy formula gives a uniform energy bound. The joint
measures
$$
 \Lambda_\varepsilon
 =[(x,y)\mapsto(x,\pi_x(y))]_\#
       (\beta\otimes\alpha_\varepsilon)
$$
have uniformly bounded densities in $L^p(\beta\otimes\sigma)$ by \eqref{prep:radial-integrated}. Separation gives weak
convergence to the analogous joint measure $\Lambda$.
Recovering the measure on directions at a fixed pin from this
joint measure is called disintegration.
Duality against continuous functions, followed by $L^{p'}$ density
\cite[theorem~6.15 and proposition~7.9]{Folland1999},
shows that $\Lambda$ has an $L^p$ density $F$, and disintegration, checked
on a countable determining family of continuous functions on $S^h$,
identifies $F(x,\cdot)$ as the radial density for $\beta$-almost every $x$,
a Borel version supplying joint measurability.

Interchanging $\alpha$ and $\beta$ gives the density in the opposite
orientation.
\end{proof}

An $L^p$ density may be large on a small set of directions, so we truncate
it, which gives a pointwise cap estimate and controls the discarded mass by
the $L^p$ norm. This is carried out in suitable projections, and the
directional bound is then lifted to the original space, the geometry of the
chosen projections guaranteeing that a small cap upstairs maps into a
small cap in at least one of the projected configurations.

\begin{lemma}[Integer caps and lifting]\label{prep:integer-caps}
Under the hypotheses of lemma~\ref{prep:radial-lp}, for
every $e>0$ and $R=2^N$ one can remove pair mass
$O(R^{-\frac{e(p-1)}{2}})$ so that the two remaining directional measures
give every cap of radius $u$ mass at most $R^eu^h$, for all
$R^{-1}\le u\le1$ and all sufficiently large $R$.

The same cap bound lifts through a fixed orthogonal projection from
a higher-dimensional space, provided the original and projected
cross supports are bounded and separated. Constants may change,
and may be absorbed by a smaller preliminary choice of $e$.
\end{lemma}
\begin{proof}
Let $F(x,w)$ and $H(y,w)$ be the jointly Borel densities from
lemma~\ref{prep:radial-lp}. Remove the rows and columns of marginal
measure zero on which the density identities fail. The exceptional
pin sets may be taken Borel by testing the identities on a countable
determining family on $S^h$. With $T=R^{\frac{e}{2}}$, also remove
the pairs satisfying $F(x,\pi_x(y))>T$ or $H(y,\pi_y(x))>T$. The
product mass discarded by the first threshold is bounded by
$$
 \iint_{\{w:F(x,w)>T\}}F(x,w)d\sigma(w)d\beta(x)
 \le T^{1-p}\iint F(x,w)^pd\sigma(w)d\beta(x).
$$
The other threshold has the same bound with the roles reversed.
For each pin with a density identity, the remaining directional
measure is bounded by $T\sigma$. At an exceptional pin its fiber is
empty. Each cap of radius $u$ therefore has mass at most
$CTu^h\le R^eu^h$ for sufficiently large $R$.

Let $P$ be the fixed projection. On all actual cross directions,
$|Pw|\ge c>0$, because projected separation is bounded below and
original distances are bounded above. The map
$$
 \Phi(w)=\frac{Pw}{|Pw|}
$$
is Lipschitz on this region and satisfies $\Phi(\pi_x(y))=\pi_{Px}(Py)$.
The image of any original cap, intersected with the actual directions, is
therefore contained in a projected cap of comparable radius, and the
pushforward of the source measure gives the asserted bound at each original
pin, with no injectivity of $P$ required and in both orientations. For
vectors $u,v$ with $|u|,|v|\ge c$, the triangle inequality gives
$$
\left|\frac{u}{|u|}-\frac{v}{|v|}\right|
\le\frac{2}{c}|u-v|.
$$
Applying this with $u=Pw$ and $v=Pw'$, and using $|P(w-w')|\le|w-w'|$,
proves the claimed Lipschitz bound for $\Phi$ on the actual
directions and hence the asserted comparison of cap radii.
\end{proof}

\subsection{A decreasing coarse filter in odd dimensions}

For an affine subspace $L$ and $r>0$, write
$L_r=\{z:\operatorname{dist}(z,L)\le r\}$. In particular,
$H_r=\{z:\operatorname{dist}(z,H)\le r\}$ when $H$ is an affine hyperplane.
The following construction gives the spatial regularity that will survive
all the way to the limiting distance measure.

The proof adapts the heavy-plate localization and exclusion argument of
Shmerkin~\cite[appendix~B]{Shmerkin2023}, in the quantitative form
developed by Ren~\cite[section~7]{Ren2023}. Here the exclusions are
unions of coarse cell pairs at successive dyadic scales, producing
decreasing pair filters with positive limiting product mass.

At each dyadic width we remove pairs near heavy plates. After an initial
range of scales, the mass first removed at a given width has a power
bound, and its sum over the remaining widths will control the distance
between the coarse measure and its limit, while the finer cap estimates
will use a separate deletion at each scale.

\begin{lemma}[Coarse exclusion of heavy plates]\label{prep:coarse-filter}
Let $k\ge2$, and let $\rho_1,\rho_2$ be bounded separated
$s$-Frostman probabilities in ${\mathbb R}^{k+1}$, with
$k-\frac{1}{2}<s<k$. Suppose $\rho=\rho_1+\rho_2$ gives zero mass
to every affine hyperplane. For every fixed $0<\omega\le1$, there
are decreasing Borel sets $\Gamma_N\subset\operatorname{supp}\rho_1\times\operatorname{supp}\rho_2$
such that
\begin{enumerate}
\item $\Gamma_\infty=\bigcap_N\Gamma_N$ has positive product mass.
For some $c_\Gamma>0$ and all sufficiently large $N$,
\begin{equation}\label{prep:coarse-tail}
 (\rho_1\otimes\rho_2)(\Gamma_N\setminus\Gamma_\infty)
 \le C2^{-c_\Gamma N}.
\end{equation}
\item $\Gamma_N$ is a union of pairs of cells in a fixed dyadic
grid of depth at most $\omega N+O(1)$. At most
$C2^{(k+1)\omega N}$ cells meet either bounded support at that depth.
\item Both directional measures of $\Gamma_N$ satisfy
$\lambda\{w:\operatorname{dist}(w,V)\le u\}\le C_\omega u^\zeta$
for $2^{-N}\le u\le1$, with some $\zeta>0$.
\end{enumerate}
\end{lemma}
\begin{proof}
Put $a=s-(k-1)>\frac{1}{2}$. Covering a bounded portion of an affine
$(k-1)$-plane by $O(r^{-(k-1)})$ balls gives
\begin{equation}\label{prep:lower-plane}
\rho(L_r)\le C r^a
\end{equation}
uniformly over such planes and $0<r\le1$. This also bounds
neighborhoods of all lower-dimensional affine subspaces.
We shall fix a sufficiently small $b>0$ and then a sufficiently large
integer $j_0$.

Call $H$ heavy at width $r$ if $\rho(H_r)\ge r^b$.
We first prove that its parameter set has an $r$-net of size
\begin{equation}\label{prep:heavy-net}
O(r^{-C_0b}),
\end{equation}
where $C_0$ is independent of sufficiently small $b$. Restrict to the
compact family of hyperplanes meeting a fixed bounded neighborhood of the
supports, and use the usual metric given by the normal direction, modulo
sign, together with the offset. Fix $L>\frac{2}{a}$ and, within $H_r$,
choose $k+1$ points successively, each at distance at least $r^{Lb}$ from
the affine span of its predecessors. The estimate~\eqref{prep:lower-plane}
bounds each forbidden set by $Cr^{Lab}\le \frac{r^b}{2}$ when $r$ is
sufficiently small, so that the product measure of admissible tuples is at
least $c r^{(k+1)b}$, and every such tuple has $k$-simplex volume at least
$c_k r^{kLb}$.

If a fixed admissible tuple belongs to the width-$r$ plates of two
hyperplanes, then elementary linear algebra bounds the parameter distance
between these hyperplanes by $Cr^{1-kLb}$. Indeed, the matrix of the $k$
difference vectors has bounded operator norm and smallest nonzero singular
value at least a constant times its $k$-volume, whereas both normals have
scalar products $O(r)$ with every difference, and the offset is then
controlled by either first vertex. Since the hyperplane parameter space has
dimension $k+1$, an $r$-separated list contains at most $Cr^{-k(k+1)Lb}$
hyperplanes compatible with this tuple, and integrating tuple incidences
over any $r$-separated list of heavy hyperplanes proves
\eqref{prep:heavy-net}, with, for example, $C_0=(k+1)(1+kL)$ after the
fixed constants have been increased. A maximal separated list is the
desired net, and its plates, enlarged to width $C_kr$, contain every
original width-$r$ heavy plate in the bounded region, as well as every grid
cell of side $r$ meeting one.

For each $j\ge j_0$, put $r_j=2^{-j}$ and delete every pair of
depth-$j$ cells which both meet some heavy width-$r_j$ plate, and let
$E_j$ be the deleted pair set, a finite union of products of Borel cells.
With $J_N=\lfloor\omega N\rfloor$, define
\begin{equation}\label{prep:gamma-definition}
\Gamma_N=(\operatorname{supp}\rho_1\times\operatorname{supp}\rho_2)
              \setminus\bigcup_{j_0\le j\le J_N}E_j.
\end{equation}
These sets decrease, nested grids give the stated partition, and the
cell count follows from the boundedness of the supports in
$\mathbb R^{k+1}$. Write $F_j=E_j\setminus\bigcup_{j_0\le i<j}E_i$ for
the pairs first removed at depth $j$; these sets are disjoint.

At a fixed width $r$, use the enlarged plates of its net. Two such plates
with normal angle greater than $\sqrt r$, if they meet the bounded region,
have their intersection within an $O(\sqrt r)$-neighborhood of an affine
$(k-1)$-plane, as one sees by solving the two normal equations, the
smallest singular value of the pair of unit normals being comparable to
their angle. The equations~\eqref{prep:lower-plane}
and~\eqref{prep:heavy-net} therefore show that the union of all these
intersections has $\rho_1$-mass at most
\begin{equation}\label{prep:plate-exception}
C r^{\frac{a}{2}-2C_0b}.
\end{equation}
Outside this exceptional set, all the listed plates containing a fixed
first endpoint have pairwise normal angle at most $\sqrt r$, so that their
bounded portions lie in one plate of width $C\sqrt r$, the common endpoint
controlling the offsets as well as the normals. Consequently all the second
endpoints lost at this scale lie in that single plate.

For $j>4j_0$, set $j'=\lfloor \frac{j}{3}\rfloor$, taking $j_0$ large
enough that $j'\ge j_0$ and $C2^{-\frac{j}{2}}<2^{-j'}$. The single plate
just obtained is then contained in a width-$2^{-j'}$ plate containing the
first endpoint. If this enclosing plate is heavy at that earlier width,
every pair under consideration was already deleted by $E_{j'}$, since the
plate contains the actual endpoints and both of their depth-$j'$ cells meet
it; and if it is light, all the newly affected second endpoints have total
$\rho_2$-mass at most $2^{-bj'}$. It follows from
\eqref{prep:plate-exception} that the new product mass lost at the scale
$j$ is at most
\begin{equation}\label{prep:late-loss}
C2^{-j(\frac{a}{2}-2C_0b)}+2^{-b\lfloor \frac{j}{3}\rfloor}.
\end{equation}
Treat the initial scales $j_0\le j\le4j_0$ together.
Their combined net has at most $C_b2^{4C_0bj_0}$ members.
Enlarge them all to width $C2^{-j_0}$ and repeat the intersection
argument at angle $2^{-\frac{j_0}{2}}$. The exceptional first-endpoint
mass is at most
$$
 C_b2^{-(\frac{a}{2}-8C_0b)j_0}.
$$
For every other first endpoint, all lost second endpoints lie in
one width-$C2^{-\frac{j_0}{2}}$ hyperplane plate.
The hypothesis of zero hyperplane masses implies
\begin{equation}\label{prep:uniform-vanishing}
\sup_H\rho(H_\varepsilon)\longrightarrow0
        \text{ as } \varepsilon\downarrow0.
\end{equation}
Indeed, a contrary sequence of hyperplanes would have a convergent
subsequence in the compact parameter space, and its shrinking plates would
force positive mass on the limiting hyperplane. The total initial loss
therefore tends to zero as $j_0\to\infty$.

Choose $b>0$ so small that $8C_0b<\frac{a}{4}$ and $kLb<\frac{1}{2}$,
so that the sum of \eqref{prep:late-loss} over $j>4j_0$ also tends to
zero as $j_0\to\infty$, and then choose $j_0$ so large that the combined
loss is less than $\frac{1}{2}$, which proves positive limiting mass.

For $J_N>4j_0$ the mass that remains to be removed is the sum of
the masses of the $F_j$ with $j>J_N$. Set
$$
 b_* = \min\left\{\frac{a}{2}-2C_0b,\frac{b}{3}\right\}>0.
$$
The equation~\eqref{prep:late-loss} gives
$(\rho_1\otimes\rho_2)(F_j)\le C2^{-b_*j}$ in this range, including
the bounded factor due to rounding $\frac{j}{3}$. Hence
$$
 (\rho_1\otimes\rho_2)(\Gamma_N\setminus\Gamma_\infty)
 \le C\sum_{j>\lfloor\omega N\rfloor}2^{-b_*j}
 \le C'2^{-\omega b_*N}.
$$
This proves \eqref{prep:coarse-tail} with $c_\Gamma=\omega b_*$.
The initial scale $j_0$ may depend on the measures through
\eqref{prep:uniform-vanishing}. Once the measures are fixed, this
choice determines where the tail estimate starts.

Finally fix a pin $x$. Directions in a linear-hyperplane $u$-neighborhood
correspond to sources in a width-$Cu$ affine plate through $x$, and at an
available comparable exclusion scale a heavy plate leaves no retained
pairs, while a light plate has source mass at most its width to the power
$b$. At radii below the last exclusion scale use the last scale instead,
the finitely many radii above the first scale costing a fixed constant.
Both fibers therefore obey
$$
 \lambda\{w:\operatorname{dist}(w,V)\le u\}
 \le C\max(u,2^{-J_N})^b
 \le C_\omega u^{\omega b}
  \text{ for } 2^{-N}\le u\le1.
$$
Dyadic rounding and fixed width constants change only $C_\omega$. As
for the last power comparison, if $u\ge2^{-J_N}$ then
$u^b\le u^{\omega b}$ because $0<\omega\le1$, and otherwise
$2^{-J_N}\le2\cdot2^{-\omega N}\le2u^\omega$ since $u\ge2^{-N}$, so that
the displayed bound holds at every allowed radius. We may take
$\zeta=\omega b$, or any smaller positive exponent.
\end{proof}

\subsection{Starting pair measures}

We can now collect the properties of the pair measures used in the
rest of the proof. The cell partitions will be used for the Fourier
estimate, the cap bounds will control the mass lost to angular tests,
and the power tail will let us compare the result at every scale with a
fixed positive measure.

\begin{proposition}[Starting pair measures]\label{prep:starting-pairs}
Let $E\subset{\mathbb R}^d$ be compact, $d\ge2$, and
$\dim_HE>S=\frac{d}{2}$. Either lemma~\ref{prep:classical-plane} gives
a nonzero measure on $\Delta(E)$ with the power shell estimate
\eqref{prep:classical-shell} in an affine subspace, or the following
construction is possible.

After a similarity there are compact separated sets $E_1,E_2\subset E$
and probabilities $\mu_i$ supported on them, each satisfying
$$
 \mu_i(B(x,u))\le C u^{S+\chi}
$$
for a fixed $\chi>0$. Their cross directions lie in a fixed
arbitrarily small spherical patch. For every sufficiently small fixed
$\omega>0$, there are decreasing Borel pair sets
$\Gamma_N\subset E_1\times E_2$ such that
\begin{enumerate}
\item With $\lambda=\mu_1\otimes\mu_2$ and
$\Gamma_\infty=\bigcap_N\Gamma_N$, one has $\lambda(\Gamma_\infty)>0$ and
$\lambda(\Gamma_N\setminus\Gamma_\infty)\le C2^{-c_\Gamma N}$ for some
$c_\Gamma>0$ and all sufficiently large $N$.
\item In even dimensions $\Gamma_N=E_1\times E_2$. If $d=2k-1$, there
are Borel partitions $(X_i)$ of $E_1$ and $(Y_j)$ of $E_2$, each with
at most $C2^{(k+1)\omega N}$ members, and a matrix $A=(A_{ij})$ with
entries zero or one such that $1_{\Gamma_N}(x,y)=A_{ij}$ on every
$X_i\times Y_j$.
\item For every $e>0$, there are $c(e)>0$, $C_e<\infty$ and
Borel $B_N\subset\Gamma_N$ with
$(\mu_1\otimes\mu_2)(B_N)\le C_e2^{-c(e)N}$, such that for
$G_N=\Gamma_N\setminus B_N$, all $2^{-N}\le u\le1$, and
all spherical caps $B(w,u)$,
$$
\begin{aligned}
 \mu_2\{y:(x,y)\in G_N,\ \pi_x(y)\in B(w,u)\}
     &\le2^{eN}u^S,\\
 \mu_1\{x:(x,y)\in G_N,\ \pi_y(x)\in B(w,u)\}
     &\le2^{eN}u^S.
\end{aligned}
$$
\end{enumerate}
The last bounds integrate the original opposite marginal, without
normalizing its individual fibers. The $B_N$ need not have a product
partition or be nested. All constants are independent of $N$, but
may depend on the fixed measures, $\omega$, and $e$.
\end{proposition}
\begin{proof}
Use lemma~\ref{prep:projections} to choose an initial Frostman
probability of exponent greater than $S$. In odd dimension
$d=2k-1$, if it charges an affine $k$-plane, its normalized positive
restriction to that plane still has this exponent, and
$$
 S=k-\frac{1}{2}\ge\frac{k+1}{2} \text{ for } k\ge2.
$$
Lemma~\ref{prep:classical-plane}, applied isometrically within the
plane, proves the first alternative, the distances in question being
original distances. Otherwise project generically to ${\mathbb R}^{k+1}$, preserving energies
above $S$ and giving zero mass to all hyperplanes.

In both parities project further, if this is needed, to
${\mathbb R}^{h+1}$, where $h=\lfloor S\rfloor$, all the target dimensions
being greater than $S$; in even dimension this is the only auxiliary
projection. Use the simultaneous potential restrictions of
lemma~\ref{prep:projections} to obtain ball exponents above $S$ in all the
spaces, and choose two distinct support locations whose smallest projected
images are distinct, which is possible because that projected law is
nonatomic. Restrict to small compact neighborhoods of these locations of
positive measure; their supports are then separated in every required
space, and their original cross directions lie in a fixed small patch.
Normalize, and decrease the common exponent to $S+\chi$ if necessary.

Choose $0<\epsilon<\min(\chi,h,1)$ small enough that both laws in
${\mathbb R}^{h+1}$ have finite $I_{h+\epsilon}$ energy. Lemma~\ref{prep:radial-lp} therefore gives the
integrated radial estimates in both orientations for these projected
laws. Denote the separated probabilities on the original sets by
$\mu_1,\mu_2$.

If $d$ is even, $h=S$. Take $\Gamma_N=E_1\times E_2$ and
apply lemma~\ref{prep:integer-caps}, lifting to the original space.
This includes $d=2$, where the smallest projection is the identity.

If $d=2k-1$, write $\mu_i'$ for the laws in ${\mathbb R}^{k+1}$, whose
exponent may be decreased to an $s\in(k-\frac{1}{2},k)$ and which give
zero mass to every affine hyperplane, and apply
lemma~\ref{prep:coarse-filter} with the prescribed $\omega$. The integer
cap bounds obtained in ${\mathbb R}^k$ from lemma~\ref{prep:integer-caps}
lift to ${\mathbb R}^{k+1}$ and persist on restriction to the coarse
filters, so that $k-1$ is available there, and
proposition~\ref{prep:thin-upgrade} then gives availability of
$k-\frac{1}{2}=S$. Lift the resulting cap bounds and the exceptional pair
sets to the original space.

Pullback preserves intersections and product masses, so it preserves the
estimate \eqref{prep:coarse-tail}. In the auxiliary space let $X_i'$ and
$Y_j'$ be the cells at the finest exclusion depth that meet the
respective supports, of which there are at most $C2^{(k+1)\omega N}$ on
either side. Each cell product is either retained or deleted, its
indicator value gives $A_{ij}$, and the preimages of the cells give the
required partitions $X_i,Y_j$. In even dimensions the one product has
zero tail at every scale.
\end{proof}

\begin{remark}[Order of parameters]\label{prep:parameters}
The Frostman gap $\chi$ and the measures are fixed first.
The parameter $\beta>0$ will be the excess level used to choose
the entry depth, $\kappa$ sets the block length as a fraction of
$N$, and $\sigma$ bounds the small power losses.
Later arguments choose $\beta,\kappa,\sigma$ using that gap, and
then choose the coarse cutoff $\omega$ sufficiently small to pay for
the product partition. Only after this choice is the plate exponent
$\zeta=\omega b$ fixed and the thin-tube upgrade invoked with its
desired loss, such as $e=\sigma$, and the resulting deletion exponent
$c(e)$ may be very small. The final regularity exponent is bounded
by the fine deletion rate, the coarse tail rate $c_\Gamma$, and half
the Fourier shell saving, each of which is positive for the fixed
parameters.
\end{remark}

\section{Mass curves and angular tests}\label{test:section}

The directional estimates must now be put into a form which can be used
inside the iteration. One test limits how many occupied cells can lie close
to a line in a given direction, while a second limits how often a scalar
projection places two points close together; these answer different
questions and both are needed. We first express their bounds in terms of
the mass curve, and then prove that imposing all the required tests removes
little of the original pair mass.

We construct the directional tests from normalized regular classes and use
the cap estimates of proposition~\ref{prep:starting-pairs} to bound their
failure mass. A parallel test counts all the occupied class cube centres
satisfying its stated ball and tube conditions, whereas a projection test
uses the full conditional probability in a specified class cube, and the
exceptional-direction estimate applies to a spherical subprobability. At
its outer use, that subprobability is an unnormalized directional fiber of
the original opposite starting law, restricted to the cap-controlled pair
set.

By a spatial multiplier we mean a bounded measurable function of an
endpoint which multiplies its measure, for instance the indicator of a
Borel set, and an angular multiplier is a function of direction used in the
same way. For each frequency shell, a specified product of the masks
constructed below will weight the original product measure restricted to
$\Gamma_N$, the smaller cap-controlled subset being used only to estimate
the original pair mass removed by these masks. The cell matrix for
$\Gamma_N$ permits the Fourier estimate to use bounded spatial
multipliers at the two endpoints. The internal positive counts restrict
products of normalized class laws by the good tests. The angular energies
use endpoint masks together with the spatial multipliers.

Let $S=\frac{d}{2}$ and $R=2^N$. We use the separated probabilities from
proposition~\ref{prep:starting-pairs}, of Frostman exponent $S+\chi$.

In the inequalities below, $x\ll_{\mathcal P}y$ means that $x>0$ is chosen
sufficiently small relative to $y$ and to the fixed parameters listed as
$\mathcal P$; it is an instruction about the order in which the choices are
made, and it does not allow the parameters to vary with $N$. From this
point on, the letter $\sigma$ denotes a small loss exponent whenever it
occurs in a power of $R$, while surface measure continues to occur as
$d\sigma(w)$ in spherical integrals, the two uses being distinguished by
their arguments.
Choose the parameters in the following order.
$$
0<\beta<\min\left(\frac{\chi}{10},\frac{1}{10}\right), 0<\kappa\ll_d\frac{\beta}{d} \text{ and } 0<\sigma\ll_{d,\beta,\kappa}\kappa.
$$
The spatial cutoff parameter $\omega$ in
proposition~\ref{prep:starting-pairs} is chosen afterwards, as small as
required by these fixed parameters. All estimates concern
sufficiently large $N$ after these choices. The constants multiplying
$\kappa$, $\sigma$, or $\frac{\sigma}{\kappa}$ in exponents are dimensional;
ordinary multiplicative constants and the rates in $o(1)$ may depend
on the fixed parameters and the starting geometry. Use the preparation
proposition with angular loss $R^\sigma$. Let $\rho$ range over the
normalized classes supplied by applying lemma~\ref{reg:regularization}
with parameter $\sigma$ to each starting probability, and write $m$ for
its mass curve. All definitions in this section are made separately for each
class.

\subsection{Mass-curve increments and conditional laws}

The notation $I(g,p)$ introduced below is a mass-curve increment between
the depths $g$ and $p$, and it should not be confused with the Riesz energy
$I_s(\mu)$, which carries a subscript and takes a measure as its argument.
We also write $\min_{[g,p]}A$ as shorthand for
$\min\{A_n:g\le n\le p,\ n\in{\mathbb Z}\}$. For integer depths set
\begin{equation}\label{test:profile-definitions}
\begin{aligned}
M(n)&=R^{-m(\frac{n}{N})},& A_n&=m\left(\frac{n}{N}\right)-\frac{Sn}{N},\\
I(g,p)&=m\left(\frac{p}{N}\right)-m\left(\frac{g}{N}\right),&
k(g,p)&=A_p-\min_{g\le n\le p}A_n.
\end{aligned}
\end{equation}
Every minimum over depths is taken over integers unless we say otherwise.
Here $M(n)$ is the model mass of a depth-$n$ cell and $I(g,p)$ records the
mass curve increment between two depths, while the residual $A_n$ subtracts
the critical slope $S$ and $k(g,p)$ is the rise of $A_p$ above its minimum
since $g$. These are the quantities which will produce the cancellations in
the two directional estimates below.

The two estimates in this lemma come from opposite uses of regularity. To
count cells one divides an upper bound for the mass in a large ball by a
lower bound for the mass of one small cell, whereas to estimate a
conditional measure one divides the mass in a small ball by the mass of its
parent, and writing $m$ as its critical linear part plus $A$ exposes
the excess remaining in both calculations.

Whenever $P$ is a dyadic cell, $c_P$ denotes its center and
$\rho_P=\frac{\rho|_P}{\rho(P)}$ denotes its full conditional
probability, provided $\rho(P)>0$. In particular, the rescaling in
the next lemma is around the center of the parent cell.
\begin{lemma}[Mass-curve bounds]\label{test:profile-bounds}
There is a dimensional constant $C_d$ such that
$$
A_0=0, |A_i-A_j|\le d\frac{|i-j|}{N} \text{ and } A_n\ge\chi \frac{n}{N}-C_d\sigma.
$$
Every occupied depth-$n$ cube has mass
$M(n)R^{\pm C_d\sigma}$.  In particular, a ball of radius
$C2^{-g}$ contains at most
\begin{equation}\label{test:centre-count}
C'R^{I(g,p)+C_d\sigma}
\end{equation}
occupied depth-$p$ cube centres, for $0\le g\le p\le N$.
Here $C$ is any fixed constant and $C'$ may depend on $C$.
For integers $0\le p\le t\le N$ and an occupied depth-$p$ cube $P$,
the full conditional probability $\rho_P=\frac{\rho|_P}{\rho(P)}$,
rescaled by $z\mapsto2^p(z-c_P)$, satisfies the $S$-ball bound down to
$H^{-1}$, where $H=2^{t-p}$, with loss
\begin{equation}\label{test:conditional-ball-loss}
K=R^{A_p-\min_{p\le n\le t}A_n+C_d\sigma}.
\end{equation}
\end{lemma}

\begin{proof}
For later use, the identity behind the conditional estimate is
$$
m\left(\frac{n}{N}\right)-m\left(\frac{p}{N}\right)
=\frac{S(n-p)}{N}+A_n-A_p.
$$
The first term produces the desired power of the relative radius, whereas the
most negative possible value of $A_n-A_p$ supplies the loss.

The first two assertions follow from $m(0)=0$ and the fact that the slopes
of $m$ lie in $[0,d]$, while the lower bound follows from
\eqref{reg:profile-lower-bound}, the fixed Frostman constant being absorbed
into $R^{C_d\sigma}$ for sufficiently large $N$. The occupied-cell mass
assertion is the estimate~\eqref{reg:cell-mass}. Since a ball of radius
$C2^{-g}$ meets $O_{d,C}(1)$ depth-$g$ cubes, whose full masses are at most
$M(g)R^{C_d\sigma}$ each, whereas each occupied depth-$p$ cube has mass at
least $M(p)R^{-C_d\sigma}$, summation gives \eqref{test:centre-count}, with
a larger $C_d$.

For the last assertion, if $p\le n\le t$, a physical ball of radius
$2^{-n}$ meets a bounded number of depth-$n$ cubes, and division by
$\rho(P)$ gives an upper bound
$$
R^{-[m(\frac{n}{N})-m(\frac{p}{N})]+C_d\sigma}
=2^{-S(n-p)}R^{A_p-A_n+C_d\sigma}.
$$
Dyadic comparison of radii now yields
\eqref{test:conditional-ball-loss}, and since only depths through $t$ are
used, no mass estimate below $H^{-1}$ is involved.
\end{proof}

The two mass normalizations remain distinct throughout: $M(n)$ is within
the stated comparison factor of the actual mass of an occupied depth-$n$
cell, whereas the probability $\rho_P$ is divided by the actual mass
$\rho(P)$.

\subsection{Directional tests}

A descriptor $\mathfrak t$ consists of a test type together with a pair of
integer depths. A parallel descriptor $[g,p]$, with $0\le g<p\le N$, has
anchor $p$ and nominal length $p-g$, while a projection descriptor $[p,t]$,
with $0\le p<t\le N$, has anchor $p$ and nominal length $t-p$, the type
being part of the identity even when the displayed depth pairs coincide.
These descriptors specify individual tests, and a later construction will
select fixed lists of them for depth intervals.

Set
$$
J_0=\left\lceil100\left(\frac{1}{\kappa}+\frac{1}{\sigma}\right)+100\right\rceil
 \text{ and }
 W_i(n)=\left(2-\frac{i}{J_0}\right)2^{-n}R^{500\sigma}
 \text{ for } 0\le i\le J_0.
$$
Here $J_0$ is the maximum test index, fixed independently of $N$ once
$\kappa,\sigma$ have been fixed, so that there are $J_0+1$ good-test
levels, the masks below being defined for $i<J_0$.

\begin{definition}[Good directions]\label{test:good-directions}
Fix a class, a descriptor $\mathfrak t$ with anchor $p$, and an occupied
depth-$p$ cube $P$, of centre $c_P$.  Its level-$i$ good set
$G_{i,\mathfrak t}(P)\subset S^{d-1}$ is defined as follows, with one
sufficiently large dimensional constant $C_d$ for all tests.
\begin{enumerate}
\item For a parallel descriptor $[g,p]$, a direction $w$ is good if
      the number of occupied depth-$p$ cube centres in
      $B(c_P,10(d+1)2^{-g})$ whose distance to $c_P+{\mathbb R} w$ is at most
      $W_i(p)$ is at most
\begin{equation}\label{test:parallel-threshold}
R^{k(g,p)+C_d\sigma}.
\end{equation}
\item For a projection descriptor $[p,t]$, let
      $\rho_P=\frac{\rho|_P}{\rho(P)}$.  A direction $w$ is good if
\begin{equation}\label{test:projection-threshold}
(\rho_P\otimes\rho_P)
      \{|(z-z')\cdot w|\le W_i(t)\}
      \le2^{-(t-p)}
           R^{A_p-\min_{p\le n\le t}A_n+C_d\sigma}.
\end{equation}
\end{enumerate}
For an unoccupied anchor cube set $G_{i,\mathfrak t}(P)=\varnothing$.
For $x$ in the root cube, write
$G_{i,\mathfrak t}(x)=G_{i,\mathfrak t}(P)$ for its unique depth-$p$
descendant $P$, and set $G_{i,\mathfrak t}(x)=\varnothing$ outside the root.
The parallel test uses all occupied class centres in its stated ball and
tube, whereas the projection test uses the full conditional class
probability in $P$, and subsequent pair restrictions do not change
either test.
\end{definition}

\begin{figure}[htbp]
\centering
\begin{tikzpicture}[>=Stealth,font=\small]
\begin{scope}
\draw[gray] (0,0) rectangle (5.4,3);
\fill[gray!20] (0,1.1) rectangle (5.4,1.9);
\draw[->,thick] (0,1.5)--(5.7,1.5) node[right] {$w$};
\foreach \x/\y in {.8/1.3,1.6/1.65,2.6/1.45,4/1.7,4.8/1.35}{\fill (\x,\y) circle (2pt);}
\foreach \x/\y in {1/2.5,2/.4,3.9/2.5,4.8/.4}{\draw (\x,\y) circle (2pt);}
\draw[<->] (3.3,1.5)--(3.3,1.9);
\node[above] at (3.3,1.95) {$W_i(p)$};
\node[below,align=center] at (2.7,-.25)
{parallel test counts centers\\near $c_P+\mathbb R w$};
\end{scope}
\begin{scope}[xshift=8cm]
\draw[gray] (0,0) rectangle (5.4,3);
\fill[gray!20] (2.1,0) rectangle (2.9,3);
\draw[->,thick] (0,-.0)--(5.7,0) node[right] {$w$};
\fill (2.25,.6) circle (2pt) node[left] {$z$};
\fill (2.7,2.35) circle (2pt) node[right] {$z'$};
\draw[dashed] (2.25,.6)--(2.25,0);
\draw[dashed] (2.7,2.35)--(2.7,0);
\draw[<->] (2.25,.2)--(2.7,.2);
\node[below,align=center] at (2.7,-.25)
{projection test counts pairs\\with $|(z-z')\cdot w|\le W_i(t)$};
\end{scope}
\end{tikzpicture}
\caption{The two directional tests measure different concentrations.
On the left, the filled points represent occupied cell centers inside
a thin tube parallel to $w$; the full test also imposes the stated
ball condition. On the right, points may be far apart in directions
perpendicular to $w$ and still have close scalar projections. This
test uses two independent samples from the full conditional law
$\rho_P$. The drawings show geometry only, not the relative masses
of the points.}
\label{fig:two-tests}
\end{figure}

\subsection{A Fourier estimate for half-dimensional directions}

The two ball bounds below use the critical exponent $S=\frac{d}{2}$, and
since their exponents sum to $d$, the annular estimate contributes a
constant at each dyadic frequency scale.
For this preview, let $\alpha$ be the spatial probability, let $K$ be the
loss in its $S$-ball bound, let $L$ be the loss in the $S$-cap bound
for the direction measure, and let $B\ge1$ be the radius of a frequency
annulus. These are the data specified in the following lemma, and it is
useful to predict the estimate before proving it. On a frequency
annulus of radius $B$ the directional measure puts at most $LB^{-S}$ mass
in a unit ball, while the spatial ball bound gives
$\int_{|u|\le4B}|\widehat\alpha(u)|^2du\lesssim KB^{d-S}$, so that
multiplying the two gives $LK B^{d-2S}=LK$. Each dyadic annulus therefore
has the same cost, and summing over the annuli produces the logarithm. The
convolution argument below makes this multiplication legitimate
\cite[lemma~6.3]{OpenAI2026}.

\begin{lemma}[Half-dimensional Fourier estimate]\label{test:half-dimensional-fourier}
Let $1\le H\le R$, $K,L\ge1$, and $S=\frac{d}{2}$.  Suppose $\alpha$ is a
probability in a fixed bounded subset of ${\mathbb R}^d$ with
$$
\alpha(B(x,r))\le Kr^S \text{ for } H^{-1}\le r\le1,
$$
and $\xi$ is a subprobability on $S^{d-1}$ with
$$
\xi(B(w,r))\le Lr^S \text{ for } R^{-1}\le r\le1.
$$
Then, with Fourier convention
$\widehat\alpha(v)=\int e^{-iz\cdot v}d\alpha(z)$,
\begin{equation}\label{test:fourier-estimate}
\int_{-H}^{H}\int_{S^{d-1}}|\widehat\alpha(sw)|^2d\xi(w)ds
 \le C LK(1+\log H).
\end{equation}
The constant depends only on $d$ and the fixed bound on the support's
diameter.  In particular it does not require ball bounds below
$H^{-1}$ or any Fourier decay.
\end{lemma}

\begin{proof}
We estimate one dyadic annulus first. The directional measure on that
annulus is in general singular, so we smooth it at unit frequency scale
before combining it with the spatial $L^2$ estimate, the bounded spatial
support of $\alpha$ permitting this smoothing at a cost of no more than a
fixed constant.

The interval $|s|\le1$ costs at most $2$.  For a dyadic $B\ge1$,
let $\lambda_B$ be the pushforward of
$$
1_{\{B\le|s|\le\min(2B,H)\}}dsd\xi(w)
$$
under $(s,w)\mapsto sw$. It is supported in the annulus $B\le|v|\le2B$ and
has total mass at most $2B$, while the preimage of a unit frequency ball
has radial length $O(1)$ and lies in at most two spherical caps of radius
$O(B^{-1})$, and hence the cap bound applies, since $B\le H\le R$. A fixed
covering handles radii differing from $B^{-1}$ by fixed constants, and the
case $B=O(1)$ uses the total mass bound. Thus
\begin{equation}\label{test:frequency-unit-ball}
\sup_v\lambda_B(B(v,1))\le C LB^{-S}.
\end{equation}
Translation of $\alpha$ changes only the phase of its transform, so
assume its support is in a fixed ball.  Choose a smooth compactly
supported spatial function equal to one on that ball.  Fourier
convolution and the weighted Cauchy--Schwarz inequality give a
nonnegative, integrable, rapidly decreasing kernel $\Phi$, independent
of $B$, for which
$$
|\widehat\alpha(v)|^2
 \le C\int_{{\mathbb R}^d}|\widehat\alpha(u)|^2\Phi(v-u)du.
$$
Indeed, if the spatial cutoff is $\chi$, then
$\chi\alpha=\alpha$ and Fourier inversion yields
$\widehat\alpha=(2\pi)^{-d}\widehat\chi*\widehat\alpha$.
Cauchy--Schwarz with weight $|\widehat\chi(v-u)|$ gives
$$
|\widehat\alpha(v)|^2
\le(2\pi)^{-2d}\|\widehat\chi\|_{L^1}
\int|\widehat\alpha(u)|^2|\widehat\chi(v-u)|du.
$$
Choose a nonnegative rapidly decreasing majorant of
$|\widehat\chi|$ to obtain the asserted kernel $\Phi$.  Integrating against
$\lambda_B$ and writing $W_B=\Phi*\lambda_B$, we obtain
$$
\int|\widehat\alpha(v)|^2d\lambda_B(v)
 \le C\int|\widehat\alpha(u)|^2W_B(u)du.
$$
Summing \eqref{test:frequency-unit-ball} over unit cubes
against the decay of $\Phi$ gives $W_B\le CLB^{-S}$.  On $|u|>4B$,
the support and total mass of $\lambda_B$ give, for every fixed $M$,
\begin{equation}\label{test:frequency-tail}
W_B(u)\le C_M B(1+|u|)^{-M}.
\end{equation}
A Gaussian majorant and its positive Fourier transform show that
$$
\begin{aligned}
\int_{|u|\le4B}|\widehat\alpha(u)|^2du
&\le CB^d\iint e^{-cB^2|z-z'|^2}d\alpha(z)d\alpha(z')\\
&\le CKB^{d-S}.
\end{aligned}
$$
The last bound follows by splitting the $z'$ integral
into the ball of radius $B^{-1}$ and the annuli of radii $2^jB^{-1}$, every
radius being at least $H^{-1}$. Up to radius one use the assumed ball
bound, and above radius one use total mass at most one, which is bounded by
$Kr^S$, the resulting sum being at most
$CKB^{-S}\sum_{j\ge0}2^{jS}e^{-c'4^j}$. The part with $|u|\le4B$ therefore
costs at most
$$
CLK B^{d-2S}=CLK.
$$
For the complement, use $|\widehat\alpha|\le1$ and \eqref{test:frequency-tail}, choosing $M>d+1$.  Its
contribution is at most $C B^{d+1-M}\le C$, which is absorbed by
$CLK$.  Summing over the $O(1+\log H)$ dyadic shells proves \eqref{test:fourier-estimate}.
\end{proof}

\subsection{Exceptional mass of the tests}

For a parallel test we average the number of centers in the tube over
directions and then apply Markov's inequality, whereas for a projection test
we turn the small scalar difference into a one-dimensional Fourier integral
and use the preceding lemma. In both cases the threshold is chosen above
the average by enough powers of $R^\sigma$ to make the exceptional
directional mass at most $R^{-\sigma}$.

\begin{lemma}[Exceptional directions]\label{test:exceptional-directions}
The constant $C_d$ in definition~\ref{test:good-directions} can be
chosen depending only on $d$ so that the following holds for all fixed
parameters and sufficiently large $N$.  If $\xi$ is any spherical
subprobability satisfying
$$
\xi(B(w,r))\le R^\sigma r^S
        \text{ for } R^{-1}\le r\le1,
$$
then every level-zero test, uniformly in the class and the occupied
anchor cube, has
\begin{equation}\label{test:exceptional-bound}
\xi(S^{d-1}\setminus G_{0,\mathfrak t}(P))\le R^{-\sigma}.
\end{equation}
\end{lemma}

\begin{proof}
First take a parallel descriptor.  The anchor centre itself costs at
most one in the count.  Divide the other centres into annuli of
distance comparable to $2^{-h}$ from $c_P$, where
$g-O_d(1)\le h\le p+O_d(1)$.  A given centre in this annulus is counted
only for directions in two caps of radius
$$
O(2^{-(p-h)}R^{500\sigma}).
$$
For radii at most one the angular mass is at most
$C R^{(1+500S)\sigma}2^{-S(p-h)}$, and for larger radii the same bound,
after adjusting the constant, follows from the subprobability property, all
the small radii used being at least a fixed multiple of $R^{-1}$, a change
of that multiple costing only a fixed cap covering. The number of centres
is bounded by $C R^{I(h,p)+C_d'\sigma}$ by \eqref{test:centre-count}, and
if an endpoint annulus has $h$ outside $[g,p]$, we replace it by the nearer
endpoint, which changes the scales only by a dimensional constant. The
identity
$$
I(h,p)-S\frac{p-h}{N}=A_p-A_h\le k(g,p)
$$
then gives expected count at most
$$
C(N+1)R^{k(g,p)+C_d''\sigma}.
$$
Here $C_d''$ is dimensional, including the factor $1+500S$. Choose
$C_d>C_d''+3$; then for fixed $\sigma$ and sufficiently large $N$, Markov's
inequality at the threshold in \eqref{test:parallel-threshold} gives
\eqref{test:exceptional-bound}.

For a projection descriptor, rescale $P$ to unit size, let $\alpha$
denote its conditional law, and set $H=2^{t-p}$ and
$\delta_A=A_p-\min_{p\le n\le t}A_n\ge0$, so that by
lemma~\ref{test:profile-bounds} the $S$-ball loss of $\alpha$ is
$K=R^{\delta_A+C_d'\sigma}$, while the rescaled tolerance of the
level-zero test is $\varepsilon=2H^{-1}R^{500\sigma}$. A constant
multiple of
$$
F_\varepsilon(u)=
 \left(\frac{\sin\left(\frac{u}{2\varepsilon}\right)}{\frac{u}{2\varepsilon}}\right)^2,
$$
where the quotient at $u=0$ is interpreted by continuity so that
$F_\varepsilon(0)=1$, majorizes $1_{[-\varepsilon,\varepsilon]}$, and
its Fourier transform is nonnegative, supported in
$[-\varepsilon^{-1},\varepsilon^{-1}]\subseteq[-H,H]$, and bounded by
$C\varepsilon$. Fourier inversion and Tonelli's theorem therefore bound
the mean collision mass by
$$
\begin{aligned}
\int (\alpha\otimes\alpha)
 \{|(z-z')\cdot w|\le\varepsilon\}d\xi(w)
&\le C\varepsilon
 \int_{-H}^H\int|\widehat\alpha(sw)|^2d\xi(w)ds\\
&\le C H^{-1}(1+\log H)
 R^{\delta_A+(501+C_d')\sigma}.
\end{aligned}
$$
The last inequality is lemma~\ref{test:half-dimensional-fourier} with
$L=R^\sigma$, and increasing the same dimensional $C_d$ so that
$C_d>504+C_d'$ absorbs the constants and $1+\log H$ and gives
\eqref{test:exceptional-bound} by Markov's inequality. This argument also
covers large tolerances; alternatively, if the threshold is at least one,
then the test holds for every direction. The choices of $C_d$ depend only
on previously fixed dimensional constants, and not on any small
parameter.
\end{proof}

\subsection{Smooth masks and products}

The width decreases as the index increases while the count or mass
threshold stays fixed, so that the sets $G_{i,\mathfrak t}$ increase with
$i$, and the considerable gap between the width exponent $500\sigma$ and
the smoothing exponent $50\sigma$ gives a uniform buffer.

Characteristic functions of the good sets would not permit repeated angular
differentiation, so we replace them by smooth functions with values between
zero and one. This requires space between the set on which a mask equals
one and the complement of its support, and the slightly narrower tube at
level $i+1$ creates that space. For projection tests we check the
same margin by bounding how much a scalar product changes when the
direction is moved. This construction is the buffered smoothing of
\cite[lemma~6.5]{OpenAI2026}.

For a function $w(x,\omega)$, we write $\operatorname{supp}_\omega w$ for its support as
a function of direction with $x$ held fixed. A tangent field is a smoothly
varying vector tangent to the sphere, and applying such a field to a
function means taking its directional derivative. In the following lemma
$T_1,\ldots,T_j$ are chosen from one fixed finite family of smooth tangent
fields spanning every tangent space, and $\|\cdot\|_\infty$ denotes the
supremum over the variables in the stated domain; the directional variable
there is written $w$, so that its partial support is denoted $\operatorname{supp}_w$.
\begin{lemma}[Masks with a directional buffer]\label{test:smooth-masks}
For a descriptor $\mathfrak t$ of nominal length $l'$, put
$$
D_{l'}=\max(1,2^{l'}R^{-50\sigma}).
$$
For all sufficiently large $N$, every $0\le i<J_0$ admits a mask
$w_{i,\mathfrak t}(x,w)\in[0,1]$, measurable in $x$ and smooth in $w$, such
that, on every occupied anchor cube,
\begin{equation}\label{test:mask-properties}
w_{i,\mathfrak t}=1\text{ on }G_{i,\mathfrak t}, \operatorname{supp}_w w_{i,\mathfrak t}\subset G_{i+1,\mathfrak t} \text{ and } \|T_1\cdots T_jw_{i,\mathfrak t}\|_\infty\le C_jD_{l'}^j.
\end{equation}
Here the $T_j$ are any members of a fixed finite collection of smooth
tangent fields spanning the tangent spaces of the sphere.  The masks
are antipodally symmetric and constant as functions of $x$ on each
anchor cube; they are zero on unoccupied cubes and outside the root.
Moreover, the
$cD_{l'}^{-1}$-neighbourhood of $G_{i,\mathfrak t}$ lies in $G_{i+1,\mathfrak t}$, with
the same assertion for its closure after decreasing the fixed $c$.
\end{lemma}

\begin{proof}
If the direction moves by $\theta$, then a point at distance $O(2^{-g})$
from an anchor line moves transversely by $O(2^{-g}\theta)$, while for a
projection test the scalar difference changes by at most $|z-z'|\theta$. We
compare each of these errors with $W_i-W_{i+1}$, the angular convolution
radius below being chosen smaller than both margins.

Consider directions $w,w'$ at distance at most $cD_{l'}^{-1}$.
In a parallel test $[g,p]$, changing direction changes the distance
of a tested centre to the line by at most
$$
C_d2^{-g}|w-w'|
 \le C_dc2^{-p}R^{50\sigma}.
$$
Here $l'=p-g$; if $D_{l'}=1$, use
$2^{l'}\le R^{50\sigma}$.  In a projection test $[p,t]$ the change
in the tested scalar displacement is at most
$$
C_d2^{-p}|w-w'|
 \le C_dc2^{-t}R^{50\sigma},
$$
by the same argument with $l'=t-p$.  By contrast, the consecutive
width gap at the relevant fine depth $n$ is
$$
W_i(n)-W_{i+1}(n)=J_0^{-1}2^{-n}R^{500\sigma}.
$$
For fixed parameters and sufficiently large $N$, each of the perturbations
above is strictly less than half this gap. Every configuration counted at
width $W_{i+1}$ in the direction $w'$ was therefore counted at width $W_i$
in the direction $w$; thus if $w$ is good at level $i$, then $w'$ is
good at level $i+1$, and the strict remaining margin also gives the
asserted inclusion for the closure of a smaller neighbourhood.

To construct the mask, take $r=c_0D_{l'}^{-1}$ with $c_0$ sufficiently
small, together with the indicator of the spherical $2r$-neighbourhood of
$G_{i,\mathfrak t}$, and convolve it on the sphere with a nonnegative
rotationally invariant smooth kernel of integral one, supported where the
chordal distance is less than $r$. Such a kernel has angular derivative
$L^1$ bounds $C_jr^{-j}$, by rescaling in a fixed finite collection of
spherical charts. The convolution equals one on $G_{i,\mathfrak t}$ and is
supported in the closure of its $3r$-neighbourhood, which is contained in
$G_{i+1,\mathfrak t}$, and the derivative bounds in
\eqref{test:mask-properties} follow. Since both the tests and the kernel
are invariant under $w\mapsto-w$, the mask is even. The construction is
made separately for each of the finitely many anchor cubes, and no spatial
differentiability is required.
\end{proof}

The iteration will carry products of masks, and differentiating a product
can only leave its support inside the intersection of the original
supports, so that a differentiated mask still certifies its good test at
the next support level. The size estimate comes from elsewhere: there are
only polynomially many descriptors, while the total number of derivatives
is fixed before $N$ grows, and therefore the product rule costs only a subpower
factor.

By the derivative budget we mean the fixed upper bound on the total order
of derivatives already assigned to factors in a product; it does not count
the later derivatives used to estimate that product.
\begin{lemma}[Products and retained supports]\label{test:product-symbols}
Consider a product with one factor for each occurrence in one or more
descriptor lists whose total size is at most $N^{C_*}$, with $C_*$ fixed
independently of $N$. A factor may be a mask $w_{i,\mathfrak t}$ for any
$0\le i<J_0$, or a derivative of order $q$ divided by $D_{l'}^q$, where
$l'$ is its nominal length. The indices may differ among factors. Assume
the total order of these prescribed derivatives is bounded independently
of $N$. A bounded number, independent of $N$, of measurable spatial
multipliers and smooth angular factors are also allowed, with fixed
numerical bounds for their sup norms and for all angular derivatives under
consideration. These bounds are common to the family of products in any
later supremum. Let $L$ be the largest nominal length, with $L=0$ for an
empty list, and set
\begin{equation}\label{test:total-symbol-scale}
D_*=R^\sigma\max(1,2^LR^{-50\sigma}).
\end{equation}
At each fixed angular derivative order $j$, the resulting product
$b$ satisfies
$$
\|T_1\cdots T_jb\|_\infty\le C_jD_*^j.
$$
For each factor with descriptor $\mathfrak t$ and index $i$, the support
of the product still imposes its $G_{i+1,\mathfrak t}$ condition, including
when that factor was already differentiated. All these assertions are
pointwise in the spatial variables. Constants may depend on the fixed
list-size exponent, prescribed derivative budget, and numerical bounds
for the extra factors, as well as on the previously fixed parameters.
\end{lemma}

\begin{proof}
An undifferentiated mask has norm at most one. Only a bounded number of
factors are prescribed derivatives, and their normalized derivatives have
bounds $C_jD_{l'}^j$ by lemma~\ref{test:smooth-masks}, uniformly over the
allowed indices. At order $j$ the product rule gives at most polynomially
many terms in $N$, by the assumed list-size bound and because $j$ is fixed,
and each term is bounded by a constant times
$\max(1,2^LR^{-50\sigma})^j$. For fixed $\sigma>0$ its polynomial
multiplicity is at most $R^{\sigma j}$ for sufficiently large $N$, which
gives \eqref{test:total-symbol-scale}, the order-zero bound having no such
multiplicity. Finally, the support of a derivative of a smooth function is
contained in the support of the function, so that none of the
differentiations removes an imposed good-direction condition.
\end{proof}

\subsection{Fixed lists under interval containment}

The individual tests and masks above are available for every descriptor,
and we now choose a fixed list for each depth interval. The list for a
larger interval has to contain the list for each subinterval, and each
chosen descriptor must have a selecting interval extending by its nominal
length on both sides of its anchor. Each interval makes its own scale
choice once, and that choice is reused whenever the interval contributes to
a containing list.

Put $q_0=\lceil\kappa N\rceil$. Partition $\{0,\ldots,N\}$ into
consecutive blocks of $q_0$ indices, with the last block possibly
shorter. In each block choose, once and for all, a depth minimizing
$A$; break ties by a fixed rule. An interval $[a,t]$ is
terminal if $t-a\le100q_0$. In a nonterminal interval select
$p$ by minimizing $A$ among the representatives of blocks wholly
contained in $[a+q_0,t-q_0]$, again using fixed tie breaking. Set
$$
l=\min(p-a,t-p).
$$
The choice is a prefix choice if $2p\le a+t$, and a
suffix choice otherwise.

For each interval $[a,t]$, form the list $\mathcal M(a,t)$ as follows.
For every nonterminal subinterval $[a',t']\subseteq[a,t]$, use its
fixed choice $p,l$ and include all parallel descriptors
$$
[g,p] \text{ with } p-l\le g<p.
$$
If the choice is a suffix choice, include also the projection
descriptor $[p,t']$. Repetitions within either type are removed.

For a membership $\mathfrak t\in\mathcal M(a,t)$, an origin means a
nonterminal subinterval $[a',t']\subseteq[a,t]$ whose fixed choice
generated $\mathfrak t$. The origin witnesses this particular membership,
and it is not part of the identity of the descriptor, so that different
memberships may use different origins. The initial entry will use one
additional parallel descriptor outside these lists, to which the individual
test and mask estimates apply without an origin.

Choosing the tests afresh after every restriction could lose a test needed
by a later reduction, and we avoid this by assigning to an interval, at the
outset, the tests arising from all of its subintervals, and therefore containment
of intervals implies containment of lists. The remaining assertions below
describe the geometry of the origin of a test, and they will permit
the two-depth payment argument later on.

\begin{lemma}[Descriptor geometry]\label{test:descriptor-properties}
For a nonterminal $[a,t]$, the union of the blocks eligible for its
choice is an integer interval $[a_*,t_*]$, called its inner interval, with
$$
a+q_0\le a_*\le a+2q_0, t-2q_0\le t_*\le t-q_0
 \text{ and }
 A_n\ge A_p \text{ for } a_*\le n\le t_*.
$$
The lists $\mathcal M(a,t)$ are nested under interval containment.
Across the scheduled lists, there are $O(\frac{1}{\kappa})$ possible anchors
and polynomially many descriptors. For every
$\mathfrak t\in\mathcal M(a,t)$ with anchor $p$ and nominal length $l'$,
there is an origin $[a',t']\subseteq[a,t]$ such that
\begin{equation}\label{test:descriptor-origin}
[p-l',p+l']\subseteq[a',t'].
\end{equation}
That origin also supplies every parallel descriptor $[g,p]$ with
$p-l'\le g<p$. In particular $l'\le\frac{t-a}{2}$.
\end{lemma}

\begin{proof}
The eligible blocks are consecutive, the first starting less than one block
length after $a+q_0$ and the last ending less than one block length before
$t-q_0$, and nonterminality ensures that eligible blocks exist. For an
index in their union, its own block representative has no larger excess
value, whereas the selected $p$ has no larger value than that
representative, which proves the first assertions.

The choice belonging to each origin is made independently of the interval
in whose list it is being included, so that an origin contributing to a
smaller interval contributes the same descriptors to every larger one,
which proves the nesting. Anchors are block representatives, of which there
are at most $1+\frac{N+1}{q_0}$, and a descriptor is specified by its type,
its anchor and one other integer depth, and hence their number is at most
$2(N+1)(1+\frac{N+1}{q_0})$.

A parallel descriptor coming from $[a',t']$ has $l'=p-g\le l$, while a
projection descriptor occurs only for a suffix choice and then has
$l'=t'-p=l\le p-a'$. These observations prove
\eqref{test:descriptor-origin}, since the list of parallel descriptors
supplied by the same origin contains all the asserted ones, and listing
each descriptor only once preserves the conclusion for every list
containing it. By definition at least one of the origins which generated a
descriptor in $\mathcal M(a,t)$ lies inside $[a,t]$, and that origin may be
chosen for this membership independently of any choice made for another
list.
\end{proof}

The scheduled lists satisfy the polynomial bound in
lemma~\ref{test:product-symbols}, as do those lists with any fixed number
of extra descriptors. Their origin properties are needed only when the
schedule is being used, whereas the individual test, exceptional-direction
and mask estimates apply to every descriptor defined above.

\section{Spherical oscillatory integrals and distance shells}
\label{osc:section}

Why do spherical Fourier integrals enter a problem about distances? At high
frequency their stationary points occur in the directions parallel to
$x-y$, and their oscillation contains $e^{\pm ir|x-y|}$, which is the
bridge from a distance-frequency band to a spherical energy. We need that
bridge in both directions, and with the directional tests still visible in
the coefficients, and the calculations below supply it, stating exactly how
many angular derivatives it uses.

By an endpoint symbol we mean the factor depending on one endpoint and on
the angular variable inside its Fourier integral, built here from masks,
their angular derivatives, and bounded multipliers; the permitted factors
and their bounds are specified in definition~\ref{osc:symbols}. This
section proves a Fourier estimate for weighted distance measures built from
separated class laws and endpoint masks, whose output is a spherical
angular energy with the same endpoint data. The proof uses global
stationary-phase operators on the sphere, so that their derivatives can be
distributed between the two endpoint symbols without creating a coefficient
depending jointly on the endpoints. Section~\ref{rec:section} will apply
this estimate to high frequencies in a positive count of pairs with nearly
equal distances.

Throughout this section, spatial and frequency depths are integers in
$[0,N]$, and in an angular energy they satisfy $0\le a\le v\le N$. The
largest nominal length of an empty descriptor list is zero, and we write
$n_*=d-1$.

\subsection{Separated cells and endpoint symbols}

The center of a dyadic cell $X$ is denoted by $c_X$.
Fix a power of two $K_d$ with
$1000(d+1)<K_d<4000(d+1)$.  Two depth-$a$ cells $X,Y$ are
called separated if
\begin{equation}\label{osc:separation}
\frac{K_d2^{-a}}{8}\le |c_X-c_Y|\le 4K_d2^{-a}.
\end{equation}
For $x\in X$, $y\in Y$ we then have $D(x,y)\asymp_d2^{-a}$, and their
directions $e(x,y)=\frac{x-y}{D(x,y)}$ lie in a cone of small fixed aperture.
Choose $0\le\psi\le1$ smooth on $S^{d-1}$, identically one on a fixed
neighborhood of these directions and zero on a fixed neighborhood of
their negatives.  Its derivatives of every fixed order are uniformly
bounded.  We use the same convention for the fixed separated outer
supports, with $a=0$; constants there may depend on the initial geometry.

The symbol $\partial_{w_v}$ denotes differentiation in the coordinate
$w_v$, and the combinations below are tangent to the sphere, so that their
action on a function defined on the sphere does not depend on an extension
off the sphere. Later $\nabla^m$ will denote the collection of derivatives
of order $m$, and the Hessian $\nabla^2D$ the matrix of second derivatives
of the distance function in its spatial variables. All the angular
derivatives below may be taken in the fixed smooth tangent fields
$$
 \Omega_{uv}=w_u\partial_{w_v}-w_v\partial_{w_u}
  \text{ for } 1\le u<v\le d.
$$
These fields span every tangent space of the sphere.  A word of length
$q$ in these fields is denoted by $\Omega^\alpha$, $|\alpha|=q$; the empty
word is the identity.  Recall that a descriptor $\mathfrak t$ of nominal
length $l_{\mathfrak t}$ has angular scale
$$
 D_{l_{\mathfrak t}}=\max(1,2^{l_{\mathfrak t}}R^{-50\sigma}).
$$

In the next definition $f(x)$ is a measurable spatial multiplier with
$|f(x)|\le1$ and $q(w)$ is a smooth angular multiplier, while
$\|q\|_{C^j(S^{d-1})}$ denotes the maximum of the suprema of all the
derivatives $\Omega^\alpha q$ with $|\alpha|\le j$, the empty derivative
included; using a sum instead would change only a constant. The numbers
$Q_j$ bound these norms uniformly over the symbols, an empty product has
value one, and a bar over a complex quantity denotes complex
conjugation.
\begin{definition}[Admissible symbols]\label{osc:symbols}
Fix an indexed descriptor list
$\mathcal A=(\mathfrak t_\nu)_{\nu\in\mathcal I}$ with
$\#\mathcal I\le N^{C_{\mathcal A}}$, a smoothing index
$i\in\{0,\ldots,J_0-1\}$, and a nonnegative integer derivative budget
$B$.  The exponent $C_{\mathcal A}$ and budget $B$ are fixed independently
of $N$ and are common to the symbol family under consideration.  At each
angular order under consideration, also fix a common numerical bound
$Q_j$ independent of $N$.  An admissible endpoint symbol has the form
\begin{equation}\label{osc:symbol-form}
b_x(w)=f(x)q(w)
       \prod_{\nu\in\mathcal I}
       D_{l_{\mathfrak t_\nu}}^{-|\alpha_\nu|}
       \Omega^{\alpha_\nu}w_{i,\mathfrak t_\nu}(x,w)
  \text{ with } \sum_{\nu\in\mathcal I}|\alpha_\nu|\le B.
\end{equation}
The index $\nu$ distinguishes occurrences, so repeated descriptors may
receive different derivative words.
Here $f$ is an arbitrary measurable spatial function with $|f|\le1$,
and $q$ is a smooth angular function satisfying
$\|q\|_{C^j(S^{d-1})}\le Q_j$ at every order under consideration.
The constants in estimates may depend on $C_{\mathcal A}$, $B$ and the
fixed $Q_j$.
Fixed scalar factors can be absorbed in those estimates.  Empty lists
and empty products are allowed, as are complex conjugates of such
symbols.  Symbols do not depend on the radial variable $r$.
\end{definition}

No spatial regularity is required in definition~\ref{osc:symbols}, and in
particular $f$ may impose an arbitrary Borel restriction inside a cell.
Every differentiated mask has the same angular support as the original
mask, or a smaller one; thus every descriptor present in a symbol still
imposes its level-$(i+1)$ good condition.

An angular derivative is expensive only when it falls on a mask which
changes on a small angular scale, and the nominal length of its descriptor
records that scale. Once the anchor depth has been reached, the
cell containing the endpoint is fixed, so that the associated mask no
longer depends on the endpoint variable and can be removed from a positive
independent angular energy by using its absolute bound.

\begin{lemma}[Derivative bounds and reached anchors]
\label{osc:symbol-bounds}
Let $L$ be the largest nominal length in the symbol's list, or in all
lists for a product of symbols, with $L=0$ when those lists are empty, and put
\begin{equation}\label{osc:derivative-scale}
D_*=R^\sigma\max(1,2^L R^{-50\sigma}).
\end{equation}
For every fixed $j$, an admissible symbol, or a product of a fixed number
of endpoint symbols, obeys
\begin{equation}\label{osc:angular-derivatives}
\sup_x\|b_x\|_{C^j(S^{d-1})}\le C_{j,B}D_*^j
\end{equation}
for sufficiently large $N$.  The constant may depend on the fixed
list-size exponents, parameters, derivative budgets and common angular
bounds, but not on the spatial multipliers.  On a depth-$a$ cell, the product of all mask factors
whose anchors have depth at most $a$ is independent of the spatial
integration variable and has uniformly bounded modulus.  We call these
anchors reached at depth $a$.
\end{lemma}

\begin{proof}
Apply lemma~\ref{test:product-symbols} to each endpoint symbol, using the
polynomial list bound, the prescribed derivative budget and the common
bounds for $q$; a fixed product of endpoint symbols then satisfies the same
estimate by the finite product rule, even when their smoothing indices
differ. For the last assertion, an anchor cell at depth at most $a$
contains the depth-$a$ cell and its mask is constant in the spatial
variable on that cell, a constancy preserved by angular differentiation.
The product has bounded modulus because only boundedly many factors were
prescribed derivatives, while all the other masks have modulus at most
one.
\end{proof}

\subsection{General angular energies}

In the next display the classes have probabilities $\rho_1,\rho_2$ and
model masses $M_1,M_2$, the cells being $X$ on the first side and $Y$ on
the second, with endpoint symbols $b_x(w)$ and $c_y(w)$. The transforms
$U_X$ and $V_Y$ are therefore functions of the same radial frequency $r$
and direction $w$, but use different endpoint measures, and when the
classes coincide the two copies retain the same mass curve. For occupied
depth-$a$ cells of one or two classes, with specified admissible endpoint
symbols, define
\begin{equation}\label{osc:transforms}
U_X(r,w)=\frac{1}{M_1(a)}\int_X e^{irw\cdot x}b_x(w)d\rho_1(x)
  \text{ and }
  V_Y(r,w)=\frac{1}{M_2(a)}\int_Y e^{irw\cdot y}c_y(w)d\rho_2(y).
\end{equation}
These transforms are normalized by the mass curve model masses rather than
by the actual cell masses. Before the endpoint factors are inserted, the
first underlying measure has mass $\frac{\rho_1(X)}{M_1(a)}=R^{O_d(\sigma)}$,
and similarly on the second side, and therefore it need not be a conditional
probability. Fix $\Psi\ge0$ in $C_c^\infty((\frac{1}{4},4))$ with
$\Psi\ge1$ on $[\frac{1}{2},2]$, and denote surface measure on the sphere
by $dw$. The word independent below refers to two separate angular
integrations, one for each endpoint, and imposes no probabilistic
assumption on the endpoint measures.

\begin{definition}[General angular energies]\label{osc:states}
For $0\le a\le v\le N$, the independent angular energy is
\begin{equation}\label{osc:H}
\mathsf H^{X,Y}(a,v)
 =2^{d(v-a)}\int\Psi\left(\frac{r}{2^v}\right)
   \left(\int_{S^{d-1}}|U_X(r,w)|^2dw\right)
   \left(\int_{S^{d-1}}|V_Y(r,w)|^2dw\right)
   \frac{dr}{2^v}.
\end{equation}
It is defined without any separation condition.  In separated geometry
we also use the form with one shared angle
\begin{equation}\label{osc:HF}
\mathsf H_F^{X,Y}(a,v)
 =2^{d(v-a)}\int\Psi\left(\frac{r}{2^v}\right)
   \left|\int_{S^{d-1}}\psi(w)
           U_X(r,w)\overline{V_Y(r,w)}dw\right|^2
   \frac{dr}{2^v}.
\end{equation}
The subscript $F$ is merely a label for the form with one shared angular
integration, not an additional function or a new parameter.
These are values for the specified class laws, cells and endpoint
symbols, including the lists present in those symbols.  We use maxima over
finite collections of cell pairs and suprema when symbols vary; the family
in a supremum will be stated.
\end{definition}

Cauchy--Schwarz in $w$, with $0\le\psi\le1$, gives at each radius
$$
 \left|\int\psi(w)U_X(r,w)\overline{V_Y(r,w)}dw\right|^2
 \le \left(\int|U_X(r,w)|^2dw\right)
      \left(\int|V_Y(r,w)|^2dw\right).
$$
The radial measure and prefactor in the two forms are the same, so
\begin{equation}\label{osc:angular-cs}
\mathsf H_F^{X,Y}(a,v)\le\mathsf H^{X,Y}(a,v)
\end{equation}
with the same endpoint symbols and their support conditions.

There is a second operation available only in the positive independent
energy. On a depth-$a$ cell, factor the reached mask factors out of its
transform as an angular function. Lemma~\ref{osc:symbol-bounds} bounds its
modulus, so that bounding that modulus inside the corresponding positive
angular integral gives, at bounded cost, the independent energy with those
factors deleted, the remaining endpoint symbols being left unchanged. The
same operation is not permitted inside the shared angular integral in
$\mathsf H_F$ before its absolute square has been taken, and the reason is
an elementary one about cancellation. From $|a(w)|\le1$ one can infer
$\int |a(w)U(w)|^2dw\le\int|U(w)|^2dw$, whereas one cannot infer
$|\int a(w)U(w)dw|\le|\int U(w)dw|$, since removing $a$ may create
cancellation which was absent before. For this reason the reached
factors are removed only after passing to the positive independent
energies. Here $a$
and $U$ are arbitrary bounded integrable functions, unrelated to the depth
variable $a$.

\subsection{The weighted shell estimate}

The scalar input below includes a distance weight which cancels the leading
stationary distance power, and the conclusion retains the specified
endpoint lists and the same cell depth.

\begin{proposition}[Weighted shell conversion]\label{osc:shell-conversion}
Let $0\le a\le v\le N$ be integer depths, and let $X,Y$ be separated
cells at depth $a$ of one or two classes.  Let
$m_x(w)$ and $n_y(w)$ be products of the undifferentiated masks at one
index $i<J_0$ for two specified indexed descriptor lists.  Assume both
lists have at most $N^{C_{\mathcal A}}$ occurrences, where one exponent
$C_{\mathcal A}$ is fixed independently of $N$ for the whole family under
consideration, as in definition~\ref{osc:symbols}.  Let $L$ be the largest
nominal length in either list, with $L=0$ when both lists are empty.
For measurable $|f|,|g|\le1$
define a possibly complex measure $\nu$ on the distance line by its integrals against bounded Borel functions $F$,
\begin{equation}\label{osc:weighted-measure}
\begin{aligned}
 \int F(s)d\nu(s)
 &=\frac{1}{M_1(a)M_2(a)}
   \int_{X\times Y} F(D(x,y))
       (D(x,y)2^a)^{-\frac{n_*}{2}}\\
 &  \cdot f(x)\overline{g(y)}
       m_x(e(x,y))n_y(e(x,y))
       d\rho_1(x)d\rho_2(y).
 \end{aligned}
\end{equation}
If
\begin{equation}\label{osc:shell-regime}
v-a\ge20\sigma N \text{ and } L\le\frac{v-a}{2}+\sigma N,
\end{equation}
then for every fixed $T>0$,
\begin{equation}\label{osc:shell-estimate}
2^{-a}\int_{2^{v-1}}^{2^{v+1}}|\widehat\nu(r)|^2dr
 \le R^{o(1)}\sup\mathsf H_F(a,v)+O(R^{-T}).
\end{equation}
The supremum has the same cells, class normalizations, lists and
smoothing index as in \eqref{osc:weighted-measure}.  Its symbols
are admissible with a derivative budget and common bounds for smooth angular factors
depending only on $T,d,C_{\mathcal A}$ and the fixed parameters.  All
these bounds are independent of $N$.
When $f=g=1$, the measure $\nu$ is positive.  The assertion also holds on
the fixed outer supports with $a=0$.
\end{proposition}

The rest of this section proves proposition~\ref{osc:shell-conversion}. We
obtain the global sphere expansions first, and then separate the remaining
smooth spatial coefficient into bounded functions of the two endpoints.

\subsection{Stationary phase with global angular coefficients}

The operator $\operatorname{Id}$ below is the identity operator on
smooth functions. Products of differential operators mean composition.
The gradient $\nabla$ and divergence $\operatorname{div}$ in this
subsection are taken in the induced Riemannian metric of the sphere.
Thus their composition is the spherical Laplacian, rather than the
ordinary Laplacian in the ambient Euclidean space.
Give the sphere Laplacian the convention
$\Delta_{S^{n_*}}=\operatorname{div}\nabla$, so its eigenvalues are
nonpositive.  Define
\begin{equation}\label{osc:L-operators}
\begin{aligned}
 c_{n_*}&=(2\pi)^{\frac{n_*}{2}}e^{-i\pi \frac{n_*}{4}},\\
 L_0&=c_{n_*}\operatorname{Id},\\
 L_j&=\frac{c_{n_*}}{(2i)^j j!}
   \prod_{q=0}^{j-1}
   \left(\Delta_{S^{n_*}}+q(q+1)
                    -\frac{n_*(n_*-2)}{4}\right) \text{ for } j\ge1.
 \end{aligned}
\end{equation}
For the opposite stationary point put
$L_{j,-1}b=\overline{L_j\overline b}$, and write $L_{j,+1}=L_j$.
These operators are defined globally, including when $d=2$.

For a fixed nonzero vector $x-y$, the phase $w\mapsto w\cdot(x-y)$ on the
unit sphere has two stationary points, in the directions
$\pm\frac{x-y}{|x-y|}$, each contributing a factor of size
$(r|x-y|)^{-\frac{d-1}{2}}$. We need more than this leading-order
description, since the coefficients must remain globally smooth functions
of direction and the error must be small after angular derivatives have
been taken, and the expansion below proves these assertions and then solves
recursively for the inverse expansion. Two scales occur in it. The
oscillation parameter $\Lambda=r|\Delta|$ measures frequency times
separation, where $\Delta$ is the nonzero spatial difference vector, while
the number $A\ge1$ in the lemma bounds the angular variation of the
amplitude $b$. Since each further stationary-phase coefficient
differentiates $b$ twice and contributes a factor $\Lambda^{-1}$, the
useful ratio is $A^2\Lambda^{-1}$, and the smoothing margins were chosen
earlier so as to make this ratio small where the expansion is applied. The
letter $A$ here is a derivative scale, unrelated to $A_n$.

The underlying analytic method is classical stationary phase
\cite[chapter~VIII]{Stein1993}. One integrates by parts away from
the stationary points and reduces the phase near each stationary point to a
quadratic form. Here we also write the coefficients as global differential
operators on the sphere and keep track of the angular derivatives in the
remainder, as in \cite[lemma~7.5]{OpenAI2026}, both of which are needed
when the expansion is used repeatedly in the reductions.

\begin{lemma}[Direct and inverse sphere expansions]
\label{osc:stationary}
Suppose $b\in C^\infty(S^{n_*})$ obeys
$\|b\|_{C^m}\le C_m A^m$ at the finitely many orders needed below, where
$A\ge1$.  Let $\Delta\ne0$, $e=\frac{\Delta}{|\Delta|}$,
$\Lambda=r|\Delta|$, and $r>0$.  A cutoff $\psi$ is assumed to equal one
on a fixed neighborhood of $e$ and zero on a fixed neighborhood of
$-e$, with uniformly bounded derivatives.  For every integer $J\ge1$,
\begin{equation}\label{osc:direct}
\begin{aligned}
 \int_{S^{n_*}}e^{irw\cdot\Delta}\psi(w)b(w)dw
 &=\Lambda^{-\frac{n_*}{2}}
   \left[e^{i\Lambda}\sum_{j=0}^{J-1}
                    \Lambda^{-j}(L_jb)(e)+E_J\right],\\
 |E_J|&\le C_J\Lambda^{-J}A^{2J+d+2},
 \end{aligned}
\end{equation}
provided $\Lambda\ge1$.  Without a cutoff the corresponding formula is
\begin{equation}\label{osc:full-sphere}
\int_{S^{n_*}}e^{irw\cdot\Delta}b(w)dw
 =\Lambda^{-\frac{n_*}{2}}
   \left[\sum_{\epsilon\in\{+1,-1\}}e^{i\epsilon\Lambda}
       \sum_{j=0}^{J-1}\Lambda^{-j}
                 (L_{j,\epsilon}b)(\epsilon e)+E_J'\right],
\end{equation}
with the same remainder bound.  Define global differential operators
\begin{equation}\label{osc:B-operators}
B_0=L_0^{-1} \text{ and } B_j=-L_0^{-1}\sum_{h=1}^j L_h B_{j-h} \text{ for } j\ge1.
\end{equation}
Each $B_j$ is a polynomial in $\Delta_{S^{n_*}}$ of differential order at
most $2j$.  If, for a fixed $c>0$,
\begin{equation}\label{osc:stationary-regime}
\Lambda\ge R^{c\sigma}, \frac{A^2}{\Lambda}\le R^{-c\sigma} \text{ and } 1\le A\le R^2,
\end{equation}
then for every prescribed $T>0$ a fixed sufficiently large $J$ gives
\begin{equation}\label{osc:inverse}
e^{i\Lambda}b(e)
 =\Lambda^{\frac{n_*}{2}}\sum_{j=0}^{J-1}\Lambda^{-j}
       \int_{S^{n_*}}e^{irw\cdot\Delta}\psi(w)(B_jb)(w)dw
       +O(R^{-T}).
\end{equation}
The extracted remainders in the equations~\eqref{osc:direct} and
\eqref{osc:full-sphere} are also $O(R^{-T})$ after increasing $J$.
All constants are uniform in $e$ in the stipulated geometry.

For $r=2^v s$ with $s$ in a fixed compact subinterval of $(0,\infty)$,
each coefficient $\Lambda^{-j}(L_{j,\epsilon}b)(\epsilon e)$ has uniformly
bounded derivatives of every fixed order in $s$ under \eqref{osc:stationary-regime}.  Every differential coefficient
retains the supports of all the mask factors in $b$.
\end{lemma}

\begin{proof}
There are three things to do. Away from the two stationary points we
integrate by parts; near each stationary point we choose coordinates making
the phase exactly quadratic and expand the amplitude; and finally we
compose the resulting differential operators to solve for the inverse
expansion. The last of these is finite algebra, but its derivative order
has to be fixed before we ask for the final error bound.

We first obtain a local expansion with a controlled remainder. Away from
$e$ and $-e$ the angular gradient of $w\cdot e$ is bounded below, so that
repeated integration by parts with this gradient bounds the integral there
by any chosen power of $\Lambda^{-1}$, at the cost of at most the
corresponding angular derivatives of $b$. Near the maximum, rotate $e$ to
the last coordinate vector and use
$$
 w(u)=\left(u\sqrt{1-\frac{|u|^2}{4}},1-\frac{|u|^2}{2}\right)
 \text{ and }
 dw=\left(1-\frac{|u|^2}{4}\right)^{\frac{n_*-2}{2}}du.
$$
The phase is exactly $w\cdot e=1-\frac{|u|^2}{2}$ and the Jacobian equals
one at $u=0$, these charts and their derivatives being chosen uniformly on
finitely many patches of possible $e$. The local integral is therefore
$e^{i\Lambda}\int e^{-i\Lambda\frac{|u|^2}{2}}a_e(u)du$, with $a_e$
compactly supported in a fixed coordinate ball and
$\|a_e\|_{C^m}\lesssim_m A^m$.

The Fourier formula for a Gaussian, obtained first with positive damping
and then by taking its limit, expresses the quadratic integral
$\int e^{-i\Lambda\frac{|u|^2}{2}}a_e(u)du$ as
$c_{n_*}\Lambda^{-\frac{n_*}{2}}$ times the integral of
$\widehat a_e(\xi)e^{i\frac{|\xi|^2}{2\Lambda}}$ with the Fourier inversion
normalization. Explicitly,
$$
 \int e^{-i\Lambda\frac{|u|^2}{2}}a_e(u)du
 =\frac{c_{n_*}\Lambda^{-\frac{n_*}{2}}}{(2\pi)^{n_*}}
       \int \widehat a_e(\xi)e^{i\frac{|\xi|^2}{2\Lambda}}d\xi.
$$
Taylor expansion through order $J-1$ has remainder at
most
$$
 C_J\Lambda^{-J}
       \int |\xi|^{2J}|\widehat a_e(\xi)|d\xi
 \le C_J\Lambda^{-J}A^{2J+d+2}.
$$
The last inequality follows by integrating by parts in the transform to an
order greater than $2J+n_*$, all the supports having fixed diameter. For
the nonstationary portion take $M=J+\lceil \frac{n_*}{2}\rceil$
integrations by parts; after extracting $\Lambda^{-\frac{n_*}{2}}$ its
bound is $C_M\Lambda^{\frac{n_*}{2}-M}A^M\le C_J\Lambda^{-J}A^{2J+d+2}$.
This proves \eqref{osc:direct}, with differential coefficients of order at
most $2j$ and leading coefficient $c_{n_*}b(e)$.

We identify these coefficients without choosing moving endpoint
coordinates. Fix for this purpose a smooth $b$ and use a polar cutoff
$\chi(w\cdot e)$ equal to one near $w=e$ and zero outside a slightly larger
fixed cap, and write $\Delta_\xi=\sum_{j=1}^d\partial_{\xi_j}^2$ for the
Euclidean Laplacian. The Fourier extension
$F(\xi)=\int e^{iw\cdot\xi}b(w)dw$ satisfies $(\Delta_\xi+1)F=0$, where by
the outgoing part we mean the contribution near the stationary maximum,
with phase $e^{i\Lambda}$. The same identity for its polar-cutoff outgoing
part has an error decreasing faster than every power of $|\xi|$, since
derivatives of the cutoff are supported away from both stationary points,
so that the preceding integration by parts applies, and does so after any
fixed number of derivatives as well. In polar coordinates the differential
operator is
$$
 \partial_\Lambda^2+\frac{n_*}{\Lambda}\partial_\Lambda
       +1+\Lambda^{-2}\Delta_{S^{n_*}}.
$$
The Gaussian expansion may be differentiated any fixed number of times in
its angular and radial parameters: differentiate the compactly supported
coordinate amplitudes and the Gaussian formula, and take more terms in
order to bound the differentiated remainder. We may therefore substitute
the local outgoing expansion
$e^{i\Lambda}\Lambda^{-\frac{n_*}{2}}
\sum_{j\ge0}\Lambda^{-j}a_j(e)$, which is an asymptotic expansion, and
we compare any fixed finite number of powers after choosing a
sufficiently small remainder.
For a smooth function $H(\Lambda,e)$, direct differentiation gives
$$
\begin{aligned}
 &\left(\partial_\Lambda^2+\frac{n_*}{\Lambda}\partial_\Lambda
           +1+\Lambda^{-2}\Delta_{S^{n_*}}\right)
     \left(e^{i\Lambda}\Lambda^{-\frac{n_*}{2}}H\right)\\
 &=e^{i\Lambda}\Lambda^{-\frac{n_*}{2}}
     \left(\partial_\Lambda^2H+2i\partial_\Lambda H
       +\Lambda^{-2}\left(\Delta_{S^{n_*}}
                    -\frac{n_*(n_*-2)}{4}\right)H\right).
\end{aligned}
$$
Substitute the finite sum for $H$ and collect the coefficient of
$\Lambda^{-j-2}$.  Uniqueness of its power coefficients
gives
$$
 2i(j+1)a_{j+1}
  =\left(\Delta_{S^{n_*}}+j(j+1)
                     -\frac{n_*(n_*-2)}{4}\right)a_j.
$$
Starting from $a_0=c_{n_*}b$, this is the
recursion~\eqref{osc:L-operators}. A cutoff equal to one near $e$ has the
same coefficients, since the difference is nonstationary, and the minimum
has the conjugate Gaussian phase, which proves \eqref{osc:full-sphere}.

For inversion, apply \eqref{osc:direct} to $B_jb$, expanding its
series through index $J-j-1$.  The equation~\eqref{osc:B-operators} cancels
every coefficient of total index $1,\ldots,J-1$, and its index-zero
coefficient is $b$.  The derivative bound for $B_jb$ at order $m$ is
$C_{j,m}A^{m+2j}$.  After multiplication by $\Lambda^{-j}$ its extracted
remainder is therefore at most
$C_J\Lambda^{-J}A^{2J+d+2}$, uniformly in $0\le j<J$.
Under \eqref{osc:stationary-regime} this is at most
$$
 C_J\left(\frac{A^2}{\Lambda}\right)^J A^{d+2}
 \le C_J R^{-c\sigma J+2(d+2)}.
$$
Taking $J$ large proves all the asserted arbitrary-power estimates.

Finally, $e$ and $b$ are independent of $s$, and each $s$ derivative of
$\Lambda^{-j}(L_{j,\epsilon}b)(\epsilon e)$ costs only a fixed constant
times $\Lambda^{-j}A^{2j}$, which is bounded. The same elementary Taylor
remainder argument accommodates any prescribed finite number of radial
derivatives, by increasing the truncation order, and since
$\Lambda\lesssim R^2$ in the present applications, even radial derivatives
of the full error, its oscillatory exponential included, can be made an
arbitrary negative power by increasing that order once more. For the
support assertion, use
$$
 \Delta_{S^{d-1}}=\sum_{u<v}\Omega_{uv}^2.
$$
Leibniz expansion of any $L_j$ or $B_j$ differentiates each original mask
factor but never removes it, and derivatives of a smooth function cannot
enlarge its support, which proves the assertion, for arbitrary measurable
spatial multipliers as well.
\end{proof}

\subsection{Separation of the remaining spatial coefficient}

After stationary phase a smooth coefficient can still depend on both
endpoints, so that in order to apply separate Fourier transforms we must
express this coefficient as a sum of products of one-endpoint functions.
Rescaling the two cells to unit size makes their derivative bounds uniform,
and a Fourier series on a fixed box then has absolutely summable
coefficients after sufficiently many integrations by parts, as in
\cite[lemma~7.6]{OpenAI2026}.

\begin{lemma}[Separated spatial expansion]\label{osc:spatial-separation}
On separated depth-$a$ cells, for every fixed $j\ge0$ there are bounded
measurable functions $f_\ell$ on $X$ and $g_\ell$ on $Y$ and constants
$c_\ell$ such that
\begin{equation}\label{osc:spatial-series}
(D(x,y)2^a)^{-j}
   =\sum_\ell c_\ell f_\ell(x)g_\ell(y)
   \text{ with } |f_\ell|,|g_\ell|\le1 \text{ and } \sum_\ell|c_\ell|\le C_{d,j}.
\end{equation}
The convergence is uniform.  The same assertion holds in the fixed
outer geometry, with a constant depending on that geometry.
\end{lemma}

\begin{proof}
Write $x=c_X+2^{-a}u$ and $y=c_Y+2^{-a}z$, and therefore both rescaled endpoint
cells lie in fixed boxes. By \eqref{osc:separation}, $|2^a(c_X-c_Y)+u-z|$
is bounded above and below by positive constants on a fixed neighborhood of
their product, so that its negative $j$th power has a smooth extension,
with uniform bounds of every fixed order, supported in a slightly larger
product box. Regard this extension as a smooth function on a fixed torus of
dimension $2d$, obtained by identifying opposite faces of a box larger than
the support, and hence the extension is smooth and periodic. Its Fourier
series, indexed by $k,l\in{\mathbb Z}^d$, has coefficients bounded by
$C_m(1+|k|+|l|)^{-m}$, and choosing $m>2d$ makes their absolute sum
bounded. Each Fourier exponential splits into one function of $u$ and one
of $z$, each of modulus one, and pulling them back proves the claim. The
proof on the fixed separated outer supports is identical.
\end{proof}

\subsection{Proof of the weighted shell estimate}

\begin{proof}[Proof of proposition~\ref{osc:shell-conversion}]
Put $q=v-a$, $s=\frac{r}{2^v}$, and
$D_{\max}=\max(1,2^L R^{-50\sigma})$.  On the shell in question,
$\Lambda=rD(x,y)\asymp2^q$ and
$$
 D_*=R^\sigma D_{\max}
 \le\max\left(R^\sigma,2^{\frac{q}{2}}R^{-48\sigma}\right)
 \text{ and }
 \frac{D_*^2}{\Lambda}
 \lesssim\max(R^{2\sigma}2^{-q},R^{-96\sigma})
 \lesssim R^{-18\sigma}.
$$
Moreover $\Lambda\gtrsim R^{20\sigma}$ and $D_*\le R^2$, and the fixed
geometric constants can be absorbed by decreasing these positive exponents
slightly, so that lemma~\ref{osc:stationary} applies uniformly to the
endpoint product $m_x(w)n_y(w)$.

We estimate first $\widehat\nu(-r)$, which has the outgoing phase
$e^{irD}$, and this suffices for $\widehat\nu(r)$ as well, since
conjugating the spatial multipliers conjugates $\nu$ and
$\widehat{\overline\nu}(-r)=\overline{\widehat\nu(r)}$. Use
\eqref{osc:inverse} pointwise for each $(x,y)$ and multiply by the weight
in \eqref{osc:weighted-measure}; the leading distance power then cancels
exactly, as follows.
\begin{equation}\label{osc:weight-cancellation}
(D2^a)^{-\frac{n_*}{2}}(rD)^{\frac{n_*}{2}}(rD)^{-j}
 = (r2^{-a})^{\frac{n_*}{2}}(r2^{-a})^{-j}(D2^a)^{-j}.
\end{equation}
The operator $B_j$ is global, and its Leibniz expansion on $m_xn_y$ is a
sum of products of angular derivatives of the individual masks, of total
order at most $2j$. Divide every differentiated mask by its own scale
$D_{l_{\mathfrak t}}$ to the corresponding order; the resulting factors are
admissible, they retain every descriptor, and their extracted scales have
product at most $D_{\max}^{2j}$, while their number is polynomial in $N$
because $j$ is fixed. Smooth angular factors, if any are produced, may be
included in either endpoint symbol.

The inverse factor pays for all these scales, since
\begin{equation}\label{osc:scale-payment}
2^{-jq}D_{\max}^{2j}
 \le[\max(2^{-q},R^{-98\sigma})]^j\le1.
\end{equation}
Apply lemma~\ref{osc:spatial-separation} to the remaining factor
$(D2^a)^{-j}$.  Its absolutely summable expansion only multiplies the
original arbitrary spatial factors by bounded endpoint functions.
After summing the finitely many inverse orders, we obtain
\begin{equation}\label{osc:shell-factorization}
\widehat\nu(-2^vs)
  =2^{n_*\frac{q}{2}}\sum_\alpha c_\alpha(s)
       \int\psi(w)U_{X,\alpha}(2^vs,w)
                       \overline{V_{Y,\alpha}(2^vs,w)}dw
       +\mathcal E(2^vs),
\end{equation}
where
$$
 \sum_\alpha\|c_\alpha\|_{L^\infty([\frac{1}{2},2])}
 \le R^{o(1)}.
$$
The factors $s^{\frac{n_*}{2}-j}$ are included in $c_\alpha(s)$ and are
uniformly bounded on this interval. The index set may be countable because
of the spatial Fourier series, whose absolute summability justifies all the
interchanges. All the endpoint transforms in
\eqref{osc:shell-factorization} have the form in \eqref{osc:transforms},
with the stated enlarged fixed budget.

The pointwise error $\mathcal E$ is smaller than any prescribed negative
power of $R$. Indeed, lemma~\ref{osc:stationary} supplies such an error
before the endpoint integrations, the rescaled distance weight is bounded,
and the normalized cell masses are
$\frac{\rho_j(X)}{M_j(a)}\le R^{O_d(\sigma)}$, so that increasing the
initial accuracy absorbs these factors and any subsequent square. The shell
length is at most $2^{N+1}=2R$, since $v\le N$, and is absorbed by the same
choice.

Apply the triangle inequality in $L^2([\frac{1}{2},2],ds)$ to \eqref{osc:shell-factorization}.  Since $\Psi\ge1$ there, the
square of each angular integral is bounded by its contribution to the
corresponding $\mathsf H_F$ energy.  The complete prefactor is
\begin{equation}\label{osc:normalization-identity}
2^{-a}2^v2^{n_*(v-a)}=2^{d(v-a)}.
\end{equation}
This is the normalization in \eqref{osc:HF}, which proves \eqref{osc:shell-estimate}, and complex conjugation gives the
positive-frequency formulation stated in the proposition.
\end{proof}

\section{The graph estimate}\label{graph:section}

After passing to angular energies we split each measure into smaller cells,
which at first produces a sum over four cell labels. Oscillation removes
most of the combinations, and the parallel tests show that each retained
combination has only a controlled number of partners; thus the sum can
be bounded by Cauchy--Schwarz. This is the source of the geometric cost
which the descent potential will have to pay.

The independent angular energy uses a common radial frequency and two
independent directions, and we now pass from a pair of depth-$a$ cells to
depth-$p$ cells. After the endpoint transforms have been split into cell
pieces, angular and radial oscillation restrict the four cell labels to a
graph of small row and column degrees, the degree bound being supplied by
the parallel tests \cite[proposition~8.1 and lemmas~8.2--8.3]{OpenAI2026}.

The cell terms remain oscillatory until the nonstationary terms have been
removed, after which we take absolute values of the retained terms and
bound them by positive independent angular energies on the smaller cells.
At the end of the section, Cauchy--Schwarz transfers this estimate to the
separated form $\mathsf H_F$ with the same endpoint symbols.

\begin{proposition}[Graph estimate for the independent angular energy]
\label{graph:estimate}
Let $\rho_1,\rho_2$ be two regularized class laws, possibly equal, and
write $M_j$, $A_{j,n}$, $I_j$ and $k_j$ for their model cell masses and associated
quantities from section~\ref{test:section}.  Let $X,Y$ be occupied
cells of depth $a$.  Put $X_1=X$, $X_2=Y$, and let
$b^{(j)}_z(w)$ be an admissible endpoint symbol on $X_j$ in the sense
of definition~\ref{osc:symbols}.  Let $L$ be the largest nominal
length of a descriptor in either symbol, with $L=0$ if both lists are
empty.  Suppose that the integer depths satisfy $0\le a<p\le v\le N$
and
\begin{equation}\label{graph:range}
p-a\le v-p+2\sigma N, L\le v-p+2\sigma N \text{ and } v-p\ge10\sigma N.
\end{equation}
At endpoint $j$, suppose that the symbol imposes the parallel test
$[g_j,p]$ at a level $\ell_j\in\{0,\ldots,J_0\}$, meaning that for
$\rho_j$-almost every $z\in X_j$ and every $w\in S^{d-1}$,
\begin{equation}\label{graph:support-test}
b^{(j)}_z(w)\ne0
  \Longrightarrow w\in G^{(j)}_{\ell_j,[g_j,p]}(z).
\end{equation}
The superscript $j$ indicates the class used in the test.  Assume
\begin{equation}\label{graph:anchor-range}
a\le g_j<p \text{ and } g_j\le a+2\sigma N.
\end{equation}
Form the general independent energy $\mathsf H^{X,Y}(a,v)$ of
definition~\ref{osc:states} using these specified symbols.  For every
fixed $T>0$,
\begin{equation}\label{graph:bound}
\begin{aligned}
 \mathsf H^{X,Y}(a,v)
 &\le
 R^{-\sum_{j=1}^2(A_{j,p}-A_{j,a})
       +\min_{j=1,2}k_j(a,p)+O_d(\sigma)+o(1)}\\
 & \cdot
 \max_{\substack{P\subset X,\ Q\subset Y\\
                   P,Q\text{ occupied depth-}p\text{ cells}}}
       \mathsf H^{P,Q}(p,v)
 +O(R^{-T}).
 \end{aligned}
\end{equation}
No separation of $X$ and $Y$ is required.  Every child in
\eqref{graph:bound} is formed by restricting the same endpoint symbols
to $P$ and $Q$.  It retains every original mask factor, including any
whose anchor is $p$.  Constants are uniform over cells and allowed
symbols with fixed list-size exponents and derivative budgets, and with
fixed common numerical bounds for their smooth angular factors at the
derivative orders used.  The $o(1)$ tends to zero as $N\to\infty$ with
these bounds, the parameters, and $T$ fixed.
\end{proposition}

We give the localization and degree arguments separately.  Set
\begin{equation}\label{graph:scales}
h=p-a, q=v-p, \delta=\min(c_0,2^{-h}R^{10\sigma}) \text{ and } \tau=2^{-p}R^{30\sigma},
\end{equation}
where $c_0>0$ is a sufficiently small dimensional constant. By a
nonnegative partition of unity we mean a family of nonnegative functions
whose sum is one, and for functions, bounded overlap refers to their
supports. Take a nonnegative smooth partition of unity $\{\eta_\nu\}$ on
the sphere, with bounded overlap, such that $\eta_\nu$ is supported in a
cap of radius $C\delta$ centered at $w_\nu$, and
\begin{equation}\label{graph:bin-derivatives}
\|\nabla^m\eta_\nu\|_\infty\le C_m\delta^{-m}.
\end{equation}
Such a partition is obtained by normalizing scaled bumps on a finite atlas,
the number of caps being $O_d(\delta^{-(d-1)})$, and no relation between
the two cap centers used below is imposed.

For occupied depth-$p$ cells, use the transforms $U_P,V_Q$ of \eqref{osc:transforms}, divided by the model masses $M_1(p)$
and $M_2(p)$, respectively.  The exact splitting identities are
\begin{equation}\label{graph:splitting}
U_X=\lambda_1\sum_{P\subset X}U_P \text{ and } V_Y=\lambda_2\sum_{Q\subset Y}V_Q,
  \text{ where } \lambda_j=\frac{M_j(p)}{M_j(a)}=R^{-I_j(a,p)}.
\end{equation}
Expand the two squares in $\mathsf H^{X,Y}(a,v)$ and insert the angular
partition in both variables.  Index each term by a row $(P,Q)$ and a
column $(P',Q')$, with
$x\in P$, $x'\in P'$, $y\in Q$, and $y'\in Q'$.  The $P$-labels belong
to $\rho_1$ and the $Q$-labels to $\rho_2$.  Put
$$
 \Delta_1=x-x', \Delta_2=y-y', d_1=c_P-c_{P'} \text{ and } d_2=c_Q-c_{Q'}.
$$
With the spatial points fixed, the term has phase
\begin{equation}\label{graph:phase}
r(w_1\cdot\Delta_1+w_2\cdot\Delta_2).
\end{equation}
Each displacement joins two points in the same depth-$a$ parent cell.
Consequently
\begin{equation}\label{graph:displacement-size}
|\Delta_j|+|d_j|\lesssim_d2^{-a}.
\end{equation}
The locations of $X$ and $Y$ relative to one another do not enter this
bound.

For an endpoint cell $P$ with measure $\rho$ and symbol $b_x$, define
$$
 E(P,\nu)=\int_P\int_{S^{d-1}}
                  \eta_\nu(w)|b_x(w)|^2dwd\rho(x).
$$
Here symbols are set to zero on any fixed spatial null set where their
stated support conditions have not been specified, and the cell--cap
incidence $(P,\nu)$ is eligible precisely when $E(P,\nu)>0$. If
$E(P,\nu)=0$, then Fubini's theorem shows that for almost every $w$ with
$\eta_\nu(w)>0$ the symbol vanishes for $\rho$-almost every $x\in P$, so
that the corresponding transform vanishes for every $r$ and all the
expanded terms involving this incidence have zero integral and may be
omitted. Positive $E(P,\nu)$ supplies a point--angle witness at which the
symbol is nonzero and the stated support conditions hold, and at each fixed
$N$ only finitely many witnesses are chosen, so that no measurable
selection and no lower bound for their masses is needed.

Fix two caps, with centers $w_1^0,w_2^0$.  Retain an edge precisely when
all its cell--cap incidences are eligible and
\begin{equation}\label{graph:retained-tests}
|P_{(w_j^0)^\perp}d_j|\le\tau
  \text{ for } j=1,2 \text{ and } |w_1^0\cdot d_1+w_2^0\cdot d_2|\le\tau.
\end{equation}
Eligibility and the edge decision concern an entire cell--cap term, and
never cut its angular support by a nonsmooth pointwise condition.

The cell decomposition produces two displacement vectors, and angular
derivatives of the phase detect their components perpendicular to the
corresponding cap directions, while a radial derivative detects the sum of
their scalar projections. If any of these quantities is large, integration
by parts in that variable makes the term negligible, and it is the
smoothing margin in the masks which allows this operation without losing
the good-direction conditions.

\begin{lemma}[Nonstationary localization]\label{graph:localization}
Under \eqref{graph:range}, replacing the expanded independent energy by its
sum over the retained edges changes it by $O(R^{-T})$, for every fixed
$T>0$.  This replacement is made in the integrated oscillatory form,
before taking absolute values of the expanded cell terms.
\end{lemma}

\begin{proof}
Fix the four spatial points in one cell term. Fubini's theorem is
applicable here, since the spatial measures are finite and all the
amplitudes are bounded with compact radial support. Replacing centers by
points, or cap centers by angles in their cap supports, changes each test
in \eqref{graph:retained-tests} by at most
\begin{equation}\label{graph:test-error}
E\lesssim_d2^{-p}+2^{-a}\delta
 \lesssim_d2^{-p}R^{10\sigma}=o(\tau).
\end{equation}
The second inequality also holds when $\delta=c_0$, since that case implies
$c_0\le2^{-h}R^{10\sigma}$, the fixed geometric constants being absorbed by
the power gap between $R^{10\sigma}$ and $R^{30\sigma}$.

Suppose first that the last inequality in
\eqref{graph:retained-tests} fails.  Throughout this entire term,
$$
 |w_1\cdot\Delta_1+w_2\cdot\Delta_2|\ge\frac{\tau}{2}
$$
for sufficiently large $N$. In the variable $u=\frac{r}{2^v}$ each radial
integration by parts gains a factor at most $C(2^v\tau)^{-1}$, and since
the symbols and the angular partition do not depend on $u$, only the fixed
smooth radial cutoff is differentiated. Since
\begin{equation}\label{graph:oscillation-size}
2^v\tau=2^qR^{30\sigma}\ge R^{40\sigma},
\end{equation}
repeating this operation gives arbitrary power decay.

Otherwise, one of the two perpendicular inequalities fails.  Suppose
it is the first.  Let $t$ be the unit tangent vector at $w_1^0$ in the
direction of $P_{(w_1^0)^\perp}d_1$, and extend it
by the smooth sphere vector field
$$
 V(w)=t-(t\cdot w)w.
$$
All derivatives of this vector field are bounded independently of the
cap.  By \eqref{graph:test-error},
$$
 |V(w_1)\cdot\Delta_1|\ge\frac{\tau}{2}
$$
throughout its support.  Applied to the phase
$\Phi=r(w_1\cdot\Delta_1+w_2\cdot\Delta_2)$, the field therefore
satisfies
$$
 |V\Phi|\gtrsim2^v\tau
 \text{ and }
  |V^m\Phi|\lesssim_m2^v2^{-a} \text{ for } m\ge1.
$$
The same reasoning applies in the second variable if its perpendicular
test fails.

We shall need the following derivative bounds in order to iterate this
integration by parts. By definition~\ref{osc:symbols},
lemma~\ref{test:product-symbols} and \eqref{graph:bin-derivatives}, the
product amplitude has angular derivatives of order $m$ bounded by
$C_m(1+D_*+\delta^{-1})^m$, where
$$
 D_*=R^\sigma\max(1,2^LR^{-50\sigma}).
$$
Differentiating $(V\Phi)^{-1}$ contributes, in addition to its initial
factor $(2^v\tau)^{-1}$, at most powers of $1+\frac{2^{-a}}{\tau}$, since
every positive-order derivative of $V\Phi$ is bounded by $C_m2^v2^{-a}$
whereas $|V\Phi|\gtrsim2^v\tau$, while the bounded divergence and
derivatives of $V$ cause only constant factors. Thus $m$ integrations by
parts give a bound by a constant $C_m$ times
\begin{equation}\label{graph:ibp-ratio}
\left(
 \frac{1+D_*+\delta^{-1}+\frac{2^{-a}}{\tau}}{2^v\tau}
 \right)^m
\end{equation}
for the angular--radial integral of the fixed spatial term.
To see what is differentiated, write $b$ for the amplitude in
one angular variable. Since $(iV\Phi)^{-1}V(e^{i\Phi})=e^{i\Phi}$,
integration by parts gives
$$
 \int e^{i\Phi}bdw
 =-\int e^{i\Phi}\operatorname{div}
          \left(\frac{bV}{iV\Phi}\right)dw.
$$
Each repetition therefore differentiates the old amplitude, the vector
field, or the reciprocal of $V\Phi$, and the bounds just proved apply to
all three. A smooth extension of the reciprocal to a slightly larger cap is
permitted, because the same lower bound holds there.

The three hypotheses in \eqref{graph:range} imply
$$
\begin{aligned}
 D_*&\le R^\sigma+2^qR^{-47\sigma},\\
 \delta^{-1}&\le c_0^{-1}+2^qR^{-8\sigma},\\
 \frac{2^{-a}}{\tau}&=2^hR^{-30\sigma}\le2^qR^{-28\sigma}.
\end{aligned}
$$
Consequently the numerator in \eqref{graph:ibp-ratio} is at most
$C2^qR^{4\sigma}$ and the ratio is at most $CR^{-26\sigma}$, and therefore there
is a fixed positive power gain at each integration by parts, and taking
sufficiently many derivatives gives $O(R^{-T'})$ for any prescribed
$T'$.

Finally, all the sums and normalizations outside the fixed spatial term
cost at most a fixed power of $R$, since there are at most $R^{4d}$ cell
quadruples and $O(R^{2(d-1)})$ cap pairs, the mass curves lie between $0$
and $d$, and the energy prefactors are at most $R^d$, while the spatial
integrations contribute at most the same kind of fixed power through their
normalizations. Choose $T'$ larger than $T$ by these fixed powers before
summing. The required derivative order depends on the fixed parameters and
on $T$, but not on $N$, and since no spatial derivative has been taken,
bounded measurable spatial multipliers are allowed. This proves the
assertion.
\end{proof}

The retained terms form a bipartite graph with the same cell labels on
both sides. For fixed caps, interchanging $(P,Q)$ and $(P',Q')$ negates
both displacements in \eqref{graph:retained-tests}, while the absolute
value tests and the four cell--cap eligibility conditions are invariant
under this interchange, so that the adjacency matrix is symmetric even
when the two endpoint laws differ, and it suffices to bound the number of
neighbors of a fixed row.

\begin{lemma}[Both degrees of the retained graph]\label{graph:degrees}
For each fixed cap pair, both the row degrees and the column degrees
of the retained graph are at most
\begin{equation}\label{graph:degree-bound}
D_{\mathrm{gr}}=
 R^{\min_{j=1,2}k_j(a,p)+O_d(\sigma)+o(1)}.
\end{equation}
This estimate does not require the two cap directions to be separated.
\end{lemma}

\begin{proof}
Throughout this count, ``occupied'' refers to the full regular class and
not to a later restriction of its mass, the parallel test having been
defined using exactly these full-class centers. Restricting a symbol can
only remove vertices or edges, so that it cannot make either degree
larger.

Fix a row $(P,Q)$, whose unknown column members $P'$ and $Q'$ have mass
curves $1$ and $2$ respectively, and consider the unknown of mass curve
$j$. Its perpendicular constraint places its center within $\tau$ of a line
through its fixed mate in the direction $w_j^0$, so group its eligible
candidate cells by their depth-$g_j$ ancestor. Since the candidate cells
lie in one depth-$a$ parent, the number of these groups is at most
$$
 2^{d(g_j-a)}\le R^{2d\sigma}.
$$
In one nonempty group choose an eligible cell $P_*$ and a point--angle
witness $(z_*,w_*)$ from its incidence.  The support condition
\eqref{graph:support-test} gives its mass curve's parallel test
$[g_j,p]$ at $z_*\in P_*$. This test is based at the center $c_*$ of the
anchor cell $P_*$, and in applications where a normalized derivative of a
mask occurs, the same implication holds by
lemma~\ref{test:smooth-masks}.

For every other candidate center $c$ in this group, the two
perpendicular constraints and $|w_*-w_j^0|\lesssim\delta$ give
\begin{equation}\label{graph:witness-tube}
|c-c_*|\le\sqrt d2^{-g_j}
\text{ and }
\operatorname{dist}(c,c_*+\mathbb Rw_*)
 \lesssim_d2\tau+2^{-g_j}\delta
 \lesssim_d2^{-p}R^{30\sigma}.
\end{equation}
For large $N$ the latter is smaller than every test width
$W_i(p)\ge2^{-p}R^{500\sigma}$, while the first inequality lies inside the
testing ball $B(c_*,10(d+1)2^{-g_j})$. The parallel test therefore bounds
the number of candidates in this group by
$$
 R^{k_j(g_j,p)+C_d\sigma}
 \le R^{k_j(a,p)+C_d\sigma}.
$$
The last inequality holds because $[g_j,p]\subseteq[a,p]$.
Summing the groups bounds the number of possibilities for this first
unknown by $R^{k_j(a,p)+O_d(\sigma)}$.

Fix such a first unknown, so that its displacement, and with it its
contribution to the scalar constraint in \eqref{graph:retained-tests}, is
now fixed. That scalar constraint confines the longitudinal coordinate of
the remaining unknown to an interval of length $O(\tau)$ in its own cap
direction, while its perpendicular constraint confines its other $d-1$
coordinates to a ball of radius $O(\tau)$, and therefore together they confine
its center to a ball of radius $O_d(\tau)$, containing at most
$$
 C_d(1+2^p\tau)^d\lesssim_d R^{30d\sigma}
$$
depth-$p$ grid centers. Since we may begin with either mass curve $j$,
taking the better of the resulting two counts proves the asserted
row bound and hence, by symmetry, the column bound, the choice of
which endpoint to count first being what produces the minimum in
\eqref{graph:degree-bound}. Once that endpoint is fixed, the scalar test
controls the remaining coordinate of the other endpoint in its own cap
direction, which is also why the bound is uniform when the cap
directions approach one another.
\end{proof}

In the graph picture a row label $\alpha=(P,Q)$ and a column label
$\beta=(P',Q')$ each specify two cells, and at fixed $r,w_1,w_2$ we set
$F_\alpha=U_P(r,w_1)V_Q(r,w_2)$ and
$G_\beta=U_{P'}(r,w_1)V_{Q'}(r,w_2)$. The number $D_{\mathrm{gr}}$ is a
common upper bound for the row and the column degrees, that is, for the
number of edges incident to any one vertex.
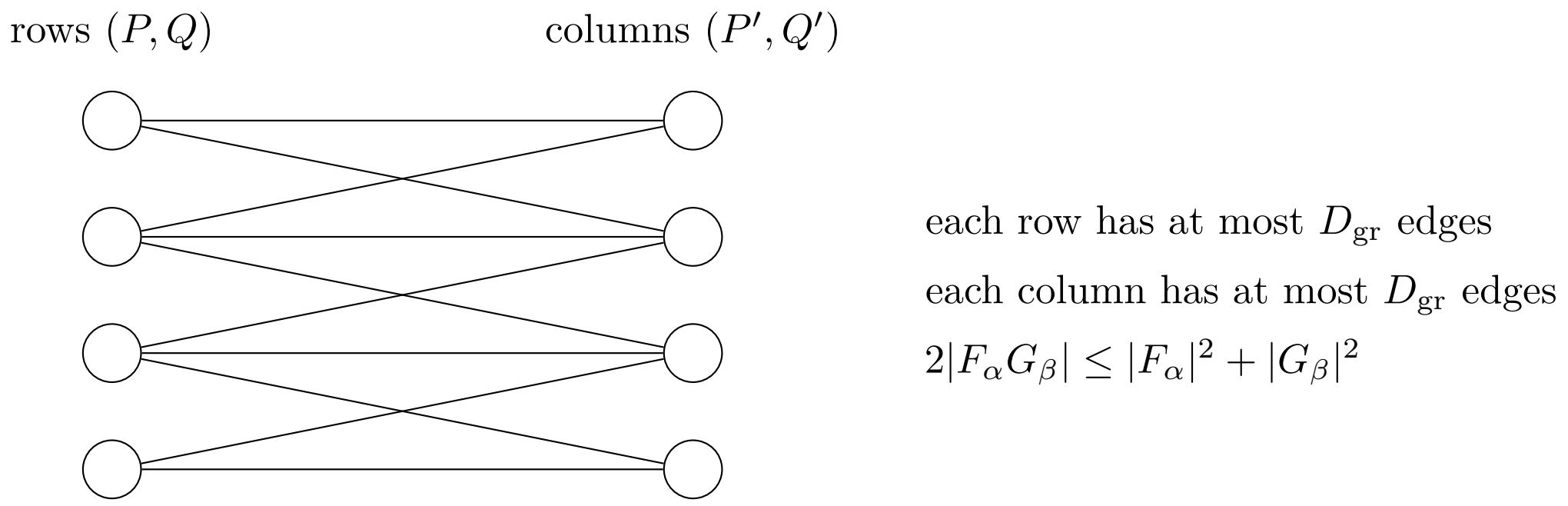
\begin{figure}[htbp]
\centering
\begin{tikzpicture}[>=Stealth,font=\small,
 vertex/.style={circle,draw,fill=white,minimum size=5mm}]
\foreach \i/\y in {1/0,2/-1,3/-2,4/-3}{
\node[vertex] (l\i) at (0,\y) {};
\node[vertex] (r\i) at (5,\y) {};}
\foreach \a/\b in {1/1,1/2,2/1,2/2,2/3,3/2,3/3,3/4,4/3,4/4}{\draw (l\a)--(r\b);}
\node[above] at (0,.45) {rows $(P,Q)$};
\node[above] at (5,.45) {columns $(P',Q')$};
\node[align=left,text width=6cm] at (10,-1.5)
{each row has at most $D_{\mathrm{gr}}$ edges\\[5pt]
each column has at most $D_{\mathrm{gr}}$ edges\\[5pt]
$2|F_\alpha G_\beta|\le|F_\alpha|^2+|G_\beta|^2$};
\end{tikzpicture}
\caption{The finite graph after nonstationary terms have been removed.
An edge records a four-cell combination satisfying the localization
conditions. The parallel tests bound the degree, and the symmetry
of the graph gives the same bound on each side. Summing the displayed
inequality over edges charges each row square and each column square
at most $D_{\mathrm{gr}}$ times. The resulting squares are exactly the
positive child energies, with all endpoint masks retained.}
\label{fig:graph-degrees}
\end{figure}

\begin{proof}[Proof of proposition~\ref{graph:estimate}]
Apply lemma~\ref{graph:localization}, and for fixed $r,w_1,w_2$ in a cap
pair take absolute values of the retained terms only. In the row--column
indexing above the amplitudes are
$$
 F_{P,Q}=U_P(r,w_1)V_Q(r,w_2) \text{ and } G_{P',Q'}=U_{P'}(r,w_1)V_{Q'}(r,w_2).
$$
The elementary inequality $2|FG|\le|F|^2+|G|^2$, together with both
degree bounds in lemma~\ref{graph:degrees}, gives the finite matrix form
of Schur's estimate \cite{Schur1911},
\begin{equation}\label{graph:schur}
\sum_{(\alpha,\beta)\text{ edge}}|F_\alpha G_\beta|
 \le\frac{D_{\mathrm{gr}}}{2}
       \left(\sum_\alpha|F_\alpha|^2+
             \sum_\beta|G_\beta|^2\right).
\end{equation}
The sums on the right may be enlarged to all the row and column labels.
Summing the nonnegative partition weights over cap pairs gives one, so that
no factor equal to the number of caps appears, and we integrate against the
same nonnegative radial cutoff
$\Psi(\frac{r}{2^v})\frac{dr}{2^v}$. Each row square then produces the
independent-angle energy of its child cell pair, and so does each column
square, whence
\eqref{graph:schur} involves only child quantities
$\mathsf H^{P,Q}(p,v)$ formed from the restricted symbols.

Let $I_j=I_j(a,p)$.  The occupied-cell estimate in
lemma~\ref{test:profile-bounds} bounds the number of row labels, and
also of column labels, by
$$
 R^{I_1+I_2+O_d(\sigma)+o(1)}.
$$
The splitting factors in \eqref{graph:splitting} appear squared,
giving $R^{-2(I_1+I_2)}$.  Finally, the parent prefactor divided by
the child prefactor is
$$
 \frac{2^{d(v-a)}}{2^{d(v-p)}}=2^{d(p-a)}.
$$
Since $S=\frac{d}{2}$ and
$I_j=A_{j,p}-A_{j,a}+S\frac{p-a}{N}$, their product is
$$
\begin{aligned}
 R^{-2(I_1+I_2)}R^{I_1+I_2}2^{d(p-a)}
 &=R^{-(I_1+I_2)+d\frac{p-a}{N}}\\
 &=R^{-\sum_{j=1}^2(A_{j,p}-A_{j,a})}.
\end{aligned}
$$
Multiplying by \eqref{graph:degree-bound} proves
\eqref{graph:bound}, including the arbitrary-power error from
lemma~\ref{graph:localization}. Since the argument has only restricted
symbols to smaller cells, every child retains all its original mask factors
and support conditions.
\end{proof}

The graph estimate was proved for the independent angular energy, where
positivity is available, whereas for the form with a common angular
variable Cauchy--Schwarz first produces two such positive energies. The
following corollary records the resulting bound, without changing the
endpoint masks or their support information.

\begin{corollary}[The separated angular form]\label{graph:separated-form}
Under the hypotheses of proposition~\ref{graph:estimate}, suppose also
that $X,Y$ and their supported laws have the separated geometry used
to define $\mathsf H_F^{X,Y}(a,v)$ in section~\ref{osc:section}.  Then
\eqref{graph:bound} holds with its left-hand side replaced by
$\mathsf H_F^{X,Y}(a,v)$ and with the same child energies.  This includes
the outer root $a=0$ with separated actual supports inside possibly
equal depth-zero cells $X,Y$.
\end{corollary}

\begin{proof}
The estimate~\eqref{osc:angular-cs} gives
$\mathsf H_F^{X,Y}(a,v)\le\mathsf H^{X,Y}(a,v)$ with the same endpoint
symbols, so that their tests and model-mass normalizations are unchanged,
and proposition~\ref{graph:estimate} supplies the child energies in
\eqref{graph:bound}.
\end{proof}

\section{Reductions for one mass curve}\label{rec:section}

All the estimates are now available, and what remains is to specify how
they may be used repeatedly without losing their hypotheses. The two kinds
of quantity involved are a positive count of almost equal distances and a
positive angular Fourier energy, and the reduction table below says which
quantity can follow which, how the depth interval changes, and what the
change costs. The support conditions recorded there are part of the
estimates themselves, and a numerical bound alone would not justify the
next step.

We now use the general angular estimates to bound positive counts of pairs
with nearly equal distances. Fix one regularized probability $\rho$, with
mass curve $m$, and write $A_n,M(n),I(a,b),k(a,b)$ for its mass curve data.
Both endpoint laws in an angular state have this mass curve, though their
cells and bounded spatial multipliers may differ.

There are two kinds of recursive state. A positive state records a count of
pairs with nearly equal distances and uses the tests scheduled by its
current interval, whereas an angular state uses the independent energy of
definition~\ref{osc:states}, together with a retained base interval
specifying which tests must still be present, and we define these data
before stating the reductions which preserve them. By an angular episode we
mean a consecutive run of angular states, ending when the argument returns
to a positive count or stops. Low and high refer to the frequency ranges
specified in the reduction, whereas near and far refer to the sizes of the
point displacements after the cells have been subdivided.
For the diagram and the table, $\mathsf D_i(a,t)$ denotes the positive
distance-collision count at depths $a,t$ and
$\mathsf H_i^{\mathcal B}(a,v)$ the independent angular energy at depths
$a,v$, their full formulas following immediately below. The integer $i$ is
the smoothing level, the retained base $\mathcal B$ is the earlier interval
whose test list must still be imposed, and the depth $p$ is selected by the
fixed intermediate-depth rule, while a far term uses separated depth-$b$
cells. If $L$ is the largest nominal length of the remaining tests, put
$L_0=\max\{L,\lceil10\sigma N\rceil\}$ and $h=v-2L_0$; this $h$ is the new
start of a near term, and it is local to the reduction discussion, distinct
from the separation length $10q_0$ used in the potential.

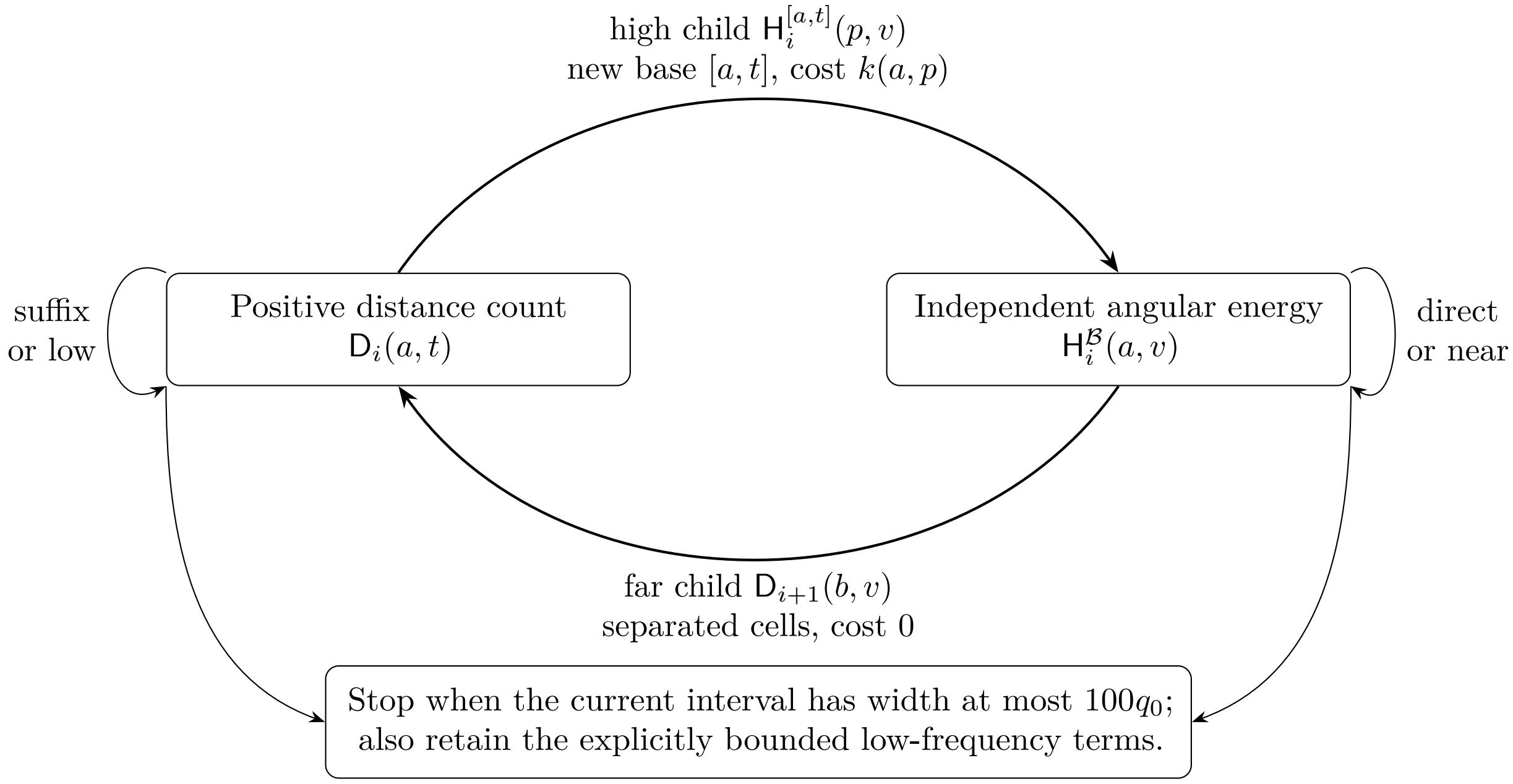
\begin{figure}[htbp]
\centering
\begin{tikzpicture}[>=Stealth,font=\small,
 box/.style={draw,rounded corners,align=center,text width=4.7cm,minimum height=1.2cm}]
\node[box] (D) at (0,0) {Positive distance count\\$\mathsf D_i(a,t)$};
\node[box] (H) at (7.7,0) {Independent angular energy\\$\mathsf H_i^{\mathcal B}(a,v)$};
\draw[->,thick] (D.north) to[out=55,in=125]
node[above,align=center] {high child $\mathsf H_i^{[a,t]}(p,v)$\\new base $[a,t]$, cost $k(a,p)$}(H.north);
\draw[->,thick] (H.south) to[out=235,in=305]
node[below,align=center] {far child $\mathsf D_{i+1}(b,v)$\\separated cells, cost $0$}(D.south);
\draw[->] (D.north west) .. controls (-3.3,1.0) and (-3.3,-1.0) ..
node[left,align=center,pos=.5] {suffix\\or low}(D.south west);
\draw[->] (H.north east) .. controls (10.8,1.0) and (10.8,-1.0) ..
node[right,align=center,pos=.5] {direct\\or near}(H.south east);
\node[box,text width=9cm] (stop) at (3.85,-4.2)
{Stop when the current interval has width at most $100q_0$;\\
also retain the explicitly bounded low-frequency terms.};
\draw[->] (D.south west) to[out=270,in=160] (stop.west);
\draw[->] (H.south east) to[out=270,in=20] (stop.east);
\end{tikzpicture}
\caption{The two quantities used in the descent. An angular energy
can use identical endpoint cells, whereas a positive distance count
requires separated cells, and a far reduction supplies that separation
before returning to a distance count. The retained base $\mathcal B$ of an
angular energy specifies its remaining tests, whereas its current
interval is the interval used in the potential; the table in
proposition~\ref{rec:reductions} gives the exact endpoints and costs.
Near and far terms are estimated separately, so that removing reached
masks from a near child does not remove them from the far term.}
\label{fig:reductions}
\end{figure}

\begin{table}[htbp]
\centering\small
\renewcommand{\arraystretch}{1.3}
\begin{tabular}{p{2.1cm}p{5.3cm}p{5.8cm}}
\toprule
Reduction & What makes progress & What must remain valid\\
\midrule
Suffix & End depth decreases by at least $q_0$. & The full conditional projection test transfers to the coarser count.\\
Low & End depth decreases by at least $\sigma N$. & The exact normalization $M(a)$ is unchanged; there is no $R^{C_d\sigma}$ loss.\\
High & Start advances to a fixed representative. & Converted symbols keep the old support tests; the new base is the old interval.\\
Direct & Start advances to a fixed representative. & The old base and index remain; reached factors are removed only from positive energies.\\
Near & Start advances to $h=v-2L_0$. & A passed anchor or the next direct step accounts for each nonterminal near step.\\
Far & The angular episode ends. & Separation and the new positive test list hold at depth $b$; the index increases.\\
\bottomrule
\end{tabular}
\caption{Why iteration is legal and finite. The potential pays
costs; the progress in this table bounds the number of steps,
including those with zero cost.}
\label{tab:progress}
\end{table}

\subsection{The two recursive states}

\begin{definition}[Positive distance count]\label{osc:positive-count}
Let $0\le a\le t\le N$ be integer depths, let
$i\in\{0,\ldots,J_0\}$, and let $X,Y$ be a separated pair of occupied
depth-$a$ cells of $\rho$, in the geometry of \eqref{osc:separation}.  Define the set of pairs satisfying the tests
$$
 \mathcal P_i^{X,Y}(a,t)=
 \{(x,y)\in X\times Y:
 e(x,y)\in G_{i,\mathfrak t}(x)\cap G_{i,\mathfrak t}(y)
 \text{ for every }\mathfrak t\in\mathcal M(a,t)\},
$$
and let
$$
 \theta_i^{X,Y}(a,t)
 =1_{\mathcal P_i^{X,Y}(a,t)}
       (\rho|_X\otimes\rho|_Y).
$$
This is an unnormalized pair restriction.  Put $C_0=100(d+1)$ and
\begin{equation}\label{osc:count}
\mathsf D_i^{X,Y}(a,t)
 =2^{t-a}M(a)^{-4}
   \int1_{\{|D(x,y)-D(x',y')|\le C_0 2^{-t}\}}
   d\theta_i^{X,Y}(a,t)(x,y)
   d\theta_i^{X,Y}(a,t)(x',y').
\end{equation}
A positive state consists of these depths, index and cells, with value
$\mathsf D_i^{X,Y}(a,t)$.  Its current interval is $[a,t]$; it has no
separate retained base.  We suppress the cells in $\mathsf D_i(a,t)$ when
they are fixed or are varied in a stated maximum.
\end{definition}

\begin{definition}[Angular energy with retained tests]\label{rec:angular-state}
Let $\mathcal B=[a_0,t_0]\subseteq[0,N]$ be an integer interval, called
the retained base, and let $[a,v]$ be the current interval, with
\begin{equation}\label{osc:base-contract}
0\le a_0\le a\le v\le t_0\le N
  \text{ and } t_0\le v+\sigma N.
\end{equation}
Write
$$
 \mathcal M_{>a}(\mathcal B)=
 \{\mathfrak t\in\mathcal M(a_0,t_0):
                 \text{the anchor of }\mathfrak t\text{ is deeper than }a\}.
$$
A recursive angular state consists of these intervals, an index
$i\in\{0,\ldots,J_0-1\}$, two occupied depth-$a$ cells $X,Y$ of
$\rho$, and two admissible endpoint symbols.  No separation of the cells
is required.  Each symbol's descriptor list consists of exactly one occurrence of each descriptor in
$\mathcal M_{>a}(\mathcal B)$, at the common index $i$.
Lemma~\ref{test:descriptor-properties} gives a common fixed polynomial
size bound for these lists throughout the recursive family.  The state also
specifies a fixed derivative budget and common numerical bounds for its
smooth angular factors, as in definition~\ref{osc:symbols}; these bounds
are independent of $N$.

The value of this state is the general energy already defined for those
data,
$$
 \mathsf H_i^{\mathcal B}(a,v)=\mathsf H^{X,Y}(a,v).
$$
The compact notation displays the current interval, base and index.  The
cells, symbols and fixed bounds remain part of the state data even when
they are not displayed.  We use a maximum only for finite choices of cells,
and a supremum whenever admissible symbols vary.  Neither operation varies
the displayed base or index, or the fixed bounds, unless this is stated.
\end{definition}

The origin of a descriptor is not a further factor in a state. For each
membership $\mathfrak t\in\mathcal M_{>a}(\mathcal B)$,
lemma~\ref{test:descriptor-properties} supplies an originating interval
inside this particular base, and if the descriptor has anchor $p$ and
length $L$, that origin contains $[p-L,p+L]$ and supplies the parallel
descriptors $[g,p]$ for $p-L\le g<p$. An origin may be chosen when it is
used; it need not be the current interval, and a descriptor listed only
once need not have one origin valid for every list containing it.

The general graph estimate returns energies with the restricted old
symbols, the reached factors included. Suppose that a one-mass-curve graph
child at depth $a'$ carries the list $\mathcal M_{>a}(\mathcal B)$ from a
prior start $a\le a'$. The removal operation after \eqref{osc:angular-cs}
deletes its factors anchored at depths at most $a'$ at bounded cost,
leaving exactly $\mathcal M_{>a'}(\mathcal B)$, and once the interval
inequalities in \eqref{osc:base-contract} hold as well, the resulting
general energy is the value of a recursive angular state.

The identity responsible for the same normalization in every dimension is
\begin{equation}\label{rec:normalization}
2^{d(b-a)}\left(\frac{M(b)}{M(a)}\right)^2
 =2^{d(b-a)}R^{-2I(a,b)}
 =R^{-2(A_b-A_a)}.
\end{equation}
The last equality is where $S=\frac{d}{2}$ is used.

\subsection{Positive collision estimates}

For a positive measure on the line, a count of close pairs is comparable to
a sum of squares of masses in short intervals. Shifting one copy changes
which nearby intervals interact, but Cauchy--Schwarz still bounds the
result by the unshifted sum of squares, and a decaying kernel is handled by
summing these shifted comparisons. Positivity is essential here, and the
argument is not available for an arbitrary signed oscillatory integral.

\begin{lemma}[Binning positive measures]\label{rec:binning}
Let $\eta,\zeta$ be finite positive measures on the line and let
$\varepsilon>0$.  Write
$$
 B_\varepsilon(\eta)
 =\sum_{j\in\mathbb Z}\eta([j\varepsilon,(j+1)\varepsilon))^2.
$$
For any fixed $C\ge1$,
\begin{equation}\label{rec:bins-lower}
B_\varepsilon(\eta)
 \le(\eta\otimes\eta)\{|s-s'|\le\varepsilon\},
\end{equation}
\begin{equation}\label{rec:bins-shift}
\sup_{u\in\mathbb R}
 (\eta\otimes\eta)\{|s-s'-u|\le C\varepsilon\}
 \le C' B_\varepsilon(\eta)
\end{equation}
and
\begin{equation}\label{rec:bins-kernel}
\iint\left(1+\frac{|s-s'|}{\varepsilon}\right)^{-2}d\eta(s)d\zeta(s')
 \le C' B_\varepsilon(\eta)^{\frac{1}{2}}
                    B_\varepsilon(\zeta)^{\frac{1}{2}}.
\end{equation}
Moreover,
\begin{equation}\label{rec:collision-fourier}
(\eta\otimes\eta)\{|s-s'|\le C\varepsilon\}
 \le C'\varepsilon
       \int_{|r|\le\varepsilon^{-1}}|\widehat\eta(r)|^2dr.
\end{equation}
\end{lemma}

\begin{proof}
Partition the line into intervals at the collision scale. If two points are
within that scale, then their intervals have indices differing by a bounded
amount, and this finite overlap proves the unshifted comparison. The other
assertions follow by the same comparison together with Cauchy--Schwarz, and
then by Fourier inversion for the smooth kernel.

Put $b_j=\eta([j\varepsilon,(j+1)\varepsilon))$ and define $c_j$
using $\zeta$.  Cauchy--Schwarz gives, for every integer $h$,
$$
 \sum_jb_jb_{j+h}\le\sum_jb_j^2
 \text{ and }
  \sum_jb_jc_{j+h}\le
       \left(\sum_jb_j^2\sum_jc_j^2\right)^{\frac{1}{2}}.
$$
An interval of fixed length $2C\varepsilon$, with any center, permits
only a bounded number of bin offsets.  This proves
\eqref{rec:bins-lower}--\eqref{rec:bins-shift}.  The kernel in
\eqref{rec:bins-kernel} is bounded on the bin pair with offset $h$ by
$C(1+|h|)^{-2}$, whose sum is finite.

For the last assertion, choose a fixed sufficiently small $b>0$ and
put $K(t)=(\frac{\sin(bt)}{bt})^2$, with its continuous value at zero.
Then $K\ge c>0$ on $[-C,C]$, and its nonnegative Fourier transform
is bounded and supported in $[-1,1]$.  Fourier inversion and Fubini,
applied to $K(\frac{s-s'}{\varepsilon})$, prove
\eqref{rec:collision-fourier}.  These arguments require no density for
either measure.
\end{proof}

\subsection{The reduction mechanism}

A positive state uses the depth $p$ selected by the fixed choice rule for
its current interval $[a,t]$. At a suffix choice, the projection test for
the full conditional law permits a direct passage to the coarser distance
count at depth $p$, whereas at a prefix choice a positive Fourier estimate
splits the count into frequency shells, the low shells giving shorter
positive states while the high shells use the weighted shell conversion and
the graph estimate to enter an angular state.

For an angular state, let $L$ be the largest length in its remaining
list, with $L=0$ for an empty list, and put
\begin{equation}\label{rec:angular-split-depth}
L_0=\max\{L,\lceil10\sigma N\rceil\} \text{ and } h=v-2L_0.
\end{equation}
This compares the longest remaining test with half the current frequency
range. If $h\le a$, an origin of a longest descriptor supplies a parallel
test for a graph step, while if $h>a$, splitting at depth $h$ leaves enough
frequency range for stationary phase on the far terms, the near terms
having only boundedly many neighboring cells. An angular episode therefore
continues through direct or near steps and returns to a positive state
through a far step.

The current interval is $[a,t]$ for $\mathsf D_i(a,t)$ and $[a,v]$ for
$\mathsf H_i^{\mathcal B}(a,v)$, and terminality refers to this interval
and not to the retained base, the state being terminal exactly when its
width is at most $100q_0$.

The reductions below follow \cite[proposition~9.4]{OpenAI2026}, with
the retained descriptors and the transition costs made explicit.

\begin{proposition}[State reductions]\label{rec:reductions}
Suppose the current state is nonterminal, its smoothing index satisfies
$i<J_0$, and $N$ is sufficiently large with all parameters fixed.  For
every prescribed $T>0$, it is bounded, apart from an error $O(R^{-T})$
from oscillatory operations, by a sum of child bounds whose number is
polynomial in $N$, and possibly a low-frequency stopping term.  For every
non-low child,
if the old start is $a$ and the new start is $a'$, its multiplier is
\begin{equation}\label{rec:transition-factor}
R^{-2(A_{a'}-A_a)+\mathrm{cost}+C_d\sigma+o(1)}.
\end{equation}
A low child has cost zero and the entire multiplier $R^{o(1)}$, with no
fixed $C_d\sigma$ loss.  A low-frequency stopping term is
$R^{C_d\sigma+o(1)}$.  The legal children are
\begin{center}
\renewcommand{\arraystretch}{1.2}
\begin{tabular}{lll}
Reduction & Parent and child & Cost \\ \hline
Suffix & $\mathsf D_i(a,t)\to\mathsf D_{i+1}(a,p)$
 & $A_p-\min_{[p,t]}A$ \\
Low & $\mathsf D_i(a,t)\to\mathsf D_{i+1}(a,v)$
 & $0$; full factor $R^{o(1)}$ \\
High & $\mathsf D_i(a,t)\to\mathsf H_i^{[a,t]}(p,v)$
 & $k(a,p)$ \\
Direct & $\mathsf H_i^{\mathcal B}(a,v)\to\mathsf H_i^{\mathcal B}(p,v)$
 & $k(a,p)$ \\
Near & $\mathsf H_i^{\mathcal B}(a,v)\to\mathsf H_i^{\mathcal B}(h,v)$
 & $0$ \\
Far & $\mathsf H_i^{\mathcal B}(a,v)\to\mathsf D_{i+1}(b,v)$
 & $0$
\end{tabular}
\end{center}

For a positive parent, $p$ is the depth selected by the fixed choice
rule for its current interval, and the suffix row is used when
$2p>a+t$. When $2p\le a+t$, the shell index is an integer $a<v\le t$,
the low row having $t-v\ge\sigma N$ and the high row
$0\le t-v<\sigma N$. A high child has base $[a,t]$ and carries
the fixed enlarged derivative budget and common angular bounds produced
by the inverse conversion; these depend on $T$ and the fixed parameters,
not on $N$.

For an angular parent, use $L_0,h$ from \eqref{rec:angular-split-depth}.  If $h\le a$, the direct row
chooses a remaining descriptor of maximal length $L$, its anchor $p$,
and an origin $[b_0,u_0]\subseteq\mathcal B$ for that membership.  This
choice has $L=L_0$, $a<p<v$,
$$
 b_0\le p-L, p+L\le u_0 \text{ and } v-a\le2L,
$$
and it supplies the parallel descriptor $[\max(a,p-L),p]$ used by the
graph estimate.  If $h>a$, the near row uses the displayed depth $h$.
The far rows range over integer depths $b$ arising from the comparable
far displacements, with $a<b\le h+O_d(1)<v$, and use separated depth-$b$
cell pairs.  Direct and near children keep the same base, index,
derivative budget and common angular bounds.  A far child has the fresh
positive list $\mathcal M(b,v)$.

Finite maxima over child cells use the stated geometry.  Whenever admissible
symbols vary, take a supremum over the symbol family just specified.  An
empty collection of child terms contributes zero.  Constants outside powers
of $R$ may depend on
the fixed parameters, derivative budgets and angular bounds.  The
constants $C_d$ in exponent losses depend only on $d$.  A terminal state
is bounded by
\begin{equation}\label{rec:terminal}
R^{C_d(\kappa+\sigma)+o(1)}.
\end{equation}
\end{proposition}

\begin{proof}
We prove every reduction, including the support conditions needed for
its iteration.

Consider first the suffix reduction, and let $\theta=\theta_i^{X,Y}(a,t)$.  Retain a pair $(P,Q)$ of
depth-$p$ subcells of the separated parent cells whenever
$\mathcal P_i^{X,Y}(a,t)\cap(P\times Q)\ne\varnothing$, and choose
one pair in that intersection as a witness.  Put $e_0=e(c_P,c_Q)$.  Separation of the parent cells gives
the uniform Hessian bound $\|\nabla^2D\|\lesssim2^a$, and hence
\begin{equation}\label{rec:suffix-taylor}
D(x,y)=e_0\cdot(x-y)+O_d(2^{a-2p})
       =e_0\cdot(x-y)+O_d(2^{-t})
        \text{ for } x\in P \text{ and } y\in Q.
\end{equation}
The last equality follows from the suffix inequality $2p>a+t$, and there is
no missing constant term in this Taylor expansion, since at the cell
centers $D(c_P,c_Q)=e_0\cdot(c_P-c_Q)$, the linear term is
$e_0\cdot((x-c_P)-(y-c_Q))$, and their sum is $e_0\cdot(x-y)$. The witness
angle differs from $e_0$ by $O_d(2^{a-p})$, and hence for $y,y'\in Q$,
changing between these two angles changes the projection of $y-y'$ by
$O_d(2^{a-2p})$.

This passage from distances to projections is also used in
\cite[proposition~4.1 and its proof]{KeletiShmerkin2019}. The useful
scale range is determined by the same requirement that the quadratic
Taylor error fit below the observation scale. Here that requirement
is $2^{a-2p}\le2^{-t}$, and we apply it to collision counts.

The witness satisfies the projection test $[p,t]$. This test is formulated
for the full conditional probability $\frac{\rho|_Q}{\rho(Q)}$ rather than
for a witness-dependent restriction, and its width absorbs every fixed
multiple of $2^{-t}$, so that the conditional projected measure in the
direction $e_0$ has self-collision mass at tolerance $C2^{-t}$ bounded
by
\begin{equation}\label{rec:suffix-B}
B=2^{-(t-p)}R^{A_p-\min_{p\le n\le t}A_n+C_d\sigma}.
\end{equation}
By \eqref{rec:bins-shift}, the same bound, up to a fixed constant,
holds uniformly for any translated collision window.  Fixing
$x,x'\in P$ and using \eqref{rec:suffix-taylor} therefore gives
$$
 B_{2^{-t}}(D_\#(\rho|_P\otimes\rho|_Q))
 \lesssim B\rho(P)^2\rho(Q)^2.
$$

Group the retained cell pairs by the bin of width $2^{-p}$ containing their
center distance, and write $W_G=\sum_{(P,Q)\in G}\rho(P)\rho(Q)$. For a
retained pair $(P,Q)$, let $v_{P,Q}$ be its vector of masses in the fine
distance bins, so that the preceding bound says
$\|v_{P,Q}\|_{\ell^2}\le C B^{\frac{1}{2}}\rho(P)\rho(Q)$. For a coarse
group $G$ the triangle inequality therefore gives
$$
\left\|\sum_{(P,Q)\in G}v_{P,Q}\right\|_{\ell^2}
\le C B^{\frac{1}{2}}\sum_{(P,Q)\in G}\rho(P)\rho(Q)
=C B^{\frac{1}{2}}W_G.
$$
Squaring gives $CBW_G^2$, an estimate which uses the sum of the pair masses
and does not introduce the number of pairs in $G$. The distances of each
group occupy a fixed enlargement of its coarse bin, and these enlargements
have bounded overlap; thus the fine collision mass of the parent is at
most $CB\sum_GW_G^2$, by lemma~\ref{rec:binning}.

Every full pair in a retained $P\times Q$ satisfies the child list
$\mathcal M(a,p)$ at level $i+1$. Indeed, its descriptors belong to the
parent list, have anchors no deeper than $p$, and have nominal lengths at
most $\frac{p-a}{2}$, whereas the spatial anchor cells of the full pair and
of the witness agree and the change of direction is $O_d(2^{a-p})$, which
lies inside the transfer neighborhood of lemma~\ref{test:smooth-masks} for
every such descriptor. Two full pairs in one coarse group have distance
difference at most $(1+4\sqrt d)2^{-p}<C_0 2^{-p}$, so that the definition
of the positive child count gives
$$
 \sum_GW_G^2
 \le 2^{-(p-a)}M(a)^4\mathsf D_{i+1}(a,p).
$$
Multiplying by the parent normalization, and substituting
\eqref{rec:suffix-B}, yields
$$
 2^{t-a}M(a)^{-4}B
 2^{-(p-a)}M(a)^4
 =R^{A_p-\min_{[p,t]}A+C_d\sigma}.
$$
This proves the suffix reduction with no change of start.

We turn to the prefix reduction at low frequencies. Let
$\theta=\theta_i^{X,Y}(a,t)$, and let $\zeta$ be the positive
distance pushforward of
$$
 M(a)^{-2}(D(x,y)2^a)^{-\frac{d-1}{2}}
 \prod_{\mathfrak t\in\mathcal M(a,t)}
 w_{i,\mathfrak t}(x,e(x,y))w_{i,\mathfrak t}(y,e(x,y))
 d\rho|_X(x)d\rho|_Y(y).
$$
On the good pairs defining $\theta$, every mask is one.  The
remaining distance weight is bounded above and below by positive
dimensional constants on the separated geometry.  Thus positivity and
\eqref{rec:collision-fourier} imply
\begin{equation}\label{rec:prefix-lowpass}
\mathsf D_i(a,t)
 \lesssim 2^{-a}\int_{|r|\le2^t}|\widehat\zeta(r)|^2dr.
\end{equation}
Since $\zeta(\mathbb R)\le R^{C_d\sigma}$, the frequencies below $2^a$
contribute $R^{C_d\sigma}$, which is a stopping term, and reality allows us
to treat only the positive frequencies above $2^a$, split into shells
$[2^{v-1},2^v]$ with $a<v\le t$.

Suppose $t-v\ge\sigma N$.  A smooth majorant of the shell has
inverse Fourier transform bounded in modulus by
$C2^v(1+2^v|s-s'|)^{-2}$.  Its contribution to
\eqref{rec:prefix-lowpass} is therefore at most
$$
 C2^{v-a}M(a)^{-4}
 \iint(1+2^v|D(x,y)-D(x',y')|)^{-2}
 d\theta'(x,y)d\theta'(x',y'),
$$
where $\theta'$ is the unweighted pair restriction to the supports of the
old masks. Every such pair satisfies $\mathcal M(a,v)$ at level $i+1$, by
list containment and mask support, and lemma~\ref{rec:binning} bounds the
displayed expression by $C\mathsf D_{i+1}(a,v)$. The normalization uses the
same exact $M(a)$ on both sides, and therefore this step incurs no
$R^{C_d\sigma}$ mass-comparison loss, while summing the shell indices costs
only $R^{o(1)}$.

For the prefix reduction at high frequencies, suppose instead that
$t-v<\sigma N$.  The largest nominal length in the
list is at most $\frac{t-a}{2}$.  Nonterminality gives
$$
 v-a>100q_0-\sigma N>20\sigma N \text{ and } L\le\frac{v-a}{2}+\sigma N.
$$
Thus proposition~\ref{osc:shell-conversion} applies to $\zeta$.
For the prefix choice $p$, the descriptor $[a,p]$ is present,
and
$$
\begin{aligned}
 p-a&\le t-p\le v-p+\sigma N,\\
 L&\le\frac{t-a}{2}\le t-p\le v-p+\sigma N,\\
 v-p&\ge\frac{t-a}{2}-\sigma N>10\sigma N.
\end{aligned}
$$
These verify the hypotheses of corollary~\ref{graph:separated-form} with
$g=a$, and its one-mass-curve specialization gives first general
independent energies on depth-$p$ cells, formed with the restricted
converted symbols and still carrying every old mask. In those positive
energies, remove the factors whose anchors are at most $p$; the remaining
list at each endpoint is then exactly $\mathcal M_{>p}([a,t])$, with
$a\le p\le v\le t\le v+\sigma N$. The resulting values are therefore
recursive states $\mathsf H_i^{[a,t]}(p,v)$, and the high contribution is
bounded by
$$
 R^{-2(A_p-A_a)+k(a,p)+C_d\sigma+o(1)}
       \sup\mathsf H_i^{[a,t]}(p,v)+O(R^{-T}).
$$
Conversion produced normalized angular derivatives of masks at index $i$,
with a fixed enlarged budget and common angular bounds, and their supports
still impose the existing level-$i+1$ tests, which are accepted by the
graph estimate. The new angular state therefore keeps the index $i$.

For the direct angular reduction, consider
$\mathsf H_i^{\mathcal B}(a,v)$ with base $\mathcal B=[a_0,t_0]$, whose
remaining anchors are strictly deeper than $a$, with
$t_0\le v+\sigma N$. If $h\le a$, nonterminality forces $L=L_0$, since
otherwise $v-a\le2\lceil10\sigma N\rceil<100q_0$. Choose a
maximal-length descriptor, with anchor $p>a$ and an origin $[b_0,u_0]$
in the base, so that
$$
 b_0\le p-L, p+L\le u_0\le v+\sigma N
  \text{ and } v-a\le2L.
$$
It follows that
$$
 L\le v-p+\sigma N, p-a\le v-p+2\sigma N \text{ and } v-p\ge\frac{v-a}{2}-\sigma N>10\sigma N.
$$
Set $g=\max(a,p-L)$.  Then
$$
 a\le g<p, p-L\le v+\sigma N-2L\le a+\sigma N
  \text{ and } g\le a+\sigma N.
$$
The same origin supplies the parallel descriptor $[g,p]$, whose anchor has
not yet been reached, so that proposition~\ref{graph:estimate} applies to
the general energy underlying $\mathsf H_i^{\mathcal B}(a,v)$, without any
separation of its two cells, and gives the direct reduction together with
the cost $k(a,p)$. The general graph child still carries the old symbols,
and we remove its reached factors only after forming its positive
independent energy; the remaining list is then exactly
$\mathcal M_{>p}(\mathcal B)$, which gives $\mathsf H_i^{\mathcal B}(p,v)$
with the same fixed symbol bounds.

We come next to the near angular reduction. If $h>a$, split each transform
into its depth-$h$ cell pieces, keeping for the moment all the remaining
factors, and expand the two squares in \eqref{osc:H}. A near term is one
for which both center displacements, within the first and the second parent
cell respectively, are at most $10(d+1)2^{-h}$. In the graph with row
$(P,Q)$ and column $(P',Q')$ these conditions give bounded row and column
degrees, since a fixed grid cell has only a dimensional number of neighbors
at this distance. Apply $2|uv|\le|u|^2+|v|^2$ to each row--column product
and integrate the resulting independent-angle squares. There are at most
$R^{2I(a,h)+C_d\sigma}$ rows or columns, and the four transform splitting
factors, the label count, and the ratio of energy prefactors give
$$
 2^{d(h-a)}R^{-4I(a,h)}R^{2I(a,h)+C_d\sigma}
 =R^{-2(A_h-A_a)+C_d\sigma},
$$
by \eqref{rec:normalization}, which is the zero-cost near estimate. Only
after each positive independent energy has been obtained do we remove the
factors anchored at depths at most $h$, which on its spatial cells are
angular factors of bounded modulus; the remaining list is then exactly
$\mathcal M_{>h}(\mathcal B)$, giving the child
$\mathsf H_i^{\mathcal B}(h,v)$ with the same fixed symbol bounds. This
removal does not apply to the far terms which follow.

We now prove the analytic estimate for the far reduction. In a remaining
term at least one point displacement has length
$\gtrsim_d2^{-h}$.  Fix the four spatial points and denote the two
displacements by $\Delta_1,\Delta_2$, with lengths $d_1,d_2$.
Put $d_+=\max(d_1,d_2)$ and $d_-=\min(d_1,d_2)$.
For the larger displacement, on the radial support,
$$
 \Lambda=r d_+\gtrsim 2^{v-h}=2^{2L_0}
 \text{ and }
 \frac{D_*^2}{\Lambda}
 \lesssim
 \begin{cases}
 R^{-18\sigma},&2^LR^{-50\sigma}\le1,\\
 R^{-98\sigma},&2^LR^{-50\sigma}>1.
 \end{cases}
$$
Also $\Lambda\ge R^{c\sigma}$ and $D_*\le R^2$ for sufficiently large $N$,
so that lemma~\ref{osc:stationary} expands the corresponding full-sphere
integral to any prescribed negative power accuracy. After the leading
stationary gain has been extracted, the coefficients and any fixed number
of derivatives in $\frac{r}{2^v}$ are bounded, and the factor
$(r d_+)^{-j}$ pays for the order-$2j$ symbol derivatives.

If $d_+>2d_-$, leave the other angular integral exact. The derivative in
$r$ of each resulting phase then has modulus at least
$d_+-d_->\frac{d_+}{2}$, uniformly in the unexpanded angle, and repeated
integration by parts in $\frac{r}{2^v}$ gives $O(R^{-T})$, including when
$d_-=0$. Leaving that integral exact is legitimate because its symbols do
not depend on $r$.

Otherwise the lengths are comparable and the second angular integral admits
the same expansion. Choose an integer $n$ with $d_+\in(2^{-n-1},2^{-n}]$
and set $b=n+\log_2K_d$, so that both lengths lie between
$\frac{K_d2^{-b}}{4}$ and $K_d2^{-b}$. Moving to their depth-$b$ cell
centers changes these bounds by at most $2\sqrt d2^{-b}$, and the choice of
$K_d$ puts both pairs in the exact separated geometry of
definition~\ref{osc:positive-count}. Since each displacement joins points
in one depth-$a$ cell, we have $d_+\le\sqrt d2^{-a}$ and hence $b>a$, while
the far lower bound gives $b\le h+O_d(1)<v$.

The distances in this part of the argument are the internal distances
$|x-x'|$ and $|y-y'|$, and not the original cross distance $|x-y|$, which
is legitimate because the recursive angular state uses two endpoint laws
with the same mass curve, and therefore both internal pairs can be estimated by
the same positive-count recursion once separation at depth $b$ has been
established. At the outer entry, where the endpoint mass curves may differ,
Cauchy--Schwarz first produces one state for each curve separately. The two
leading stationary gains multiply to $O_d(2^{-(d-1)(v-b)})$, and radial
integration, using the coefficient derivative bounds just established,
bounds each remaining kernel by
\begin{equation}\label{rec:far-kernel}
C2^{-(d-1)(v-b)}
       (1+2^v|d_1-d_2|)^{-2}.
\end{equation}
For opposite stationary signs the phase involves $d_1-d_2$, while for equal
signs it involves $d_1+d_2$, whose modulus is no smaller, which proves
\eqref{rec:far-kernel} for every choice of signs.

It remains to check the support and the normalization in the far reduction.
The interval $[b,v]$ is contained in the base, so that $\mathcal M(b,v)$ is
a sublist of its list, every anchor of this child list is strictly deeper
than $b>a$, and no such factor has been removed, in particular none at
depth $h$. The far expansion used the symbols before the removal made for
the near child, and since each coefficient differentiates the original
masks only angularly, it retains all of their supports. For example, the
amplitude for the internal pair $(x,x')$ is $b_x(w)\overline{b_{x'}(w)}$,
whose stationary coefficient is evaluated at $w=\pm\frac{x-x'}{|x-x'|}$, so
that a nonzero coefficient imposes the level-$i+1$ tests of
$\mathcal M(b,v)$ at both spatial endpoints $x$ and $x'$. The same holds
for the other internal pair $(y,y')$, antipodal symmetry of the tests
taking care of the negative stationary direction, while the expansion
remainders are already arbitrary-power errors and require no support
assertion.

The far terms were obtained by splitting the transforms at depth
$h$.  Once their four cell transforms are opened, the four splitting
factors cancel their model denominators exactly, giving
$$
 \left(\frac{M(h)}{M(a)}\right)^4M(h)^{-4}=M(a)^{-4}.
$$
For each comparable length band defining $b$, the stationary estimate is
now a nonnegative kernel bound with the child-support indicators, and on
this upper bound we may discard the indicator selecting the far $h$-cell
labels while keeping the comparable length band. Summing the disjoint
$h$-cell partition then reassembles the raw four-point measure on
$X^2\times Y^2$, which we partition into depth-$b$ cells, and consequently
no count of $h$-cell labels remains, whether $b$ is smaller or larger than
$h$.

For a fixed quadruple of separated depth-$b$ cells $(P,P',Q,Q')$,
put
$$
 \eta_1=D_\#\theta_{i+1}^{P,P'}(b,v) \text{ and } \eta_2=D_\#\theta_{i+1}^{Q,Q'}(b,v).
$$
These are the unnormalized positive distance measures for the two
internal pairs, each restricted by the child tests at both endpoints.
By lemma~\ref{rec:binning} and the exact normalization of
$\mathsf D_{i+1}$,
$$
 \iint(1+2^v|s-s'|)^{-2}d\eta_1(s)d\eta_2(s')
 \lesssim M(b)^4 2^{-(v-b)}
              \max\mathsf D_{i+1}(b,v).
$$
Indeed, each bin-square sum is bounded by its unshifted positive collision
count, and the geometric mean of the two counts is at most their maximum,
while once the first cell of either internal pair is fixed, its partner has
only a dimensional number of possible labels, since their center distance
is $O_d(K_d2^{-b})$. There are thus at most $R^{2I(a,b)+C_d\sigma}$
quadruple labels. The two stationary gains and the scalar collision
normalization combine to $2^{-d(v-b)}$, converting the parent frequency
prefactor into $2^{d(b-a)}$, while the four depth-$b$ masses and the two
free first-cell labels leave $R^{-2I(a,b)+C_d\sigma}$. Explicitly,
$$
\begin{aligned}
 &2^{d(v-a)}M(a)^{-4}
   2^{-(d-1)(v-b)}
   M(b)^4 2^{-(v-b)}R^{2I(a,b)+C_d\sigma}\\
 &
 =2^{d(b-a)}R^{-2I(a,b)+C_d\sigma}
 =R^{-2(A_b-A_a)+C_d\sigma}.
\end{aligned}
$$
Summing the possible scales $b$, the coefficient terms and the shell
indices introduces only $R^{o(1)}$, which proves the zero-cost far
reduction, with one increase of the smoothing index.

It remains to describe the stopping rule. For a terminal interval the
width is at most $100q_0$.
The definitions of the two state values, bounded symbol
moduli, and the occupied-cell bounds give respectively
$$
 \mathsf H_i^{\mathcal B}(a,v)\lesssim2^{d(v-a)}R^{C_d\sigma} \text{ and } \mathsf D_i(a,t)\lesssim2^{t-a}R^{C_d\sigma}.
$$
Since $\frac{q_0}{N}=\kappa+O(\frac{1}{N})$, both satisfy
\eqref{rec:terminal}, and the low-frequency stopping term was already
estimated after \eqref{rec:prefix-lowpass}.
\end{proof}

\section{Paying for every reduction with one potential}\label{ledger:section}

The reductions of proposition~\ref{rec:reductions} have three transitions
with potentially positive cost.  For a suffix step $[a,t]\to[a,p]$,
if a depth to the left of the choice rule's inner interval minimizes $A$
on both intervals, the minimum does not increase and cannot pay a
positive cost $A_p-\min_{[p,t]}A$.  We instead use a pair of well-separated
depths.  In a suffix child, separation places the later depth in the
inner interval, so it can be replaced by a minimizer on $[p,t]$
while the earlier depth remains.  A high prefix or direct reduction can
replace the earlier depth.

Write $h=10q_0$. The separation between the two depths is needed on the
long intervals where a reduction has positive cost. At a short interval we
can charge for the missing separation at the rate at which $A$ can
change in one step.

\begin{definition}[Two-depth potential]\label{ledger:potential}
For an integer interval $J=[a,t]$, let
\begin{equation}\label{ledger:potential-formula}
 \mathcal V_h(J)=\min_{a\le i\le j\le t}
       \left\{A_i+A_j+\frac{d}{N}(h-(j-i))_+\right\},
\end{equation}
where $u_+=\max(u,0)$.
\end{definition}

Suppose $t-a\ge h$ and a candidate pair has gap $g<h$. In $J$ there
is room to move the two indices outward by a total of $h-g$ unit steps.
At each step the sum of the two excesses increases by at most
$\frac{d}{N}$, by lemma~\ref{test:profile-bounds}, and the penalty
decreases by exactly $\frac{d}{N}$. The resulting pair has gap $h$
and its objective is at most that of the original pair. Hence on such
intervals the definition agrees with OpenAI's separated minimum
\cite[definition 10.1]{OpenAI2026},
$$
 \mathcal V_h(J)=\min_{\substack{a\le i<j\le t\\j-i\ge h}}(A_i+A_j).
$$
For intervals of every width the choice $i=j=a$ gives
\begin{equation}\label{ledger:terminal-budget}
 \mathcal V_h([a,t])-2A_a\le\frac{dh}{N}.
\end{equation}

Each state uses its current interval, $[a,t]$ for
$\mathsf D_i(a,t)$ and $[a,v]$ for $\mathsf H_i^{\mathcal B}(a,v)$.
The retained base $\mathcal B$ determines the angular state's
descriptor list. The origin of a descriptor determines the interval on
which its anchor minimizes $A$. In the direct graph reduction below,
these intervals play different roles.

Restricting $J$ to a subinterval removes candidates from the minimum
in \eqref{ledger:potential-formula}. Thus $\mathcal V_h$ increases under
restriction. To bound a positive cost we take a minimizing pair in the
child and exhibit a cheaper pair in the parent. The difference between
their sums is the amount available to pay that cost.

\begin{lemma}[Payment of the costs]\label{ledger:payment}
Every transition of proposition~\ref{rec:reductions}, with parent
interval $J$ and child interval $J'$, satisfies
\begin{equation}\label{ledger:payment-inequality}
 \mathrm{cost}\le \mathcal V_h(J')-\mathcal V_h(J).
\end{equation}
\end{lemma}

\begin{proof}
A zero cost transition follows from $J'\subseteq J$, including a
transition that stops at a child of width less than $h$. In the three
remaining cases we check that the child is long enough for the separated
minimum, choose a minimizing child pair of gap at least $h$, and replace
one of its depths.

For a suffix step we have $J=[a,t]$, $J'=[a,p]$, and
$\mathrm{cost}=A_p-\min_{[p,t]}A$. Since $2p>a+t$, the child width is
more than half the nonterminal parent width, and hence greater than
$10q_0$. Let $(i,j)$ be admissible in $[a,p]$. Then $j\ge a+10q_0$,
whereas the left endpoint of the origin's inner interval is at most
$a+2q_0$, and since $j\le p$, the index $j$ belongs to that inner
interval, so that $A_j\ge A_p$. Choose $r\in[p,t]$ attaining
$\min_{[p,t]}A$. The pair $(i,r)$ is admissible in the parent because
$r\ge j$, and
$$
 A_i+A_r\le A_i+A_j-
          (A_p-\min_{[p,t]}A),
$$
which is \eqref{ledger:payment-inequality}.

For a high prefix step, $J=[a,t]$, $J'=[p,v]$, $t-v<\sigma N$, and the
cost is $k(a,p)$. The prefix inequality gives
$$
 v-p\ge\frac{t-a}{2}-\sigma N>10q_0.
$$
For an admissible child pair $(i,j)$ we have
$p\le i\le v-10q_0\le t-10q_0$, so that $i$ is in the origin's inner
interval and $A_i\ge A_p$. Choose $r\in[a,p]$ attaining $\min_{[a,p]}A$.
The pair $(r,j)$ stays in $[a,t]$ and remains admissible because
$r\le i$, and its sum is at most $A_i+A_j-k(a,p)$, as required.

For a direct graph reduction, $J=[a,v]$ and $J'=[p,v]$. Let $L$ be the
maximal remaining descriptor length chosen by the reduction, and let
$[b_0,u_0]$ be its origin, so that
\begin{equation}\label{ledger:origin-geometry}
b_0\le p-L, p+L\le u_0,
  v\le a+2L.
\end{equation}
The reduction gives $v-p\ge\frac{v-a}{2}-\sigma N$, which
nonterminality and $\sigma\ll\kappa$ make greater than $10q_0$. If
$k(a,p)=0$, containment suffices. Otherwise choose $i_0\in[a,p]$ with
$$
 A_{i_0}=\min_{[a,p]}A=A_p-k(a,p)<A_p.
$$
The origin's choice rule makes $A_p$ no larger than every value in its
inner interval, and since $i_0\le p$, the index $i_0$ must lie strictly
to the left of that inner interval, whether or not it lies in the origin
itself. Consequently
$$
 a\le i_0<b_0+2q_0,
  v\le a+2L<b_0+2q_0+2L\le u_0+2q_0.
$$
For any admissible child pair $(i,j)$,
$$
 p\le i\le v-10q_0<u_0-8q_0.
$$
The origin's inner interval begins no later than $p$ and ends no earlier
than $u_0-2q_0$, so that $i$ lies in it and $A_i\ge A_p$. Replacing $i$
by $i_0$ gives a pair $(i_0,j)$ in the actual parent $[a,v]$, not merely
in the origin, which remains admissible because $i_0\le p\le i$ and whose
sum is smaller by at least $k(a,p)$. This proves the direct case.
\end{proof}

\begin{figure}[htbp]
\centering
\begin{tikzpicture}[>=Stealth,font=\small]
\draw[thick] (0,0)--(12,0);
\draw[very thick] (0,-0.12)--(7,-0.12);
\foreach \x/\name in {0/a,2/i,5/j,7/p,10/r,12/t}{
\draw (\x,0.1)--(\x,-0.2);
\node[below] at (\x,-0.25) {$\name$};}
\draw[->,thick] (5,0.3) to[bend left=25] node[above]
{replace $j$ by $r$} (10,0.3);
\node[align=center] at (6,-1.5)
{$A_j\ge A_p$, $A_r=\min_{[p,t]}A$\\
$r-i\ge j-i\ge10q_0$};
\end{tikzpicture}
\caption{For a suffix reduction, replacing the later depth of a child
pair by $r$ gives an admissible parent pair and reduces its sum by at
least $A_p-A_r$. This is the cost of the reduction.}
\label{fig:potential-replacement}
\end{figure}
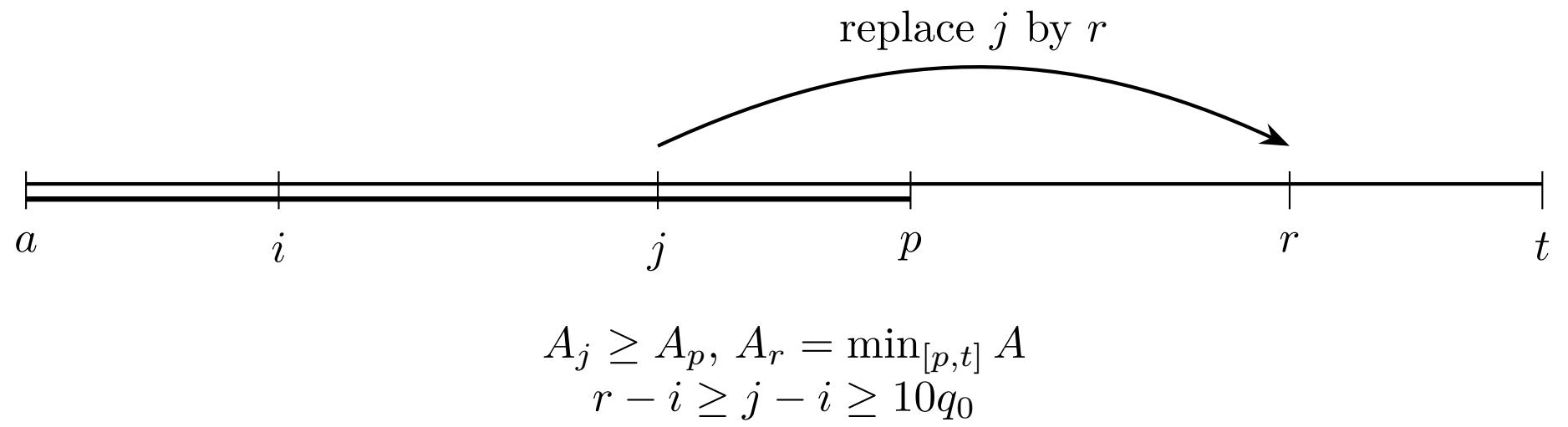

\subsection{The descent and its stopping time}

Since $|A_n|\le\frac{d}{2}$, the equation~\eqref{ledger:terminal-budget}
gives $\mathcal V_h(J)\le d+\frac{dh}{N}<4d$ for sufficiently small
fixed $\kappa$ and large $N$. Thus $\Phi(J)=4d-\mathcal V_h(J)$ is a
nonnegative decreasing potential on every interval. If $J_0',\ldots,J_m'$
are the successive intervals of a branch and $c_l$ is the cost of its
$l$th transition, lemma~\ref{ledger:payment} gives
$$
 \sum_{l=0}^{m-1}c_l\le\Phi(J_0')-\Phi(J_m').
$$
The inequality holds across a change from a distance count to an angular
energy as well as across refinements of either state.

To bound the length of the branch we also have to count the transitions
of cost zero. Since the two ends of the interval move in only one
direction and the block representatives are fixed in advance, we can
count the steps which move the endpoints and charge the remaining near
steps to representatives just passed or to the next direct step, and the
same count bounds the number of smoothing levels used.

\begin{lemma}[Finite descent and smoothing budget]\label{ledger:termination}
Every branch of the state reductions has at most $\frac{20}{\kappa}$
non-low transitions and at most $\frac{1}{\sigma}$ low transitions, for
sufficiently small fixed parameters and sufficiently large $N$.
Along such a branch the smoothing index increases by at most
$\frac{20}{\kappa}+\frac{1}{\sigma}$.  If the initial index is zero, every index
required by the iteration, including the next support level used by
a mask, is less than
$$
 J_0=\left\lceil100\left(\frac{1}{\kappa}+\frac{1}{\sigma}\right)+100\right\rceil.
$$
\end{lemma}

\begin{proof}
Along a branch the start depth never decreases and the end depth never
increases. Let $K$ be the number of fixed block representatives, so that
$K\le\frac{1}{\kappa}+2$ for large $N$. Every high or direct step strictly
advances the start to a representative, and since such a representative can
never be passed twice, these steps number at most $K$ together. Each suffix
step lowers the end by at least $q_0$, so that there are at most
$\frac{N}{q_0}\le\frac{1}{\kappa}$ suffix steps, whereas low steps lower the
end by at least $\sigma N$, and therefore there are at most
$\frac{1}{\sigma}$.

An angular episode starts at the root or immediately after a high step and
ends at a far step or a stopping term, direct and near steps remaining
inside that episode. There are at most $K+1$ episodes, and hence at most
$K+1$ far steps.

Consider a near step whose child is nonterminal. That child has width
$2L_0>100q_0$, so that the floor in the formula for $L_0$ is inactive and
$L=L_0$. Remove the factors whose anchors have now been reached. If the
maximum remaining nominal length decreases, then at least one formerly
maximal anchor has just been passed, and we charge the near step to one
such anchor, these charges being distinct along the branch. If instead the
maximum stays equal to $L$, then recomputing $v-2L$ at the child gives
exactly its start depth, so that the very next step is direct, and we
charge this near step to that next direct step. There is finally at most
one near step with terminal child, and therefore near steps number at most
$2K+1$, and the combined non-low count is at most $5K+2$, which is less
than $\frac{20}{\kappa}$ for sufficiently small $\kappa$.

Only suffix, low and far steps increase the smoothing index, and each
increases it by one, so that if $i_{\mathrm{in}}$ is the initial index, the
preceding counts give
$$
 i\le i_{\mathrm{in}}+\frac{20}{\kappa}+\frac{1}{\sigma}
$$
at every state on the branch. High, direct and near steps use the
level-$i+1$ supports of their index-$i$ symbols without increasing the
index, and for $i_{\mathrm{in}}=0$ even this additional support level
remains below $J_0$, since
$\frac{20}{\kappa}+\frac{1}{\sigma}+1<J_0$. Every mask and support level
required for the iteration from index zero is therefore available.
\end{proof}

The next table records the exact correspondence with the abstract descent
lemma, and it is at the same time a list of the facts which must survive
each reduction, since a proof of payment alone would leave the stopping
terms, the number of steps, and the accumulated analytic errors
unaccounted for.

\begin{table}[htbp]
\centering\small
\renewcommand{\arraystretch}{1.25}
\begin{tabular}{p{3.05cm}p{5.05cm}p{5.05cm}}
\toprule
Data in lemma~\ref{bik:tree} & Choice here & Justification\\
\midrule
$X_v$, $u_v$, $V_v$ & Current count or energy, $2A_a$,
$\mathcal V_h([a,t])$ & proposition~\ref{rec:reductions}
and definition~\ref{ledger:potential}\\
Edge cost $c_{vw}$ & The cost of the selected reduction &
lemma~\ref{ledger:payment}\\
Marked edges & All non-low transitions & At most
$\lceil\frac{20}{\kappa}\rceil$ by lemma~\ref{ledger:termination}\\
Stopping allowance $\tau$ & A dimensional constant times
$\kappa+\sigma$ & Terminal estimates and
$\mathcal V_h([a,t])-2A_a\le\frac{dh}{N}$\\
Finite records & Depth bounded in terms of $\kappa,\sigma$ and
polynomially many scale choices & lemma~\ref{ledger:termination}
and proposition~\ref{rec:reductions}\\
Error exponent $T$ & Fixed larger than $B_0+K\varepsilon+1$ &
The arbitrary-order oscillatory estimates, after fixing all parameters\\
\bottomrule
\end{tabular}
\vskip.125in 
\caption{How the analytic estimates meet the hypotheses of the
descent lemma. The constants $B_0,K,\varepsilon$ are those of
lemma~\ref{bik:tree}, with $K$ the marked-edge count and
$\varepsilon$ the loss allowed at each marked edge.}
\label{tab:descent-hypotheses}
\end{table}
\begin{theorem}[Single-mass-curve angular bound]\label{ledger:bound}
Let $\mathsf H_0^{[c,N]}(c,N)$ be a recursive angular state of
nonterminal width, as in definition~\ref{rec:angular-state}.
Fix its initial derivative budget $B$ and its common numerical angular
bounds $Q_j$, independently of $N$, as in
definition~\ref{osc:symbols}.
Then
\begin{equation}\label{ledger:bound-equation}
\mathsf H_0^{[c,N]}(c,N)
 \le R^{2A_c-\mathcal V_h([c,N])
                   +C_d(\kappa+\frac{\sigma}{\kappa})+o(1)}.
\end{equation}
This bound holds uniformly over its two occupied cells, including
two copies of the same cell, and over its admissible bounded
spatial multipliers and symbols with these fixed bounds.  The $o(1)$
term may depend on $B$ and the $Q_j$.
\end{theorem}

\begin{proof}
We apply lemma~\ref{bik:tree}. At a vertex whose current interval is
$J=[a,t]$, let $X$ be the positive distance count or angular energy
specified by that vertex. Set
$$
u=2A_a, V=\mathcal V_h(J) \text{ and } Y=R^{\mathcal V_h(J)-2A_a}X.
$$
For a non-low child with starting depth $a'$, the coefficient in
proposition~\ref{rec:reductions} is
$$
R^{2A_a-2A_{a'}+c+C_d\sigma+o(1)}.
$$
After normalization its exponent becomes
$$
c+\mathcal V_h(J)
-\mathcal V_h(J')+C_d\sigma+o(1)
\le C_d\sigma+o(1).
$$
For a low child the start does not change, its cost is zero, and the whole
coefficient is only $R^{o(1)}$, and hence its normalized coefficient is also
at most $R^{o(1)}$. This verifies the transition hypothesis of the general
descent lemma, with the non-low edges marked, and the distinction
is numerically essential. There may be on the order of $\frac{1}{\sigma}$
low steps, and charging $C_d\sigma$ at each would give a total exponent of
order one, which could exceed the small gain chosen at entry, whereas their
actual losses are $o(1)$, and a fixed number of such losses, depending on
$\sigma$, still tends to zero as $N$ tends to infinity. The fixed power
losses occur only on the order of $\frac{1}{\kappa}$ times, which yields
the permitted exponent $C_d\frac{\sigma}{\kappa}$.

At a terminal vertex, proposition~\ref{rec:reductions} gives
$X\le R^{C_d(\kappa+\sigma)+o(1)}$. The equation~\eqref{ledger:terminal-budget} gives
$$
V-u\le\frac{dh}{N}=O_d(\kappa)+O_d\left(\frac{1}{N}\right).
$$
Hence
$$
Y\le R^{C_d(\kappa+\sigma)+o(1)}.
$$
The same conclusion applies to a low-frequency stopping term at a
nonterminal vertex, since its unnormalized value is at most
$R^{C_d\sigma+o(1)}$. Thus the stopping allowance in
lemma~\ref{bik:tree} can be taken to be
$\tau=C_d(\kappa+\sigma)$, after enlarging the dimensional constant.

Lemma~\ref{ledger:termination} bounds the total number of edges on a
branch by a constant depending on $\kappa,\sigma$, of which at most
$K=\lceil\frac{20}{\kappa}\rceil$ are non-low, and each reduction
has polynomially many scale choices. Maxima over cell choices and
suprema over admissible symbols are bounded uniformly at the next
stage rather than summed, so that for the actual sums the number of
records is polynomial in $N$ to a fixed power and contributes only
$R^{o(1)}$; this is the finite-record chaining step.

All $A_n$ and all values of $\mathcal V_h$ are bounded by a
dimensional constant. Consequently $|V-u|\le B_0$ with $B_0$ depending
only on $d$. Choose the local oscillatory accuracy $T$ so that
$$
T>B_0+C_d\sigma\left\lceil\frac{20}{\kappa}\right\rceil+10.
$$
This choice is made after the fixed parameters and before $N$ tends to
infinity. The normalized errors then satisfy the error hypothesis of
lemma~\ref{bik:tree}. We obtain
$$
Y_{\mathrm{root}}
\le R^{C_d(\kappa+\sigma)+C_d\sigma\lceil\frac{20}{\kappa}\rceil+o(1)}
\le R^{C_d(\kappa+\frac{\sigma}{\kappa})+o(1)}.
$$
Undoing the normalization at the root proves
\eqref{ledger:bound-equation}.

It remains to specify the derivative bounds required by the oscillatory
estimates. Every passage from a distance count to an angular energy
starts with undifferentiated masks, and once the inverse stationary-phase
order has been chosen for the fixed accuracy $T$, this order together
with the initial derivative order in the theorem gives one finite bound
for the angular derivatives retained in the iteration. The graph estimate
only restricts these symbols to smaller cells, removing factors whose
depths have been reached introduces no additional derivatives, and a far
reduction may use a longer stationary-phase expansion, but its child is a
positive distance count and retains none of these derivatives. We can
therefore choose the direct expansion and integration orders after
fixing the persistent derivative bound, and their constants depend on
the fixed parameters and the initial angular bounds, but not on $N$.
Thus the quantities used in every application of the descent lemma have
the common bounds required by the analytic reductions, and no new
derivative allowance is needed at any generation.
\end{proof}

\subsection{What has changed in the iteration}

The local state reductions and the replacement of a depth in a minimizing
pair are the arguments of OpenAI
\cite[proposition 9.4 and lemma 10.2]{OpenAI2026}, and the difference
lies in the accounting, as explained in section~\ref{sec:changes}. The
formula \eqref{ledger:potential-formula} keeps their payment exact on
long intervals and controls every stopping interval by the same
minimization, so that after multiplication by
$R^{\mathcal V_h(J)-2A_a}$ the only fixed losses left on a branch occur
at the marked transitions. This is the role of the descent principle of
Burstein, Iosevich, and Krause \cite{BursteinIosevichKrause2026} in the
iteration, and the choice of its initial interval will supply the saving
needed for the Fourier estimate.

\section{Entry, shell summation, and the limiting distance measure}
\label{asm:section}

The two mass curves may enter the descent at different heights, and we
first choose one depth at which their combined saving exceeds the cost
of entry. This gives a Fourier estimate for each pair of regular classes
with arbitrary bounded spatial multipliers at its endpoints, and we use
those multipliers to put back the coarse pair filter before summing with
the original class masses. The geometric estimate on late exclusions then
compares the resulting measure on each shell to one fixed measure on the
original distance set in $\mathbb R^d$.

\subsection{The entry depth}

Let $S=\frac{d}{2}$.  Fix the starting probabilities $\mu_1,\mu_2$ and the number
$\chi>0$ supplied by proposition~\ref{prep:starting-pairs}.
Make one common similarity placing $E$ in $(\frac{1}{4},\frac{3}{4})^d$, and retain
the notation for the transformed sets and measures.  This multiplies
every distance by one fixed positive constant and preserves the
desired positive-measure conclusion.  Choose one small common
translation of the dyadic grid whose boundaries have zero mass for
both starting laws.  Its depth-zero cube $Q_0$ still contains $E$ in
its interior.  Use this same grid for both class decompositions.
The two actual supports remain separated by a fixed positive
distance, and their cross directions remain in the prepared patch.
At scale
$R=2^N$, apply lemma~\ref{reg:regularization} with parameter $\sigma$ to
both probabilities.  Write $\mathcal C_{j,N}$ for the retained classes on
side $j$.  For $C\in\mathcal C_{j,N}$, let
$$
 a_{j,C}=\mu_j(C) \text{ and } \rho_{j,C}=a_{j,C}^{-1}\mu_j|_C.
$$
The classes on each side are disjoint, their number is $R^{o(1)}$, and
their discarded original mass is at most $R^{-\sigma+o(1)}$, all the $o(1)$
terms below being uniform over these classes with the fixed parameters held
constant. Every normalized class is supported in the same root cube $Q_0$
and assigns it mass one, since every mass curve has $m_j(0)=0$, so that its
root normalization is exactly $M_j(0)=R^{-m_j(0)}=1$.

For a pair of classes, suppress their labels and write
$$
 A_{j,n}=m_j\left(\frac{n}{N}\right)-\frac{Sn}{N} \text{ and } k_j(g,p)=A_{j,p}-\min_{g\le n\le p}A_{j,n}.
$$
Recall from lemma~\ref{test:profile-bounds} that
\begin{equation}\label{asm:profile-input}
A_{j,0}=0, |A_{j,n}-A_{j,n'}|\le d\frac{|n-n'|}{N} \text{ and } A_{j,n}\ge\chi \frac{n}{N}-C_d\sigma.
\end{equation}
Here and in what follows, constants denoted by $C_d$ may be enlarged a
finite number of times and depend only on $d$.

We seek a common depth $c$ for which one entry cost $k_j(0,c)$ is small
while both excess mass curves remain above a positive level on $[c,N]$. The
graph estimate contributes the smaller cost, and after Cauchy--Schwarz has
separated the mass curves, the two single-mass-curve estimates must supply
a larger potential saving. The choice by a last crossing of a fixed level
comes from \cite[lemma~11.1]{OpenAI2026}.

Put $h=10q_0$. For a nonterminal interval $[c,N]$, denote by
$\mathcal V_{h,j}([c,N])$ the potential in
definition~\ref{ledger:potential} for mass curve $j$. On an interval
of this width it equals the minimum sum of two excess values at depths
separated by at least $h$. The net saving at entry is
\begin{equation}\label{asm:net-saving}
\mathfrak g(c)
=\frac{\mathcal V_{h,1}([c,N])+\mathcal V_{h,2}([c,N])}{2}
 -\min_{j=1,2}k_j(0,c).
\end{equation}
The normalization explains this difference. The graph step
contributes the exponent $-A_{1,c}-A_{2,c}+\min_j k_j(0,c)$, while the
single-curve descent bound contributes
$2A_{j,c}-\mathcal V_{h,j}([c,N])$ on the side $j$, and Cauchy--Schwarz takes
half the sum of these two exponents. The terms $A_{1,c}$ and $A_{2,c}$
therefore cancel, and the exponent which remains is $-\mathfrak g(c)$,
before the small analytic losses. The estimate~\eqref{asm:entry-cancellation}
will check this with the full energy normalizations in place.

This is the useful comparison with the reciprocal estimates of Iosevich and
Yavicoli \cite{IosevichYavicoli2026}, in that we retain the separate
endpoint excesses until one common depth allows their contributions to be
combined. The two depths within each potential and the two endpoints of the
distance have different roles here, the first paying for the descent of one
curve, the second being combined only in this outer calculation.

The last crossing of a level is more useful here than the first, since at
the first crossing one of the mass curves could fall below the level again,
whereas at the last crossing of $\min(A_1,A_2)$ both curves stay above it
afterwards while at least one is still at or below the level at the
crossing itself. These are the two properties needed in order to
compare the entry cost with the later saving.

\begin{lemma}[Entry depth]\label{asm:entry}
Fix $0<\beta<\min(\frac{\chi}{10},\frac{1}{10})$.  If $\sigma$ is sufficiently small
in terms of $d,\chi,\beta$, then, for all sufficiently large $N$, the
integer
\begin{equation}\label{asm:entry-definition}
c=\max\{n\in\{0,\ldots,N\}:
                  \min(A_{1,n},A_{2,n})\le\beta\}
\end{equation}
is defined and satisfies
\begin{equation}\label{asm:entry-range}
 \left\lfloor\frac{\beta N}{2d}\right\rfloor\le c<\frac{N}{4},
\end{equation}
\begin{equation}\label{asm:entry-lower}
 A_{j,n}\ge\beta-\frac{d}{N}
        \text{ for } j=1,2 \text{ and } c\le n\le N,
\end{equation}
and
\begin{equation}\label{asm:entry-cost}
 \min_{j=1,2} k_j(0,c)\le\beta+C_d\sigma.
\end{equation}
If $\kappa$ is sufficiently small and $q_0=\lceil\kappa N\rceil$, the
interval $[c,N]$ is nonterminal.  Its potentials from
definition~\ref{ledger:potential} obey
\begin{equation}\label{asm:entry-potential}
\mathcal V_{h,j}([c,N])
 =\min_{\substack{c\le u<v\le N\\v-u\ge h}}
                         (A_{j,u}+A_{j,v})
 \ge2\beta-\frac{2d}{N}.
\end{equation}
\end{lemma}

\begin{proof}
The index $0$ qualifies in \eqref{asm:entry-definition}. For
$n\le\frac{\beta N}{2d}$, the Lipschitz bound in
\eqref{asm:profile-input} gives $A_{j,n}\le\frac{\beta}{2}$, which proves
the lower bound in \eqref{asm:entry-range}, while for $n\ge \frac{N}{4}$
the lower mass curve bound gives $A_{j,n}\ge\frac{\chi}{4}-C_d\sigma>\beta$
when $\sigma$ is sufficiently small, which proves its upper bound.

By maximality of $c$, both $A_{j,c+1}$ exceed $\beta$, so that
$A_{j,c}>\beta-\frac{d}{N}$ and both mass curves exceed $\beta$ at every
later index, whereas for at least one side $j_*$ we also have
$A_{j_*,c}\le\beta$. Since $\min_{0\le n\le c}A_{j_*,n}\ge-C_d\sigma$, this
gives \eqref{asm:entry-cost}. The width of $[c,N]$ exceeds $\frac{3N}{4}$
and therefore exceeds $100q_0$ for sufficiently small fixed $\kappa$ and
large $N$, and every pair entering its potential consists of two indices
covered by \eqref{asm:entry-lower}, which proves
\eqref{asm:entry-potential}.
\end{proof}

The potential does not by itself produce the final negative exponent. At
entry each mass curve contributes a potential of at least
$2\beta-O(\frac{1}{N})$, and the equations~\eqref{asm:entry-cost}
and~\eqref{asm:entry-potential} give the useful conclusion
$$
\mathfrak g(c)\ge\beta-C_d\sigma-\frac{2d}{N}.
$$
The lemma therefore produces a positive net saving, with room for the
analytic losses when $\sigma$ is small and $N$ is large. Cauchy--Schwarz
later takes the geometric mean of the two resulting bounds; thus their
combined saving is at least $2\beta-O(\frac{1}{N})$, whereas the graph
entry costs at most $\beta+O_d(\sigma)$, since it uses the smaller of the
two entry costs. The difference leaves the positive margin $\beta$, before
the small analytic losses, and the
estimate~\eqref{asm:entry-cancellation} performs this computation with
every normalization present.

\subsection{Positive approximants and restoration of mass}

Write
$$
 D(x,y)=|x-y|, e(x,y)=\frac{x-y}{|x-y|} \text{ and } W(x,y)=D(x,y)^{-\frac{d-1}{2}}.
$$
The fixed separated supports give constants $0<d_0<d_1<\infty$ with
\begin{equation}\label{asm:distance-range}
d_0\le D(x,y)\le d_1
  \text{ on } \operatorname{supp}\mu_1\times\operatorname{supp}\mu_2.
\end{equation}
In particular $W$ is bounded above and below there by positive constants.

A subscript or superscript $j$ on a descriptor list or mask
indicates that it is constructed using the class at endpoint $j$.
For each pair of classes $(C_1,C_2)$, use its own entry index $c$ from \eqref{asm:entry-definition}.  At endpoint $j$ take the extra
parallel descriptor $[0,c]$, together with every descriptor in
$\mathcal M_j(c,N)$.  Let
\begin{equation}\label{asm:endpoint-mask}
b_{j;C_1,C_2}(z,w)
 =w^{(j)}_{0,[0,c]}(z,w)
                  \prod_{\mathfrak t\in\mathcal M_j(c,N)}
                         w^{(j)}_{0,\mathfrak t}(z,w).
\end{equation}
By lemma~\ref{test:smooth-masks}, every factor
lies in $[0,1]$ and equals one when its level-zero condition holds.
Antipodal symmetry allows both endpoints to be evaluated at $e(x,y)$.
Define the Borel function $\Theta_N$ on the fixed separated product
$E_1\times E_2$ by
$$
 \Theta_N(x,y)=
 b_{1;C_1,C_2}(x,e(x,y))
 b_{2;C_1,C_2}(y,e(x,y))
  \text{ if } (x,y)\in C_1\times C_2,
$$
and set it equal to zero whenever either endpoint belongs to no retained
class.
In particular $0\le\Theta_N\le1$.

Let $\Gamma_N$ be the decreasing pair filters of
proposition~\ref{prep:starting-pairs}.  With $D_\#$ denoting pushforward
under the scalar map $D$, define
\begin{equation}\label{asm:measures}
\alpha_N=D_\#(W1_{\Gamma_N}
                          (\mu_1\otimes\mu_2))
 \text{ and }
 \tau_N=D_\#(W1_{\Gamma_N}\Theta_N
                          (\mu_1\otimes\mu_2)).
\end{equation}
These measures use the original class masses.  Only $\tau_N$ depends on
the shell's class decomposition and angular masks.

We now return to the original measures. A class of small original mass may
have become a probability during the internal estimates, but it must carry
its original weight when the classes are recombined. The
exceptional-direction bound is integrated against the original opposite
fiber, which is a subprobability, and the deletion estimate therefore
survives recombination without an inverse class-mass factor.

\begin{lemma}[Restoration of mass]\label{asm:retained-measures}
There is $\eta>0$, depending on the fixed parameters and starting
measures, such that
\begin{equation}\label{asm:tv-approximation}
0\le\tau_N\le\alpha_N \text{ and } \|\alpha_N-\tau_N\|_{\mathrm{TV}}\lesssim2^{-\eta N}.
\end{equation}
Moreover $\alpha_N$ decreases to the nonzero finite positive measure
\begin{equation}\label{asm:limit-measure}
\alpha_\infty=D_\#(W1_{\Gamma_\infty}
                               (\mu_1\otimes\mu_2)) \text{ with } \Gamma_\infty=\bigcap_N\Gamma_N.
\end{equation}
All these measures are carried by the compact set $\Delta(E)$ of original
distances. For some $e>0$ and all sufficiently large $N$,
\begin{equation}\label{asm:fixed-tv}
 \|\tau_N-\alpha_\infty\|_{\mathrm{TV}}\lesssim 2^{-eN}.
\end{equation}
\end{lemma}

\begin{proof}
Apply the cap preparation in proposition~\ref{prep:starting-pairs} with
loss $R^\sigma$.  It gives $\mathcal G_N\subset\Gamma_N$ such that
$$
 (\mu_1\otimes\mu_2)(\Gamma_N\setminus\mathcal G_N)
       \lesssim R^{-c_*}
  \text{ for some } c_*>0.
$$
For $\mu_1$-almost every first endpoint $x$, the angular subprobability
$$
 \xi_{N,x}=e(x,\cdot)_\#
             (1_{\mathcal G_N}(x,\cdot)\mu_2)
$$
has the required $S$-cap bound with loss $R^\sigma$ down to scale
$R^{-1}$.  The change from $\pi_x(y)$ to $e(x,y)=-\pi_x(y)$ preserves
cap bounds, since the antipodal map sends caps to caps.  The analogous
assertion holds for $\mu_2$-almost every second endpoint with the
endpoints reversed.
These are bounds for restrictions of the original opposite probability.

Fix $C_1,C_2$ and one descriptor $\mathfrak t$ used at the first endpoint.
For $\mu_1$-almost every $x\in C_1$, the anchor cell determined by
$\mathfrak t$ and $x$ is occupied and $\xi_{N,x}$ has the cap bound above.  Let
$F_{C_1,C_2,\mathfrak t,x}$ be the directions failing this level-zero test.
For such $x$, lemma~\ref{test:exceptional-directions} applies directly
to $\xi_{N,x}$ and gives
\begin{equation}\label{asm:original-fiber-estimate}
\mu_2\{y\in C_2:(x,y)\in\mathcal G_N,
       \ e(x,y)\in F_{C_1,C_2,\mathfrak t,x}\}\le R^{-\sigma}.
\end{equation}
Restricting the set of $y$ to $C_2$ only decreases the left side, and no
normalization by $\mu_2(C_2)$ occurs, while integrating in $x$ over $C_1$
costs at most $a_{1,C_1}R^{-\sigma}$. Sum first over $C_1$, whose original
masses have sum at most one, and then over $C_2$ and the descriptors, the
latter two counts being $R^{o(1)}$ by lemma~\ref{reg:regularization} and
lemma~\ref{test:descriptor-properties}. The first-endpoint failure mass is
therefore $R^{-\sigma+o(1)}$, and the same argument in the other
orientation gives this bound for second-endpoint failures. No sum over
anchor cells occurs, since a point and a descriptor already specify one.

The product in \eqref{asm:endpoint-mask} equals one if all
its level-zero tests hold.  Hence $1-\Theta_N$ on retained class pairs is
bounded by the indicator of the union of these failures.  Including
the original mass outside retained class pairs and
$\Gamma_N\setminus\mathcal G_N$ gives
$$
 \int_{\Gamma_N}(1-\Theta_N)d(\mu_1\otimes\mu_2)
 \lesssim R^{-c_*}+R^{-\sigma+o(1)}.
$$
The bounded weight $W$ proves \eqref{asm:tv-approximation}, for example
with any fixed $0<\eta<\frac{\min(c_*,\sigma)}{2}$ and sufficiently large
$N$. The set $\mathcal G_N$ was used only in this estimate, and it is not a
multiplier in \eqref{asm:measures}.

The filters $\Gamma_N$ decrease and their intersection has positive
product mass, so that monotone convergence gives
\eqref{asm:limit-measure}, the lower bound for $W$ implies
$\alpha_\infty({\mathbb R})>0$, and the distance images are contained in
$\Delta(E)$ since the prepared supports are compact subsets of $E$.

In odd dimensions \eqref{prep:coarse-tail} gives
$(\mu_1\otimes\mu_2)(\Gamma_N\setminus\Gamma_\infty)\lesssim 2^{-c_\Gamma N}$,
and since the measure $\alpha_N-\alpha_\infty$ is positive,
$$
 \|\alpha_N-\alpha_\infty\|_{\mathrm{TV}}
 =\int_{\Gamma_N\setminus\Gamma_\infty} Wd(\mu_1\otimes\mu_2)
 \lesssim 2^{-c_\Gamma N}.
$$
Combining this with \eqref{asm:tv-approximation} proves
\eqref{asm:fixed-tv} with $e=\min(\eta,c_\Gamma)$. In even dimensions
$\Gamma_N=E_1\times E_2$ at every scale and we may take $e=\eta$.
\end{proof}

\subsection{The shell estimate}

\begin{proposition}[Decay of the approximating distance measures]
\label{asm:shell-decay}
The fixed parameters can be chosen, in the order
$\beta,\kappa,\sigma,\omega$, and therefore some $\gamma>0$ satisfies
\begin{equation}\label{asm:shell-bound}
\int_{2^{N-1}\le |r|\le2^{N+1}}
                   |\widehat{\tau_N}(r)|^2dr
       \lesssim2^{-\gamma N}
\end{equation}
for all sufficiently large $N$.
\end{proposition}

\begin{proof}
Fix one normalized class pair and real spatial multipliers $f(x)$ and
$g(y)$ with $|f|,|g|\le1$. Write $\nu_{f,g}$ for the pushforward by $D$
of the signed measure
$$
 W(x,y)f(x)g(y)
 b_{1;C_1,C_2}(x,e(x,y))b_{2;C_1,C_2}(y,e(x,y))
             d\rho_{1,C_1}(x)d\rho_{2,C_2}(y).
$$
The shell conversion permits both multipliers, including their signs,
and no angular derivative falls on either of them.
At depth zero take $X=Y=Q_0$, each with its
own normalized class law; these are occupied root cells, with
$M_1(0)=M_2(0)=1$, and their supported laws lie in the fixed separated sets
$E_1,E_2$, so that the outer-support geometry of
section~\ref{osc:section} applies. The cell-center separation in
\eqref{osc:separation} is the internal-cell convention, whereas the outer
root uses the separation of the supported laws; in particular all the root
spatial diameters in the graph argument are $O_d(1)$, while the actual
cross distances have a fixed positive lower bound. By
lemma~\ref{test:descriptor-properties}, the largest nominal length of an
endpoint descriptor is at most
$$
 L\le\max\left(c,\frac{N-c}{2}\right)\le\frac{N}{2}.
$$
The root shell conversion, proposition~\ref{osc:shell-conversion},
therefore applies with $a=0$ and $v=N$. Its weight is $W$ itself, and
there is no internal normalization of the distance weight at this root. For
any prescribed fixed $T>0$ it bounds the squared Fourier norm on the
positive shell by $R^{o(1)}\sup\mathsf H_F(0,N)+O(R^{-T})$, the supremum
retaining the root cells, the class normalizations, the lists and the
smoothing index zero, while its admissible symbols have the fixed
derivative budget and the common numerical angular bounds supplied by that
proposition for $T$.

Apply corollary~\ref{graph:separated-form} at $p=c$ to each energy in this
supremum, the required parallel descriptor being the extra $[0,c]$ on each
side, with $g=0$. All the scale conditions hold, since $c>0$ for large $N$,
$c<\frac{N}{4}$, $L\le \frac{N}{2}<N-c$ and $N-c\ge10\sigma N$, and since
$A_{j,0}=0$, the graph factor is
\begin{equation}\label{asm:graph-factor}
R^{-A_{1,c}-A_{2,c}+\min_j k_j(0,c)+O_d(\sigma)+o(1)}.
\end{equation}
The graph estimate produces a general two-mass-curve energy
$\mathsf H(c,N)$ on a pair of occupied $c$-cells, still carrying the
same endpoint symbols.

On these cells the extra $[0,c]$ factor is spatially constant within
each endpoint cell, so that it, and any of its normalized differentiated
versions produced by shell conversion, can be removed at bounded
cost from the positive independent angular energies, the remaining
descriptor lists being exactly those with base $[c,N]$ and with reached
anchors omitted. This removal is made in the general positive energy
$\mathsf H$, where the two angular factors are nonnegative.

To separate the two mass curves, write their nonnegative angular
energies as $E_1(r),E_2(r)$ and use the common measure
$d\lambda_N(r)=\Psi(\frac{r}{2^N})\frac{dr}{2^N}$.  Cauchy--Schwarz gives the exact
normalization
\begin{equation}\label{asm:two-mass-curve-cs}
\begin{aligned}
 2^{d(N-c)}\int E_1(r)E_2(r)d\lambda_N(r)
 &\le
 \left(2^{d(N-c)}\int E_1(r)^2d\lambda_N(r)\right)^{\frac{1}{2}}\\
 & \cdot
 \left(2^{d(N-c)}\int E_2(r)^2d\lambda_N(r)\right)^{\frac{1}{2}}.
 \end{aligned}
\end{equation}
Each parenthesis is the independent angular energy formed by duplicating
one endpoint law, cell and symbol, and with the unreached descriptors of
$\mathcal M_j(c,N)$ it is the single-mass-curve recursive state
$\mathsf H_0^{[c,N]}(c,N)$ for the mass curve $j$, a value which we denote
by $\mathsf H_j(c,N)$. Its symbols have the fixed derivative budget and the
common numerical angular bounds supplied by shell conversion, both of which
are preserved by restriction to the $c$-cells and by duplication. These
states impose no separation condition, so that theorem~\ref{ledger:bound}
applies to each of them, including the case of two identical cells, and
yields
$$
 \mathsf H_j(c,N)
 \le R^{2A_{j,c}-\mathcal V_{h,j}([c,N])
                         +O_d(\kappa+\frac{\sigma}{\kappa})+o(1)}.
$$
The mass curve terms in \eqref{asm:graph-factor} and the
geometric mean of these bounds combine to
\begin{equation}\label{asm:entry-cancellation}
\begin{aligned}
 &-A_{1,c}-A_{2,c}+\min_j k_j(0,c)
       +\frac{1}{2}\sum_{j=1}^2
                    (2A_{j,c}-\mathcal V_{h,j}([c,N]))\\
 &
 =\min_j k_j(0,c)-\frac{1}{2}
                     (\mathcal V_{h,1}([c,N])+\mathcal V_{h,2}([c,N]))\\
 &=-\mathfrak g(c)\le-\beta+C_d\sigma+\frac{2d}{N},
 \end{aligned}
\end{equation}
where the last inequality uses the cost and potential bounds in
lemma~\ref{asm:entry}.  Including the stated losses proves uniformly
over these spatial multipliers that
\begin{equation}\label{asm:component-shell}
\int_{2^{N-1}\le r\le2^{N+1}}|\widehat{\nu_{f,g}}(r)|^2dr
 \lesssim
 R^{-\beta+C_d(\kappa+\frac{\sigma}{\kappa})+o(1)}+R^{-T}.
\end{equation}
Here $T$ may be any fixed sufficiently large number, once the expansion and
integration-by-parts orders have been chosen, and the term $O_d(\sigma)$
has been included in $O_d(\frac{\sigma}{\kappa})$, with $\kappa<1$.

Put $B_N^2$ equal to the right side of \eqref{asm:component-shell},
including its implicit constant. We now put back the factor
$1_{\Gamma_N}$, which in even dimensions equals one, so that the bound
already applies there. Suppose then that $d=2k-1$. By
proposition~\ref{prep:starting-pairs}, the two supports have partitions
into at most $M\lesssim R^{(k+1)\omega}$ cells each; denote the first
partition by $(X_i)$ and the second by $(Y_j)$, and write $A_{ij}=1$ if
$X_i\times Y_j$ is retained and $A_{ij}=0$ otherwise. The angular masks
and the two regular classes remain fixed in the following calculation.

Assign independent uniform signs $\varepsilon_i$ to the cells $X_i$ and
put $f_\varepsilon(x)=\varepsilon_i$ on $X_i$. For each column put
$$
 S_j=\{i:A_{ij}=1\}, m_j=\#S_j \text{ and } Z_j=\sum_{i\in S_j}\varepsilon_i.
$$
If $m_j>0$, define on $Y_j$
$$
 g_\varepsilon(y)=\frac{m_j}{\mathbb E|Z_j|}\operatorname{sgn}(Z_j) \text{ with } \operatorname{sgn}(0)=0.
$$
Take $g_\varepsilon=0$ on an empty column. For $i\in S_j$, symmetry
makes $\mathbb E[\varepsilon_i\operatorname{sgn}(Z_j)]$ independent of
$i$, and its sum over $i\in S_j$ is $\mathbb E|Z_j|$, whereas for
$i\notin S_j$ the expectation vanishes by independence. Thus, pointwise
on the prepared product,
\begin{equation}\label{asm:signed-filter}
 \mathbb E[f_\varepsilon(x)g_\varepsilon(y)]=1_{\Gamma_N}(x,y).
\end{equation}
The size of the second multiplier costs only a square root, since
$\mathbb E Z_j^2=m_j$ and $\mathbb E Z_j^4=3m_j^2-2m_j$ for $m_j>0$,
and H\"older's inequality
$\mathbb E Z_j^2\le(\mathbb E|Z_j|)^{\frac{2}{3}}
(\mathbb E Z_j^4)^{\frac{1}{3}}$ then gives
$\mathbb E|Z_j|\ge\sqrt{\frac{m_j}{3}}$, so that
$\|g_\varepsilon\|_\infty\le\sqrt{3M}$ for every choice of the signs.
The lower bound for $\mathbb E|Z_j|$ is a nonsharp case of the
Khintchine inequality \cite{Haagerup1981}.

Write $\nu_{C_1,C_2}^{\Gamma_N}$ for the normalized class pair measure with
its angular masks and the coarse spatial filter inserted. The expectation
in \eqref{asm:signed-filter} is a finite average, so linearity under
integration gives
$\nu_{C_1,C_2}^{\Gamma_N}=\mathbb E\nu_{f_\varepsilon,g_\varepsilon}$.
Apply \eqref{asm:component-shell} to $f_\varepsilon$ and
$\frac{g_\varepsilon}{\sqrt{3M}}$. For $I_N=[2^{N-1},2^{N+1}]$, the triangle
inequality gives
$$
 \|\widehat{\nu_{C_1,C_2}^{\Gamma_N}}\|_{L^2(I_N)}
 \le \mathbb E\|\widehat{\nu_{f_\varepsilon,g_\varepsilon}}\|_{L^2(I_N)}
 \le \sqrt{3M}B_N .
$$
The angular masks, their derivative estimates, and the entry depth
belong to the original class pair and are independent of the signs.

It remains to restore the original class weights, for which the exact
identity is
$$
 \tau_N=\sum_{C_1,C_2}a_{1,C_1}a_{2,C_2}\nu_{C_1,C_2}^{\Gamma_N}
 \text{ with }
 \sum_{C_1,C_2}a_{1,C_1}a_{2,C_2}
 =\left(\sum_{C_1}a_{1,C_1}\right)\left(\sum_{C_2}a_{2,C_2}\right)\le1.
$$
The triangle inequality therefore introduces no further loss. Put
$\ell_d=0$ for even $d$ and $\ell_d=k+1$ for $d=2k-1$. Squaring the
norm bound yields
\begin{equation}\label{asm:assembled-shell}
 \int_{2^{N-1}}^{2^{N+1}}|\widehat{\tau_N}(r)|^2dr
 \lesssim R^{-\beta+C_d(\kappa+\frac{\sigma}{\kappa})+\ell_d\omega+o(1)}
                    +R^{-T+\ell_d\omega}.
\end{equation}

\begin{table}[htbp]
\centering\small
\renewcommand{\arraystretch}{1.25}
\begin{tabular}{p{3.2cm}p{4.7cm}p{5.15cm}}
\toprule
Contribution & Exponent or bound & Reason\\
\midrule
Entry cost & $\min_j k_j(0,c)\le\beta+C_d\sigma$ & One curve is at or below $\beta$ at the last crossing.\\
Descent saving & $\frac{\mathcal V_{h,1}+\mathcal V_{h,2}}{2}\ge2\beta-\frac{2d}{N}$ & Both curves stay above $\beta-\frac{d}{N}$ from entry onward.\\
Stopping and reductions & $C_d(\kappa+\frac{\sigma}{\kappa})$ & At most $\frac{20}{\kappa}$ transitions carry a fixed power loss.\\
Coarse filter & $\ell_d\omega$ & In odd dimensions the squared norm costs at most $3M$. In even dimensions the filter is the full product.\\
Class weights & Total at most one & Each normalized class receives its original mass.\\
\bottomrule
\end{tabular}
\vskip.125in 
\caption{The contributions to the squared Fourier norm in
\eqref{asm:assembled-shell}. The last crossing supplies a saving of
$\beta$ before the analytic and spatial losses.}
\label{tab:final-exponent}
\end{table}

Choose $\beta$ as in lemma~\ref{asm:entry}, then $\kappa$ so small that
the entry interval is nonterminal and $C_d\kappa<\frac{\beta}{16}$, then
$\sigma$ so that the mass curve and transition hypotheses hold and
$C_d\frac{\sigma}{\kappa}<\frac{\beta}{16}$, and finally $\omega>0$
within the range of proposition~\ref{prep:starting-pairs} so that
\begin{equation}\label{asm:omega-choice}
 \ell_d\omega<\frac{\beta}{8}.
\end{equation}
The geometric preparation then gives the strictly positive exponent $e$
in \eqref{asm:fixed-tv}, and the reconstruction below permits any such
$e$.

Fix $T>10+\ell_d\omega$ and apply theorem~\ref{ledger:bound} with local
oscillatory accuracy sufficiently large to leave the error $R^{-T}$ at
the root. Its proof fixes the inverse expansion order, and hence the
persistent symbol derivative budget, followed by the direct expansion
and integration orders, while the extra spatial multipliers have
uniformly bounded modulus after the normalization above and take no
angular derivatives, so that all these choices are fixed before $N$
tends to infinity. For sufficiently large $N$ the right side of
\eqref{asm:assembled-shell} is $O(R^{-\frac{\beta}{2}})$, and since
$\tau_N$ is a positive measure,
$\widehat{\tau_N}(-r)=\overline{\widehat{\tau_N}(r)}$, which proves
\eqref{asm:shell-bound} on both halves of the shell.
\end{proof}

\subsection{Frequency bands of the fixed distance measure}

The estimate \eqref{asm:fixed-tv} holds in total variation, whereas
\eqref{asm:shell-bound} holds in a Fourier $L^2$ norm, and passing the
total variation error directly to the latter norm would cost the square
root of the shell length. Since the target measure has compact support,
we can instead use both estimates in physical space, and the same
observation gives regularity of its density and pointwise decay of its
Fourier transform.

\begin{lemma}[Frequency bands of a fixed compactly supported measure]
\label{asm:fixed-shell}
Let $\alpha$ be a finite positive measure supported in a compact interval
$I\subset\mathbb R$. Suppose there are finite measures $\tau_N$ and
numbers $e,\gamma>0$ such that for every sufficiently large $N$,
\begin{equation}\label{asm:criterion-hypotheses}
 \|\tau_N-\alpha\|_{\mathrm{TV}}\lesssim 2^{-eN}
 \text{ and }
 \int_{2^{N-1}\le |r|\le 2^{N+1}}|\widehat{\tau_N}(r)|^2dr
       \lesssim 2^{-\gamma N}.
\end{equation}
Put $\delta=\min\{e,\frac{\gamma}{2}\}$. Then $\alpha=fdt$ with
$f\in B^\delta_{1,\infty}(\mathbb R)$ and
\begin{equation}\label{asm:pointwise-decay}
 |\widehat\alpha(r)|\lesssim (1+|r|)^{-\delta}.
\end{equation}
Moreover $f\in L^p(\mathbb R)$ for $1\le p<\infty$ satisfying
$1-\frac{1}{p}<\delta$.
\end{lemma}

\begin{proof}
Let $\vartheta$, $Q_N$ and $P_N$ be the cutoff, the Fourier multiplier
and the band cutoff defined in the introduction, so that the symbol of
$P_N$ is supported on $2^{N-1}\le |r|\le2^{N+1}$ and its convolution
kernel is $2^N\psi(2^N\cdot)$ for a fixed Schwartz function $\psi$,
whose $L^1$ norm is therefore independent of $N$. By
\eqref{asm:criterion-hypotheses} and Plancherel,
$$
 \|P_N(\alpha-\tau_N)\|_1\lesssim 2^{-eN} \text{ and } \|P_N\tau_N\|_2\lesssim 2^{-\frac{\gamma N}{2}}.
$$
On the fixed interval $I^*=I+[-1,1]$, Cauchy--Schwarz gives
$\|P_N\tau_N\|_{L^1(I^*)}\le |I^*|^{\frac{1}{2}}\|P_N\tau_N\|_2$.
On its complement we use the compact support of $\alpha$ directly.
For any fixed $L>0$, integration of the convolution kernel gives
$$
 \int_{\mathbb R\setminus I^*}|P_N\alpha(t)|dt
 \le\|\alpha\|_{\mathrm{TV}}\int_{|u|\ge 2^N}|\psi(u)|du
 \lesssim_L 2^{-LN}.
$$
Combining the two regions, we obtain
\begin{equation}\label{asm:band-l1}
 \|P_N\alpha\|_{L^1(\mathbb R)}\lesssim 2^{-\delta N}.
\end{equation}

The low frequency term $Q_0\alpha$ is integrable, and the series
$\sum_{N\ge1}P_N\alpha$ converges absolutely in $L^1$ by
\eqref{asm:band-l1}, since the finitely many initial terms are integrable
as well. The partial sums after adding $Q_0\alpha$ are $Q_N\alpha$, whose
kernels form an approximate identity, so that $Q_N\alpha\to\alpha$ in
distributions and the $L^1$ limit is a density $f$ of $\alpha$. The
Besov condition
$$
 \|Q_0f\|_1+\sup_{N\ge1}2^{\delta N}\|P_N f\|_1<\infty
$$
follows from the same estimate.

For large $|r|$, the band symbols at $r$ sum to one, and only a bounded
number of them can be nonzero, each with $2^N$ comparable to $|r|$.
Since $|\widehat{P_N\alpha}(r)|\le\|P_N\alpha\|_1$, summing at these
indices proves \eqref{asm:pointwise-decay}, while at bounded frequencies
the total mass of $\alpha$ suffices. Finally insert a smooth multiplier
equal to one on the band of $P_N$; Young's inequality with its rescaled
convolution kernel gives the Bernstein inequality
\cite[lemma~2.1]{BahouriCheminDanchin2011},
$$
 \|P_N f\|_p\lesssim_p
 2^{N(1-\frac{1}{p})}\|P_N f\|_1,
$$
whose right side is summable in the stated range of $p$, and since the
low frequency term is in $L^p$ by convolution with a Schwartz function,
summing proves the assertion.
\end{proof}

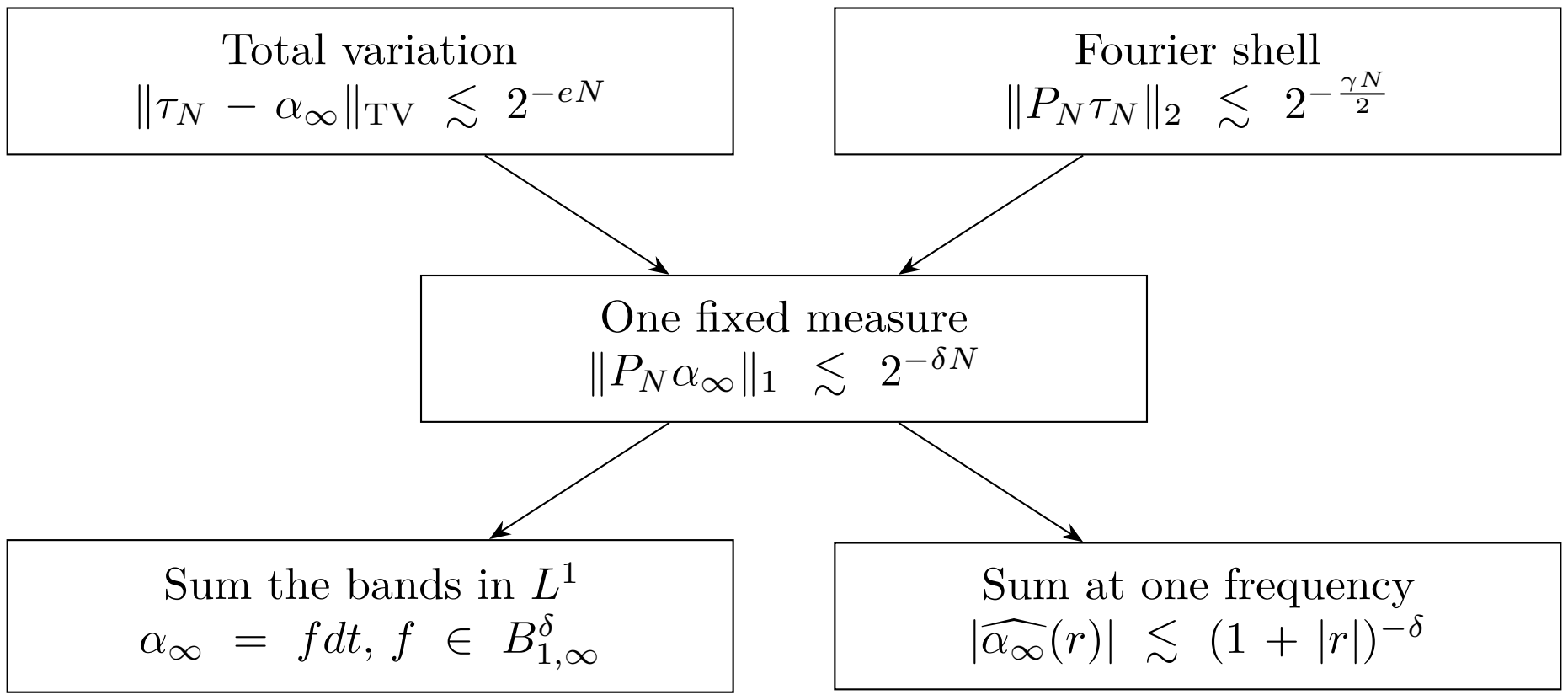
\begin{figure}[htbp]
\centering
\begin{tikzpicture}[>=Stealth,font=\small,
 box/.style={draw,align=center,text width=5.35cm,minimum height=1.15cm,
 inner sep=5pt}]
\node[box] (mass) at (0,0)
 {Total variation\\
 $\|\tau_N-\alpha_\infty\|_{\mathrm{TV}}\lesssim 2^{-eN}$};
\node[box] (shell) at (6.5,0)
 {Fourier shell\\
 $\|P_N\tau_N\|_2\lesssim 2^{-\frac{\gamma N}{2}}$};
\node[box] (band) at (3.25,-2.1)
 {One fixed measure\\
 $\|P_N\alpha_\infty\|_1\lesssim 2^{-\delta N}$};
\node[box] (density) at (0,-4.2)
 {Sum the bands in $L^1$\\
 $\alpha_\infty=fdt$, $f\in B^\delta_{1,\infty}$};
\node[box] (decay) at (6.5,-4.2)
 {Sum at one frequency\\
 $|\widehat{\alpha_\infty}(r)|\lesssim(1+|r|)^{-\delta}$};
\draw[->] (mass)--(band);
\draw[->] (shell)--(band);
\draw[->] (band)--(density);
\draw[->] (band)--(decay);
\end{tikzpicture}
\caption{On a fixed interval containing the distance support, the two
upper estimates control the $L^1$ norm of each band. Outside a unit
enlargement of the support, the Schwartz kernel gives rapid decay.
The same band estimate yields a density and Fourier decay.}
\label{fig:reconstruction}
\end{figure}

\subsection{Completion of the proof}

\begin{proposition}\label{asm:completion}
For the prepared pair of proposition~\ref{prep:starting-pairs},
the positive measure $\alpha_\infty$ in \eqref{asm:limit-measure}
has a density in $B^\delta_{1,\infty}(\mathbb R)$ for some $\delta>0$
and satisfies
$|\widehat{\alpha_\infty}(r)|\lesssim(1+|r|)^{-\delta}$.
Its mass is positive and its support lies in $\Delta(E)$.
\end{proposition}

\begin{proof}
Choose the parameters as in proposition~\ref{asm:shell-decay}.
Lemma~\ref{asm:retained-measures} gives \eqref{asm:fixed-tv} and the mass
and support assertions. By \eqref{asm:distance-range}, the support lies
in the fixed interval $[d_0,d_1]$. The remaining hypothesis of
lemma~\ref{asm:fixed-shell} is \eqref{asm:shell-bound}. Applying the
lemma with $\delta=\min\{e,\frac{\gamma}{2}\}$ gives the conclusion.
\end{proof}

\begin{proof}[Proof of theorems~\ref{thm:main} and~\ref{thm:fourier-distance}]
Choose a Frostman probability on $E$ with exponent strictly greater than
$\frac{d}{2}$. Whenever proposition~\ref{prep:starting-pairs} supplies a
prepared pair, proposition~\ref{asm:completion} gives the desired
measure on the original distances. If the exceptional affine-plane case
occurs instead, the measure $\tau$ constructed in
lemma~\ref{prep:classical-plane} is supported in a compact interval, and
in the dimension $n$ of the plane its Frostman exponent $s$ satisfies
$s>\frac{n+1}{2}$, so that the shell estimate there has right side
$O(R^{\frac{n+1}{2}-s}+R^{-1})$, on both frequency halves by
conjugation. Using $\alpha=\tau$ and $\tau_N=\tau$ in
lemma~\ref{asm:fixed-shell}, with $\gamma=\min\{s-\frac{n+1}{2},1\}>0$,
again gives a nonzero measure with the asserted regularity on the
original distances.

In either case divide the resulting measure by its mass. Its density is
integrable and is supported on $\Delta(E)$, so that this set has positive
Lebesgue measure, and its Fourier transform is $O((1+|r|)^{-\delta})$,
which gives $\dim_F\Delta(E)\ge\min\{1,2\delta\}>0$. All distances used
in these measures are Euclidean distances in the original space.
\end{proof}

\section{From Falconer to separated discrete distances}\label{discrete:section}

What does the theorem say about a large finite set? A finite set has
Hausdorff dimension zero, no matter how many points we put into it, which
is not a promising start. We can, however, replace each point by a small
disk, apply the continuous argument, and then shrink the disks back to
their centers, the size of the disks determining how far apart the
resulting distances will be.

This passage from Falconer to Erd\H{o}s goes back to Iosevich and
\L aba~\cite{IosevichLaba2005}. The energy formulation for large finite
sets was developed by Iosevich, Rudnev, and
Uriarte-Tuero~\cite{IosevichRudnevUriarteTuero2014}. We use that
formulation here, keeping track of separation between the distances.
In the plane it recovers the separated-distance conclusion of Guth,
Solomon, and Wang~\cite[theorem~1.4]{GuthSolomonWang2019} under a more
flexible hypothesis on the points.

The scale of the thickening matters, since every disk has an interval of
distances, however small its radius, and to get information about the
original points we need a lower bound for the length of the thickened
distance set which is uniform in the number of points. This is the
hypothesis in \cite[theorem~3.6]{IosevichRudnevUriarteTuero2014}, and
the quantitative planar theorem of Cheng, Fang, Li, Qu, and Xiao gives
such a bound with an explicit dependence on the energy constant.

\subsection{A uniform length estimate in the plane}

Recall that the $t$-energy of a probability measure $\mu$ on
${\mathbb R}^2$ is
$$
 I_t(\mu)=\iint |x-y|^{-t}d\mu(x)d\mu(y).
$$
\begin{lemma}[Uniform positive length in the plane]
\label{discrete:uniform-length}
Suppose $1<t<2$, $A\ge1$, and $E\subset[-1,2]^2$ is compact and
supports a probability measure $\mu$ with $I_t(\mu)\le A$. For each
$u>\frac{1}{t-1}$ there is a constant $c_{t,u}>0$ such that
$$
 {\mathcal L}^1(\Delta(E))\ge c_{t,u}A^{-u}.
$$
\end{lemma}

\begin{proof}
The potential $U(x)=\int |x-y|^{-t}d\mu(y)$ has $\mu$-average at most
$A$, so that $U\le4A$ on a set of mass at least $\frac{3}{4}$. Choose a
compact subset $K\subset E$ of this set with $\mu(K)\ge\frac{1}{2}$ and
put $\sigma=\frac{\mu|_K}{\mu(K)}$. If $B(z,r)$ meets $K$, take
$x\in B(z,r)\cap K$, and then
$$
 \sigma(B(z,r))\le2\mu(B(x,2r))
 \le 2(2r)^t U(x)\le C_t A r^t.
$$
Let $D=\operatorname{diam} K$, so that $0<D\le3\sqrt{2}$, since a measure
of finite $t$-energy has no atoms. Applying the ball bound at radius $2D$
shows that we can enlarge $C_t$ so that $C_tAD^t\ge1$, and therefore
$\sigma$ satisfies the diameter-normalized bound
$\sigma(B(z,r))\le C_tAD^t(\frac{r}{D})^t$ for $0<r\le D$. Corollary~1.3
of Cheng, Fang, Li, Qu, and Xiao~\cite{ChengFangLiQuXiao2026}, applied
with $\varepsilon=t-1$, gives
$$
 {\mathcal L}^1(\Delta(E))\ge {\mathcal L}^1(\Delta(K))
 \ge c_{t,u}D(C_t A D^t)^{-u}
 \ge c'_{t,u}A^{-u}.
$$
The last inequality uses the upper bound for $D$ and $1-tu<0$.
\end{proof}

\subsection{Thickening the points and counting the distances}

For a finite set $P\subset[0,1]^2$ of $N$ points and $0<s<2$,
write
$$
 \mathcal E_s(P)=\frac{1}{N^2}
       \sum_{\substack{x,y\in P\\x\ne y}}|x-y|^{-s}.
$$
A family is $s$-adaptable if this energy is uniformly bounded and if
distinct points are at least $aN^{-\frac{1}{s}}$ apart for some fixed
$a>0$. Our bound uses the averaged energy alone. For a finite set
$D\subset{\mathbb R}$, let $\mathcal P_r(D)$ be the greatest cardinality
of a subset whose elements are pairwise at least $r$ apart.

For the more general case of unequal masses, let $w_x\ge0$ for $x\in P$,
with $\sum_{x\in P}w_x=1$, and put
$$
 \mathcal E_s(P,w)=\sum_{\substack{x,y\in P\\x\ne y}}
            w_x w_y |x-y|^{-s}
 \text{ and }
 M=\left(\sum_{x\in P}w_x^2\right)^{-1}.
$$
For equal weights $w_x=\frac{1}{N}$ we have $M=N$ and
$\mathcal E_s(P,w)=\mathcal E_s(P)$. The diagonal term in the energy of
the thickened measure explains the appearance of $M$.

\begin{theorem}[Separated distances from a discrete energy bound]
\label{discrete:main}
Suppose $1<s<2$, $A\ge1$, and $P\subset[0,1]^2$ is a finite set with
weights as above such that $\mathcal E_s(P,w)\le A$. For every $1<t\le s$
and every $u>\frac{1}{t-1}$,
\begin{equation}\label{discrete:packing}
 \mathcal P_{M^{-\frac{1}{t}}}(\Delta(P))
 \ge c_{t,u}\left(1+A^{\frac{t}{s}}\right)^{-u} M^{\frac{1}{t}} .
\end{equation}
Consequently, for every $\varepsilon>0$,
\begin{equation}\label{discrete:nearly-linear}
 \mathcal P_{M^{-1}}(\Delta(P))\ge c_{s,A,\varepsilon}M^{1-\varepsilon}.
\end{equation}
In particular an unweighted $N$-point set with $\mathcal E_s(P)\le A$
determines at least $c_{s,A,\varepsilon}N^{1-\varepsilon}$ distances
separated by $N^{-1}$.
\end{theorem}

\begin{proof}
Fix $1<t\le s$ and initially let $0<\delta\le1$. Replace each point by
the uniform probability on its closed disk of radius $\delta$ and form
their weighted average
$$
 d\mu_{P,\delta}(x)=\sum_{z\in P}w_z
        \frac{1_{B(z,\delta)}(x)}{\pi\delta^2}dx .
$$
It is supported on $E_{P,\delta}=\bigcup_{z\in P}\overline B(z,\delta)
\subset[-1,2]^2$. In the energy calculation we separate equal centers
from unequal centers as in
\cite[theorem~2.7]{IosevichRudnevUriarteTuero2014}, allowing the disks
to overlap. This gives
\begin{equation}\label{discrete:thickened-energy}
 I_t(\mu_{P,\delta})\le C_t\left(
         \frac{1}{M\delta^t}+\mathcal E_t(P,w)\right).
\end{equation}
Indeed, scaling two copies of one disk to the unit disk gives a contribution
$C_t\delta^{-t}w_z^2$ from the center $z$, with $C_t$ finite since $t<2$.
For distinct centers at distance $h\ge4\delta$, all distances between their
disks are at least $\frac{h}{2}$. For $0<h<4\delta$ and $x$ in either disk,
the other disk is contained in $B(x,6\delta)$. The average of $|x-y|^{-t}$
on that disk is at most
$$
 \frac{1}{\pi\delta^2}\int_{B(x,6\delta)}|x-y|^{-t}dy
 \le C_t\delta^{-t}\le C_t4^th^{-t}.
$$
This proves \eqref{discrete:thickened-energy} in all cases.

For $t<s$, H\"older's inequality gives
$\mathcal E_t(P,w)\le\mathcal E_s(P,w)^{\frac{t}{s}}$, since the sum of the
off-diagonal weights is at most one. At $t=s$ we use the given bound
directly. With $\delta=M^{-\frac{1}{t}}$ the diagonal term equals one, and
lemma~\ref{discrete:uniform-length} gives
$$
 {\mathcal L}^1(\Delta(E_{P,\delta}))\ge
 c_{t,u}\left(1+A^{\frac{t}{s}}\right)^{-u}=:c_0 .
$$
Choose a maximal $\delta$-separated subset $D_1$ of $\Delta(P)$ and put
$H=\#D_1$. Its radius-$\delta$ intervals cover $\Delta(P)$, and moving
both endpoints by at most $\delta$ changes their distance by at most
$2\delta$, so that
$$
 \Delta(E_{P,\delta})\subset
 \bigcup_{a\in D_1}[a-3\delta,a+3\delta]
 \text{ and therefore }
 c_0\le6\delta H.
$$
This proves \eqref{discrete:packing}, after changing the constant.
For \eqref{discrete:nearly-linear}, take $t\in(1,s)$ close enough to
one that $\frac{1}{t}>1-\varepsilon$ and fix any $u>\frac{1}{t-1}$.
The distances just obtained are at least $M^{-\frac{1}{t}}\ge M^{-1}$
apart.
\end{proof}

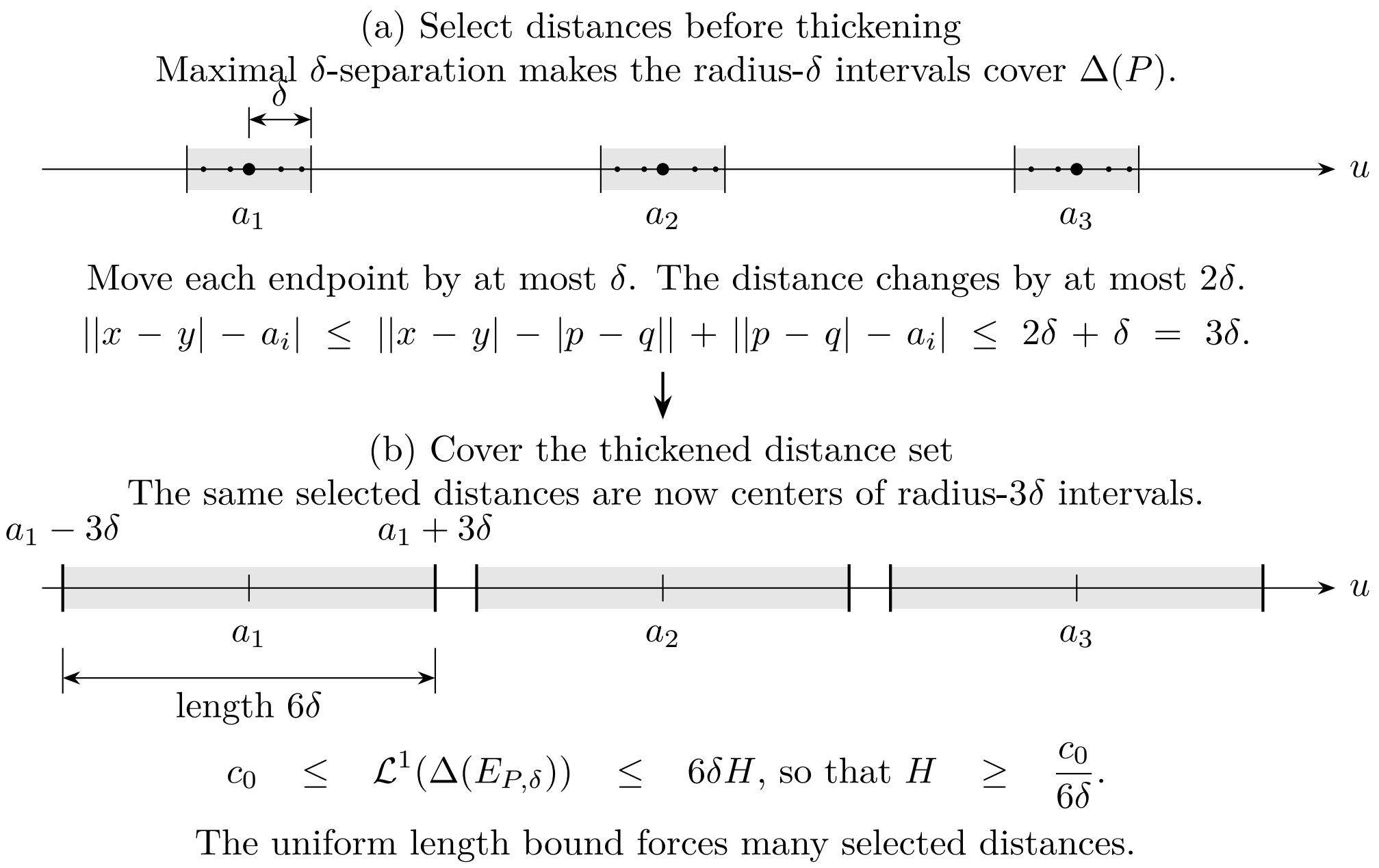
\begin{figure}[!htbp]
\centering
\begin{tikzpicture}[x=1cm,y=1cm,>=Stealth,font=\small]
\node[align=center,text width=13cm] at (6,5.65)
 {(a) Select distances before thickening\\
  Maximal $\delta$-separation makes the radius-$\delta$ intervals cover $\Delta(P)$.};
\foreach \a in {2,6,10} {
 \fill[gray!20] ({\a-.6},4.3) rectangle ({\a+.6},4.7);
 \draw[thin] ({\a-.6},4.27)--({\a-.6},4.73);
 \draw[thin] ({\a+.6},4.27)--({\a+.6},4.73);
 \foreach \b in {-.44,-.18,.31,.51} {
  \fill ({\a+\b},4.5) circle[radius=.027];
 }
 \fill (\a,4.5) circle[radius=.06];
}
\draw[->] (0,4.5)--(12.5,4.5) node[right] {$u$};
\node[below] at (2,4.26) {$a_1$};
\node[below] at (6,4.26) {$a_2$};
\node[below] at (10,4.26) {$a_3$};
\draw[<->] (2,4.98)--node[above] {$\delta$}(2.6,4.98);
\draw[thin] (2,4.8)--(2,5.1);
\draw[thin] (2.6,4.8)--(2.6,5.1);
\node[align=center,text width=13cm] at (6,3.15)
 {Move each endpoint by at most $\delta$. The distance changes by at most $2\delta$.\\[3pt]
  $\displaystyle||x-y|-a_i|
       \le||x-y|-|p-q||+||p-q|-a_i|
       \le2\delta+\delta=3\delta$.};
\draw[->,thick] (6,2.54)--(6,2.08);
\node[align=center,text width=13cm] at (6,1.6)
 {(b) Cover the thickened distance set\\
  The same selected distances are now centers of radius-$3\delta$ intervals.};
\foreach \a in {2,6,10} {
 \fill[gray!20] ({\a-1.8},.25) rectangle ({\a+1.8},.65);
 \draw[thick] ({\a-1.8},.22)--({\a-1.8},.68);
 \draw[thick] ({\a+1.8},.22)--({\a+1.8},.68);
 \draw (\a,.32)--(\a,.58);
}
\draw[->] (0,.45)--(12.5,.45) node[right] {$u$};
\node[above] at (.2,.73) {$a_1-3\delta$};
\node[above] at (3.8,.73) {$a_1+3\delta$};
\node[below] at (2,.22) {$a_1$};
\node[below] at (6,.22) {$a_2$};
\node[below] at (10,.22) {$a_3$};
\draw[thin] (.2,-.13)--(.2,-.57);
\draw[thin] (3.8,-.13)--(3.8,-.57);
\draw[<->] (.2,-.42)--node[below] {length $6\delta$}(3.8,-.42);
\node[align=center,text width=13cm] at (6,-1.62)
 {$\displaystyle c_0\le{\mathcal L}^1(\Delta(E_{P,\delta}))
       \le6\delta H$, so that $\displaystyle H\ge\frac{c_0}{6\delta}$.\\[5pt]
  The uniform length bound forces many selected distances.};
\end{tikzpicture}
\caption{The continuous-to-discrete conversion. The drawing shows
three selected distances, with schematic positions on the distance
axis $u$; the larger dots are the selected $a_i$ and the smaller dots
represent other distances in $\Delta(P)$. Maximality supplies the
first covering, and the triangle inequality supplies the second.
In general there are $H$ covering intervals, which may overlap,
but whose total length is still at most $6\delta H$, while
lemma~\ref{discrete:uniform-length} supplies $c_0$ independently
of the weights and their number when $\delta=M^{-\frac{1}{t}}$.}
\label{fig:discrete-conversion}
\end{figure}

Figure~\ref{fig:discrete-conversion} records the last step of the proof,
in which $D_1=\{a_1,\ldots,a_H\}$ denotes the selected distances. The upper
intervals cover the distances between the centers, whereas the lower
intervals cover all the distances between their thickened disks.

The loss $\varepsilon$ comes from taking $t$ close to one, and both the
admissible power $u>\frac{1}{t-1}$ and the constant from the planar
length theorem depend on this choice.

\subsection{Comparison with Guth, Solomon, and Wang}

Guth, Solomon, and Wang assume that each ball of radius
$N^{-\frac{1}{2}}$ contains at most $K$ points, with $K$ fixed
independently of $N$, and their theorem produces, for every
$\varepsilon>0$, at least $c_{K,\varepsilon}N^{1-\varepsilon}$ distances
separated by $N^{-1}$. For equal weights, \eqref{discrete:nearly-linear} has
the same exponent and the same separation scale, while its hypothesis is an
averaged energy bound at any one exponent $s>1$.

The bounded occupancy hypothesis enters the argument as follows. Choose a
maximal subset $Q\subset P$ with minimum separation
$r=N^{-\frac{1}{2}}$. The radius-$r$ balls centered at $Q$ cover $P$, and
each contains at most $K$ points, so that, with $M=\#Q$, we have
$\frac{N}{K}\le M\le N$. Disjoint disks of radius $\frac{r}{2}$ around the
points of $Q$ show that, for $u\ge r$,
$$
 \#(Q\cap B(x,u))\le C\left(\frac{u}{r}\right)^2=CNu^2.
$$
For a fixed $q\in Q$, split the other points into the annuli
$2^jr\le|q-q'|<2^{j+1}r$, with $j\ge0$ an integer, stopping when $2^jr$
exceeds the diameter of the square. For $0<s<2$ this yields
$$
 \sum_{q'\in Q\setminus\{q\}}|q-q'|^{-s}
 \le C_sN\sum_{\substack{j\ge0\\2^jr\le2}}(2^jr)^{2-s}
 \le C_s'N.
$$
The last sum is bounded because $2-s>0$ and its radii form a
geometric progression ending at a fixed scale. Summing in $q$ and
dividing by $M^2$ proves $\mathcal E_s(Q)\le C_s'K$.
The point separation also satisfies
$r\ge K^{-\frac{1}{s}}M^{-\frac{1}{s}}$. Thus $Q$ satisfies both
requirements of $s$-adaptability, uniformly for each fixed $0<s<2$.
Theorem~\ref{discrete:main}, applied to $Q$, therefore recovers
the separated-distance conclusion for $P$, since $M$ is a fixed
proportion of $N$ and $M^{-1}\ge N^{-1}$.

The passage to $Q$ is not a formality. Bounded occupancy allows a few
points to be arbitrarily close, and those pairs could make the energy of
the original set large, whereas the selection removes this difficulty at a
bounded cost. We have therefore recovered the distance conclusion of Guth,
Solomon, and Wang by the continuous argument, their incidence estimates for
tubes being separate results.

To see how much room the energy hypothesis leaves, consider the
rectangular grids
$$
 P_L=\left\{\left(\frac{i}{L^2},\frac{j}{L}\right)
           :1\le i\le L^2,\ 1\le j\le L\right\}
 \text{ with } N=L^3,
$$
where $L\ge2$ is an integer. These grids have horizontal spacing $L^{-2}$
and vertical spacing $L^{-1}$, and we shall check that they are uniformly
$s$-adaptable for every fixed $1<s\le\frac{3}{2}$. The concentration along
a row makes the bounded occupancy condition fail. In
the figure below, $X\asymp Y$ means that $X$ lies between fixed positive
multiples of $Y$, independently of $L$.

\begin{figure}[!htbp]
\centering
\begin{tikzpicture}[x=1cm,y=1cm,>=Stealth,font=\small]
\node[align=center,text width=6cm] at (2.5,5.75)
 {(a) A square grid\\$64$ points, spacing $\frac{1}{8}$};
\node[align=center,text width=6cm] at (9.5,5.75)
 {(b) The rectangular grid $P_4$\\$64$ points in $4$ rows};
\draw[gray!60] (0,0) rectangle (5,5);
\draw[gray!60] (7,0) rectangle (12,5);
\fill[gray!15] (2.5,2.5) circle[radius=.625];
\fill[gray!15] (9.5,2.5) circle[radius=.625];
\foreach \i in {1,...,8} {
 \foreach \j in {1,...,8} {
  \fill ({5*\i/8},{5*\j/8}) circle[radius=.024];
 }
}
\foreach \i in {1,...,16} {
 \foreach \j in {1,...,4} {
  \fill ({7+5*\i/16},{5*\j/4}) circle[radius=.024];
 }
}
\draw[thick] (2.5,2.5) circle[radius=.625];
\draw[thick] (9.5,2.5) circle[radius=.625];
\draw[->] (2.5,2.5)--(2.94,2.94);
\draw[->] (9.5,2.5)--(9.94,2.94);
\node[anchor=west,fill=white,inner sep=2pt] at (2.9,3.15)
 {$N^{-\frac{1}{2}}$};
\node[anchor=west,fill=white,inner sep=2pt] at (9.9,3.15)
 {$N^{-\frac{1}{2}}$};
\node[anchor=north east] at (0,0) {$0$};
\node[anchor=north east] at (7,0) {$0$};
\node[anchor=south west] at (5,5) {$(1,1)$};
\node[anchor=south west] at (12,5) {$(1,1)$};
\draw[<->] (12.45,1.25)--node[right] {$L^{-1}$}(12.45,2.5);
\draw[thin] (12.1,1.25)--(12.58,1.25);
\draw[thin] (12.1,2.5)--(12.58,2.5);
\node[align=center,text width=13.2cm] at (6,-.75)
 {Both circles have radius $N^{-\frac{1}{2}}=\frac{1}{8}$.
  The lower drawing enlarges one row of a general $P_L$.};
\draw[thick] (2,-2.1)--(10,-2.1);
\foreach \u in {2,3,...,10} {\fill (\u,-2.1) circle[radius=.038];}
\draw[dashed] (1.8,-1.9)--(1.8,-2.35);
\draw[dashed] (10.2,-1.9)--(10.2,-2.35);
\draw[<->] (1.8,-1.65)--node[above] {$2L^{-\frac{3}{2}}$}(10.2,-1.65);
\draw[<->] (5,-2.65)--node[below] {$L^{-2}$}(6,-2.65);
\draw[thin] (5,-2.2)--(5,-2.8);
\draw[thin] (6,-2.2)--(6,-2.8);
\node[align=center,text width=12cm] at (6,-3.65)
 {Points in the row segment $\displaystyle\asymp
       \frac{L^{-\frac{3}{2}}}{L^{-2}}=L^{\frac{1}{2}}=N^{\frac{1}{6}}$};
\end{tikzpicture}
\caption{The two upper panels use the same unit square and the same
number of points. In the rectangular family, a ball of radius
$N^{-\frac{1}{2}}=L^{-\frac{3}{2}}$ centered at an interior point of a row
meets no other row and contains order $L^{\frac{1}{2}}$ points. The
enlarged row is schematic, and its labels indicate the general scales
rather than the value $L=4$ used above. This growing occupancy is
compatible with a bounded averaged $s$-energy when $1<s\le\frac{3}{2}$.}
\label{fig:discrete-rectangles}
\end{figure}
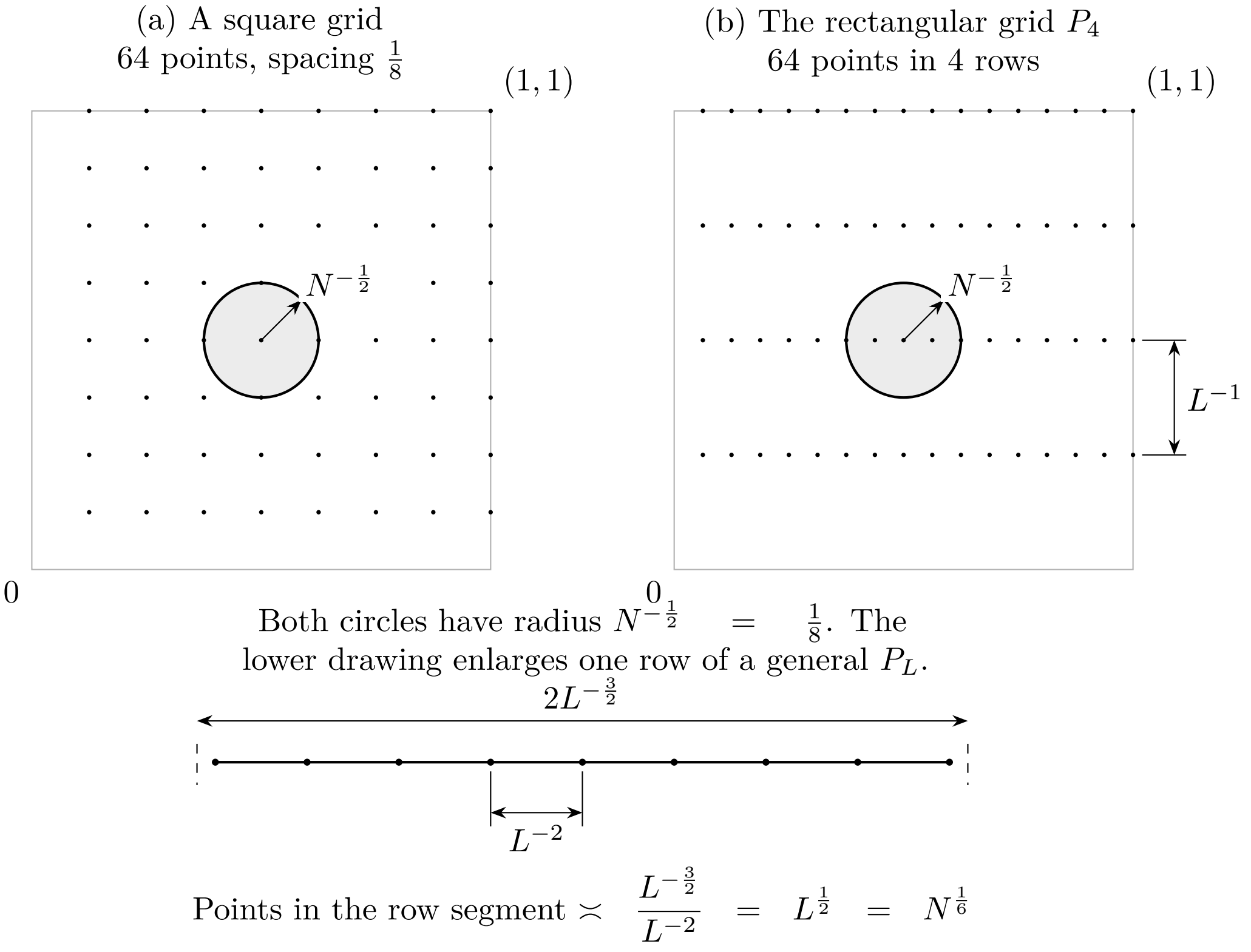

Fix a point of $P_L$. For the points in its own row, the sum of
inverse $s$th powers of distances is at most
$$
 2L^{2s}\sum_{h=1}^{\infty}h^{-s}\le C_sL^{2s},
$$
because $s>1$. For a row at vertical distance $a=\frac{k}{L}$,
where $1\le k\le L$, the corresponding sum is at most
$$
 \sum_{h\in{\mathbb Z}}\left(\left(\frac{h}{L^2}\right)^2+a^2\right)^{-\frac{s}{2}}
 \le a^{-s}+2L^2\int_0^\infty(u^2+a^2)^{-\frac{s}{2}}du
 \le C_sL^2a^{1-s}.
$$
The first inequality compares a decreasing function with its integral,
while in the last inequality we substitute $u=av$, the resulting integral
in $v$ being finite for $s>1$ and $a^{-s}\le L^2a^{1-s}$ because
$a\ge\frac{1}{L}$. Since there are at most two rows at each such vertical
distance, summing their contributions gives
$$
 C_sL^{s+1}\sum_{k=1}^{L}k^{1-s}\le C_s'L^3,
$$
where we used $s<2$. Summing over the $N=L^3$ choices of the fixed
point and dividing by $N^2$ proves
\begin{equation}\label{discrete:rectangular-energy}
 \mathcal E_s(P_L)\le C_s(1+L^{2s-3}).
\end{equation}
This is uniformly bounded for $1<s\le\frac{3}{2}$, and the minimum point
separation is $L^{-2}\ge N^{-\frac{1}{s}}$ throughout this range, so that
both requirements of $s$-adaptability hold, with separation constant
$a=1$.

On the other hand, a ball of radius $N^{-\frac{1}{2}}=L^{-\frac{3}{2}}$
centered well inside a row contains a number of points comparable to
$L^{\frac{1}{2}}=N^{\frac{1}{6}}$, so that no fixed occupancy constant
works for this family. The energy hypothesis therefore permits
configurations outside the homogeneous class, while respecting the minimum
separation required by adaptability.

The thickening calculation extends to balls in any dimension, with
$t$ below the ambient dimension, and its continuous input is a uniform
length bound over the thickened measures. In the odd dimensional
construction of section~\ref{prep:section}, however, the first coarse
scale is chosen using \eqref{prep:uniform-vanishing}, and the later tail
estimate \eqref{prep:coarse-tail} holds only after that choice, so that a
uniform length bound by this route would also require a uniform
treatment of the initial scales.

\section*{Source and changes}

This paper is based on OpenAI's manuscript \cite{OpenAI2026}, repository
revision \texttt{adc7f12}. The geometric state reductions and the
replacement of one depth in a minimizing pair come from proposition~9.4
and lemma~10.2 of that manuscript, and its equation~(5.16) estimates the
mass first discarded at a late coarse scale; summing that estimate gives
\eqref{prep:coarse-tail}, which allows the Fourier estimates to be
applied to one fixed measure in lemma~\ref{asm:fixed-shell}. The
modifications of the argument are the ones described in
section~\ref{sec:changes}: the single potential of
definition~\ref{ledger:potential}, used on every interval and combined
with the chaining lemma~\ref{bik:tree}; the two minimizations of
lemma~\ref{prep:contact}, which give explicit margins for the radial
improvement of lemma~\ref{prep:one-scale} at each finite scale;
lemma~\ref{reg:pruning-multiplication}, which records the full
conditional estimate used in both radial arguments; the sign identity
\eqref{asm:signed-filter}, which reduces the cost of summing the coarse
spatial pieces; the fixed-measure argument leading to
theorem~\ref{thm:fourier-distance}; and the separated-distance estimates
of section~\ref{discrete:section}, derived from a quantitative planar
length bound. The references for the remaining comparisons are given
where they enter the proof.

The retained source material is distributed under the Apache License,
Version 2.0.

\end{document}